\documentclass[a4paper]{article}

\usepackage[utf8]{inputenc}    
\usepackage[T1]{fontenc}
\usepackage[english]{babel}
\usepackage[margin=3.5cm]{geometry}

\usepackage{amsthm}
\usepackage{amsmath}
\usepackage{amssymb}

\usepackage{amsfonts}
\usepackage{mathrsfs}
\usepackage{soul}
\usepackage{graphicx}
\usepackage{enumitem} \setlist[enumerate]{label=(\arabic*)}\usepackage{dsfont}
\usepackage{placeins}
\usepackage{lipsum}
\usepackage{float}
\usepackage{fancyhdr}
\usepackage{tikz-cd}
\usepackage{extarrows}
\usepackage{comment}
\usepackage{adjustbox}
\usepackage{csquotes}
\usepackage{mathtools}
\usepackage{hyperref}
\usepackage{cleveref}
\usepackage{stmaryrd}
\usepackage{tabularx}

\newtheorem{lemma}{Lemma}[subsection]
\newtheorem{defi}[lemma]{Definition}
\newtheorem{prop}[lemma]{Proposition}
\newtheorem{thm}[lemma]{Theorem}
\newtheorem{ex}[lemma]{Example}
\newtheorem{coro}[lemma]{Corollary}
\newtheorem{rk}[lemma]{Remark}
\newtheorem{conj}[lemma]{Conjecture}
\newtheorem{notation}[lemma]{Notation}

\newtheorem{propx}{Proposition}
\newtheorem{thmx}[propx]{Theorem}

\newtheorem*{mainthm}{Main Theorem}

\newtheorem*{defi*}{Definition}
\newtheorem*{prop*}{Proposition}
\newtheorem*{thm*}{Theorem}
\newtheorem*{ex*}{Example}
\newtheorem*{coro*}{Corollary}
\newtheorem*{rk*}{Remark}

\excludecomment{q}
\excludecomment{imprecis}
\excludecomment{todo}

\newcommand \com {\textcolor{orange}}
\newcommand \research {\textcolor{purple}}
\newcommand \comdef {\textcolor{black}}

\renewcommand\to{\longrightarrow}

\newcommand{\N}{\mathbb N}
\newcommand{\Z}{\mathbb Z}

\newcommand{\F}{\mathbb F}
\newcommand{\G}{\mathbb G}

\newcommand{\pro}{\mathbb P}
\newcommand{\Spaces}{\mathcal S}

\newcommand{\X}{\mathcal X}

\newcommand{\g}{\mathfrak g}
\newcommand{\p}{\mathfrak p}

\newcommand{\T}{\mathbb T}
\newcommand{\Fo}{\mathcal F}

\newcommand{\cat}{\mathcal C}
\newcommand{\Cat}{\mathrm{Cat}}
\newcommand{\diag}{\mathrm{diag}}
\newcommand{\FMonoid}{\mathrm{FMonoid}}

\newcommand{\iso}{\xlongrightarrow{\simeq}}

\newcommand{\indp}{\ind_P^G}
\newcommand{\indpf}{\ind_{P_f}^G}

\newcommand{\overl}[1][P]{\overline{#1}}
\newcommand{\la}{\lambda}
\newcommand{\thet}{\theta}
\newcommand{\al}{\alpha}
\newcommand{\M}{\mathcal{M}}
\newcommand{\Nr}{\mathcal{N}}
\newcommand{\Y}{\mathcal{Y}}

\newcommand{\s}{s_{\alpha} \cdot}
\newcommand{\PP}{\mathcal{P}}

\newcommand{\f}{\underline{f}}

\newcommand{\calg}{\mathrm{CAlg}}
\newcommand{\acts}{\curvearrowright}
\newcommand{\Def}{\operatorname{Def}}
\newcommand{\Mod}{\mathrm{Mod}}
\newcommand{\soc}{\mathrm{soc }}

\newcommand{\Hom}{\mathrm{Hom}}
\newcommand{\Dist}{\mathrm{Dist}}

\DeclareMathOperator{\spec}{Spec} 
\DeclareMathOperator{\Spec}{Spec} 
\DeclareMathOperator{\ind}{ind}
\DeclareMathOperator{\res}{res}
\DeclareMathOperator{\Aut}{Aut}
\DeclareMathOperator{\FMP}{FMP}
\DeclareMathOperator{\LieAlg}{LieAlg}

\newcommand{\pliealgd}[2]{\mathrm{LieAlgd}^{\pi}_{#1/#2}}
\newcommand{\pliealg}[1]{\mathrm{LieAlg}^{\pi}_{#1}}
\newcommand{\modstk}{\mathrm{ModuliStk}}

\newcommand{\prstk}{\mathrm{PrStk}}
\newcommand{\prstklaft}{\mathrm{PrStk}^{\mathrm{laft}}}
\newcommand{\prstkdef}{\mathrm{PrStk}^{\mathrm{def}}}
\newcommand{\prstklaftdef}{\mathrm{PrStk}^{\mathrm{laft-def}}}
\newcommand{\art}{\mathrm{art}}

\newcommand{\laftdef}{\mathrm{laft-def}}

\newcommand{\Grp}{\mathrm{Grp}}

\newcommand{\Grpdlaft}{\mathrm{Grpd}^{\mathrm{laft}}}

\newcommand{\Grpd}{\mathrm{Grpd}}
\newcommand{\FGrpd}{\mathrm{FGrpd}}
\newcommand{\FGrp}{\mathrm{FGrp}}
\newcommand{\infAut}{\widehat{\mathrm{Aut}}}
\newcommand{\sqz}{\mathrm{sqz}}
\newcommand{\QC}[1]{\mathrm{QC^{\vee}_{#1}}}
\newcommand{\laft}{\mathrm{laft}}

\newcommand{\Sch}{\mathrm{Sch}}
\newcommand{\Schaff}{\mathrm{Sch}^{\mathrm{aff}}}
\newcommand{\Schaffart}{\mathrm{Sch}^{\mathrm{aff, art}}}
\newcommand{\Schaffnil}{\mathrm{Sch}^{\mathrm{aff, nil}}}

\newcommand{\Schaffaft}{\mathrm{Sch}^{\mathrm{aff, aft}}}

\newcommand{\Schaffop}[1]{(\mathrm{Sch}^{\mathrm{aff}})^{\mathrm{op}}}
\newcommand{\Schaffartop}[1]{(\mathrm{Sch}^{\mathrm{aff, art}}_{#1})^{\mathrm{op}}}
\newcommand{\Schaffnilop}[1]{(\mathrm{Sch}^{\mathrm{aff, nil}}_{#1})^{\mathrm{op}}}

\newcommand{\Monoid}{\mathrm{Monoid}}
\newcommand{\End}{\mathrm{\underline{End}}}
\newcommand{\Act}{\mathrm{Act}}

\newcommand{\coind}{\mathrm{coind}}
\newcommand{\rad}{\mathrm{rad}}
\newcommand{\an}{\mathrm{an}}

\newcommand{\cl}{\mathrm{cl}}
\newcommand{\CRing}{\mathrm{CRing}}

\newcommand{\point}{\mathbin{\scalebox{0.7}{$\bullet$}}}

\newcommand{\sto}{\rightarrow}
\newcommand{\fib}{\mathrm{fib}}
\newcommand{\fibdef}{\mathrm{FibDef}}

\newcommand{\Map}{\mathrm{Map}}

\newcommand{\Stk}{\mathrm{Stk}}
\newcommand{\Aff}{\mathrm{Aff}}
\newcommand{\dR}{\mathrm{dR}}
\newcommand{\Sp}{\mathrm{Sp}}
\newcommand{\id}{\mathrm{id}}
\newcommand{\sch}{\mathrm{Sch}}

\newcommand{\qc}{\mathrm{QC}}

\newcommand{\triv}{\mathrm{triv}}
\newcommand{\ob}{\mathrm{ob}}

\newcommand{\qfrak}{\mathfrak{q}}
\newcommand{\<}{\langle}
\renewcommand{\>}{\rangle}
\newcommand{\Ext}{\mathrm{Ext}}
\newcommand{\kernel}{\mathrm{ker}}
\newcommand{\im}{\mathrm{im}}
\newcommand{\der}{\mathrm{der}}

\newcommand{\Set}{\mathrm{Set}}
\newcommand{\Stab}{\mathrm{Stab}}
\newcommand{\Defcl}{\overline{\Def}}
\newcommand{\act}{\mathrm{act}}

\newcommand{\SL}{\mathrm{SL}}
\newcommand{\Grass}{\mathrm{Gr}}
\newcommand{\Fl}{\mathrm{Fl}}
\newcommand{\supp}{\mathrm{supp}}

\newcommand{\grp}{\mathrm{Grp}}

\newcommand{\nil}{\mathrm{nil}}
\newcommand{\red}{\mathrm{red}}
\newcommand{\contra}{\mathrm{op}}

\DeclareMathOperator*{\colim}{colim}

\newcommand{\deff}{\mathrm{def}}

\renewcommand{\Schaffop}[1]{(\mathrm{Sch}^{\mathrm{aff}}_{#1})^{\mathrm{op}}}

\renewcommand \com [1] {\textcolor{orange}{}}
\renewcommand \research [1] {\textcolor{purple}{}}

\title{Equivariant deformation theory and deformations of homogeneous varieties}

\author{Noé Sotto}
\date{}

\begin{document}

\maketitle
{\let\thefootnote\relax\footnotetext{Laboratoire de Mathématiques d'Orsay, Université Paris-Saclay, France. \\ \textit{Email address}: \href{mailto:sottonoe@gmail.com}{\texttt{sottonoe@gmail.com}}}}

\begin{abstract}

We show that, in characteristic $p$, almost all projective homogeneous varieties under a reductive group and with non-reduced stabilizer do not lift to characteristic $p^2$. To this end, first we develop a general framework for equivariant deformation theory, leading to practical criteria ensuring that the deformations of a scheme come with deformations of a group scheme action. Second, we establish a formula to compute the cohomology of equivariant vector bundles over projective homogeneous varieties with non-reduced stabilizer, enabling us to determine the $G$-module structure of the first cohomology groups of their tangent bundles.

\end{abstract}

\tableofcontents

\addtocontents{toc}{\protect\setcounter{tocdepth}{1}}

%\newpage 

\section{Introduction}

When one studies deformations of an algebraic variety, there is a series of results which guarantee that, in various situations, deformations are \textit{unobstructed}. Let us mention, among others, the Bogomolov-Tian-Todorov theorem \cite{bogomolov1974decomposition,tian1987smoothness,todorov1989weil} which states that, in characteristic $0$, deformations of Calabi--Yau varieties are unobstructed, and was recently extended to mixed characteristic under some conditions in \cite{brantner2025deformationsliftscalabiyauvarieties}; or results by Bott \cite{bott1957homogeneous} and Demazure \cite{Demazure1977} implying that deformations of projective homogeneous varieties under a reductive group are unobstructed; or the fact that Frobenius split varieties in characteristic $p$ lift to characteristic $p^2$ \cite{zdanowicz2016liftabilitysingularitiesfrobeniusmorphism}.

On the other hand, there has been a long-standing effort to find examples of varieties (as ``nice'' as possible, e.g. projective, smooth, Calabi--Yau or Fano...) having obstructed deformations (as obstructed as possible). In \cite{Serre1961}, Serre  gave the first example of a smooth projective variety which does not lift to characteristic $0$. This quest culminated in \cite{vakil2006murphy}, where Vakil proved that any deformation behaviour can happen in general, even for smooth projective varieties (for instance, smooth projective surfaces lifting to $\Z/p^n\Z$ but not to $\Z/p^{n+1}\Z$). However, his constructions are of general type, not Calabi--Yau or Fano.

It can already be useful to prove that a variety lifts to characteristic $p^2$: Deligne and Illusie \cite{deligne1987relevements} proved that, over a field of characteristic $p$, if a variety lifts to $W_2(k)$ and has dimension less than $p$, it has good properties like vanishing of the Hodge-de Rham spectral sequence and Kodaira vanishing. Raynaud \cite{raynaud1978contre} gave the first example of a (smooth projective) variety which violates Kodaira vanishing, and, by \cite{deligne1987relevements}, that does not lift to $W_2(k)$. Hirokado \cite{Hirokado1999} gave the first example of an obstructed Calabi--Yau variety, and Ekedahl \cite{ekedahl2004nonliftable} proved that Hirokado's variety does not lift to $W_2(k)$. In \cite[Example 4]{HaboushLauritzen1993}, Haboush and Lauritzen gave an example of a $6$-fold violating Kodaira vanishing, hence which does not lift to $W_2(k)$ when $p \geq 7$. This variety arises as an \emph{unseparated flag variety}, that is, a projective homogeneous variety whose stabilizer is non-reduced (this situation was not considered by Demazure). Totaro \cite{totaro2017failure} gave another example of an unseparated flag variety violating Kodaira vanishing, and which is moreover Fano. In this example, the dimension condition $\dim X < p$ is not verified so Deligne-Illusie's result cannot be used. Zdanowicz \cite{zdanowicz2018arithmetically} proved by other methods that Totaro's variety does not lift to $W_2(k)$. Moreover, Achinger and Zdanowicz \cite{achinger2017nonliftablecalabiyauvarietiescharacteristic} exhibited an obstructed Calabi--Yau variety as a hypersurface of this unseparated flag variety.

In this work, we prove that in fact almost all unseparated flag varieties do not lift to $W_2(k)$.

\begin{notation}[\Cref{definition parabolic}]
    Let $G$ be a reductive group over an algebraically closed field, let $T$ be a maximal torus of $G$ and $B$ be a Borel subgroup of $G$ containing $T$.
    A (smooth or non-smooth) parabolic subgroup $P$ of $G$ admits a unique decomposition as a product of Frobenius kernels of smooth parabolic subgroups $P = P^0_{r_0} P^1_{r_1} \dots P^l_{r_l}$, 
    where $P^l \subset ... \subset P^1 \subset P^0 = G$ and $0 \leq r_0 < r_1 < \dots < r_{l-1} < r_l = +\infty$. 
    
    The integer $l$ is called the width of $P$. When the Dynkin diagram of $G$ is connected, we say that $P$ is special if the subset of simple roots which generates $P^1$ has cardinality $|S|-1$ and forms a connected subgraph of the Dynkin diagram of $G$.
\end{notation}

\begin{mainthm}[\Cref{conclusion finale}]
    Let $k$ be an algebraically closed field of characteristic $p \geq 5$. Let $G$ be a reductive group with connected Dynkin diagram over $k$, let $P$ be a parabolic subgroup of $G$, and let $X = G/P$. Deformations of $X$ exhibit a dichotomous behaviour:
    \begin{enumerate}
    
    \item If $P$ has width at most $1$, then the variety $X$ has unobstructed deformations.

    \item If $P$ has width at least $2$ and is not special, then the variety $X$ does not lift to $W_2(k)$.
\end{enumerate}
\end{mainthm}

We notice that these unseparated flag varieties lift to $W_2(k)$ if and only if they are Frobenius split. Moreover, their canonical bundle has an explicit expression \cite[3.7]{HaboushLauritzen1993}: it would be interesting to determine which of them are Fano, and which of them give rise to obstructed Calabi--Yau hypersurfaces.
We have not been able to prove that varieties with special stabilizer of width at least $2$ do not lift to $W_2(k)$, but we conjecture that it still holds true.

\paragraph{Strategy of proof.}

Our strategy is inspired by the methods of \cite{zdanowicz2018arithmetically}. Deformations of an algebro-geometric object $X$ can be captured by a deformation functor $\Def_X$ defined over local Artinian rings with residue field $k$, that we call \emph{Artinian extensions} of $k$ (note that when $k$ has characteristic $p$, Artinian extensions need not be $k$-algebras --- most notably $W_2(k)$).
Consider a reductive group $G_k$ over an algebraically closed field $k$. It admits a lift $G_{\Z}$ to $\Z$. For any Artinian extensions $A$ of $k$, let us write $G_A = G_{\Z} \otimes_{\Z} A$. Let $X$ be a projective $G_k$-homogeneous variety with a $k$-point, that is, $X = G_k/P$ where $P$ is a parabolic subgroup of $G_k$. We relate deformations of the variety $X$ to deformations of the stabilizer $P$ of the action of $G_k$ on $X$ through the chain of (classical) deformation functors
\[\begin{tikzcd}
	{\Defcl_{P \sto G_k}^{\Grp}} & {\Defcl_{G_k \acts X}^*} & {\Defcl_{G_k \acts X}} & {\Defcl_{X}}
	\arrow["Q", from=1-1, to=1-2]
	\arrow["p", from=1-2, to=1-3]
	\arrow["U_G", from=1-3, to=1-4]
\end{tikzcd}\]

where \begin{enumerate}
    \item $\Defcl_X$ denotes the functor of flat deformations of $X$ as a scheme.
    \item $\Defcl_{G_k \acts X}$ denotes the functor of flat deformations of the action of $G_k$ on $X$, and the map $U_G$ denotes the natural map which forgets the group action.
    \item $\Defcl_{G_k \acts X}^*$ denotes the functor of \textit{pointed} flat deformations of the action of $G_k$ on $X$ and the map $p$ denotes the natural map which forgets the pointing.
    \item $\Defcl_{P \sto G_k}^{\grp}$ denotes the functor of flat deformations of $P$ \textit{as a subgroup scheme of $G_k$}, and the map $Q$ sends a deformation $\PP$ of $P$ over an Artinian extension $A$, to the quotient homogeneous variety $G_A/\PP$.
    \end{enumerate}

We want to transfer the study of deformations of $X$ to the study of deformations of the subgroup $P$ of $G_k$. For this, we need \textit{unobstructed} maps (also often called ``smooth''), that is, maps which transfer obstructions for one side to another. It is the relative version of the notion of unobstructed deformation functor.

\begin{defi}[\Cref{defi unobstructed map}]
    A map of deformation functors $f \colon F \to G$ is said to be unobstructed if for every surjection $A \to B$ of Artinian extensions of $k$, the induced map $$F(A) \to F(B) \times_{G(B)} G(A)$$ is surjective. 
\end{defi}

We show without much difficulty that the map $Q \colon \Defcl_{P \sto G_k}^{\grp} \to \Defcl_{G_k \acts X}^*$ is in fact an isomorphism, and that the map $p \colon \Defcl_{G_k \acts X}^* \to \Defcl_{G_k \acts X}$ is unobstructed.

\begin{prop}[\Cref{section transitive case}]
    \begin{enumerate}
    \item The quotient map $Q \colon \Defcl_{P \sto G_k}^{\grp} \to \Defcl_{G_k \acts X}^*$ is an isomorphism.

    \item The forgetful map $p \colon \Defcl_{G_k \acts X}^* \to \Defcl_{G_k \acts X}$ is unobstructed.
    \end{enumerate}
\end{prop}

Suppose for a moment that the forgetful map $U_G \colon \Defcl_{G_k \acts X} \to \Defcl_{X}$ is unobstructed. We also show that when $P$ is a non-reduced parabolic subgroup, $P$ does not lift to $W_2(k)$ as a subgroup scheme of $G_k$.

\begin{prop}[\Cref{no lift}]
    A parabolic subgroup $P$ of $G_k$ lifts to $W_2(k)$ as a subgroup scheme of $G_{W_2(k)}$ if and only if $P$ is a reduced parabolic subgroup.
\end{prop}

Finally, by transferring obstructions through unobstructed maps of deformation functors, we get that the homogeneous variety $X = G_k/P$ does not lift to $W_2(k)$.

Now, the hardest part is to know when the forgetful map $U_G \colon \Defcl_{G_k \acts X} \to \Defcl_X$ is unobstructed. In \Cref{section equivariant deformation theory}, we study the map $U_G$: this is what we call equivariant deformation theory. We develop a very general framework to address this question, allowing us, in the end, to obtain practical criteria which ensure that the map $U_G$ is unobstructed. In \Cref{section homogeneous varieties}, we carry out a general study of the cohomology of vector bundles over unseparated flag varieties, allowing us, in the end, to check that those criteria are verified for most unseparated flag varieties. Each of these studies is of independent interest and we think they could find many other applications. \Cref{section equivariant deformation theory} and \Cref{section homogeneous varieties} are completely independent from one another. Notice that \Cref{section equivariant deformation theory} is written in the language of derived algebraic geometry and $\infty$-categories, while \Cref{section homogeneous varieties} is written in the language of classical algebraic geometry.

\subsection{Equivariant deformation theory} 

Let $G$ be a group scheme over $W(k)$ and denote $G_A = G \otimes_{W(k)} A$ for any $W(k)$-algebra $A$.
Consider a $k$-scheme $X$ endowed with an action of $G_k$. When one studies deformations of $X$, one naturally wonders if the $G_k$-action on $X$ admits a deformation when $X$ admits a deformation. An analogous method has been for instance used by Serre \cite{Serre1961} to construct his non-liftable variety, or by Ekedahl \cite{ekedahl2004nonliftable} to prove that Hirokado's variety does not lift to $W_2(k)$. But in both examples, the authors study the deformations of the action of an abstract group, not of a group scheme which itself deforms. On the other hand, Zdanowicz \cite{zdanowicz2018arithmetically} studied the linear version of this question, that is, deformations of a group scheme action on vector bundles over $X$, now taking into account deformations of the group scheme itself.

Linear equivariant deformations can be studied as deformations of vector bundles over the quotient stack $X/G$, making it possible to use the general theory of deformations developed since \cite{Illusie1971CotangentI} and further in ``deformations of objects'' of Lurie \cite[16.5]{lurie2018spectral}. It allowed Zdanowicz to get criteria which ensure that any deformation of a vector bundle over $X$ induces a deformation of the action of $G$.

In the non-linear case, for an abstract group $G$, Rim \cite{rim1980equivariant} gave criteria to ensure existence of a semi-universal $G$-equivariant deformation of a variety in characteristic $0$, in terms of vanishing of $G$-cohomology groups of certain non-explicit $G$-modules. Schröer and Takayama \cite{schroer2018equivariant} extended Rim's results to mixed characteristic deformations of varieties and gave more explicit conditions. Doan \cite{doan2023semiprorepresentabilityformalmoduliproblems} also extended some of Rim's results about equivariant semi-universal deformations in the setting of (derived) formal moduli problems in characteristic $0$.

Our work can be seen either as an extension of Rim's, Schröer-Takayama's, Doan's works when taking into account mixed characteristic deformations of the group scheme $G$ or as a non-linear extension of Zdanowicz's work. We give very tractable criteria which ensure that any deformation of $X$ to an Artinian extension $A$ of $k$ admits an action of $G_A$ which is a deformation of the action of $G_k$ on $X$.

\subsubsection{Equivariant deformation theory: statements of results}

As we explained above, wondering if deformations of $X$ come with deformations of the $G_k$-action amounts to wondering if the forgetful map of deformation functors $$U_G \colon \Defcl_{G_k \acts X} \to \Defcl_X$$ is unobstructed. 
Now, Manetti-Fantechi \cite{fantechi1998obstruction} showed that unobstructedness (or ``smoothness'') can be read off from the induced maps on ``tangent spaces'' and ``obstruction spaces'' when deformation functors satisfy Schlessinger's conditions. But in the classical setting, it is hard to build canonical and functorial obstruction spaces associated to deformation functors. The modern way to circumvent this issue is to consider derived deformation functors, called formal moduli problems after the work of Lurie \cite{lurie2011formal}. Those are functors from simplicial (or ``animated'') Artinian extensions of $k$ to spaces (or $\infty$-groupoids) satisfying a derived version of Schlessinger's conditions. They carry a canonical and functorial obstruction theory, in the form of an associated tangent complex (a chain complex). The classical tangent space and obstruction space correspond to the first homotopy groups of this tangent complex. In this setting, we can reformulate Fantechi-Manetti conditions as follows.

\begin{prop}[\Cref{criterion smooth}] \label{intro criterion unobstructed}
    Let $F \colon \M \to \Nr$ be a map of formal moduli problems. Let $Q$ be the fiber of the associated map of tangent complexes $T_{\M} \to T_{\Nr}$ in $\Mod_k$. If $\pi_{-1}Q = 0$, then the map $F$ is unobstructed.
\end{prop}

The classical deformation functors $\Defcl_X$ and $\Defcl_{G_k \acts X}$ have natural enhancements to formal moduli problems, that we denote $\Def^{\sch}_X$ and $\Def^{\sch}_{G_k \acts X}$.
Then, we need to study the map of tangent complexes associated to the map of formal moduli problems $$U_G \colon \Def_{G_k \acts X} \to \Def_X.$$

It has long been known that the tangent complex of the formal moduli problem $\Def^{\sch}_X$ identifies with the shifted complex of (derived) global vector fields $R\Gamma(X, T_X)[1]$. We show that the tangent complex of the equivariant formal moduli problem identifies with the space of (derived) equivariant global vector fields.

\begin{thmx}[\Cref{equivariant deformations as deformations over BG}] \label{intro tangent equivariant}
    Let $G$ be a smooth group scheme over $W(k)$. Let $X$ be a proper derived scheme locally of finite presentation over $k$ endowed with a $G_k$-action. The action of $G_k$ on $X$ induces an action of $G_k$ on $R\Gamma(X, T_X)$ and there is an equivalence in $\Mod_k$
    $$T_{\Def^{\sch}_{G_k \acts X}} \simeq R\Gamma(X, T_X)^{G_k}[1],$$
    where $(-)^{G_k}$ denotes the derived fixed points functor. Moreover, the tangent map of the forgetful transformation $U_G$ identifies with the natural map $$R\Gamma(X, T_X)^{G_k}[1] \to R\Gamma(X, T_X)[1].$$
\end{thmx}

\begin{rk} \label{rk base change}
    In fact, we prove a stronger statement. We show that the action of $G_k$ on $X$ even induces a natural action of $G$ on the formal moduli problem $\Def_X$, in the sense that there is a geometric object $\M$ with a map $\M \to BG$  such that the pullback along the point $\Spec(k) \to BG$ identifies with $\Def^{\sch}_X$ (hence, $\M$ is a kind of $G$-torsor). The object $\M$ identifies with deformations of $X/G$ over the base $BG$ and under the base $BG_k$, $\Def_{X/G}^{BG_k/BG}$ (see \cite{gaitsgory2019study} and \cite{brantnermagidsonnuiten2025} for the theory of formal moduli stacks under/over a base). From another point of view, it means that the action of $G_k$ on $R\Gamma(X, T_X)$ is compatible with the Lie algebra structure. 
    Moreover, the formal moduli problem $\Def^{\sch}_{G_k \acts X}$ then identifies with the pushforward of $\M$ along the map $BG \to \Spec(k)$. In other words, the formal moduli problem $\Def^{\sch}_{G_k \acts X}$ directly identifies with the $G_k$-fixed points of the formal moduli problem $\Def^{\sch}_X$: $$\Def^{\sch}_{G_k \acts X} \simeq (\Def^{\sch}_X)^{G_k}$$
\end{rk}

Using the criterion of \Cref{intro criterion unobstructed}, we deduce from \Cref{intro tangent equivariant} criteria which ensure that the deformations of $X$ come with a deformation of the action of $G$.

\begin{thmx}[\Cref{criterion}]
    Let $G$ be a smooth group scheme over $W(k)$, let $X$ be a smooth \comdef{proper} scheme over $k$ endowed with an action of $G_k$. The $G_k$-action induces an action on $T_X$ and on $R\Gamma(X, T_X)$. If the following conditions are verified:
    \begin{enumerate}
    \item $H^1(G_k, k) = 0$

    \item $H^2(G_k, H^0(X, T_X)) = 0$ 
    \item $H^1(X, T_X)$ is a trivial $G_k$-module,
    \end{enumerate}
    then the forgetful map $U_G \colon \Def^{\sch}_{G_k \acts X} \to \Def^{\sch}_X$ is unobstructed.
\end{thmx}

\subsubsection{Equivariant deformation theory: strategy of proof}

Let $G$ be a smooth group scheme over $W(k)$. We place the study of the map $$U_G \colon \Def^{\sch}_{G_k \acts X} \to \Def^{\sch}_X$$ in a more general framework. The data of a $G_k$-action on $X$ corresponds to the data of the classifying map (of Artin stacks) $X/G_k \to BG_k$ in the etale topology.
Hence, we will more generally study deformations of maps $Y \to W'_k$ with fixed deformations of the target ($W'_k$ plays the role of $BG_k$). $W'$ (which plays the role of $BG$) will be a $W(k)$-prestack and we write $W'_A = W' \otimes_{W(k)} A$ for all $W(k)$-algebra $A$. In the case of schemes, this is a classical deformation problem, studied for instance in \cite[3.4.2]{sernesi2007deformations}, where Sernesi computes its associated classical ``tangent space'' and ``obstruction space''. Though, to be able to apply this setting to classifying stacks, which are Artin stacks, we choose to place ourselves in a very general setting (the most general possible in some sense): (locally almost of finite type) prestacks having deformation theory (see \Cref{defi deformation theory}), a category which encompasses Artin stacks.

\begin{defi}[\Cref{defi deformation of maps mixed}]
    Let $Y \in \prstklaftdef_{W'_k}$. 
    Let $\Def^{k/W(k)}_{Y \to W'_k}$ be the functor which sends an Artinian extension $A$ of $k$ to the space of cartesian squares

\[\begin{tikzcd}
	Y & W'_k \\
	\Y & W'_A
	\arrow[from=1-1, to=1-2]
	\arrow[from=1-1, to=2-1]
	\arrow["\lrcorner"{anchor=center, pos=0.125}, draw=none, from=1-1, to=2-2]
	\arrow[from=1-2, to=2-2]
	\arrow[from=2-1, to=2-2]
\end{tikzcd}\]
    It is called the functor of deformations of the map $Y \to W'_k$ with fixed deformations of the target.
\end{defi}

We extend Sernesi's results by computing the tangent complex of this formal moduli problem.

\begin{prop}[\Cref{pushforward mixed deformations}] \label{intro pushforward mixed deformations}
    Let $Y \in \prstklaftdef_{W'_k}$. Write $f \colon Y \to W'_k$ and $g \colon W'_k \to \Spec(k)$. Then the functor $\Def^{k/W(k)}_{Y \sto W'_k}$ is a formal moduli problem and its tangent complex is given by  $$T_{\Def^{k/W(k)}_{Y \sto W'_k}} \simeq g_* \circ f_* T_{Y \sto W'_k}[1].$$ 
\end{prop}

By applying \Cref{intro pushforward mixed deformations} to the case where the map $Y \to W'_k$ is the classifying map $X/G_k \to BG_k$, we recover \Cref{intro tangent equivariant}.

\begin{rk}
    In fact, we prove the following stronger statement. We show that the formal moduli problem (under $\Spec(k)$) $\Def^{k/W(k)}_{Y \sto W'_k}$ identifies with the pushforward along $g$ of the deformation functor (under $W'_k$) $\Def_Y^{W'_k/W'}$:$$\Def^{k/W(k)}_{Y \sto W'_k} \simeq g_* \Def^{W'_k/W'}_Y.$$ 
Moreover, if one considers a point $\Spec(k) \to W'_k$ and writes the fiber $X = Y \times_{W'_k} \Spec(k)$, then the formal moduli problem (under $\Spec(k)$) $\Def_X$ identifies with the pullback along the point of the deformation functor (under $W'_k$) $\Def_Y^{W'_k/W'}$: $$\Def^{k/W(k)}_X \simeq e^* \Def_Y^{W'_k/W'}.$$
    Applying that to the case where the map $Y \to W'_k$ is the classifying map $X/G_k \to BG_k$, we recover \Cref{rk base change}
\end{rk}

\subsection{Equivariant vector bundles over homogeneous varieties with non-reduced parabolic stabilizer: cohomology and irreducibility}

Let $G$ be a reductive group over an algebraically closed field $k$, and let $X$ be a projective $G$-homogeneous variety with a $k$-point. The stabilizer of the action is a parabolic subgroup of $G$. It is usually supposed that parabolic subgroups are smooth: in this case, $X$ is called a (generalized) flag variety, and its properties have been thoroughly studied. In any characteristic, cohomology of equivariant vector bundles can be computed on the complete flag variety $G/B$, where equivariant vector bundles always admit a filtration with line bundles. 

In characteristic $0$, the Borel-Weil-Bott theorem gives the cohomology of line bundles and states that simple $G$-modules correspond to global sections of line bundles associated to dominant weights. Moreover, all $G$-modules are semi-simple and their higher $G$-cohomology vanishes.
In positive characteristic, the situation is more complicated, as explained in the reference book of Jantzen \cite{jantzen2003representations}. In general, $G$-modules are not semi-simple, their higher cohomology does not vanish, and the Borel-Weil-Bott theorem fails. However, simple $G$-modules are still classified by their highest weight, Kempf vanishing holds, and some partial results about cohomology of line bundles exist.

However, the stabilizer of the action of $G$ can also be a non-reduced parabolic subgroup of $G$. In that case, $X = G/P$ is called an \emph{unseparated flag variety}. Non-reduced parabolic subgroups were first classified by Wenzel \cite{Wenzel1993} and unseparated flag varieties have been studied in more depth by Haboush and Lauritzen \cite{HaboushLauritzen1993}. In that case it is not possible to transfer the study of cohomology from $X = G/P$ to $G/B$ anymore, vector bundles do not admit a filtration with line bundles and Kempf vanishing fails in general (\cite{HaboushLauritzen1993}, Example 3). Thus, the cohomology of vector bundles over $X$ has remained quite mysterious.

In this work, we establish a formula to compute the cohomology of equivariant vector bundles over unseparated flag varieties. We establish a filtration of the tangent bundle of $X$ whose factors are Frobenius twists of tangent bundles of usual flag varieties. Then, we show that the $G$-module structure of $H^0(X, T_X)$ and $H^1(X, T_X)$ is very simple. However, we show for instance that in general $H^2(X, T_X) \neq 0$, in contrast with the case of usual flag varieties (\cite{bott1957homogeneous} in characteristic $0$, \cite{Demazure1977} in characteristic $p$).

Non-reduced parabolic subgroups also give rise to non-reduced reductive subgroups, as a kind of Levi subgroups. We show that, as for classical reductive groups and parabolic subgroups, simple representations of non-reduced parabolic subgroups or of non-reduced reductive subgroups are classified by appropriate dominant weights. We also prove that the Steinberg decomposition holds for those non-reduced reductive groups.

\subsubsection{Homogeneous varieties with non-reduced parabolic stabilizer : statement of results}

The (homotopy) pullback square of quotient stacks%
\[\begin{tikzcd}
	{G/P} & {*} \\
	BP & BG
	\arrow[from=1-1, to=1-2]
	\arrow[from=1-1, to=2-1]
	\arrow["\lrcorner"{anchor=center, pos=0.125}, draw=none, from=1-1, to=2-2]
	\arrow[from=1-2, to=2-2]
	\arrow[from=2-1, to=2-2]
\end{tikzcd}\]
shows by the base change formula that an equivariant vector bundle $M$ over $X = G/P$ corresponds to a $P$-module, and that the cohomology groups $H^i(G/P, M)$ are naturally endowed with a $G$-module structure given by the induction $R^i\ind_P^G(M)$.

Let $P$ be a parabolic subgroup of width $l$. Write $P = P^0_{r_0}P^1_{r_1}\dots P^l_{r_l}$ with $r_0 < \dots < r_{l-1} < r_l$, and where $P^k$ is a reduced parabolic subgroup for all $0 \leq k \leq l$ associated to the subset of simple roots $I_k \subset S$.

\begin{defi}[\Cref{defi Frobenius}]
    Let $G$ be an algebraic group. Denote by $G_r$ the kernel of the $r$-th iterated Frobenius endomorphism of $G$. We will write $F^r_G \colon G \to G/G_r$. Thus, if $M$ is a $G/G_r$-module, the $G$-module $(F^r_G)^*(M)$ denotes the Frobenius twist of $M$.
\end{defi}

We begin by establishing a formula which expresses the cohomology of a Frobenius-twisted vector bundle as the composite of inductions between the parabolic subgroups $P^k$ and Frobenius twists. Notice that induction from $P^{k+1}$ to $P^k$ can be computed at each step simply as the induction from the Borel subgroup $B$ to $P^k$ (after restricting to $B$). Hence it is enough to know the $B$-module structure at each step of the composite.

\begin{thmx}[\Cref{induction}] \label{intro induction formula}
    Let $M$ be a $P^l$-module. Let $r \geq r_{l-1}$ (with the convention $r_{-1} = 0$). 
    Then there is an isomorphism of $G$-modules $$\ind_{P}^{G}  (F^r_{P_f})^* (M) \simeq (F_{P^0}^{r_0})^*  \ind_{P^1}^{P^0}  (F_{P^1}^{r_1-r_0})^*  \ind_{P^2}^{P^1}  \dots  (F_{P^{l-1}}^{r_{l-1} - r_{l-2}})^*  \ind_{P^{l}}^{P^{l-1}}  (F_{P^{l}}^{r - r_{l-1}})^* (M)$$
\end{thmx}

Whereas in the reduced case, the tangent bundle can be studied thanks to a filtration with $1$-dimensional factors (see \cite{bott1957homogeneous} and \cite{Demazure1977}), that is not possible in the non-reduced case. We decompose the map $X \to \Spec(k)$ into a sequence of maps corresponding to ``width-$1$'' situations. In a ``width-$1$'' situation, unseparated flag varieties are simply given by the Frobenius. Hence, we get a filtration of the tangent bundle of $X$ by Frobenius twists of tangent spaces of classical quotient varieties.

\begin{defi}
    For all $0 \leq k \leq l$, denote $$\overl^k = P^0_{r_0}\dots P^{k-1}_{r_{k-1}}.$$
    There is a filtration $P = \overl^l \subset \overl^{l-1} \subset \dots \subset \overl^1 \subset \overl^0 = G$ of parabolic subgroups.
\end{defi}

\begin{thmx}[\Cref{filtration}]
    Consider the chain $$X = G/P \to G/\overl^{l-1} \to \dots \to G/\overl^1 \to G/\overl^0 = *.$$ The sequence $(T_{G/P \sto G/\overl^k})_{0 \leq k \leq l}$ forms a filtration of $T_X$ whose factors are 
    $$\res_P^{\overl^{k+1}} (F^{r_k}_{\overl^{k+1}})^* (\frak{p}^k/\frak{p}^{k+1})$$ for $0 \leq k \leq l-1$.
\end{thmx}

The cohomology of the factors of this filtration can be computed using \Cref{intro induction formula}. With some additional work, we deduce the following result about the $G$-module structure of the first cohomology groups of the tangent bundle.

\begin{thmx}[\Cref{conditions fulfilled}]
    Suppose that the Dynkin diagram of $G$ is connected and that $P$ is a non-special parabolic subgroup of $G$. Let $\thet$ be the highest root of $G$.
Then the following assertions hold:
\begin{enumerate}
\item $H^1(X, T_X)$ is a trivial $G$-module
\item for all $i \geq 2$, $H^i(G, H^0(X, T_X)) \simeq H^i(G, (F_G^{r_0})^*H_G^0(\thet)).$  
\item if $r_0 = 0$, then for all $i \geq 2$, $H^i(G, H^0(X, T_X)) = 0$.
\end{enumerate}
\end{thmx}

Write $G^k$ for the Levi subgroup of $P^k$ for all $0 \leq k \leq l$. To any parabolic subgroup $P$ is associated its Levi subgroup $G_f = G^0_{r_0} \dots G^l_{r_l}$, which is a ``non-reduced reductive'' subgroup of $G$, in the sense that it contains no connected unipotent normal subgroup. There is also a filtration $G_f = \overl[G]^l \subset \overl[G]^{l-1} \subset \dots \subset \overl[G]^1 \subset \overl[G]^0 = G$. We prove that such non-reduced reductive subgroups have two properties analogous to the reduced case: correspondence between simple representations and $G_f$-dominant weights, and Steinberg decomposition.

\begin{defi}[\Cref{defi dominant weights}]
    Let $I \subset S$ and $r \geq 1$.
    Let $X_+^I(T)$ be the set of $I$-dominant weights $$X_+^I(T) := \{\la \in X(T), \forall \alpha \in I, \<\la, \al^{\vee}\> \geq 0 \}.$$
    
    Let $X_r^I(T)$ be the subset of $p^r$-restricted $I$-dominant weights
    $$X_r^I(T) := \{\la \in X(T), \forall \alpha \in I, 0 \leq \<\la, \al^{\vee}\> < p^r \} \subset X_+^I(T).$$
\end{defi}

\begin{thmx}[\Cref{BW theorem non-inf} + \Cref{steinberg}]
    \begin{enumerate}
    \item Up to isomorphism, simple $G_f$-modules are in bijection with the set $X^{I_l}_+(T)$. In one direction, it sends a simple $G_f$-module $M$ to its highest weight, and in the other direction it sends a $G_f$-dominant weight $\la$ to the simple $G_f$-module $$L_{G_f}(\la) := \soc_{G_f} \ind_{B_f}^{G_f}(\la).$$
    \item Let $\la = \la_0 + \la_1 p + \dots + \la_n p^n \in X(T)$. For all $r \in \llbracket 0, n \rrbracket$, denote by $k_r$ the unique integer such that $r_{k_r-1} \leq r < r_{k_r}$ (with convention $r_{-1} = 0$) and suppose that $\la_r \in X_1^{I_{k_r}}(T)$.
    Then there is an isomorphism of $G_f$-modules
    $$L_{G_f}(\la_0+\dots +\la_n p^n) \simeq L_{G_f}(\la_0) \otimes (F^1_{G_f})^* L_{G_f/(G_f)_1}(\la_1) \otimes \dots \otimes (F^n_{G_f})^* L_{G_f/(G_f)_n}(\la_n)$$ 
    which can be reformulated as 
    $$L_{G_f}(\la_0+\dots +\la_n p^n) \simeq \res_{G_f}^{G^0} L_{G^0}(\la_0) \otimes \res_{G_f}^{\overl[G]^{k_1}} (F^1_{\overl[G]^{k_1}})^* L_{G^{k_1}}(\la_1) \otimes \dots \otimes \res_{G_f}^{\overl[G]^{k_n}} (F^n_{\overl[G]^{k_n}})^* L_{G^{k_n}}(\la_n).$$ 
    \end{enumerate}
\end{thmx}

\paragraph{Acknowledgements.} I am deeply grateful to my PhD advisor Lukas Brantner for his constant guidance throughout this project, and to my co-advisor Yonatan Harpaz for his insightful suggestions. I thank Piotr Achinger for his careful reading of the manuscript and for helping fill a gap in one of the arguments. I have also benefited greatly from numerous discussions with Hugo Pourcelot.

\addtocontents{toc}{\protect\setcounter{tocdepth}{2}}

\section{Formal moduli stacks, base change and unobstructedness: building equivariant deformation theory} \label{section equivariant deformation theory}

Let $k$ be a field, let $\Lambda$ be a local Noetherian ring with residue field $k$ (for instance $\Lambda = W(k)$ when $k$ has characteristic $p$). Let $G$ be a group scheme over $\Lambda$. Suppose that its reduction to $k$, $G_k$, acts on a variety $X$. When $X$ admits a deformation $\X$, it is natural to wonder if the action of $G_k$ admits a deformation to an action on $\X$. Authors such as \cite{rim1980equivariant} and \cite{schroer2018equivariant} have stated conditions (on some cohomology groups) under which such an equivariant deformation is possible. However, Rim works in characteristic $0$ with trivial deformations of $G_k$, and Schröer-Takayama work in mixed characteristic but only with an abstract group $G$. 
One of our motivations was to extend their results to the case of the mixed characteristic deformations of a group scheme action, which we finally do in \Cref{criterion} (in a derived algebraic geometry setting).

In the setting of deformation functors, the above question can be reformulated as: when is the forgetful transformation from the $G_k$-equivariant deformations of $X$ to the deformations of $X$ $$U_G \colon \Def^{\sch}_{G_k \acts X} \to \Def^{\sch}_X$$ an unobstructed map of deformation functors?

In \Cref{section unobstructed maps}, we discuss a folklore result giving a criterion for a map of classical deformation functors to be unobstructed using its ``tangent space'' and its ``obstruction space''. But to define properly a functorial ``obstruction space'', we need to work in the derived setting of formal moduli problems, that come with a fully functorial tangent complex. We reformulate the folklore result as a vanishing condition on one homotopy group of the fiber of the tangent map.
To use this criterion, we need to compute the tangent map of the forgetful map of formal moduli problems $$U_G \colon \Def^{\sch}_{G_k \acts X} \to \Def^{\sch}_X.$$ In \Cref{section deformations over a base}, we show that this map is induced by the natural transformation which sends, for each Artinian extension $A$ of $k$, the classifying map of a deformation $\X/G_A \to BG_A$, which fits in the pullback diagram
\[\begin{tikzcd}
	\X & {\Spec(A)} \\
	{\X/G_A} & BG_A
	\arrow[from=1-1, to=1-2]
	\arrow[from=1-1, to=2-1]
	\arrow["\lrcorner"{anchor=center, pos=0.125}, draw=none, from=1-1, to=2-2]
	\arrow[from=1-2, to=2-2]
	\arrow[from=2-1, to=2-2]
\end{tikzcd}\]
to its fiber $\X \to \Spec(A)$.

With this point of view, the above question naturally generalizes to a larger setting: let $Y \to K'$ be a map of general algebro-geometric spaces with a section $\Spec(k) \to K'$, let $X = Y \times_{K'} \Spec(k)$ denote the associated fiber. One can define a natural transformation $$\fibdef_Y^{k/\Lambda} \colon \Def^{k/\Lambda}_{Y \to K'} \to \Def^{k/\Lambda}_X$$ from the formal moduli problem of the deformations of the map $Y \to K'$ to the formal moduli problem of the deformations its fiber $X$. Then, what is the tangent map associated to the transformation $\fibdef$? We ultimately answer this question in \Cref{fiber of deformation tangent equi} and \Cref{fiber of deformation tangent mixed}. Applied to the case of a $G$-action, we deduce in \Cref{equivariant deformations as deformations over BG} that the tangent map of the forgetful transformation $U$ is, in a sense, the forgetful map from (derived) equivariant global vector fields of $X$ $$R\Gamma(X, T_X)^{G_k} \to R\Gamma(X, T_X).$$ 

In \Cref{section unobstructed equivariant deformations}, we study this tangent map in detail to extract very applicable criteria that ensure unobstructedness of the forgetful map $U_G \colon \Def^{\sch}_{G \acts X} \to \Def^{\sch}_X$, analogous to Schröer-Takayama's results in \cite{schroer2018equivariant} for abstract group actions or to Zdanowicz's results in \cite{zdanowicz2018arithmetically} for deformations of equivariant vector bundles.

Finally, in \Cref{section transitive case}, we focus on the case of a transitive $G$-action, that is, the case where $X$ is a $G$-homogeneous variety. In that situation, pointed equivariant deformations $\Defcl_{G_k \acts X}^*$ can be reformulated as deformations of the stabilizer of this action as a subgroup scheme of $G$, that we write $\Defcl^{\Grp}_{P \to G_k}$. We can then write a chain of maps of formal moduli problems 

\[\begin{tikzcd}
	{\Defcl^{\Grp}_{P \sto G_k}} & {\Defcl_{G_k \acts X}^*} & {\Defcl_{G_k \acts X}} & {\Defcl_{X}}
	\arrow["\simeq", from=1-1, to=1-2]
	\arrow[from=1-2, to=1-3]
	\arrow[from=1-3, to=1-4]
\end{tikzcd}\]
whose composite is unobstructed when $X$ is smooth, and when the conditions of \Cref{section unobstructed equivariant deformations} are satisfied. In that case, it implies that the deformations of a homogeneous $X$ are completely determined by the deformations of the stabilizer $P$ of the action of $G_k$ on $X$, which are generally simpler to study.

\subsection{Formal moduli stacks, base change and equivariant deformations} \label{section deformations over a base}

The main motivation of this section is to study abstractly the formal moduli problem of equivariant deformations of a variety. Let $X$ be a variety endowed with an action of an algebraic group $G_k$, over a field $k$. We call equivariant deformations of $X$ the deformations of $X$ together with a deformation of the action of $G_k$. This question goes back at least to \cite{pinkham1974deformations} and \cite{rim1980equivariant}.

In some sense, this question is the non-linear version of the study of deformations of $G$-equivariant vector bundles over $X$. But this latter case is simpler because it can be reformulated as the study of deformations of vector bundles over the quotient stack $X/G$, and falls under the general theory of ``deformations of objects'' of \cite[16.5]{lurie2018spectral} for instance.

It is known that deformations of $X$ are encoded in the derived global vector fields $R\Gamma(X, T_X)$ in $\Mod_k$, endowed additionally with a Lie algebra structure. When $X$ is endowed with a $G$-action, this action induces an action on the (derived) global vector fields $R\Gamma(X, T_X)$. In other words, $R\Gamma(X, T_X)$ lifts to an object of $\qc_{BG}$, where $BG$ denotes the etale classifying stack of $G$. Now one may naturally wonder if this $G$-module structure is compatible with the Lie algebra structure on $R\Gamma(X, T_X)$. If it is the case, then this object would be a natural candidate to encode equivariant deformations of $X$. This is what we conclude in the final result of this section, \Cref{equivariant deformations as deformations over BG}. 

This question naturally pushes us in the setting of formal moduli problems \textit{over} (or rather \textit{under}) a base, developed by \cite{gaitsgory2019study} (here the base will be $BG$). Whereas usual formal moduli problems correspond to infinitesimal thickenings of the field $K = \Spec(k)$, formal moduli problems under a base $K'$ correspond to infinitesimal thickenings of $K'$. At least in characteristic $0$, they are encoded in Lie algebras over $K'$ (more precisely, Lie algebra objects in ind-coherent or pro-coherent sheaves over $K'$). 

The correspondence between a $G_k$-action on $X$ and the classifying map $X/G_k \to BG_k$ gives a correspondence between equivariant deformations of $X$ and deformations of the classifying map $X/G_k \to BG_k$. Thus, we will focus on the study of deformations of maps $Y \to K'$ where $K'$ is a base which can be an Artin stack (this formal moduli problem is denoted $\Def^{k/k}_{Y \to K'}$). In \Cref{tangent deformation of map}, we show that $\Def^{k/k}_{Y \to K'}$ arises as a kind of pushforward of the formal moduli problem $\Def_Y^{K'/K'}$ of deformations of $Y$ over $K'$, which is a formal moduli problem under $K'$. Applied to the case of the map $X/G_k \to BG_k$, it says that the formal moduli problem of equivariant deformations of $X$ is the pushforward of the formal moduli problem $\Def^{BG_k/BG_k}_{X/G_k}$ of deformations of $X/G_k$ over $BG_k$. As a consequence, we obtain that the tangent Lie algebra which encodes equivariant deformations of $X$ corresponds to the (derived) $G_k$-equivariant global vector fields, $R\Gamma(X, T_X)^{G_k}$.

In this picture, three complications arise. First, we need to give a precise meaning to the pushforward of a formal moduli problem under a base. In \Cref{section base change adjunction}, we show the existence of a natural base change adjunction.

Second, to be able to study deformations of quotient stacks $X/G$, we need to prove that deformations of Artin stacks indeed form a formal moduli problem. This result was known only for Deligne-Mumford stacks (see \cite{lurie2018spectral}). In fact we prove it even for prestacks having deformation theory, a setting which encompasses in particular Artin stacks. Our method consists in expressing this functor as the delooping of the infinitesimal automorphism prestack of $X$, which is known to be a formal moduli problem when $X$ admits deformation theory. We explain that in \Cref{section split deformations}.

Third, we are interested in deformations in mixed characteristic, that is, deformations which are infinitesimal thickenings of $k$ (\textit{under} $k$), but which are not necessarily $k$-algebras (for instance the rings of truncated Witt vectors $W_n(k)$). Thus, we need to consider formal moduli problems which live \textit{under} a base $K$, but \textit{over} a different base $W$. This extension is explained in \Cref{section mixed deformations}.

In this section, $W$ will denote a base prestack.

\subsubsection{Preliminaries}

Let us begin by recalling the constructions of formal moduli stacks and formal moduli problems under a base, following \cite{gaitsgory2019study} and \cite[Appendix B]{brantnermagidsonnuiten2025}. Intuitively, a formal moduli stack under a base is a prestack with some finiteness conditions, having a kind of descent property called deformation theory, and which is an infinitesimal thickening of the base prestack.

First, the objects that we study are simply prestacks in the most general sense.

\begin{defi}[{\cite[Definition B.1]{brantnermagidsonnuiten2025}}]
    A prestack is an accessible functor from animated rings to spaces $$X \colon \calg^{\an} \to \Spaces.$$
    We will write $\prstk$ for the $\infty$-category of prestacks.
Given a base prestack $W$, we let 
    \begin{enumerate}
    \item $\prstk_W$ be the $\infty$-category of $W$-prestacks, that is, prestacks $X$ endowed with a map $X \to W$.

    \item $\Schaff_{/W}$ be the $\infty$-category of affine derived schemes with a map to $W$. 
    \end{enumerate}

    The $\infty$-category of $W$-prestacks can then also be identified with the $\infty$-category of accessible functors $\Schaffop{/W} \to \Spaces$.
\end{defi}

One needs to impose two finiteness conditions on those prestacks to be able to manipulate them: convergence and local almost finite presentation (laft).

\begin{defi}
A map of prestacks $Y \to X$ is said to be convergent if for all $A \in \calg^{\an}$, the following square is cartesian 
\[\begin{tikzcd}
	{Y(A)} & {\lim Y(\tau_{\leq n}A)} \\
	{X(A)} & {\lim X(\tau_{\leq n}A)}
	\arrow[from=1-1, to=1-2]
	\arrow[from=1-1, to=2-1]
	\arrow["\lrcorner"{anchor=center, pos=0.125}, draw=none, from=1-1, to=2-2]
	\arrow[from=1-2, to=2-2]
	\arrow[from=2-1, to=2-2]
\end{tikzcd}\]
\end{defi}

\begin{defi}
    A map of prestacks $Y \to X$ is said to be locally almost finitely presented (laft) if it is convergent and for any point $x \colon \Spec(A) \to X$ and each $n \geq 0$, the base change $Y_x$ has the property that its restriction to $n$-truncated animated $A$-algebras preserves filtered colimits.
\end{defi}

\begin{notation}
    Write $\Schaffaft_{/W}$ for the full subcategory of $\Schaff_{/W}$ of affine schemes almost of finite presentation over $W$.
\end{notation}

Then, to study infinitesimal deformations of prestacks, one restricts to prestacks verifying a kind of descent condition, called ``having deformation theory'', which is quite simple. It is the same condition as the condition for a functor to be a formal moduli problem in Lurie's work (in \cite{lurie2011formal} for instance).

\begin{defi}
    A map of animated rings $A \to B$ is said to be a nilpotent extension if the map $\pi_0(A) \to \pi_0(B)$ is surjective and its kernel is a nilpotent ideal.
\end{defi}

\begin{defi} \label{defi deformation theory}
    Let $W$ be a prestack. A $W$-prestack $\M$ is said to have deformation theory (relative to $W$) if for all points $w \colon \Spec(A) \to W$, the fiber $\M_w$ is convergent and preserves pullbacks along nilpotent extensions of animated $A$-algebras.
\end{defi}

\begin{rk} \label{pushout affine def}
    A pushout (in affine schemes) along a nilpotent extension of affine schemes is also a pushout in the category of prestacks having deformation theory $\prstkdef$. Indeed, a square is a pushout in a category $\cat$ if and only if for all $X \in \cat$, the square becomes a pullback after applying $\Map(-, X)$. The category of $\prstkdef$ can be viewed as the larger category with this property.
\end{rk}

\begin{defi}
    Let $W$ be a prestack.
    Write $\prstklaftdef_{/W}$ for the full subcategory of $\prstk_{/W}$ of $W$-prestacks $\M$ which admit deformation theory and such that the map $\M \to W$ is laft.

    Let $K \to W$ be a map of prestacks.
    Write $\prstklaftdef_{K/W}$ for the full subcategory of $\prstk_{K/W}$ of $W$-prestacks $\M$ which admit deformation theory and such that the map $K \to \M$ is laft.
\end{defi}

One formalizes the notion of ``infinitesimal thickening'' with the notion of nil-isomorphism.

\begin{defi}
    Let $W$ be a prestack. A laft map of $W$-prestacks $f \colon X \to Y$ is said to be a nil-isomorphism if for every reduced $A \in \calg^{\an}_W$, $X(A) \to Y(A)$ is an isomorphism. 
\end{defi}

Finally, the $\infty$-category of formal moduli stacks over $W$ and under $K$ is the full subcategory of $\prstklaftdef_{K/W}$ of prestacks which are nil-isomorphic to $K$.

\begin{defi}[\cite{gaitsgory2019study,brantnermagidsonnuiten2025}]
    Let $W$ be a prestack.
A formal moduli stack under a $W$-prestack $K$ is a map of $W$-prestack $K \to \M$ such that:
\begin{enumerate}
\item $\M$ has deformation theory
\item $K \to \M$ is a (laft) nil-isomorphism
\end{enumerate}
Write $\modstk_{K/W} \subset \prstklaftdef_{K/W} \subset (\prstk_W)_{K/}$ for the full subcategory of formal moduli stacks under $K$.
\end{defi}

In this situation, two bases come into play: formal moduli stacks live \textit{over} a base $W$ and \textit{under} a base $K$. That is why we will call this situation a \textit{mixed bases} setting. We need such a setting to deal with mixed characteristic deformations, when $W = \Spec(W(k))$ and $K = \Spec(k)$, hence our choice of notation. Notice also that $W$ is not an object of $\modstk_{K/W}$ in general, and that $K$ is an object of $\modstk_{K/W}$ if and only if $K$ admits deformation theory.

When $W = K$, we will talk about a \textit{split} setting. When $K$ admits deformation theory, then there is an equivalence $$\modstk_{K/K} \simeq (\modstk_{K/W})_{/K}.$$

\begin{notation}
    Let $K \to W$ be a map of prestacks. Write $\Schaffnil_{K/W}$ for the full subcategory of $\modstk_{K/W}$ whose objects are affine schemes.
\end{notation}

\paragraph{Formal completion.}

There is a canonical operation which makes a map into a nil-isomorphism: formal completion.

\begin{defi}
    \begin{enumerate}
    \item Let $A$ be an animated ring.
    The reduction of $A$ is the quotient of $\pi_0(A)$ by its nilradical.

    \item Let $X$ be a prestack. The de Rham prestack of $X$, denoted $X_{\dR}$, is defined for all animated rings $A$ as $$X_{\dR}(A) = X(A_{\red})$$
    \end{enumerate}
\end{defi}

\begin{defi}
    Let $K \to \M$ be a laft map of prestacks.
    The formal completion of $K \to \M$, denoted $\M_K ^{\wedge}$ or $(K \sto \M)^{\wedge}$, is the prestack defined as $$\M_K^{\wedge} = K_{\dR} \times_{\M_{\dR}} \M.$$
\end{defi}

\begin{prop}[{\cite[Chapter 5, 2.3.3]{gaitsgory2019study}}] \label{formal completion niliso}
    If $K \to \M$ is a laft map of prestacks and $\M$ admits deformation theory, then $K \to M_K^{\wedge}$ is a nil-isomorphism.
    Thus, formal completion determines a functor $\prstklaftdef_{K/W} \to \modstk_{K/W}$.
\end{prop}

\begin{rk} \label{niliso formal completion}
    A laft map of prestacks $K \to \M$ is a nil-isomorphism if and only if the induced map $K_{\dR} \to \M_{\dR}$ is an equivalence. In that case the natural map $\M_K^{\wedge} \to \M$ is an equivalence.
\end{rk}

\begin{prop} \label{formal completion right adjoint to forgetful functor}
Let $K$ be a $W$-prestack. 
    Formal completion defines a functor $$(-)_K^{\wedge} \colon \prstklaftdef_{K/W} \to \modstk_{K/W}$$ which is right adjoint to the inclusion $U \colon \modstk_{K/W} \to \prstklaftdef_{K/W}$.
\end{prop}

\begin{proof}
    The counit of the adjunction is given by the transformation $$\Nr_K^{\wedge} \to \Nr$$ which is natural in $\Nr \in \prstklaftdef_{K/W}$. For $\M \in \modstk_{K/W}$, \Cref{niliso formal completion} gives a natural equivalence $\M \iso \M_K^{\wedge}$: this is the unit of the adjunction. Moreover, if $\Nr \in \prstklaftdef_{K/W}$, the map $(\Nr_K^{\wedge})_K^{\wedge} \to \Nr_K^{\wedge}$ coming from the counit $\Nr_K^{\wedge} \to \Nr$ and the map $(\Nr_K^{\wedge})_K^{\wedge} \to \Nr_K^{\wedge}$ coming from the inverse of the unit evaluated at $\Nr_K^{\wedge}$ are equivalent. Then (it is a purely categorical result for a retraction with such properties), if $\M \in \modstk_{K/W}$ and $\Nr \in \prstklaftdef_{K/W}$, the unit and the counit induce a natural equivalence of spaces $$\Map_W(\M, \Nr_K^{\wedge}) \simeq \Map_W(\M, \Nr).$$

\end{proof}

\begin{lemma} \label{restriction adjunction}
    Let $L \colon \begin{tikzcd}
    {C} & {D}
	\arrow[shift left, from=1-1, to=1-2]
	\arrow[shift left, from=1-2, to=1-1]
\end{tikzcd} \colon R$ be an adjunction of $\infty$-categories.
Let $C'$ and $D'$ be full subcategories of $C$ and $D$ respectively, such that the inclusion functor $U_C \colon C' \to C$  admits a right adjoint, and such that $L$ restricts to a functor $L' \colon C' \to D'$. Denote by $F$ the right adjoint of $U_C$.
Then $F \circ R \circ U_D$ is right adjoint to $L'$.
\end{lemma}

\begin{proof}
    The unit is given by the composite $\id \to (FR)(LU_C) \iso FRU_DL$ (the first map comes from the unit of the adjunction between $FR$ and $LU_C$ and the second from the fact that $L'$ is the restriction of $L$). 
    For the counit, consider the composite $$U_DL'FRU_D \to LU_CFRU_D \to LRU_D \to U_D$$ (the first map comes from the fact that $L'$ is the restriction of $L$, the second from the counit of the adjunction between $F$ and $U_C$, the third from the counit of the adjunction between $L$ and $R$). The functor $U_D$ is fully faithful so this composite induces a map $L'FRU_D \to \id$.

    Now let $X \in C'$ and $Y \in D'$. The maps we constructed above induce natural equivalences $$\Map_{C'}(X, FRU(Y)) \simeq \Map_D(LUX, UY) \simeq \Map_{D}(UL'X, UY) \simeq \Map_{D'}(L'X, Y).$$
\end{proof}

\begin{prop}
    Let $K$ be a $W$-prestack with deformation theory. 
    Let $$U_K^W \colon \modstk_{K/K} \to \modstk_{K/W}$$ be the functor given by postcomposition with the map $K \to W$.
    Let $$\res_K^W \colon \modstk_{K/W} \to \modstk_{K/K}$$ be the functor which naturally sends $\M \in \modstk_{K/W}$ to the formal completion $$(K \to \M \times_W K)^{\wedge}.$$
Those functors form an adjunction
    $$U_K^W \colon\begin{tikzcd}
	\modstk_{K/K} & \modstk_{K/W}
	\arrow[shift left, from=1-1, to=1-2]
	\arrow[shift left, from=1-2, to=1-1]
\end{tikzcd} \colon \res_K^W.$$ 
\end{prop}

\begin{proof}
    Write the structure map $f \colon K \to W$. It induces a base change adjunction $$U_K^W \colon\begin{tikzcd}
	\prstk_{/K} & \prstk_{/W}
	\arrow[shift left, from=1-1, to=1-2]
	\arrow[shift left, from=1-2, to=1-1]
\end{tikzcd} \colon f^*.$$
and a pointed base change adjunction $$U_K^W \colon\begin{tikzcd}
	\prstk_{K/K} & \prstk_{K/W}
	\arrow[shift left, from=1-1, to=1-2]
	\arrow[shift left, from=1-2, to=1-1]
\end{tikzcd} \colon f^*.$$
where the pointing $K \to \M \times_W K$ is built from the maps $K \to \M$ and $\id_K$.
The forgetful functor $U_K^W \colon \prstk_{K/K} \to \prstk_{K/W}$ clearly restricts to the functor $U_K^W \colon \modstk_{K/K} \to \modstk_{K/W}$.
But by \Cref{formal completion right adjoint to forgetful functor}, formal completion is right adjoint to the inclusion functor, so by \Cref{restriction adjunction} the functor $U_K^W \colon \modstk_{K/K} \to \modstk_{K/W}$ has a right adjoint given by $\res_K^W \simeq (K \to \M \times_W K)^{\wedge}.$
\end{proof}

\paragraph{The tangent complex.}

The usual $\infty$-category of quasi-coherent sheaves is not very well adapted to deformation theory. It is better to consider the category of ind-coherent sheaves, following \cite{gaitsgory2019study} or, dually, of pro-coherent sheaves, following \cite{brantnermagidsonnuiten2025}. In the presence of Serre duality, they are Serre dual. Here we choose to use the setting of pro-coherent sheaves.

\begin{defi}[{\cite[A.23, A.26]{brantnermagidsonnuiten2025}}]
    Let $A$ be an animated ring. We will write $\QC{A}$ for the $\infty$-category of exact functors $F \colon \Mod_A \to \Sp$ satisfying the two following conditions: 
    \begin{enumerate}
    \item they are convergent, in the sense that they preserve limits of Postnikov towers in $\Mod_A$
    \item they are almost finitely presented, in the sense that for each $n \in \N$ and each filtered diagram of n-coconnective $A$-modules $(M_{\alpha})$, the map $\colim_{\alpha} F(M_{\alpha}) \to F(\colim_{\alpha} M_{\alpha})$ is an equivalence. 
    \end{enumerate}
    Equivalently, $\QC{A}$ is the $\infty$-category of left exact functors $\Mod_A^{<\infty} \to \Spaces$ that are almost finitely presented.
\end{defi}

\begin{defi}[{\cite[B.15]{brantnermagidsonnuiten2025}}]
    Let $K$ be a prestack. Define the stable presentable $\infty$-category $\QC{K}$ as 
    $$\QC{K} = \lim_{S \in \Aff_{/K}} \QC{O(S)}$$
    Every map $f \colon K' \to K$ between prestacks induces an adjoint pair 
    $$f^* \colon \begin{tikzcd}
	\QC{K} & \QC{K'}
	\arrow[shift left, from=1-1, to=1-2]
	\arrow[shift left, from=1-2, to=1-1]
\end{tikzcd} \colon f_*.$$
    Moreover, there is a natural transformation $\qc_{K} \to \QC{K}$ which induces an equivalence on eventually connective objects.
\end{defi}

\begin{notation}
    Let $Y$ be a prestack and $X \in \prstklaftdef_Y$. Following \cite[B.17, B.18, B.19]{brantnermagidsonnuiten2025}, the map $X \to Y$ admits a tangent complex, denoted $T_{X \rightarrow Y} \in \QC{X}$.
\end{notation}

\begin{notation}
    Let $W$ be a prestack and $K \in \prstklaftdef_W$. If $\M \in \prstklaftdef_{K/W}$, we will denote $$T_{\M} := \fib (T_{K \sto \M} \to T_{K \sto W}) [1].$$
\end{notation}

The tangent complex has the following important properties.

\begin{prop}[{\cite[B.28, B.33]{brantnermagidsonnuiten2025}}] \label{properties tangent complex}
    Let $W$ be a prestack.
    \begin{enumerate}
    \item If $e \colon K \to \M$ and $\M \to W$ are laft maps where $K$ and $\M$ have deformation theory, then there is a natural map $T_{K \sto W} \to e^*T_{\M \sto W}$ that fits into a fiber sequence $$T_{K \sto \M} \to T_{K \sto W} \to e^*T_{\M \sto W}.$$
    \item Let $K$ be a prestack with deformation theory. The natural map $T_{K \sto \M_K^{\wedge}} \to T_{K \sto \M}$ is an equivalence (the tangent complex can be computed on the formal completion). 
    \item  Let $K \in \prstklaftdef_{/W}$.
    There is a pair of adjoint functors
    $$\sqz_K \colon \begin{tikzcd}
	\QC{K} & \prstklaftdef_{K/W}
	\arrow[shift left, from=1-1, to=1-2]
	\arrow[shift left, from=1-2, to=1-1]
\end{tikzcd} \colon T_{(-)}$$
If $K = \Spec(A)$ is an affine scheme, square zero extension on pro-coherent modules is defined for $F \in \QC{A}$ as $$\sqz_A(F) := \colim_{M \in \mathcal{F}^{\contra}} \Spec(A \oplus M),$$ where $\mathcal{F}$ denotes the unstraightening of $F \colon \Mod_A^{<\infty} \to \Spaces$, and can be extended to any prestack $K$ from the affine case.
    \end{enumerate}
    
\end{prop}

\begin{rk}
    In view of \Cref{properties tangent complex} (1), if $\M \in \prstklaftdef_{K/W}$ and the map $\M \to W$ is laft, then $T_{\M} \simeq e^* T_{\M \sto W}$ where $e \colon K \to \M$ denotes the structure map. 
\end{rk}

\begin{lemma} \label{tangent restriction}

    Let $W$ be a prestack, let $K \in \prstklaftdef$.
    Then the following diagram commutes 

\[\begin{tikzcd}
	{\modstk_{K/W}} & {\modstk_{K/K}} \\
	{(\QC{K})_{/T_{K \to W}}} & {\QC{K}}
	\arrow["\res_K^W", from=1-1, to=1-2]
	\arrow["{T_{K \sto \point}}"', from=1-1, to=2-1]
	\arrow["{T_{K \sto \point}}", from=1-2, to=2-2]
	\arrow["\fib"', from=2-1, to=2-2]
\end{tikzcd}\]

In other words, for $\M \in \modstk_{K/W}$, there is an equivalence $$T_{\M} \simeq T_{K \sto \res_K^W \M}[1].$$

\end{lemma}

\begin{proof}
    Let $\M \in \modstk_{K/W}$, let us write the map $i \colon K \to \M \times_K W$ and the projection $p \colon \M \times_W K \to K$. By \Cref{properties tangent complex}, there is a natural equivalence $$T_{K \sto \res_K^W \M} \iso T_{K \sto \M \times_W K}.$$ Now consider the fiber sequence $$T_{K \sto \M \times_W K} \to T_{K \sto \M} \to i^*T_{\M \times_W K \sto \M}$$ By \Cref{properties tangent complex} there is a natural equivalence $$T_{\M \times_W K \sto \M} \simeq p^* T_{K \sto W},$$ hence a natural equivalence $i^*T_{\M \times_W K \sto \M} \simeq i^*p^*T_{K \sto W} \simeq T_{K \sto W}$. Hence there is a natural equivalence $$T_{K \sto \res_K^W \M} \simeq \fib(T_{K \sto \M} \sto T_{K \sto W}).$$ 

\end{proof}

The tangent functor retains a large amount of data about a formal moduli stack: in particular, it is conservative.

\begin{prop}[] \label{tangent conservative}
    Let $W$ be a prestack, let $K$ be a $W$-prestack with deformation theory. Then \begin{enumerate}
    \item the tangent functor $T_{K \rightarrow \point} \colon \modstk_{K/W} \to \QC{K}$ is conservative
    \item if $W$ and $K$ are affine and coherent, then the restriction functor $$\res_K^W \colon \modstk_{K/W} \to \modstk_{K/K}$$ is conservative. Equivalently, the tangent functor $\M \mapsto T_{\M}$ is conservative on formal moduli stacks.
    \end{enumerate}
\end{prop}

\begin{proof}
    \begin{enumerate}
        \item This is \cite[Proposition 5.17]{brantnermagidsonnuiten2025}.
        \item If $W$ and $K$ are affine and coherent, then formal moduli stacks are completely determined by their restriction to nilpotent extensions $\Schaffart_{K/W}$ by \Cref{equivalence modstk FMP}. But \cite[Proposition 4.3]{brantnermagidsonnuiten2025} shows that all nilpotent extensions can be decomposed as a composition of square-zero extensions by $K$-algebras. Hence, the restriction functor $\res_K^W$ is conservative. 
    \end{enumerate}
\end{proof}

\paragraph{Reminders on groupoid objects, group objects and group actions.}

We refer to \cite[6.1.2]{lurie2009higher}, \cite[4.7.1]{lurie2016higher} and \cite{nikolaus2015principal} for the following definitions.

\begin{defi}
    A groupoid object in an $\infty$-category $\mathcal{C}$ is a simplicial diagram $U_{\bullet} \colon \Delta^{op} \to \mathcal{C}$ such that for all $[n] \in \Delta$ and all partitions $[n] = S \cup S'$ into subsets $S$ and $S'$ such that $S \cap S'$ consists of a single element $m \in [n]$, the square 
\[\begin{tikzcd}
	{U_n = U_{\bullet}([n])} & {U_{\bullet}(S)} \\
	{U_{\bullet}(S')} & {U_{\bullet}(\{m\}) = U_0}
	\arrow[from=1-1, to=1-2]
	\arrow[from=1-1, to=2-1]
	\arrow["\lrcorner"{anchor=center, pos=0.125}, draw=none, from=1-1, to=2-2]
	\arrow[from=1-2, to=2-2]
	\arrow[from=2-1, to=2-2]
\end{tikzcd}\]
    is cartesian.
    We will write $\Grpd_{\cat}$ for the full subcategory of the category of simplicial objects in $\cat$ of groupoid objects.
    If $X \in \mathcal{C}$, we will write $\Grpd_{\mathcal{C}}(X)$ for the full subcategory of the category of simplicial objects in $\mathcal{C}$ of groupoid objects $U_{\bullet}$ such that $U_0 \simeq X$. It is called the $\infty$-category of groupoids acting on $X$.
\end{defi}

\begin{defi}
    A group object in an $\infty$-category $\mathcal{C}$ is a groupoid object $G_{\bullet}$ such that $G_0 \simeq *$ where $*$ is the terminal object of $\mathcal{C}$.
    We will write $\Grp_{\mathcal{C}}$ for the $\infty$-category of group objects in $\mathcal{C}$.
\end{defi}

\begin{rk}
    The $\infty$-category of group objects is equivalent to a slice category of the category of groupoids: $\Grp_{\mathcal{C}} \simeq \Grpd_{\mathcal{C}}(*)_{/diag_{*}}$ where $*$ is the terminal object of $\cat$ and $diag_{*}$ is the trivial groupoid.
\end{rk}

\begin{defi}[{\cite[6.1.2]{lurie2009higher}}]
    Let $G_{\bullet}$ be a group object in an $\infty$-category $\mathcal{C}$. An action of $G_{\bullet}$ on an object $X \in \mathcal{C}$ is a groupoid object $U_{\bullet}$ with an equivalence $U_0 \simeq X$ and a morphism of simplicial diagrams $U_{\bullet} \to G_{\bullet}$ such that for every morphism $[m] \to [n]$ in $\Delta$ sending $0 \mapsto 0$, the square 
\[\begin{tikzcd}
	{U_n} & {U_m} \\
	{G_n} & {G_m}
	\arrow[from=1-1, to=1-2]
	\arrow[from=1-1, to=2-1]
	\arrow["\lrcorner"{anchor=center, pos=0.125}, draw=none, from=1-1, to=2-2]
	\arrow[from=1-2, to=2-2]
	\arrow[from=2-1, to=2-2]
\end{tikzcd}\]
    is cartesian.
\end{defi}

\begin{notation}
    \begin{enumerate}
    \item We will denote by $\Act_{\mathcal{C}}(G_{\bullet})$ the full subcategory of $\Grpd_{\mathcal{C}, /G_{\bullet}}$ of groupoids over $G_{\bullet}$ which encode an action of $G_{\bullet}$. The projection $\Act_{\mathcal{C}}(G_{\bullet}) \to \mathcal{C}$ which sends $U_{\bullet}$ to $U_0$ is a right fibration.

    \item We will denote by $(\Grp_{\cat})_{\acts X}$ the full subcategory of $\Grpd_{\cat}^{\Delta^1} \times_{\Grpd_{\cat}} \Grp_{\cat}$ of maps from a groupoid to a group which encode a group action on $X$. The projection $$(\Grp_{\cat})_{\acts X} \to \Grp_{\cat}$$ is a right fibration.

    \item We will denote by $\Act_{\mathcal{C}}(G_{\bullet}, X)$ the space of $G$-actions on $X$, defined as the fiber over $X$ of the right fibration $\Act_{\mathcal{C}}(G_{\bullet}) \to \mathcal{C}$. It is also equivalent to the fiber over $G$ of the right fibration $(\Grp_{\cat})_{\acts X} \to \Grp_{\cat}$. 
    \end{enumerate}
    
\end{notation}

\paragraph{Equivalence between formal moduli stacks and formal groupoids.}
    
In \cite{gaitsgory2019study}, Gaitsgory and Rozenblyum show that the usual adjunction between \v{C}ech nerve and quotient can be adapted to the formal setting, with the notable remark that the quotient differs from the quotient computed in prestacks, to force the quotient to have deformation theory. It induces an equivalence between mixed formal moduli stacks and formal groupoids. This equivalence itself induces an equivalence between split formal moduli stacks and formal groups.

Let us begin  by recalling what are formal groupoids.

\begin{notation} \label{notation groups short}
    Let $K \to W$ be a map of prestacks.
    To lighten the notation, we will write $\Grpd_{W}$ instead of $\Grpd_{\prstk_{/W}}$, $\Grpd_{W}(K)$ instead of $\Grpd_{\prstk_{/W}}(K)$, $\Grp_{K}$ instead of $\Grp_{\prstk_{/K}}$, $\Act_K(G)$ instead of $\Act_{\prstk_{/K}}(G_{\bullet})$ and $\Act_K(G, X)$ instead of $\Act_{\prstk_{/K}}(G_{\bullet}, X)$.
    We will also write $\Grpdlaft_{W}$ instead of $\Grpd(\prstklaft_{/W})$.
\end{notation}

\begin{defi}
    Let $\modstk_{/K}$ be the full subcategory of $\prstk_{/K}$ of prestacks $\M$ which admit deformation theory and such that the map $\M \to K$ is a nil-isomorphism.
\end{defi}

\begin{defi}[{\cite[2.2.1]{gaitsgory2019study}}]

    Let $W$ be a prestack, let $K \in \prstklaft_{/W}$.
    \begin{enumerate}
    \item let $\FGrpd_{W}(K)$ be the full subcategory of $\Grpdlaft_{W}(K)$ of groupoids $U_{\bullet}$ with $U_1 \in \modstk_{/K}$ through the map $U_1 \to U_0 \simeq K$.

    \item let $\FGrp_{K} := \Grp(\modstk_{/K})$, called the $\infty$-category of formal groups over $K$.
    \end{enumerate}
\end{defi}

\begin{rk}
    The category of formal groups can be obtained as a slice of the category of formal groupoids: $$\FGrp_{K} \simeq \FGrpd_{W}(K)_{/\diag_K}$$
    Hence, one gets an adjunction $$U \colon \begin{tikzcd}
	\FGrp_K & \FGrpd_W(K)
	\arrow[shift left, from=1-1, to=1-2]
	\arrow[shift left, from=1-2, to=1-1]
\end{tikzcd} \colon \fib.$$
\end{rk}

\begin{thm}[{\cite[Chapter 5, 2.3]{gaitsgory2019study}}] \label{equivalence formal moduli stacks formal groupoids}
    Let $W$ be a prestack, and $K \in \prstklaftdef_{/W}$.
    Denote by $\Omega \colon \prstk_{K/W} \to \Grpd_{{/W}}(K)$ the \v{C}ech nerve construction which sends $\M \in \prstk_{K/W}$ to the groupoid $K \times_{\M} K$.
    The functor $\Omega$ restricts to a functor $\modstk_{K/W} \to \FGrpd_{/W}(K)$ which is an equivalence of categories. We denote its inverse functor $B^{\deff}$, called the quotient of a formal groupoid.
    Moreover, this equivalence induces an equivalence in the split situation between formal groups and split formal moduli stacks (by taking sliced categories), as in the following commutative diagram:

\[\begin{tikzcd}
	{\FGrpd_W(K)} & {\modstk_{K/W}} \\
	{\FGrp_K} & {\modstk_{K/K}}
	\arrow["B^{\deff}", shift left, from=1-1, to=1-2]
	\arrow["\fib"', shift right, from=1-1, to=2-1]
	\arrow["{\Omega}", shift left, from=1-2, to=1-1]
	\arrow["\res"', shift right, from=1-2, to=2-2]
	\arrow["U"', shift right, from=2-1, to=1-1]
	\arrow["B^{\deff}", shift left, from=2-1, to=2-2]
	\arrow["U"', shift right, from=2-2, to=1-2]
	\arrow["{\Omega}", shift left, from=2-2, to=2-1]
\end{tikzcd}\]
This equivalence is compatible with the tangent functor: if $G \in \FGrp_K$ is a formal group, then there is an equivalence $$T_{G \sto K} \simeq T_{K \sto B^{\deff}G}.$$
\end{thm}

\begin{rk}
    The functor $\Omega \colon \modstk_{K/W} \to \FGrpd_{\prstk_{/W}}(K)$ can be computed in prestacks because the forgetful functor $\modstk_{K/W} \to \prstk_{K/W}$ commutes with finite limits. However, it does not commute with colimits (the issue is the property of having deformation theory).
    In particular, the quotient functor $B^{\deff} \colon \FGrpd_{\prstk_{/W}}(K) \to \modstk_{K/W}$ differs from the usual quotient computed as a colimit in prestacks, which does not admit deformation theory in general: it is a kind of sheafification which forces deformation theory.
\end{rk}

\paragraph{Relation with formal moduli problems.}

There is another notion of ``infinitesimal prestacks'', called formal moduli problem, following \cite{lurie2011formal}: functors defined only on nilpotent extensions of the base, sending the base to a point and having deformation theory. 
When $W$ and $K$ are affine and coherent, it can be proven that the restriction of formal moduli stacks to nilpotent extensions of $K$ induces an equivalence between formal moduli stacks and formal moduli problems. In other words, a formal moduli stack in completely determined by its value on nilpotent extensions of the base.
However, it has not been proven that this generalizes to any bases $W$ and $K$.

\begin{defi}
    Let $K \to W$ be map of coherent affine schemes. An Artinian extension of $K$ is an affine scheme $Z \in \Schaff_{K/W}$ such that the map $K \to Z$ is a nilpotent extension, almost of finite presentation and has coconnective fiber.
    Write $\Schaffart_{K/W}$ for the full subcategory of $\Schaff_{K/W}$ of Artinian extensions of $K$.
\end{defi}

\begin{rk}
    A map in $\Schaffart_{K/W}$ is a nilpotent extension if and only if it is surjective on $\pi_0$. 
\end{rk}

\begin{defi}
    Let $K \to W$ be a map of coherent affine schemes. A formal moduli problem is a functor $F \colon \Schaffartop{K/W} \to \Spaces$ which satisfies the following two conditions: 
    \begin{enumerate}
        \item $F(K) \simeq *$
        \item the functor $F$ sends pushouts along nilpotent extensions in $\Schaffart_{K/W}$ to pullbacks.
    \end{enumerate}
    Write $\FMP_{K/W}$ for the $\infty$-category of formal moduli problems.
\end{defi}

\begin{prop}[{\cite[Theorem 5.3]{brantnermagidsonnuiten2025}}] \label{equivalence modstk FMP}
    Suppose that $K \to W$ is a map of coherent affine schemes. Then the adjunction between $\modstk_{K/W}$ and $\FMP_{K/W}$ given by left Kan extension and restriction is an equivalence. 
\end{prop}

\paragraph{Relation with classical deformation functors.}

The natural classical deformation functor associated to a formal moduli problems is indeed a deformation functor in the sense of \cite{schlessinger1968functors}. Here we suppose that $K = \Spec(k)$ and $W = \Spec(\Lambda)$ where $\Lambda$ is a local Noetherian ring with residue field $k$.

\begin{notation}
    We write $\CRing^{\art}_{K/W}$ for the ordinary category of classical algebras which are local Artinian $\Lambda$-algebras with residue field $k$.
\end{notation}

\begin{defi}
    A (classical) deformation functor is a functor $F \colon \CRing^{\art}_{K/W} \to \Set$  such that: 
    \begin{enumerate}
        \item $F(K) \simeq *$
        \item for every pullback square
\[\begin{tikzcd}
	{A'} & A \\
	{B'} & B
	\arrow[from=1-1, to=1-2]
	\arrow[from=1-1, to=2-1]
	\arrow["\lrcorner"{anchor=center, pos=0.125}, draw=none, from=1-1, to=2-2]
	\arrow[from=1-2, to=2-2]
	\arrow[from=2-1, to=2-2]
\end{tikzcd}\]
in $\CRing^{\art}_{K/}$ with $A \to B$ surjective, the following induced map is surjective:
$$\phi \colon F(A') \to F(A) \times_{F(B)} F(B').$$ 
        \item if in (2) we take $B' = k[\epsilon]/(\epsilon^2)$ and $B = k$, then the map $\phi$ is bijective.
    \end{enumerate}
\end{defi}

\begin{prop}[{\cite[Lemma 3.7]{brantner2025deformationsliftscalabiyauvarieties}}]
    If $F \colon \Schaffartop{K/W} \to \Spaces$ is a formal moduli problem, then the composite $$\CRing^{\art}_{K/W} \to \Schaffartop{K/W} \to \Spaces \to \Set$$ is a deformation functor in the classical sense.
\end{prop}

\begin{notation}
    We will denote the deformation functor associated to a formal moduli problem $F$ as above by $F^{\cl}$.
\end{notation}

\subsubsection{Base change adjunction for formal moduli stacks} \label{section base change adjunction}

If $f \colon K' \to K$ is a map of prestacks, we have the base change adjunction for pro-coherent sheaves $$f^*
 \colon \begin{tikzcd}
    \QC{K} & \QC{K'}
    \arrow[shift left, from=1-1, to=1-2]
    \arrow[shift left, from=1-2, to=1-1]
\end{tikzcd} \colon f_*.$$
What kind of non-linear analogue of this base change adjunction is there? Consider an $\infty$-category $\cat$ with pullbacks, an object $c \in \cat$, and the base change functor $\fib_c \colon \cat \to \cat_{/c}$. The base change functor admits a left adjoint, the forgetful functor. If $\cat$ is a topos, the base change functor admits a right adjoint, which should be seen as a relative ``global section'' functor, or as Weil restriction.

The base change adjunction for pro-coherent sheaves lifts to an adjunction for a category of algebras over an operad, in particular for Lie algebras over a base. In characteristic $0$, the work of \cite{gaitsgory2019study} shows that the tangent functor induces an equivalence between split formal moduli stacks and Lie algebras, hence gives a base change adjunction for split formal moduli stacks. In characteristic $p$, partition Lie algebras should be used instead of naive Lie algebras. The work of \cite{brantnermagidsonnuiten2025} shows that, over a base which is a qcqs scheme, the tangent functor induces an equivalence between mixed formal moduli stacks and partition Lie algebroids, and between split formal moduli stacks and partition Lie algebras. Those partition Lie algebroids come with a global section functor which works as the pushforward in a base change adjunction. However, those works do not imply the existence of base change adjunction in characteristic $p$ for a general base (we need at least to encompass Artin stacks for our work), and they do not give at all a base change adjunction between mixed formal moduli stacks. 

In this section, we establish a general base change adjunction for split formal moduli stacks  $$f^*
 \colon \begin{tikzcd}
    \modstk_{K/K} & \modstk_{K'/K'}
    \arrow[shift left, from=1-1, to=1-2]
    \arrow[shift left, from=1-2, to=1-1]
\end{tikzcd} \colon f_*^{\nil}.$$
when $K$ is affine, and the natural generalization to mixed formal moduli stacks (\Cref{base change adjunction}): in the setting of a pullback diagram where $K \to W$ is a map of coherent affine schemes, and $W'$ and $K'$ are laft prestacks
\[\begin{tikzcd}
	{K'} & K \\
	{W'} & W
	\arrow[from=1-1, to=1-2]
	\arrow[from=1-1, to=2-1]
	\arrow["\lrcorner"{anchor=center, pos=0.125}, draw=none, from=1-1, to=2-2]
	\arrow[from=1-2, to=2-2]
	\arrow["f"', from=2-1, to=2-2]
\end{tikzcd}\]
we prove a base change adjunction $$f^*
 \colon \begin{tikzcd}
    \modstk_{K/W} & \modstk_{K'/W'}
    \arrow[shift left, from=1-1, to=1-2]
    \arrow[shift left, from=1-2, to=1-1]
\end{tikzcd} \colon f_*.$$
After taking tangent complexes, this adjunction commutes with the linear adjunction on pro-coherent sheaves.

\paragraph{The picture in characteristic $0$ over a general base.}

Let $K$ be a laft prestack over a field of characteristic $0$. We recall the following results for reference. 

\begin{notation}
    We will write $\LieAlg_K$ for the $\infty$-category of Lie algebras objects in $\QC{K}$:
    $$\LieAlg_K := \LieAlg(\QC{K})$$
\end{notation}

\begin{thm}[\cite{gaitsgory2019study}] \label{equivalence split Lie algebras}
    The tangent complex of a formal moduli stack admits a Lie algebra structure in $\QC{K}$. Endowed with this Lie algebra structure, the tangent functor becomes an equivalence, as in the following commutative diagram: 

\[\begin{tikzcd}
	{\modstk_{K/K}} & {} & {\LieAlg_K} \\
	& {\QC{K}}
	\arrow["\simeq", from=1-1, to=1-3]
	\arrow["{T_{K \sto \point}}"', from=1-1, to=2-2]
	\arrow["U", from=1-3, to=2-2]
\end{tikzcd}\] 
\end{thm}

\begin{prop} \label{adjunction Lie algebras}
    The base change adjunction $$f^*
 \colon \begin{tikzcd}
    \QC{K} & \QC{K'}
    \arrow[shift left, from=1-1, to=1-2]
    \arrow[shift left, from=1-2, to=1-1]
\end{tikzcd} \colon f_*.$$
induces an adjunction between Lie algebras categories:

\[\begin{tikzcd}
	{\LieAlg_{K'}} & {\LieAlg_K} \\
	{\QC{K'}} & {\QC{K}}
	\arrow["{f_*}", shift left, from=1-1, to=1-2]
	\arrow["U", from=1-1, to=2-1]
	\arrow["{f^*}", shift left, from=1-2, to=1-1]
	\arrow["U", from=1-2, to=2-2]
	\arrow[shift left, from=2-1, to=2-2]
    \arrow[shift left, from=2-2, to=2-1]
\end{tikzcd}\]

Through the equivalence between split formal moduli stacks and Lie algebras from \Cref{equivalence split Lie algebras}, it gives an adjunction $$f^*
 \colon \begin{tikzcd}
    \modstk_{K/K} & \modstk_{K'/K'}
    \arrow[shift left, from=1-1, to=1-2]
    \arrow[shift left, from=1-2, to=1-1]
\end{tikzcd} \colon f_*^{\nil}.$$
\end{prop}

\begin{proof}
    The pullback is monoidal, hence its adjoint is lax monoidal, hence the adjunction between modules restricts to an adjunction between Lie algebras.
\end{proof}

In fact, the pushforward defined this way even admits a lift to a pushforward functor $\modstk_{K/W} \to \modstk_{K/K}$. 

\begin{prop}[{\cite[7.3.4]{gaitsgory2019study}}] \label{global section functor}
    Let $K' \in \prstklaftdef_K$.
    There exists an adjunction
$$\begin{tikzcd}
    (\modstk_{K/K})_{/\Def_{K'}^{K/K}} & \modstk_{K'/K}
    \arrow[shift left, from=1-1 , to=1-2]
    \arrow[shift left, from=1-2, to=1-1]
\end{tikzcd}.$$

such that the following diagram is commutative:
\[\begin{tikzcd}
	{\modstk_{K'/K}} & {(\modstk_{K/K})_{/\Def_{K'}^{K/K}}} \\
	{\modstk_{K'/K'}} & {\modstk_{K/K}}
	\arrow[shift left, from=1-1, to=1-2]
	\arrow["\res"', shift right, from=1-1, to=2-1]
	\arrow[shift left, from=1-2, to=1-1]
	\arrow["\fib"', shift right, from=1-2, to=2-2]
	\arrow["U"', shift right, from=2-1, to=1-1]
	\arrow["{f_*}", shift left, from=2-1, to=2-2]
	\arrow["\triv"', shift right, from=2-2, to=1-2]
	\arrow["{f^*}", shift left, from=2-2, to=2-1]
\end{tikzcd}\]
where the functor $\triv$, which sends $\M \in \modstk_{K/K}$ to the trivial map $\M \to \Def_{K'}^{K/K}$, is right adjoint to the functor $\fib$.

\end{prop}

\paragraph{A geometric interpretation when the base is a qcqs scheme: global sections of Lie algebroids.}

In positive characteristic, Lie algebras are not the right structure anymore to control formal moduli stacks. They must be replaced by partition Lie algebra, as developed in \cite{brantner2023deformation}, where Brantner-Matthew built the equivalence $T \colon \modstk_k \iso \LieAlg^{\pi}_k$. 

However, the theory of partition Lie algebras is still recent. For instance, partition Lie algebras over a base $X$ have not been studied out of the case where $X$ is affine (\cite{brantner2024pdoperadsexplicitpartition}), as a monad over $\QC{X}$.  So there is no obvious base change adjunction for partition Lie algebras over a base.

In \cite{brantnermagidsonnuiten2025}, when the base is a qcqs scheme, the authors define partition Lie algebroids and a global section functors $\Gamma \colon \pliealgd{X}{R} \to \pliealg{R}$, and they prove that the tangent functor induces an equivalence between mixed formal moduli stacks and partition Lie algebroids. Hence, when the base is a qcqs scheme it gives a global section functor analogous to the one in \Cref{global section functor}, which reduces to a pushforward functor for split formal moduli stacks after taking sliced categories.

However, this work is limited to the case where the base is a qcqs scheme (for instance it does not work over $BG$ which is an Artin stack).

\paragraph{A general base change adjunction for formal moduli stacks.}

We do not address the question of generalizing the equivalence between formal moduli stacks and Lie algebroids to a general base here. We show that the base change adjunction for prestacks restricts to an adjunction for formal moduli stacks. 

\begin{defi}
    Let $f \colon W' \to W$ be a map of prestacks
    Write $f^* \colon \prstk^{\laft}_{/W} \to \prstk^{\laft}_{/W'}$ for the base change functor which sends $\M \in \prstk^{\laft}_{/W}$ to $$f^*(\M) = \M \times_W W'.$$
    Write $f^{\laft}_* \colon \prstk^{\laft}_{/W'} \to \prstk^{\laft}_{/W}$ for the Weil restriction functor which sends $\Nr \in \prstk^{\laft}_{/W'}$ to the prestacks sending $Z \in \Schaffaft_{/W}$ to 
    $$f^{\laft}_*(\Nr)(Z) = \Map_{W'}(Z \times_W W', \Nr).$$
    Those functors form an adjunction $$f^* \colon \begin{tikzcd}
    {\prstk^{\laft}_{/W}} & {\prstk^{\laft}_{/W'}}
	\arrow[shift left, from=1-1, to=1-2]
	\arrow[shift left, from=1-2, to=1-1]
\end{tikzcd} \colon f^{\laft}_*.$$

    Now let $K$ be a laft $W$-prestack, and write the fiber $K' := W' \times_W K$. The above adjunction induces an adjunction $f^* \colon \begin{tikzcd}
    {(\prstk^{\laft}_{/W})_{K/}} & {(\prstk^{\laft}_{/W'})_{K'/}}
	\arrow[shift left, from=1-1, to=1-2]
	\arrow[shift left, from=1-2, to=1-1]
\end{tikzcd} \colon f^{\laft}_*$, in the following way: 
    \begin{itemize}
    \item if $\M \in (\prstk^{\laft}_{/W})_{K/}$, $f^*$ sends $\M$ to the natural map $K' \to \M \times_W W'$ through the diagram 
\[\begin{tikzcd}
	{K'} \\
	K & {\M \times_W W'} & W' \\
	& \M & W
	\arrow[bend right = 10, from=1-1, to=2-1]
	\arrow[from=1-1, to=2-2]
	\arrow[bend left = 20, from=1-1, to=2-3]
	\arrow[bend right = 10, from=2-1, to=3-2]
	\arrow[from=2-2, to=2-3]
	\arrow[from=2-2, to=3-2]
	\arrow["\lrcorner"{anchor=center, pos=0.125}, draw=none, from=2-2, to=3-3]
	\arrow[from=2-3, to=3-3]
	\arrow[from=3-2, to=3-3]
\end{tikzcd}\]

   \item if $\Nr \in (\prstk^{\laft}_{/W'})_{K'/}$ with a pointing map $e : K' \to \Nr$, $f_*$ sends $\Nr$ to the natural map $K \to f_*\Nr$ naturally defined for each $Z \in \Schaffaft_{/W}$ as the composite $$\Map_W(Z, K) \to \Map_W(Z \times_{W} W', K') \to \Map_W(Z \times_W W', \Nr),$$ where the first map is given by pullback along $W' \to W$ and the second map is postcomposition by $e$.
\end{itemize}
\end{defi}

First, we restrict this adjunction to an adjunction for $\prstklaftdef_{K/W}$.

\begin{lemma} \label{pushout universal}
    In $\prstkdef$, base change preserves pushouts along nilpotent extensions. More precisely, if $W' \to W$ is a map in $\prstkdef$, if there is a pushout in $\Schaff$
    \[\begin{tikzcd}
	{Z_{01}} & {Z_1} \\
	{Z_0} & Z
	\arrow[from=1-1, to=1-2]
	\arrow[from=1-1, to=2-1]
	\arrow[from=1-2, to=2-2]
	\arrow[from=2-1, to=2-2]
	\arrow["\lrcorner"{anchor=center, pos=0.125, rotate=180}, draw=none, from=2-2, to=1-1]
\end{tikzcd}\]
where $Z_{01} \to Z_0$ and $Z_{01} \to Z_1$ are nilpotent extensions.
Then the square 
\[\begin{tikzcd}
	{Z_{01} \times_W W'} & {Z_1  \times_W W'} \\
	{Z_0 \times_W W'} & {Z  \times_W W'}
	\arrow[from=1-1, to=1-2]
	\arrow[from=1-1, to=2-1]
	\arrow[from=1-2, to=2-2]
	\arrow[from=2-1, to=2-2]
\end{tikzcd}\]
is a pushout in $\prstkdef$.
\end{lemma}

\begin{proof}
Consider a pushout diagram in $\Schaff_{/W}$%
\[\begin{tikzcd}
	{Z_{01}} & {Z_1} \\
	{Z_0} & Z
	\arrow[from=1-1, to=1-2]
	\arrow[from=1-1, to=2-1]
	\arrow[from=1-2, to=2-2]
	\arrow[from=2-1, to=2-2]
	\arrow["\lrcorner"{anchor=center, pos=0.125, rotate=180}, draw=none, from=2-2, to=1-1]
\end{tikzcd}\]
such that $Z_{01} \to Z_0$ and $Z_{01} \to Z_1$ are nilpotent extensions. This is also a pushout in the category $\prstkdef_{/W}$ by \Cref{pushout affine def}. Let us prove that the square 
\[\begin{tikzcd}
	{Z_{01} \times_W W'} & {Z_1 \times_W W'} \\
	{Z_0 \times_W W'} & {Z \times_W W'}
	\arrow[from=1-1, to=1-2]
	\arrow[from=1-1, to=2-1]
	\arrow[from=1-2, to=2-2]
	\arrow[from=2-1, to=2-2]
	\arrow["\lrcorner"{anchor=center, pos=0.125, rotate=180}, draw=none, from=2-2, to=1-1]
\end{tikzcd}\]
is also a pushout in $\prstkdef_{/W}$.
We can write $W = \colim W^{\al}$ as a colimit of affine schemes, and for each $\al$, we can write $Z \times_W W^{\al}$ as a colimit of affine schemes (then $Z_0 \times_W W^{\al}$, $Z_1 \times_W W^{\al}$, and $Z_{01} \times_W W^{\al}$ are all affine). Hence we can reduce to the case where $W$ is affine. Then, by writing $W'$ as a colimit of affine schemes, we can reduce to the case where $W'$ is affine. In this case, the square 
\[\begin{tikzcd}
	{Z_{01} \times_W W'} & {Z_1 \times_W W'} \\
	{Z_0 \times_W W'} & {Z \times_W W'}
	\arrow[from=1-1, to=1-2]
	\arrow[from=1-1, to=2-1]
	\arrow[from=1-2, to=2-2]
	\arrow[from=2-1, to=2-2]
\end{tikzcd}\]
is a pushout in affine schemes, hence a pushout in $\prstkdef$ by \Cref{pushout affine def}.
\end{proof}

\begin{lemma} \label{adjunction def theory}
    Let $W$ be a prestack. Let $f \colon W' \to W$ and $K \to W$ be laft maps of prestacks. Write the fiber $K' := W' \times_W K$.    The base change adjunction for prestacks restricts to an adjunction
    $$f^* \colon \begin{tikzcd}
    {\prstklaftdef_{K/W}} & {\prstklaftdef_{K'/W'}}
    \arrow[shift left, from=1-1, to=1-2]
    \arrow[shift left, from=1-2, to=1-1]
\end{tikzcd} \colon f^{\laft}_*.$$
\end{lemma}

\begin{proof}
    Let us first prove that the pullback $f^*$ restricts to a functor $$f^* \colon \prstklaftdef_{K/W} \to \prstklaftdef_{K'/W'}.$$ Let $\M \in \prstklaftdef_{K/W}$. By definition, given that $\M$ has deformation theory, $\M \times_W W'$ has deformation theory. Moreover, given that the map $W' \to W$ is laft, the map $\M \times_W W' \to \M$ is laft. But the fiber $K' \to K$ is also laft and the map $K \to \M$ is laft by assumption, so the map $K' \to \M$ is laft. Hence the map $K' \to \M \times_W W'$ is laft.

    Now let us prove that the pushforward $f^{\laft}_*$ restricts to a functor $$f^{\laft}_* \colon \prstklaftdef_{K'/W'} \to \prstklaftdef_{K/W}.$$ Let $\Nr \in \prstklaftdef_{K'/W'}$. By \cite[Chapter 1, 7.2.7]{gaitsgory2019studyI}, the property of infinitesimal cohesiveness can be checked on affine schemes almost of finite type. Consider a pushout of affine schemes almost of finite presentation 
    \[\begin{tikzcd}
	{Z_{01}} & {Z_1} \\
	{Z_0} & Z
	\arrow[from=1-1, to=1-2]
	\arrow[from=1-1, to=2-1]
	\arrow[from=1-2, to=2-2]
	\arrow[from=2-1, to=2-2]
	\arrow["\lrcorner"{anchor=center, pos=0.125, rotate=180}, draw=none, from=2-2, to=1-1]
    \end{tikzcd}\]
    where $Z_{01} \to Z_0$ and $Z_{01} \to Z_1$ are nilpotent extensions. By \Cref{pushout universal} adapted to laft prestacks and affine schemes almost of finite presentation, the square 
    \[\begin{tikzcd}
	{Z_{01} \times_W W'} & {Z_1 \times_W W'} \\
	{Z_0 \times_W W'} & {Z \times_W W'}
	\arrow[from=1-1, to=1-2]
	\arrow[from=1-1, to=2-1]
	\arrow[from=1-2, to=2-2]
	\arrow[from=2-1, to=2-2]
\end{tikzcd}\]
    is a pushout in $\prstklaftdef_{/W'}$. Hence the square
    \[\begin{tikzcd}
	{\Map_{W'}(Z \times_W W', \Nr)} & {\Map_{W'}(Z_1 \times_W W', \Nr)} \\
	{\Map_{W'}(Z_0 \times_W W', \Nr)} & {\Map_{W'}(Z_{01} \times_W W', \Nr)}
	\arrow[from=1-1, to=1-2]
	\arrow[from=1-1, to=2-1]
	\arrow[from=1-2, to=2-2]
	\arrow[from=2-1, to=2-2]
\end{tikzcd}\]
is a pullback square. Thus, the prestack $f^{\laft}_* \Nr$ is infinitesimally cohesive. 

    One proves convergence also by reducing to the affine case, as in \cite[Proposition 2.2.9]{PortaSala_SimpsonsShapes_2020} for instance.
\end{proof}

\begin{rk} \label{colimits are universal}
    \Cref{adjunction def theory} shows a stronger result than \Cref{pushout universal}: all colimits that exist in $\prstklaftdef$ are preserved under base change (that is, colimits are universal in $\prstklaftdef$). 
\end{rk}

Now we restrict again to an adjunction between formal moduli stacks, using the fact that the inclusion $\modstk_{K/W} \to \prstklaftdef_{K/W}$ admits a right adjoint, the formal completion functor.

\begin{rk}
    The forgetful functor $\prstklaftdef_{K/W}\to \prstk_{K/W}$ does not admit a right adjoint. Colimits in the category of prestacks with deformation theory cannot be computed in the category of prestacks.
\end{rk}

\begin{prop} \label{base change adjunction}
    Let $W$ be a prestack, let $f \colon W' \to W$ and $K \to W$ be laft maps of prestacks. Write the fiber $K' := W' \times_W K$.
    The pullback functor for prestacks restricts to a functor $$f^* \colon \modstk_{K/W} \to \modstk_{K'/W'}.$$
    Write $f_*^{\nil} \colon \modstk_{K'/W'} \to \modstk_{K/W}$ for the composite $f^{\nil}_* = (-)_K^{\wedge} \circ f^{\laft}_* \circ U$.
    Then $f_*^{\nil}$ is right adjoint to $f^*$.

\end{prop}

\begin{proof}
    Let $\M \in \modstk_{K/W}$. We have to prove that the map $K' \to f^* \M$ is a nil-isomorphism.
    Let us show that the induced map $(K')_{\dR} \to (\M \times_W W')_{\dR}$ is an equivalence. The de Rham functor is a right adjoint so it commutes with limits, so the natural map $(\M \times_W W')_{\dR} \to \M_{\dR} \times_{W_{\dR}} W'_{\dR}$ is an equivalence. But $K \to \M$ is a nil-isomorphism so $K_{\dR} \to \M_{\dR}$ is an equivalence hence the natural composite $$(\M \times_W W')_{\dR} \to \M_{\dR} \times_{W_{\dR}} W'_{\dR} \to K_{\dR} \times_{W_{\dR}} W'_{\dR} \to (K \times_W W')_{\dR} \simeq (K')_{\dR}$$ is an equivalence.

    The adjunction between $f^*$ and $f_*^{\nil}$ comes from \Cref{restriction adjunction}, given that the inclusion $\modstk_{K/W} \to \prstklaftdef_{K/W}$ admits a right adjoint, the formal completion functor.

\end{proof}

\begin{rk}
    In the split case, that is, when $W = K$ and $W' = K'$, it becomes an adjunction $$f^* \colon \begin{tikzcd}
    {\modstk_{K/K}} & {\modstk_{K'/K'}}
    \arrow[shift left, from=1-1, to=1-2]
    \arrow[shift left, from=1-2, to=1-1] 
    \end{tikzcd}
    \colon f_*^{\nil}.$$
\end{rk}

The pullback functor commutes not only with the inclusion $\modstk_{K/W} \to \prstklaftdef_{K/W}$ but also with its adjoint, the formal completion functor $\prstklaftdef_{K/W} \to \modstk_{K/W}$.

\begin{prop} \label{formal completion commutes with pullback}
    Let $f \colon K' \to K$ be a laft map of prestacks.
    The pullback functor commutes with formal completion.
\[\begin{tikzcd}
	{\modstk_{K/K}} & {\modstk_{K'/K'}} \\
	{\prstklaftdef_{K/K}} & {\prstklaftdef_{K'/K'}}
	\arrow["{f^*}", from=1-1, to=1-2]
	\arrow["(\cdot)^{\wedge}"', from=2-1, to=1-1]
	\arrow["(\cdot)^{\wedge}", from=2-2, to=1-2]
	\arrow["{f^*}"', from=2-1, to=2-2]
\end{tikzcd}\]
\end{prop}

\begin{proof}
    Let $\M \in \prstklaftdef_{K/K}$. The natural map $\M_K^{\wedge} \to \M$ induces a map $f^*\M_K^{\wedge} \to f^*\M$. The pullback $f^*(\M_K^{\wedge}) \in \modstk_{K'/K'}$ so by adjunction between formal completion and the forgetful functor, it induces a natural map $f^*(\M_K^{\wedge}) \to (f^*\M)_{K'}^{\wedge}$ in $\modstk_{K'/K'}$.

    Now consider the following diagram, where all the squares are commutative except the front square:

\[\begin{tikzcd}
	{\modstk_{K/K}} && {\modstk_{K'/K'}} \\
	& {\QC{K}} && {\QC{K'}} \\
	{\prstklaftdef_{K/K}} && {\prstklaftdef_{K'/K'}}
	\arrow["{f^*}", from=1-1, to=1-3]
	\arrow["{T_{K \sto \point}}", dashed, from=1-1, to=2-2]
	\arrow["{T_{K \sto \point}}", from=1-3, to=2-4]
	\arrow["{f^*}"{pos=0.3}, dashed, from=2-2, to=2-4]
	\arrow["{(\cdot)^{\wedge}}"', from=3-1, to=1-1]
	\arrow["{T_{K \sto \point}}", dashed, from=3-1, to=2-2]
	\arrow["{f^*}"', from=3-1, to=3-3]
	\arrow["{(\cdot)^{\wedge}}"'{pos=0.3}, from=3-3, to=1-3]
	\arrow["{T_{K \sto \point}}"', from=3-3, to=2-4]
\end{tikzcd}\]
    It shows that the associated tangent map $T_{K' \sto f^*(\M_K^{\wedge})} \to T_{K' \to (f^*\M)_{K'}^{\wedge}}$ is an equivalence. By conservativity of the tangent functor for formal moduli stacks, it shows that the natural map $f^*(\M_K^{\wedge}) \to (f^*\M)_{K'}^{\wedge}$ in $\modstk_{K'/K'}$ is an equivalence.

\end{proof}

Finally, as expected, we show that the base change adjunction for formal moduli stacks lifts the base change adjunction for pro-coherent sheaves, in the sense that it commutes with the tangent functor.

\begin{prop} \label{tangent commutes with adjunction}

Let $f \colon K' \to K$ be a laft map of prestacks.
    The pullback functor and the pushforward functor commute with the tangent functor 

\[\begin{tikzcd}
	{\prstklaftdef_{K/K}} & {\prstklaftdef_{K'/K'}} \\
	{\QC{K}} & {\QC{K'}}
	\arrow["{f^*}", from=1-1, to=1-2]
	\arrow["{T_{K \sto \point}}"', from=1-1, to=2-1]
	\arrow["{T_{K' \sto \point}}", from=1-2, to=2-2]
	\arrow["{f^*}"', from=2-1, to=2-2]
\end{tikzcd}\]

\[\begin{tikzcd}
	{\prstklaftdef_{K/K}} & {\prstklaftdef_{K'/K'}} \\
	{\QC{K}} & {\QC{K'}}
	\arrow["{f^{\laft}_*}"', from=1-2, to=1-1]
	\arrow["{T_{K \sto \point}}"', from=1-1, to=2-1]
	\arrow["{T_{K' \sto \point}}", from=1-2, to=2-2]
	\arrow["{f_*^{\laft}}", from=2-2, to=2-1]
\end{tikzcd}\]

\end{prop}

\begin{proof}
    If $\M \in \prstklaftdef_{K'/W'}$, we have the following pullback diagram 
\[\begin{tikzcd}
	{K'} & {\M \times_K K'} \\
	K & \M
	\arrow[from=1-1, to=1-2]
	\arrow["p"', from=1-1, to=2-1]
	\arrow["\lrcorner"{anchor=center, pos=0.125}, draw=none, from=1-1, to=2-2]
	\arrow[from=1-2, to=2-2]
	\arrow[from=2-1, to=2-2]
\end{tikzcd}\]
We know by \Cref{properties tangent complex} that there is a natural equivalence $T_{K' \sto \M \times_K K' } \simeq p^* T_{K \sto \M}$, hence the commutation of $f^*$ with the tangent functor.
To prove the commutation of $f^{\laft}_*$ with the tangent functor, it suffices to have the commutation of the adjoint square

\[\begin{tikzcd}
	{\prstklaftdef_{K/W}} & {\prstklaftdef_{K'/W'}} \\
	{\QC{K}} & {\QC{K'}}
	\arrow["{f^*}", from=1-1, to=1-2]
	\arrow["\sqz_{K}", from=2-1, to=1-1]
	\arrow["{f^*}"', from=2-1, to=2-2]
	\arrow["\sqz_{K'}"', from=2-2, to=1-2]
\end{tikzcd}\]
We can reduce to the case where $K$ and $K'$ are affine schemes, and commutation of square zero extension and pullback is true in that case, see \cite[p. 88, (30)]{brantnermagidsonnuiten2025}.

\end{proof}

\subsubsection{Split deformations of prestacks with deformation theory} \label{section split deformations}

There are two usual deformation functors. First, if $X$ is a scheme over a base $K$ (usually a field), one can consider $\Def^{K/K}_X$, the functor which sends an Artinian extension $Z$ of $K$ to the space of deformations of $X$ over $Z$. Second, if $Y \to K'$ is a map of schemes over $K$, one can consider the functor of deformations of the map $Y \to K'$ leaving the target fixed, denoted $\Def^{K/K}_{Y \to K'}$ (see \cite[3.4.2]{sernesi2007deformations} for instance). It sends an Artinian extension $Z$ of $K$ to the space of maps $\Y \to K' \times_K Z$ over $K'$ which restrict to $Y \to K'$ over $K$.

In the split case, Lurie proved that if $X$ is a Deligne-Mumford stack, then the functor of deformations of $X$ among Deligne-Mumford stacks is a formal moduli problem. We show that if $X$ is a prestack with deformation theory, then the functor of deformations of $X$ among \emph{prestacks with deformation theory}, denoted $\Def_X^{K/K}$, is a formal moduli problem. In some sense, having deformation theory is the minimal assumption necessary to get a formal moduli problem. 

To show this, we use the principle that, in some sense, the right object to study is in fact the looping of the deformations of $X$. We show that the functor of deformations of $X$ is equivalent to the delooping of the formal automorphism prestack, which is a formal moduli stack (\Cref{equivalence def aut}) $$\Def_X^{K/K} \simeq B^{\deff}\infAut_K(X)$$
This point of view shows also that, as expected, the tangent complex of the deformation functor is $$T_{\Def_X^{K/K}} \simeq R\Gamma(X, T_{X \sto K})[1]$$ 

Then, we remark that the functor of deformations of a map leaving the target fixed, $\Def_{Y \to K'}^{K/K}$, is simply the pushfoward (constructed in \Cref{section base change adjunction}), of $\Def_Y^{K'/K'}$. Hence, it is a formal moduli problem. Moreover, writing $g \colon K' \to K$ and $f \colon Y \to K'$, its tangent complex is $$T_{\Def_{Y \to K'}^{K/K}} \simeq g_*f_*T_{Y \sto K'}[1]$$
This is the content of \Cref{tangent deformation of map}. It generalizes the results of \cite[3.4.2]{sernesi2007deformations}, where Sernesi computes the first and second cohomology groups in the case of schemes deformations.

Finally, if one consider a section $e \colon K \to K'$ of the map $f\colon K' \to K$ (think about a point in $K'$), and the fiber $X \to K$ of the map $Y \to K'$
\[\begin{tikzcd}
	X & K \\
	Y & {K'}
	\arrow[from=1-1, to=1-2]
	\arrow[from=1-1, to=2-1]
	\arrow["\lrcorner"{anchor=center, pos=0.125}, draw=none, from=1-1, to=2-2]
	\arrow[from=1-2, to=2-2]
	\arrow[from=2-1, to=2-2]
\end{tikzcd}\]
then one can consider the transformation from deformations of the map $Y \to K'$ to deformations of the fiber $X$, that is, a map $$\fibdef_Y^{K/K} \colon \Def_{Y \to K'}^{K/K} \to \Def_X^{K/K}.$$
In \Cref{fiber of deformation tangent equi}, we show that this map corresponds to the natural map $$f^{\nil}_* \Def_Y^{K'/K'} \to e^* \Def_Y^{K'/K'}$$ built from the counit of the base change adjunction, enabling us to compute easily the map induced on tangent complexes.

\paragraph{Deformation functors and ``fiber transformation'': definitions.}

Consider the moduli functor $$\M_K^{\laftdef} \colon (\mathrm{Sch}^{\mathrm{aff}}_{/K})^{\mathrm{op}} \to \widehat{\Cat_{\infty}}$$ which sends $Z \in \Schaff_{/K}$ to the $\infty$-category $\prstklaftdef_{/Z}$.
A choice of prestack $X \in \prstklaftdef_{/K}$ induces a map $e_X \colon K \to \M^{\laftdef}_{K}$.
Ideally, we would like to define the functor of deformations of $X$, $\Def_X^{K/K}$, as the formal completion of the map $e_X \colon K \to \M_{K}^{\laftdef}$, and then say that $\M_{K}^{\laftdef}$ has deformation theory and deduce by \Cref{formal completion niliso} that $\Def_X^{K/K}$ is a formal moduli stack. However, the moduli functor $\M_{K}^{\laftdef}$ is ``too big'': it is not even a prestack. Remark still that in \cite[16.3]{lurie2018spectral}, Lurie shows a ``clutching'' result analogous to proving that the moduli functor of Deligne-Mumford stacks has deformation theory.

In this work, we will give a definition ``by hand'' of the deformation functor, and then prove that, up to delooping, it corresponds to the formal completion of a much smaller geometric object: the infinitesimal automorphism prestack of $X$.

\begin{defi}
    Let $K$ be a prestack and let $X \in \prstklaftdef_{/K}$. Let $$\Def_X^{K/K} \colon \Schaffnilop{K/K} \to \widehat{\Cat_{\infty}}$$ be the functor which sends $Z \in \Schaffnil_{K/K}$ to  $$\Def_X^{K/K}(Z) := \prstklaftdef_{/Z} \times_{\prstklaftdef_{/K}} \{X\}$$
In other words, if $Z \in \Schaffnil_{K/K}$, $\Def_X^{K/K}(Z)$ identifies with the $\infty$-category of pullback diagrams 

\[\begin{tikzcd}
	X & K \\
	\X & Z
	\arrow[from=1-1, to=1-2]
	\arrow[from=1-1, to=2-1]
	\arrow["\lrcorner"{anchor=center, pos=0.125}, draw=none, from=1-1, to=2-2]
	\arrow[from=1-2, to=2-2]
	\arrow[from=2-1, to=2-2]
\end{tikzcd}\]
\end{defi}

Now let us define the functor of deformations of a map leaving the target fixed, studied for instance in \cite[3.4.2]{sernesi2007deformations} for schemes.

\begin{defi} \label{definition deformation with fixed target}
    Let $K' \to K$ be a laft map of prestacks, let $Y \in \prstklaftdef_{K'}$.
    Let $\Def_{Y \to K'}^{K/K} \colon \Schaffnilop{K/K} \to \widehat{\Cat_{\infty}}$ be the functor which sends $Z \in \Schaffnilop{K/K}$ to $$\Def_{Y \to K'}^{K/K} := \prstklaftdef_{K' \times_K Z} \times_{\prstklaftdef_{K'}} \{Y\}.$$
    where $K' \to K' \times_K Z$ is the natural map built from $\id_{K'}$ and the composite $K' \to K \to Z$.
    In other words, for each $Z \in \Schaffnil_{K/K}$, $\Def_{Y \to K'}^{K/K}(Z)$ identifies with the $\infty$-category of pullback diagrams 
\[\begin{tikzcd}
	Y & K' \\
	\Y & {K' \times_K Z}
	\arrow[from=1-1, to=1-2]
	\arrow[from=1-1, to=2-1]
	\arrow["\lrcorner"{anchor=center, pos=0.125}, draw=none, from=1-1, to=2-2]
	\arrow[from=1-2, to=2-2]
	\arrow[from=2-1, to=2-2]
\end{tikzcd}\]
    We will call it the functor of deformations of the map $Y \to K'$ leaving the target fixed.
\end{defi}

\begin{rk}
    We will see later that the functors $\Def_X^{K/K}$ and $\Def_{Y \to K'}^{K/K}$ land in spaces in reality.
\end{rk}

Next, when we have a section $K \to K'$, we define the natural transformation from deformations of the map $Y \to K'$ to deformations of the fiber $X = Y \times_{K'} K$.

\begin{defi} \label{def fiber of deformations}
    Let $K' \to K$ be a laft map of prestacks, with a section $K \to K'$, let $Y \in \prstklaftdef_{K'}$. Denote $X = Y \times_{K'} K$.
    We define a natural transformation  $$\fibdef_Y^{K/K} \colon \Def_{Y \to K'}^{K/K} \to \Def_X^{K/K}$$ as follows.
    Let $Z \in \Schaffnil_{K/K}$. Let $K' \to K' \times_K Z$ be the natural map built from $\id_{K'}$ and the composite $K' \to K \to Z$. Let $Z \to K' \times_K Z$ be the natural map built from the composite $Z \to K \to K'$ and the identity $\id_Z$. Those maps sit in a pullback diagram

\[\begin{tikzcd}
	K & {K'} \\
	Z & {K' \times_K Z}
	\arrow[from=1-1, to=1-2]
	\arrow[from=1-1, to=2-1]
	\arrow["\lrcorner"{anchor=center, pos=0.125}, draw=none, from=1-1, to=2-2]
	\arrow[from=1-2, to=2-2]
	\arrow[from=2-1, to=2-2]
\end{tikzcd}\]
    Thus, the pullback along $Z \to K' \times_K Z$ induces a natural map $$\fibdef_Y^{K/K}(Z) \colon \prstk_{/(K' \times_K Z)} \times_{\prstk_{/K'}} \{Y\} \to \prstk_{/Z} \times_{\prstk_{/K}} \{X\}$$
    Those natural maps $\fibdef_Y^{K/K}(Z)$ give a natural transformation $$\fibdef_Y^{K/K} \colon \Def_{Y \to K'}^{K/K} \to \Def_X^{K/K}.$$

    In other words, in view of the diagram where all the squares are pullback squares
    
\[\begin{tikzcd}
	& X && K \\
	Y && {K'} \\
	& {\Y \times_{K' \times_K Z} Z} && Z \\
	\Y && {K' \times_K Z}
	\arrow[from=1-2, to=1-4]
	\arrow[dashed, from=1-2, to=3-2]
	\arrow[from=1-4, to=3-4]
	\arrow[from=2-1, to=1-2]
	\arrow[from=2-1, to=2-3]
	\arrow[from=2-1, to=4-1]
	\arrow[from=2-3, to=1-4]
	\arrow[from=2-3, to=4-3]
	\arrow[dashed, from=3-2, to=3-4]
	\arrow[dashed, from=4-1, to=3-2]
	\arrow[from=4-1, to=4-3]
	\arrow[from=4-3, to=3-4]
\end{tikzcd}\]
    then $\fibdef_Y^{K/K}(Z) \colon \Def_{Y \to K'}^{K/K}(Z) \to \Def_X^{K/K}(Z)$ is the map which sends the front pullback square

\[\begin{tikzcd}
	Y & K' \\
	\Y & {K' \times_K Z}
	\arrow[from=1-1, to=1-2]
	\arrow[from=1-1, to=2-1]
	\arrow["\lrcorner"{anchor=center, pos=0.125}, draw=none, from=1-1, to=2-2]
	\arrow[from=1-2, to=2-2]
	\arrow[from=2-1, to=2-2]
\end{tikzcd}\]
    to the back pullback square
\[\begin{tikzcd}
	X & K \\
	{\Y \times_{K' \times_K Z} Z} & Z
	\arrow[from=1-1, to=1-2]
	\arrow[from=1-1, to=2-1]
	\arrow["\lrcorner"{anchor=center, pos=0.125}, draw=none, from=1-1, to=2-2]
	\arrow[from=1-2, to=2-2]
	\arrow[from=2-1, to=2-2]
\end{tikzcd}\]
\end{defi}

\paragraph{Delooping the infinitesimal automorphisms formal group.}

Here the general philosophy is that the loops of the infinitesimal deformations of $X$ are the infinitesimal automorphisms of $X$. 
More precisely, if $Z \in \Schaffnil_{K/K}$ is an infinitesimal extension of $K$, $\Def_X(Z)$ is the space of pullback diagrams
\[\begin{tikzcd}
	X & K \\
	\X & Z
	\arrow[from=1-1, to=1-2]
	\arrow[from=1-1, to=2-1]
	\arrow["\lrcorner"{anchor=center, pos=0.125}, draw=none, from=1-1, to=2-2]
	\arrow[from=1-2, to=2-2]
	\arrow[from=2-1, to=2-2]
\end{tikzcd}\]
Now suppose that $Z = B^{\deff}\widehat{G}$ where $\widehat{G}$ is a formal group. The square becomes 
\[\begin{tikzcd}
	X &  K \\
	\X & B^{\deff}\widehat{G}
	\arrow[from=1-1, to=1-2]
	\arrow[from=1-1, to=2-1]
	\arrow["\lrcorner"{anchor=center, pos=0.125}, draw=none, from=1-1, to=2-2]
	\arrow[from=1-2, to=2-2]
	\arrow[from=2-1, to=2-2]
\end{tikzcd}\]
Observe that this is the square which encodes the data of an action of $\widehat{G}$ on $X$, or equivalently, as we will explain more formally, a map $\widehat{G} \to \infAut_K(X)$, where $\infAut_K(X)$ denotes the formal automorphism prestack of $X$. But by adjunction a map $\Omega Z \to \infAut_K(X)$ corresponds to an element of $B^{\deff} \infAut_K(X)(Z)$.

We begin by defining the formal automorphism prestack, which appears as the formal completion of the full automorphism prestack or, equivalently, of the full endomorphism prestack. The (full) endomorphism prestack is an object more tractable than the moduli functor $\M_{K}^{\laftdef}$ and, especially, it is easy to prove that it admits deformation theory.

\begin{defi}
    Let $K$ be a \comdef{laft} prestack, let $X$ be a laft $K$-prestack.
    The \comdef{laft} endomorphism prestack of $X$, denoted $\End^{\laft}_K(X)$ is the functor which sends $Z \in \Schaffaft_{/K}$ to the space $$\End^{\laft}_K(X)(Z) = \Map_{K}(Z \times_K Y, Y) \simeq \Map_{Z}(Z \times_K Y, Z \times_K Y).$$ 
    By \cite[Chapter 2, I.7.6]{gaitsgory2019studyI}, this defines a laft prestack.
    Composition makes it into a monoid, hence $$\End^{\laft}_K(X) \in \Monoid(\prstk^{\laft}_{/K}).$$
    Moreover, $\End^{\laft}_K(X)$ is pointed by $K$ via the prestack map $K \to \End^{\laft}_K(X)$ naturally defined for each $Z \in \Schaffaft_{/K}$ as the map $\Map(Z, K) \to \Map_K(Z \times_K X, X)$ induced by the identity of $X$.
\end{defi}

\begin{prop} \label{aut def theory}
    Let $K$ be a \comdef{laft} prestack, let $X \in \prstklaftdef_{/K}$.
    Then the \comdef{laft} endomorphism prestack $\End^{\laft}_K(X)$ admits deformation theory.
\end{prop}

\begin{proof}
    The proof is the same as in \Cref{pushout affine def} and \Cref{adjunction def theory}.
\end{proof}

The formal automorphism prestack can be defined as the formal completion of the full endomorphism prestack because, intuitively, the full automorphism prestack is an open in the full endomorphism prestack, hence locally around the identity, they coincide. In other words, formal monoids and formal groups are the same.

\begin{defi}[{\cite[Chapter 7, 6.2]{gaitsgory2019study}}]
    Let $K$ be a \comdef{laft} prestack, let $X \in \prstklaftdef_{K/K}$.
    We define the infinitesimal automorphism prestack of $X$, denoted $\widehat{\Aut}_K(X)$, as the formal completion of the natural map $K \to \End^{\laft}_K(X)$.
    By \Cref{formal completion niliso},  $$\widehat{\Aut}_K(X) \in \Monoid(\modstk_{K/K}).$$
\end{defi}

\begin{lemma}[{\cite[Chapter 5, 1.6.2]{gaitsgory2019study}}] \label{formal groups are equivalent to formal monoids}
    Let $K$ be a prestack.
    The fully faithful functor $\Grp(\modstk_{/K}) \subset \Monoid(\modstk_{/K})$ is an equivalence.
\end{lemma}

\begin{coro}
    Let $K$ be a prestack, let $X \in \prstklaftdef_{K/K}$.
    Then $$\widehat{\Aut}_K(X) \in \FGrp_K.$$
\end{coro}

The endomorphism prestack is a well known object and its tangent complex is known to be the derived global vector fields over $X$. Hence, it is also the tangent complex of the formal automorphism prestack.

\begin{prop} \label{tangent endomorphisms}
    Let $K$ be a \comdef{laft} prestack, let $X \in \prstklaftdef_{K/K}$. Denote the structure map $f \colon X \to K$. Then there is an equivalence in $\QC{K}$

    $$T_{\End^{\laft}_K(X)} \simeq f_*T_{X}.$$
\end{prop}

\begin{proof}
    Let $M \in \QC{K}$. By adjunction between tangent functor and square zero extension, there is a natural equivalence $$\Map_K(M, T_{\End_K(X)}) \simeq \Map_K(\sqz_K(M), \End_K(X)).$$ But there is a natural equivalence $$\Map_K(\sqz_K(M), \End_K(X)) = \Map_K(f^* \sqz_K(M), X).$$ By \Cref{adjunction def theory}, there is a natural equivalence $$\Map_K(f^* \sqz_K(M), X) \simeq \Map_K(\sqz_K(M), f_*^{\laft} X).$$
    Again by adjunction between square zero and tangent, there is a natural equivalence $$\Map_K(\sqz_K(M), f_*^{\laft} X) \simeq \Map_K(M, T_{f_*^{\laft} X}).$$ Now by \Cref{tangent commutes with adjunction}, we know that there is an equivalence $$T_{f_*^{\laft} X} \simeq f_* T_X.$$ Finally we get a natural equivalence $$\Map_K(M, T_{\End_K(X)}) \simeq \Map_K(M, f_* T_X),$$ hence the result. 
\end{proof}

\begin{coro}
    Let $K$ be a \comdef{laft} prestack, let $X \in \prstklaftdef_{K/K}$. Denote the structure map $f \colon X \to K$ and the pointing $e \colon K \to \infAut_K(X)$. Then there is an equivalence in $\QC{K}$:

    $$T_{\infAut_K(X)} \simeq f_*T_{X}.$$
\end{coro}

\begin{proof}
    By \Cref{properties tangent complex}, $T_{\infAut_K(X)} \simeq T_{\End^{\laft}_K(X)}$ so it is a direct consequence of \Cref{tangent endomorphisms}.
\end{proof}

Moreover, the endomorphism prestack of a pullback is the pullback of the endomorphism prestack.

\begin{prop} \label{pullback aut}
Let $K' \to K$ be a map of laft prestacks. Let $Y \in \prstklaftdef_{K'/K'}$. Write the fiber $X = Y \times_{K'} K$. Then there is an equivalence of prestacks $$f^* \End^{\laft}_{K'}(Y) \simeq \End^{\laft}_{K}(X).$$
    As a consequence, there is an equivalence of formal moduli stacks  $$f^* \infAut_{K'}(Y) \simeq \infAut_{K}(X).$$
\end{prop}

\begin{proof}
    Let $Z \in \Schaffaft_{/K}$. The pullback $f^*$ for prestacks is left adjoint to the forgetful functor, that is, it identifies with the restriction functor. Thus there is a natural equivalence $f^*\End^{\laft}_{K'}(Y)(Z) \simeq \Map_Z(Z \times_{K'} Y, Z \times_{K'} Y)$.
    Then we have the following obvious natural equivalences:
    $$\Map_Z(Z \times_{K'} Y, Z \times_{K'} Y) \simeq \Map_{Z}(Z \times_{K} K \times_{K'} Y, Z \times_{K} K \times_{K'} Y) \simeq \Map_{Z}(Z \times_{K} X, Z \times_{K} X)$$
    Hence there is a natural equivalence $f^* \End^{\laft}_{K'}(Y) \simeq \End^{\laft}_{K}(X)$.
    Then we use commutation of the pullback with formal completion from \Cref{formal completion commutes with pullback} to get the result for the infinitesimal automorphism prestack.
\end{proof}

Now, a map from an infinitesimal group $\widehat{G} \to \infAut_K(X)$ corresponds to an action of $\widehat{G}$ on $X$. It comes from the fact that the endomorphism prestack is an endomorphism object in the sense of Lurie. If $G$ is a formal group, we will call $G$-action an action of $G$ once forgetting that $G$ is formal. Remember also notation for group actions from \Cref{notation groups short}.

\begin{prop}[{\cite[4.7.1]{lurie2016higher}}] \label{automorphism stack terminal}
    Let $X \to K$ be a map of prestacks. There is an equivalence of right fibrations over $\Monoid^{\laft}_K$ 
    
    $$(\Monoid^{\laft}_{K})_{\acts X} \simeq (\Monoid^{\laft}_{K})_{/ \End_K(X)}$$

    Let $G \in \Monoid^{\laft}_K$. The above equivalence induces a natural equivalence on fibers  

    $$\Act_K(G, X) \simeq \Map_{\Monoid^{\laft}_K}(G, \End_K(X)).$$
\end{prop}

\begin{coro} 
    Let $K$ be a laft prestack, let $X \in \prstklaftdef_{/K}$. The equivalence of \Cref{automorphism stack terminal} restricts to an equivalence on formal objects $$(\FGrp_K)_{\acts X} \simeq (\FGrp_K)_{/ \infAut_K(X)}$$

    Let $G \in \FGrp_K$. The above equivalence induces a natural equivalence on fibers

    $$\Act_K(G, X) \simeq \Map_{\FGrp_K}(G, \widehat{\Aut}_K(X)).$$
\end{coro}

\begin{proof}
    By \Cref{formal groups are equivalent to formal monoids}, the forgetful map $$\Map_{\FGrp_K}(G, \widehat{\Aut}_K(X)) \to \Map_{\FMonoid_K}(G, \widehat{\Aut}_K(X))$$ is an equivalence.    
    Then, the natural map $$\Map_{\FMonoid_K}(G, \widehat{\Aut}_K(X)) \to \Map_{\Monoid^{\laft}_K}(G, \End_K(X))$$ is an equivalence because formal completion is right adjoint to the forgetful functor by \Cref{formal completion right adjoint to forgetful functor}.
    Using the second assertion of \Cref{automorphism stack terminal}, we get the second assertion of the proposition.
    We also get an equivalence from $(\FGrp_K)_{/\infAut_K(X)}$ to the full subcategory of $(\Monoid_K)_{/\End_K(X)}$ whose objects are formal monoids. Then, using the first assertion of \Cref{automorphism stack terminal}, we get the first assertion of the proposition.
\end{proof}

Moreover, it is classical that, in the setting of an adjunction between \v{C}ech nerve and quotient, an action of a group $G$ on $X$ corresponds to a map $X \to BG$.

\begin{lemma} \label{action is a slice over BG}
    Let $G \in \FGrp_K$. Let $X \in \prstklaftdef_{/K}$. The quotient functor from \Cref{equivalence formal moduli stacks formal groupoids} induces a map $$B^{\deff}_G \colon \Act_K(G, X) \to \prstklaftdef_{/B^{\deff}G} \times_{\prstklaftdef_{/K}} \{X\}.$$ and the \v{C}ech nerve functor induces a map $$\Omega_G \colon \prstklaftdef_{/B^{\deff}G} \times_{\prstklaftdef_{/K}} \{X\} \to \Act_K(G, X).$$
    Those functors form an equivalence of categories $$\Act_K(G, X) \simeq \prstklaftdef_{/B^{\deff}G} \times_{\prstklaftdef_{/K}} \{X\}$$ which is natural in $G$.
    
\end{lemma}

\begin{proof}
    We follow the usual argument, from \cite[4.28]{khan2023lecturesalgebraicstacks} for instance.

    If $U \to G$ is a map from a groupoid $U$ to $G$ which encodes an action of $G$ on $X$, $U$ is a formal groupoid so one can consider the associated quotient map $B^{\deff}U \to B^{\deff}G$ where $B^{\deff}U$ is an object of $\prstklaftdef_K$ by \Cref{equivalence formal moduli stacks formal groupoids}. But by definition of a group action, $U \times_G \diag_K \simeq \diag_X$ so by \Cref{colimits are universal} there is an equivalence $B^{\deff}U \times_{B^{\deff}G} B^{\deff}(\diag_K) \simeq B^{\deff}(\diag_X)$ hence an equivalence $B^{\deff}U \times_{B^{\deff}G} K \simeq X$. Thus the map $B^{\deff}_G$ of the lemma is well defined.

    Reciprocally, consider a pullback diagram in $\prstklaftdef_K$
\[\begin{tikzcd}
	X & K \\
	\X & B^{\deff}G
	\arrow[from=1-1, to=1-2]
	\arrow[from=1-1, to=2-1]
	\arrow["\lrcorner"{anchor=center, pos=0.125}, draw=none, from=1-1, to=2-2]
	\arrow[from=1-2, to=2-2]
	\arrow[from=2-1, to=2-2]
\end{tikzcd}\]
    The \v{C}ech nerve $X \times_{\X} X$  is equivalent to the base change $\X \times_{B^{\deff}G} G_{\bullet}$. This groupoid encodes an action of $G$ on $X$, so the map $\Omega_G$ of the lemma is well defined.
    Moreover, the map $K \to B^{\deff}G$ is a nil-isomorphism so the map $X \to \X$ is a nil-isomorphism. $X$ is laft and $\X$ is laft so $X \to \X$ is laft. Finally $\X$ has deformation theory so  $\X \in \modstk_{X/K}$. Hence, by \Cref{equivalence formal moduli stacks formal groupoids},  $B^{\deff} \Omega \X \simeq \X$.
    Finally, if a groupoid $U$ encodes an action of $G$ on $X$, $X \times_{B^{\deff}U} X \simeq U$ by \Cref{equivalence formal moduli stacks formal groupoids}.
    Thus, the maps $B^{\deff}_G$ and $\Omega_G$ form an equivalence.
    This equivalence is natural in $G$ because $\Omega_G$ is given by the formula $\Omega_G \X = (\X \times_{B^{\deff}G} \X) \times_{\X} (\X \times_{B^{\deff}G} \X)$ which is natural in $G$.
\end{proof}

Finally, we can deduce that the deformation functor of $X$ is the delooping of the formal automorphisms prestack of $X$, up to the potential discrepancy between formal moduli stacks and formal moduli problems when the base is not affine. Thus, their tangent complex is the same (up to a shift).

\begin{coro} \label{equivalence def aut}
    Let $X \in \prstklaftdef_{/K}$. Denote the structure map $f \colon X \to K$. 
    There is an equivalence of functors
    $$\Def_X^{K/K} \simeq B^{\deff}\infAut_K(X)_{|\Schaffnilop{K/K}}$$

    Moreover, if $Z \in \modstk_{K/K}$, there is an equivalence $$\Map_{\modstk_{K/K}}(Z, \Def_X^{K/K}) \simeq \prstklaftdef_{/Z} \times_{\prstklaftdef_{/K}} \{X\}$$ 

    As a consequence, $\Def_X^{K/K} \in \FMP_{K/K}$.
    Moreover, if $K$ is affine and coherent, there is an equivalence in $\QC{K}$:
    $$T_{\Def_X^{K/K}} \simeq f_*T_{X \sto K}[1].$$
\end{coro}

\begin{rk}
    In particular, \Cref{equivalence def aut} shows that for all $Z \in \modstk_{K/K}$, the $\infty$-category $\prstklaftdef_{/Z} \times_{\prstklaftdef_{/K}} \{X\}$ is in fact equivalent to a space.
\end{rk}

\begin{proof}
    Let $Z \in \modstk_{K/K}$. 

    By the equivalence between split moduli stacks and formal groups from \Cref{equivalence formal moduli stacks formal groupoids}, there is a natural equivalence $$\Map_{\modstk_{K/K}}(Z, B^{\deff} \infAut_K(X)) \simeq \Map_{\FGrpd_K}(\Omega Z, \widehat{\Aut}_K(X)).$$
    By \Cref{automorphism stack terminal}, there is a natural equivalence $$\Map_{\FGrpd_K}(\Omega Z, \widehat{\Aut}_K(X)) \simeq \Act_K(\Omega Z, X).$$
    By \Cref{action is a slice over BG} there is a natural equivalence $$Act_K(\Omega Z, X) \simeq \prstklaftdef_{/B^{\deff}\Omega Z} \times_{\prstklaftdef_{/K}} \{X\}.$$
    But $\prstklaftdef_{/B^{\deff}\Omega Z} \times_{\prstklaftdef_{/K}} \{X\} \simeq \prstklaftdef_{/Z} \times_{\prstklaftdef_{/K}} \{X\}$ and when $Z \in \Schaffnil$, $\prstklaftdef_{/Z} \times_{\prstklaftdef_{/K}} \{X\} = \Def_X^{K/K}(Z)$ by definition. 
    Thus we get equivalences $$B^{\deff}\infAut_K(X)(Z) \simeq \Def_X^{K/K}(Z)$$ that are natural in $Z \in \Schaffnil_{K/K}$, hence an equivalence of functors $$B^{\deff}\infAut_K(X)_{|\Schaffnilop{K/K}} \simeq \Def_X^{K/K}.$$

    When $K$ is affine and coherent, $\Def_X^{K/K}$ can be identified with $B^{\deff}\infAut_K(X) \in \modstk_{K/K}$ by \Cref{equivalence modstk FMP} and, denoting the pointing $e \colon K \to \infAut_K(X)$, we have the equivalences $$T_{K \sto \Def_X^{K/K}} \simeq T_{K \sto B^{\deff}\infAut_K(X)} \simeq e^*T_{\infAut_K(X) \sto K} \simeq f_*T_{X \to K}.$$ 
\end{proof}

 Observe now that the functor of deformations of a map $Y \to K'$ over $K$ appears as the pushforward of the functor of deformations of $Y$ over $K'$, again up to the potential discrepancy between formal moduli stacks and formal moduli problems when the base is not affine. Thus, we get its tangent complex without more computations.

\begin{prop} \label{tangent deformation of map}
    Let $g \colon K' \to K$ be a  map of prestacks. Let $Y \in \prstklaftdef_{/K'}$. Write the structure map $f \colon Y \to K'$. There is an equivalence of functors

    $$\Def_{Y \to K'}^{K/K} \simeq g^{\nil}_* (B^{\deff} \infAut_{K'}(Y))_{|\Schaffnilop{K/K}}$$

    In particular, if $K$ is affine and coherent, then $\Def_{Y \to K'}^{K/K} \in \FMP_{K/K}$.
    Moreover, there is an equivalence in $\QC{K}$:  $$T_{\Def_{Y \to K'}^{K/K}} \simeq g_* \circ f_*(T_{Y \to K'})[1].$$
\end{prop}

\begin{proof}
    Let $Z \in \Schaffnil_{K/K}$. By the base change adjunction of \Cref{base change adjunction}, there is a natural equivalence

    $$f_* \Def_Y^{K'/K'}(Z) \simeq \Map_{\modstk_{K'/K'}}(Z \times_K K', \Def_Y^{K'/K'}).$$
    But by \Cref{equivalence def aut}, we know that there is a natural equivalence  $$\Map_{\modstk_{K'/K'}}(Z \times_K K', \Def_Y^{K'/K'}) \simeq \Def_{Y \to K'}^{K/K}(Z),$$ 
    hence the expected equivalence of functors.

    Now suppose that $K$ is affine. Then $T_{K \sto \Def_{Y \to K'}^{K/K}} \simeq T_{K \sto g_* (B^{\deff}\infAut_{K'}(Y))}$.
    But by \Cref{tangent commutes with adjunction}, there is an equivalence $$T_{K \sto g_* (B^{\deff}\infAut_{K'}(Y))} \simeq g_* T_{K' \sto B^{\deff}\infAut_{K'}(Y)}.$$
    Moreover, by \Cref{equivalence def aut}, there is an equivalence $$T_{K' \sto B^{\deff}\infAut_{K'}(Y)} \simeq f_* T_{Y \to K'}.$$
    Hence $$T_{K \sto \Def_{Y \to K'}^{K/K}} \simeq g_* \circ f_*(T_{Y \to K'}).$$ 
\end{proof}

Finally, consider a section map $e \colon K \to K'$. Using the object $\Def_Y^{K'/K'}$ which lives in $\FMP_{K'/K'}$, we can reformulate the transformation $\fibdef_Y^{K/K}$ from deformations of a map $Y \to K'$ to deformations of its fiber $X = Y \times_{K'} K$ over $K$: it corresponds simply to the map between the pushforward $f_* \Def_Y^{K'/K'}$ and the pullback $e^* \Def_Y^{K'/K'}$. It is thus easy to deduce what is the corresponding induced map at the level of tangent complexes.

\begin{notation} \label{transformation map from adjunction}
    Let $f \colon K' \to K$ be a laft map of prestacks, with a section $e \colon K \to K'$.
    \begin{enumerate}

    \item if $\M \in \modstk_{K'/K'}$, we will write $F_{\M} \colon f_*^{\nil} \M \to e^* \M$ for the map in $\modstk_{K/K}$ which comes from the counit map $f^*f_*^{\nil} \M \to \M$ to which we apply the functor $e^*$.

    \item if $M \in \QC{K'}$, we will also write $F_{M} \colon f_* M \to e^* M$ for the map in $\QC{K}$ which comes from the counit map $f^*f_* M \to M$ to which we apply the functor $e^*$.
    \end{enumerate}
    By \Cref{base change adjunction}, the tangent map associated to $F_{\M}$ identifies with $F_{T_{\M}}$.
\end{notation}

\begin{prop} \label{fiber of deformation tangent equi}
    Let $K$ be a coherent and affine scheme, let $K'$ be a prestack with a laft map $g \colon K' \to K$, with a section $e \colon K \to K'$, and let $Y \in \prstklaftdef_{K'}$. Write $f \colon Y \to K'$ and denote $X = Y \times_{K'} K$.
    Remember from \Cref{def fiber of deformations} the map $$\fibdef_Y^{K/K} \colon \Def_{Y \to K'}^{K/K} \to \Def_X^{K/K}.$$
    Then, up to the equivalences of \Cref{pullback aut} and \Cref{tangent deformation of map}, the map $\fibdef_Y^{K/K}$ identifies with the map $F_{B^{\deff}\infAut_{K'}(Y)}$ from \Cref{transformation map from adjunction}, in the sense that the following diagram is commutative:
    
\[\begin{tikzcd}
	{\Def_{Y \to K'}^{K/K}} & {\Def_X^{K/K}} \\
	{g^{\nil}_*B^{\deff}\infAut_{K'}(Y)} & {e^*B^{\deff}\infAut_{K'}(Y)}
	\arrow["{\fibdef_Y^{K/K}}", from=1-1, to=1-2]
	\arrow["\simeq"', from=1-1, to=2-1]
	\arrow["\simeq", from=1-2, to=2-2]
	\arrow["{F_{B^{\deff}\infAut_{K'}(Y)}}"'{yshift=-0.8ex}, from=2-1, to=2-2]
\end{tikzcd}\]
    
    As a consequence, in $\QC{K}$ the map $T(\fibdef_Y^{K/K})$ identifies with the map $F_{f_*T_{Y \to K'}}$, in the sense that the following diagram is commutative: 

\[\begin{tikzcd}
	{T_{\Def_{Y \to K'}^{K/K}}} & {T_{\Def_X^{K/K}}} \\
	{g_* \circ f_*(T_{Y \to K'})[1]} & {e^* \circ f_*(T_{Y \to K'})[1]}
	\arrow["{T(\fibdef_Y^{K/K})}", from=1-1, to=1-2]
	\arrow["{\simeq }"', from=1-1, to=2-1]
	\arrow["\simeq", from=1-2, to=2-2]
	\arrow["{F_{f_*T_{Y \to K'}}}"'{yshift=-0.8ex}, from=2-1, to=2-2]
\end{tikzcd}\]
    
\end{prop}

\begin{proof}
    Let us show that the map of formal groups $$\Omega (\fibdef_Y^{K/K}) \colon \Omega \Def_{Y \to K'}^{K/K} \to \Omega \Def_X^{K/K}$$ identifies with the map $$F_{\infAut_{K'}(Y)} \colon g^{\nil}_* \infAut_{K'}(Y) \to \infAut_{K}(X).$$
Let $Z \in \Schaffart_{K/K}$. Loops can be computed in prestacks so it suffices to prove a natural identification after evaluating at $Z$.
Remember that $\fibdef_Y^{K/K} \colon \Def_{Y \to K'}^{K/K} \to \Def_X^{K/K}$ is the natural transformation which sends the front square to the back square in the following diagram where all the squares are pullback squares: 
\[\begin{tikzcd}
	& X && K \\
	Y && {K'} \\
	& {\Y \times_{K' \times_K Z} Z} && Z \\
	\Y && {K' \times_K Z}
	\arrow[from=1-2, to=1-4]
	\arrow[dashed, from=1-2, to=3-2]
	\arrow[from=1-4, to=3-4]
	\arrow[from=2-1, to=1-2]
	\arrow[from=2-1, to=2-3]
	\arrow[from=2-1, to=4-1]
	\arrow[from=2-3, to=1-4]
	\arrow[from=2-3, to=4-3]
	\arrow[dashed, from=3-2, to=3-4]
	\arrow[dashed, from=4-1, to=3-2]
	\arrow[from=4-1, to=4-3]
	\arrow[from=4-3, to=3-4]
\end{tikzcd}\]
    Let $\Fo \colon f^{\nil}_*\End_{K'}(Y) \to \End_{K}(X)$ be the natural transformation which, for all $Z \in \Schaffart_{K/K}$, is defined as the natural map $$\Map_{K'}(Z \times_K K' \times_{K'} Y, Y) \to \Map_{K}(Z \times_K X, X)$$ induced by the pullback along $K \to K'$.
    $\Omega(\fibdef)$ identifies with the formal completion of $\Fo$, $\widehat{\Fo} \colon f_* \infAut_{K'}(Y) \to \infAut_{K}(X)$.

    Now let us prove that the map $\Fo$ identifies with the map $F_{\End_{K'}(Y)}$.
    Let $Z \in \Schaffart_{K'/K'}$. The counit $f^*f_* \colon \End_{K'}(Y) \to \End_{K'}(Y)$ identifies with the natural transformation $$\Map_{K'}(Z \times_K K' \times_{K'} Y, Y) \to \Map_{K'}(Z \times_{K'} Y, Y)$$ induced by the pullback along $Z \to Z \times_K K'$.
    Now if one applies the pullback $e^*$, that is, restriction along the section $e \colon K \to K'$, one gets for all $Z \in \Schaffnil_{K/K}$ a pullback diagram 

\[\begin{tikzcd}
	Z & K \\
	{Z \times_K K'} & {K'}
	\arrow[from=1-1, to=1-2]
	\arrow[from=1-1, to=2-1]
	\arrow["e", from=1-2, to=2-2]
	\arrow[from=2-1, to=2-2]
\end{tikzcd}\]
    Hence, the natural transformation  $e^*f^*f_* \colon \Map_{K'}(Z \times_K K' \times_{K'} Y, Y) \to \Map_{K}(Z \times_K X, X)$ is the map induced by the pullback along $K \to K'$, that is, identifies with the map $\Fo$.
\end{proof}

\subsubsection{Mixed deformations of prestacks with deformation theory} \label{section mixed deformations}

When one looks at the definition of the functor of deformations $\Def_Y^{K/K}$, one realizes that the fact that affine schemes live \textit{over} $K$ does not come into play in the definition. In fact, there is an obvious extension of this functor to affine schemes with a map from $K$ which is a nil-isomorphism (\textit{under} $K$), but defined \textit{over} any base $W$. The typical example is the case where $K = \Spec(k)$ and $W = \Spec(W(k))$.
In this section, we essentially show that the results of \Cref{section split deformations} hold in the mixed-bases situation.
By \Cref{tangent conservative}, the restriction functor from mixed formal moduli stacks to split formal moduli stacks is conservative, and the tangent complex can be computed on the restriction. Hence, we will deduce most of the properties in the mixed case from the equivalent results in the split case.

However, we cannot deduce the fact that the functor of mixed deformations is a formal moduli stack from the equivalent for the functor of split deformations. In the split case, we proved that fact by looping the problem and proving that the infinitesimal automorphism prestack was a formal group. In the mixed case, it could be an option to follow a similar method: we would get an ``infinitesimal automorphism groupoid'', which is however less tractable than  the infinitesimal automorphism group. We choose not to explore this path here, but rather to prove by hand that the functor of mixed deformation is a formal moduli stack. A counterpart of this choice is that we manage to prove that result only when the bases are affine, and thus that formal moduli stacks correspond to formal moduli problems.

Let us begin by a lemma which guarantees that the deformation functors that we will define land in $\infty$-groupoid rather than in $\infty$-categories.

\begin{lemma} \label{fiber conservative}
    Let $K' \to Z$ be a nil-isomorphism of prestacks. Then the fiber functor $$\prstklaftdef_{Z} \to \prstklaftdef_{K'}$$ is conservative.
\end{lemma}

\begin{proof}
    Let $G \colon \Y_1 \to \Y_2$ be a map in $\prstklaftdef_{Z}$. Denote $Y_1 = \Y_1 \times_Z K'$ and $Y_2 = \Y_2 \times_Z K'$. Denote the structure maps $f_1 \colon Y_1 \to K'$ and $f_2 \colon Y_2 \to K'$. Suppose that the induced map $g \colon Y_1 \to Y_2$ is an equivalence. Given that the map $K' \to Z$ is a nil-isomorphism, the natural maps $Y_1 \to \Y_1$ and $Y_2 \to \Y_2$ are also nil-isomorphisms, and the composite $Y_1 \to Y_2 \to \Y_2$ is a nil-isomorphism. Hence $\Y_1$ and $\Y_2$ are formal moduli stacks under $Y_1$. By conservativity of the tangent functor on formal moduli stacks (\Cref{tangent conservative}), to prove that the map $\Y_1 \to \Y_2$ is an equivalence, it suffices to prove that the induced tangent map $T(G) \colon T_{Y_1 \sto \Y_1} \to T_{Y_1 \sto \Y_2}$ is an equivalence. But by \Cref{properties tangent complex}, there are equivalences $T_{Y_1 \sto \Y_1} \simeq f_1^*T_{K' \sto Z}$ and $T_{Y_1 \sto \Y_2} \simeq g^*T_{Y_2 \sto \Y_2} \simeq g^*f_2^*T_{K' \sto Z}$. Under those equivalences, the map $T(G)$ identifies with the map $f_1^*T_{K' \sto Z} \to g^*f_2^*T_{K' \sto Z}$ which is an equivalence because the map $g \colon Y_1 \to Y_2$ is an equivalence. This concludes the proof.
\end{proof}

Now, as in the split case, we can define the functors of deformation and the ``fiber of deformations'' map. Those are the natural extension of the definitions in the split case.

\begin{defi}
    Let $K \to W$ be a map of prestacks and let $X \in \prstklaftdef_K$.
    Let $$\Def_X^{K/W} \colon \Schaffnilop{K/W} \to \hat{\Spaces}$$ be the functor that sends $Z \in \Schaffnil_{K/W}$ to $$\Def_X^{K/W}(Z) = \prstklaftdef_{/Z} \times_{\prstklaftdef_{/K}} \{X\}.$$
    In other words, if $Z \in \Schaffnil_{K/W}$, $\Def_X^{K/W}(Z)$ identifies with the $\infty$-groupoid of pullback diagrams

\[\begin{tikzcd}
	Y & K \\
	\Y & Z
	\arrow[from=1-1, to=1-2]
	\arrow[from=1-1, to=2-1]
	\arrow["\lrcorner"{anchor=center, pos=0.125}, draw=none, from=1-1, to=2-2]
	\arrow[from=1-2, to=2-2]
	\arrow[from=2-1, to=2-2]
\end{tikzcd}\]
\end{defi}

\begin{defi} \label{defi deformation of maps mixed}
    Let $W' \to W$ be a laft map of prestacks, let $K \in \prstklaftdef_{K/W}$. Write the fiber $K' = K \times_W W'$, and let $Y \in \prstklaftdef_{K'}$. An affine scheme $Z \in \Schaffnil_{K/W}$ sits in a pullback diagram
\[\begin{tikzcd}
	{K'} & {W'\times_W Z} & {W'} \\
	K & Z & W
	\arrow[from=1-1, to=1-2]
	\arrow[from=1-1, to=2-1]
	\arrow["\lrcorner"{anchor=center, pos=0.125}, draw=none, from=1-1, to=2-2]
	\arrow[from=1-2, to=1-3]
	\arrow[from=1-2, to=2-2]
	\arrow["\lrcorner"{anchor=center, pos=0.125}, draw=none, from=1-2, to=2-3]
	\arrow[from=1-3, to=2-3]
	\arrow[from=2-1, to=2-2]
	\arrow[from=2-2, to=2-3]
\end{tikzcd}\]
    Using \Cref{fiber conservative}, let $\Def_{Y \to K'}^{K/W} \colon \Schaffnilop{K/W} \to \hat{\Spaces}$ be the functor which sends $Z \in \Schaffnil_{K/W}$ to $$\prstklaftdef_{W' \times_W Z} \times_{\prstklaftdef_{K'}} \{Y\},$$
    In other words, for each $Z \in \Schaffnil_{K/W}$, $\Def_{Y \to K'}^{K/W}(Z)$ identifies with the $\infty$-groupoid of pullback diagrams
\[\begin{tikzcd}
	Y & K' \\
	\Y & {W' \times_W Z}
	\arrow[from=1-1, to=1-2]
	\arrow[from=1-1, to=2-1]
	\arrow["\lrcorner"{anchor=center, pos=0.125}, draw=none, from=1-1, to=2-2]
	\arrow[from=1-2, to=2-2]
	\arrow[from=2-1, to=2-2]
\end{tikzcd}\]
\end{defi}

\begin{defi} \label{def fiber of deformations mixed}
    Let $W' \to W$ be a laft map of prestacks with a section $W \to W'$, let $K \in \prstklaftdef_{/W}$. Write the fiber $K' = W' \times _W K$ and let $Y \in \prstklaftdef_{K'}$. Denote $X = Y \times_{K'} K$.
    An affine scheme $Z \in \Schaffnil_{K/W}$ sits in a pullback diagram

\[\begin{tikzcd}
	K & {K'} \\
	Z & {W' \times_W Z}
	\arrow[from=1-1, to=1-2]
	\arrow[from=1-1, to=2-1]
	\arrow["\lrcorner"{anchor=center, pos=0.125}, draw=none, from=1-1, to=2-2]
	\arrow[from=1-2, to=2-2]
	\arrow[from=2-1, to=2-2]
\end{tikzcd}\]
    Thus, the pullback along $Z \to W' \times_W Z$ induces a natural map $$\fibdef_Y^{K/W}(Z) \colon \prstk_{/W' \times_W Z} \times_{\prstk_{/K'}} \{Y\} \to \prstk_{/Z} \times_{\prstk_{/K}} \{X\}$$
    Those natural maps $\fibdef_Y^{K/W}(Z)$ give a natural transformation $$\fibdef_Y^{K/W} \colon \Def_{Y \to K'}^{K/W} \to \Def_X^{K/W}.$$
    In other words, $\fibdef_Y^{K/W}(Z) \colon \Def_{Y \to K'}^{K/W}(Z) \to \Def_X^{K/W}(Z)$ is the map which sends a pullback diagram

\[\begin{tikzcd}
	Y & K' \\
	\Y & {W' \times_W Z}
	\arrow[from=1-1, to=1-2]
	\arrow[from=1-1, to=2-1]
	\arrow["\lrcorner"{anchor=center, pos=0.125}, draw=none, from=1-1, to=2-2]
	\arrow[from=1-2, to=2-2]
	\arrow[from=2-1, to=2-2]
\end{tikzcd}\]
    to the pullback diagram 
\[\begin{tikzcd}
	X & K \\
	{\Y \times_{W' \times_W Z}} & Z
	\arrow[from=1-1, to=1-2]
	\arrow[from=1-1, to=2-1]
	\arrow["\lrcorner"{anchor=center, pos=0.125}, draw=none, from=1-1, to=2-2]
	\arrow[from=1-2, to=2-2]
	\arrow[from=2-1, to=2-2]
\end{tikzcd}\]
\end{defi}

\begin{rk}  \label{restriction to split case}
    The above constructions, over mixed bases, extend the constructions in the split case from \Cref{section split deformations}. Let $U$ denote the forgetful functor $\Schaffnilop{K/K} \to \Schaffnilop{K/W}$ induced by the map $K \to W$. There are natural equivalences of functors
    \begin{enumerate}
        \item $\Def_{X}^{K/W} \circ U \simeq \Def_{X}^{K/K}$

        \item $\Def_{Y \to K'}^{K/W} \circ U \simeq \Def_{Y \to K'}^{K/K}$

        \item The restriction of the natural transformation $\fibdef_Y^{K/W}$ to $\Schaffnilop{K/K}$ identifies with $\fibdef_Y^{K/K}$.
    \end{enumerate}
\end{rk}

In order to prove that the deformation functor we just defined are formal moduli problem, we need a ``clutching'' result for laft prestacks with deformation theory, analogous to the clutching result for Deligne-Mumford stacks by Lurie (\cite[Theorem 16.3.0.1]{lurie2018spectral}. In general, such clutching result comes from three important results for pushouts along ``nice'' maps: existence of pushouts, universality of pushouts, and descent over pushouts. We have proven a result analogous to universality of pushouts in \Cref{pushout affine def}. Now we prove that the pushouts along a square zero extension and an affine nil-monomorphism exist and are again square zero extensions. It will be useful to prove descent over ``nice'' pushouts by reducing to the case of affine schemes.

\begin{lemma} \label{pushout sqz}
    Let $f \colon X \to Y$ be a map in $\prstklaftdef$, let $M \in \qc_{X}$ and let $X' = \sqz_X M$. Suppose that the map $f$ is an affine monomorphism, and a nil-isomorphism. Then the diagram
    \[\begin{tikzcd}
	{X} & {X'} \\
	{Y} & {\sqz_Y(f_*M)}
	\arrow[from=1-1, to=1-2]
	\arrow["f"', from=1-1, to=2-1]
	\arrow[from=1-2, to=2-2]
	\arrow[from=2-1, to=2-2]
\end{tikzcd}\]
is a pushout diagram in $\prstklaftdef$.
\end{lemma}

\begin{proof}
    The prestacks $Y$ and $X'$ are formal moduli stacks under $X$. By \cite[Chapter 5, Corollary 2.3.6]{gaitsgory2019study}, the pushout $Y' := Y \cup_X X'$ exists in the category of formal moduli stacks under $X$, and the tangent functor commutes with the pushout. Hence, $T_{X \sto Y} \simeq T_{X \sto X'} \oplus T_{X' \sto Y}$ (1).  By conservativity of the tangent functor on formal moduli stacks (\Cref{tangent conservative}), to prove that the map $Y' \to \sqz_{Y}(f_*M)$ is an equivalence, it suffices to prove that the induced tangent map $T_{Y \sto Y'} \to T_{Y \sto \sqz_{Y}(f_*M)}$ is an equivalence. Considering the fiber sequence $$T_{X \sto Y'} \to T_{X \sto \sqz_Y(f_*M)} \to f^*i^*T_{Y' \sto \sqz_Y(f_*M)},$$ it suffices to prove that $f^*i^*T_{Y' \sto \sqz_Y(f_*M)} = 0$. Consider the fiber sequence $$T_{Y \sto Y'} \to T_{Y \sto \sqz_Y(f_*M)} \to i^*T_{Y' \sto \sqz_Y(f_*M)}.$$ Composing with the pullback $f^*$, we get a fiber sequence $$f^*T_{Y \sto Y'} \to f^*T_{Y \sto \sqz_Y(f_*M)} \to f^*i^*T_{Y' \sto \sqz_Y(f_*M)}.$$ Hence it suffices to prove that the map $\phi : f^*T_{Y \sto Y'} \to f^*T_{Y \sto \sqz_Y(f_*M)}$ is an equivalence. Now the fiber sequence $$T_{X \sto Y} \to T_{X \sto Y'} \to f^*T_{Y \sto Y'}$$ and the equivalence (1) imply the equivalence $f^*T_{Y \sto Y'} \simeq T_{X \sto X'}$. Given the equivalence $T_{X \sto X'} \simeq M^{\vee}[-1]$ in $\QC{X}$, we get the equivalence $f^*T_{Y \sto Y'} \simeq M^{\vee}[-1]$. Moreover, there is the equivalence $$T_{Y \sto \sqz_Y(f_*M)} \simeq (f_*M)^{\vee}[-1]$$ in $\QC{Y}$. But by \cite[B.15(3)]{brantnermagidsonnuiten2025}, pullback commutes with dualizing, hence the equivalence $$f^*T_{Y \sto \sqz_Y(f_*M)} \simeq (f^*f_*M)^{\vee}[-1].$$ Under those equivalences, the map $\phi$ identifies with the dual map of the shift of the counit map $M \to f^*f_*M$. Now the map $f$ is a monomorphism so there is a pullback diagram 
\[\begin{tikzcd}
	X & X \\
	X & Y
	\arrow["\id", from=1-1, to=1-2]
	\arrow["\id"', from=1-1, to=2-1]
	\arrow["\lrcorner"{anchor=center, pos=0.125}, draw=none, from=1-1, to=2-2]
	\arrow["f", from=1-2, to=2-2]
	\arrow["f"', from=2-1, to=2-2]
\end{tikzcd}\] 
    Now $f$ is affine so it is in particular schematic and quasi-compact, so by \cite[Chapter 3, 2.2.2]{gaitsgory2019studyI}, the base change morphism $f^*f_*M \to M$ is an equivalence. Hence the map $\phi$ is an equivalence, which concludes the proof.

\end{proof}

\begin{prop}
    Let $K \to W$ be a map of coherent affine schemes. \begin{enumerate}
    \item Let $X \in \prstklaftdef_{K/W}$. Then $\Def_X^{K/W} \in \FMP_{K/W}$.
    \item Let $W' \to W$ be a laft map of prestacks. Write $K' = W' \times_W K$ and let $Y \in \prstklaftdef_{K'/W'}$. Then $\Def_{Y \to K'}^{K/W} \in \FMP_{K/W}.$
    \end{enumerate}
\end{prop}

\begin{proof}
    (1) is a consequence of (2) by taking $W' = W$. Let us prove (2).
    First, it is clear that $\Def_{Y \to K'}^{K/W}(K) \simeq *$. Then by \Cref{equivalence modstk FMP} it suffices to prove that $\Def_{Y \to K'}^{K/W}$ is infinitesimally cohesive on $\Schaffart_{K/W}$.
    Consider a pushout diagram in $\Schaffart_{K/W}$
    \[\begin{tikzcd}
	{Z_{01}} & {Z_1} \\
	{Z_0} & Z
	\arrow[from=1-1, to=1-2]
	\arrow[from=1-1, to=2-1]
	\arrow[from=1-2, to=2-2]
	\arrow[from=2-1, to=2-2]
	\arrow["\lrcorner"{anchor=center, pos=0.125, rotate=180}, draw=none, from=2-2, to=1-1]
\end{tikzcd}\]
where the map $Z_{01} \to Z_0$ is a square zero extension and $Z_{01} \to Z_1$ is a nilclosed embedding. Denote $X_{01} = W' \times_W Z_{01}$, $X_0 = W' \times_W Z_0$, $X_1 = W' \times_W Z_1$ and $X = W' \times_W Z$. As in \Cref{pushout affine def}, the square 
 \[\begin{tikzcd}
	{X_{01}} & {X_1} \\
	{X_0} & X
	\arrow[from=1-1, to=1-2]
	\arrow[from=1-1, to=2-1]
	\arrow[from=1-2, to=2-2]
	\arrow[from=2-1, to=2-2]
\end{tikzcd}\]
is a pushout diagram in $\prstklaftdef$. Moreover the map $X_{01} \to X_0$ is a square zero extension and the map $X_{01} \to X_1$ is an affine monomorphism, and a nil-isomorphism. We adapt the proof of \cite[Theorem 16.3.3.1 and Theorem 16.3.0.1]{lurie2018spectral} for $\prstklaftdef$ instead of Deligne-Mumford stacks. It suffices to prove that the square 
 \[\begin{tikzcd}
	{\prstklaftdef_X} & {\prstklaftdef_{X_1}} \\
	{\prstklaftdef_{X_0}} & {\prstklaftdef_{X_{01}}}
	\arrow[from=1-1, to=1-2]
	\arrow[from=1-1, to=2-1]
	\arrow[from=1-2, to=2-2]
	\arrow[from=2-1, to=2-2]
\end{tikzcd}\]
is a pullback square.
Write the canonical map $$G \colon \prstklaftdef_X \to \prstklaftdef_{X_0} \times_{\prstklaftdef_{X_{01}}} \prstklaftdef_{X_1}.$$ Objects in the target can be identified with diagrams 
\[\begin{tikzcd}
	{Y_0} & {Y_{01}} & {Y_1} \\
	{X_0} & {X_{01}} & {X_1}
	\arrow[from=1-1, to=2-1]
	\arrow[from=1-2, to=1-1]
	\arrow[from=1-2, to=1-3]
	\arrow["\lrcorner"{anchor=center, pos=0.125, rotate=-90}, draw=none, from=1-2, to=2-1]
	\arrow[from=1-2, to=2-2]
	\arrow["\lrcorner"{anchor=center, pos=0.125}, draw=none, from=1-2, to=2-3]
	\arrow[from=1-3, to=2-3]
	\arrow[from=2-2, to=2-1]
	\arrow[from=2-2, to=2-3]
\end{tikzcd}\]
in $\prstklaftdef$ where both squares are pullbacks squares. The map $Y_{01} \to Y_0$ is a square zero extension and the map $Y_{01} \to Y_1$ is an affine monomorphism, and a nil-isomorphism. By \Cref{pushout sqz}, the pushout $Y = Y_0 \cup_{Y_{01}} Y_1$ exists in $\prstklaftdef$ . This pushout construction determines a functor $$F \colon \prstklaftdef_{X_0} \times_{\prstklaftdef_{X_{01}}} \prstklaftdef_{X_1} \to \prstklaftdef_X$$ which is left adjoint to $G$. By \Cref{pushout universal},  base change preserves pushouts so the counit map $F \circ G \to \id$ is an equivalence. Let us prove that the unit $\id \to G \circ F$ is an equivalence. Consider a diagram 
\[\begin{tikzcd}
	& {Y_{01}} && {Y_1} \\
	{Y_0} && Y \\
	& {X_{01}} && {X_1} \\
	{X_0} && X
	\arrow[from=1-2, to=1-4]
	\arrow[from=1-2, to=2-1]
	\arrow[dashed, from=1-2, to=3-2]
	\arrow[from=1-4, to=2-3]
	\arrow[from=1-4, to=3-4]
	\arrow[from=2-1, to=2-3]
	\arrow[from=2-1, to=4-1]
	\arrow[from=2-3, to=4-3]
	\arrow[dashed, from=3-2, to=3-4]
	\arrow[dashed, from=3-2, to=4-1]
	\arrow[from=3-4, to=4-3]
	\arrow[from=4-1, to=4-3]
\end{tikzcd}\]
where the bottom and top squares are pushouts, and the back and left squares are pullbacks. We need to prove that the front and right squares are pullbacks. The maps $Y_{0} \to Y$ and $Y_1 \to Y$ are affine. Hence, we can reduce to the case where $Y$ and $X$ are affine schemes, which implies that $Y_0$, $Y_1$, $Y_{01}$, $X_0$, $X_1$ and $X_{01}$ are affine schemes as well. The result then follows from the analogous result for affine schemes.

Finally, we need to prove that for all $Z \in \Schaffart_{K/W}$, $\Def_{Y \to K'}^{K/W}(Z)$ is essentially small. By \Cref{equivalence def aut}, we know that it is true for all $Z \in \Schaffart_{K/K}$. But by \cite[Proposition 4.3]{brantnermagidsonnuiten2025}, all nilpotent extensions $Z \in \Schaffart_{K/W}$ can be decomposed as a composition of square zero extensions by $K$-algebras. But $\Def_{Y \to K'}^{K/W}$ is infinitesimally cohesive so by induction $\Def_{Y \to K'}^{K/W}(Z)$ is also essentially small.
\end{proof}

\begin{rk}
    In \Cref{section split deformations}, we proved that the functor $\Def_X^{K/K}$ was a formal moduli stack by showing that the looping $\Omega_K \Def_X^{K/K}$ identified with the infinitesimal automorphism prestack $\infAut_K(X)$, which is a formal group (here the pointing map $K(Z) \to \Def_X^{K/K}(Z)$ sends a $K$-algebra $Z$ to the trivial deformation of $X$ over $Z$). In the mixed case, another approach to show that the functor $\Def_X^{K/W}$ is a formal moduli stack would consist in proving that the looping $\Omega_K \Def_X^{K/W}$ is a formal groupoid. In that case the pointing $K(Z) \to \Def_X^{K/W}(Z)$ (a map over $W$) is the map which sends a structure of $K$-algebra over an affine scheme $Z$ to the trivial deformation of $X$ over $Z$ along this structure of $K$-algebra. In fact, it is a more natural pointing, even in the split case: indeed, even in the split case, the definition of the space $\Def_X^{K/K}(Z)$ does not depend on the $K$-algebra structure of $K$, so it is natural to take into account all possible $K$-algebra structures on $Z$ to define the pointing map.
\end{rk}

Now, computations of the tangent complexes in the mixed-bases case come directly from computations in the split case.

\begin{prop} \label{pushforward mixed deformations}
    Let $K \to W$ be a map of coherent affine schemes. \begin{enumerate}
    \item Let $X \in \prstklaftdef_{K/W}$, and write the structure map $f \colon X \to K$. Then there is an equivalence in $\QC{K}$ $$T_{\Def_X^{K/W}} \simeq f_* T_{X \sto K}[1].$$
    \item Let $W' \to W$ be a laft map of prestacks. Write $K' = W' \times_W K$ and let $Y \in \prstklaftdef_{K'/W'}$. Write the structure maps $f \colon Y \to K'$ and $g \colon K' \to K$. Then there is an equivalence in $\QC{K}$ $$T_{\Def_{Y \to K'}^{K/W}} \simeq g_* \circ f_* T_{Y \sto K'}[1].$$
    \end{enumerate}
\end{prop}

\begin{proof}
    The tangent complex can be computed on the split restriction by \Cref{tangent restriction}, so by \Cref{restriction to split case} the computation of the tangent complex here comes from the computation in the split case from \Cref{equivalence def aut}.
\end{proof}

\begin{prop} \label{fiber of deformation tangent mixed} 
Let $K \to W$ be a map of coherent affine schemes. Let $g \colon W' \to W$ be a laft map of prestacks with a section $e \colon W \to W'$. Write $K' = K \times_W W'$ and let $Y \in \prstklaftdef_{/K'}$ and $X = Y \times_{K'} K$.
    Remember from \Cref{def fiber of deformations mixed} the natural transformation $\fibdef_Y^{K/W} \colon \Def_{Y \to K'}^{K/W} \to \Def_X^{K/W}$.
    The map $T(\fibdef_Y^{K/W})$ identifies with the map $F_{f_*T_{Y \to K'}}$ from \Cref{transformation map from adjunction}, in the sense that the following diagram is commutative:

\[\begin{tikzcd}
	{T_{\Def_{Y \sto K'}^{K/W}}} & {T_{\Def_X^{K/W}}} \\
	{g_* \circ f_*(T_{Y \sto K'})[1]} & {e^* \circ f_*(T_{Y \sto K'})[1]}
	\arrow["{T(\fibdef_Y^{K/W})}", from=1-1, to=1-2]
	\arrow["{\simeq }"', from=1-1, to=2-1]
	\arrow["\simeq", from=1-2, to=2-2]
	\arrow["{F_{f_*T_{Y \sto K'}}}", from=2-1, to=2-2]
\end{tikzcd}\] 

\end{prop}

\begin{proof}
    By \Cref{tangent restriction}, the tangent functor can be computed on the split restriction. Thus, the result comes from \Cref{restriction to split case} and \Cref{fiber of deformation tangent equi}.   
\end{proof}

\subsubsection{Relation with usual deformation functors}

In this section, we take $K = \Spec(k)$ and $W = \Spec(\Lambda)$, a local Noetherian ring with residue field $k$, and $X$ a derived scheme locally of finite presentation over $k$. We show that the formal moduli problem of deformations of $X$ among prestacks with deformation theory that we constructed in \Cref{section split deformations} and \Cref{section mixed deformations}, $\Def_X^{K/W}$, is equivalent to the usual formal moduli problem of deformations of $X$ among derived schemes, constructed for instance in \cite[3.5, 3.6]{brantner2025deformationsliftscalabiyauvarieties}. It is analogous to the result in \cite[3.6]{brantner2025deformationsliftscalabiyauvarieties} stating that deformations of a derived scheme among Deligne-Mumford stacks are nothing more than deformations among derived schemes.

\begin{notation}
    For any affine scheme $Z$, write $\Sch_{/Z}$ for the $\infty$-category of derived schemes over $Z$.
\end{notation}

\begin{defi}[{\cite[3.6]{brantner2025deformationsliftscalabiyauvarieties}}]
    Let $f \colon X \to K$ be a derived scheme over $K$.
    Let $$\Def_X^{\sch} \colon \Schaffartop{K/W} \to \Spaces$$ be the functor which sends $Z \in \Schaffart_{K/W}$ to $$\Def_X^{\sch}(Z) := \Sch_Z \times_{\sch_K} \{X\}.$$
    This defines a formal moduli problem $\Def_X^{\sch} \in \FMP_{K/W}$.
    Moreover, when $X$ is locally of finite presentation, there is an equivalence in $\qc_K$ $$T_{\Def_X^{\sch}} \simeq f_*T_{X}[1].$$
\end{defi}

\begin{rk}
    If $X$ is locally of finite presentation, its tangent complex $T_X$ is perfect, so its pushforward $f_*T_X$ can be considered equivalently in the quasi-coherent or in the pro-coherent setting.
\end{rk}

\begin{prop} \label{deformations are the same}
    Let $X$ be a \comdef{proper} derived scheme locally of finite presentation over $K$. The natural forgetful transformation in $\FMP_{K/W}$ 

    $$\Def_X^{\sch} \to \Def_X^{K/W}$$
    is an equivalence.
\end{prop}

\begin{proof}
    The forgetful transformation is given naturally for $Z \in \Schaffart_{K/W}$ by the map $$\Def_X^{\sch}(Z) \to \Def_X^{K/W}(Z)$$ which sends a scheme deformation of $X$ over $Z$ to itself (a deformation of a scheme locally of finite presentation stays locally of finite presentation by \cite[7.4.3.18]{lurie2016higher}, \cite[09ZW]{stacks-project}, \cite[2.7.3.2]{lurie2018spectral}, and a scheme locally of finite presentation is a fortiori a laft prestack with deformation theory). Let $n \in \Z$. The map $$\pi_0 \Def_X^{\sch}(\sqz_K K[n]) \to \pi_0 \Def_X^{K/W}(\sqz_K K[n])$$ is injective. By adjunction, it means that the map $$\pi_{-n} T_{\Def_X^{\sch}} \to \pi_{-n} T_{\Def_X^{K/W}}$$ is injective. But by \Cref{pushforward mixed deformations}, there is an equivalence $$T_{\Def_X^{K/W}} \simeq f_*T_X[1]$$ in $\qc_K$. But the scheme $X$ is proper and locally of finite presentation so, by \cite[Theorem 3.2.2]{LurieDAGXII}, $\pi_{-n} f_*T_X$ is finite-dimensional. Hence, the injective map $\pi_{-n} T_{\Def_X^{\sch}} \to \pi_n T_{\Def_X^{K/W}}$ is an isomorphism. It means that the tangent map $$T_{\Def_X^{\sch}} \to T_{\Def_X^{K/W}}$$ is an equivalence. By conservativity of the tangent functor on formal moduli problems (\Cref{tangent conservative}), the map $$\Def_X^{\Sch} \to \Def_X^{K/W}$$ is an equivalence.
\end{proof}

Moreover, analogously to \cite{brantner2025deformationsliftscalabiyauvarieties}, if $X$ is a \comdef{proper} classical scheme locally of finite presentation, the classical deformation functor associated to the functor of deformations among prestacks with deformation theory is simply the classical deformation functor of $X$.

\begin{coro} \label{classicalisation of deformations}
    Let $X$ be a proper classical scheme locally of finite presentation.
    The classical deformation functor $(\Def_X^{K/W})^{\cl}$ associated to the formal moduli problem $\Def_X^{K/W}$ is equivalent to the functor of classical deformations of $X$, $\overl[\Def]_X$.
\end{coro}

\begin{proof}
    By \cite[3.35]{brantner2025deformationsliftscalabiyauvarieties}, the classical deformation functor $(\Def_X^{\sch})^{\cl}$ associated to the formal moduli problem of scheme deformations of $X$ is isomorphic to $\overl[\Def]_X$, hence the result by \Cref{deformations are the same}. 
\end{proof}

\subsubsection{Equivariant deformations through deformations over the classifying stack} \label{section equivariant deformations}

Let $K = \Spec(k)$ and $W = \Spec(\Lambda)$, a local Noetherian ring with residue field $k$. Let $G$ be a \comdef{smooth} group scheme over $W$ such that $G_K$ acts on a derived scheme $X$ over $K$. One can wonder if $G$ acts on deformations of $X$, for instance in the seminal work of \cite{rim1980equivariant}. 

Observing that the action of $G$ on $X$ induces an action of $G_K$ on the tangent bundle $T_X$ and on the derived global vector fields $R\Gamma(X, T_X)$, one can wonder if this action is, in some sense, compatible with the Lie algebra structure on $R\Gamma(X, T_X)$. Up to the equivalence between Lie algebras and formal moduli problems, the question in fact becomes: does the $G_K$-action on $X$ induces, in some sense, an action on the formal moduli stack $\Def_X^{K/W}$ itself.

We give a precise answer to this question. The object $\Def_X^{K/W}$ in $\FMP_{K/W}$ admits a lift to the category $\FMP_{BG_K/BG}$ of formal moduli problems over $BG$ and under $BG_K$ (the classifying stacks in the etale topology). The lift is the functor of deformations of $X/G$ relative to $BG$, $\Def_{X/G}^{BG_K/BG}$.

Now, in order to study the lifts of the $G_K$-action to deformations of $X$, one can define the equivariant deformation functor which sends an Artinian extension $Z$ of $K$ to the space of $G_Z$-actions on schemes over $Z$ which reduce to the action of $G_K$ on $X$. Denote it by $\Def_{G_K \acts X}^{\sch}$.
We show that the equivariant deformation functor $\Def^{\sch}_{G_K \acts X}$ appears simply as the ``$G$-fixed points'' of the deformation functor $\Def_X^{K/W}$, equipped with its $G$-equivariant structure: 
$$\Def^{\sch}_{G_K \acts X} \simeq (\Def_{X/G_K}^{BG_K/BG})^G,$$
where the (derived) $G$-fixed points should be understood as the pushforward from $\FMP_{BG_K/BG}$ to $\FMP_{K/W}$.
Hence, the tangent complex of the equivariant deformation functor is the vector space of $G$-equivariant global vector fields: 
$$T_{\Def^{\sch}_{G_K \acts X}} \simeq R\Gamma(X, T_X)^{G_K}[1]$$
where the (derived) $G$-fixed points should be understood as the pushforward from $\QC{BG_K}$ to the point.

The main tool that we use is the correspondence between $G_K$-action on $X$ and quotient maps $X/G_K \to BG_K$, where $X/G_K$ and $BG_K$ denote the quotients in the etale topology (different from the formal quotient $B^{\deff}$ of \cite{gaitsgory2019study} and \Cref{section split deformations}!). It allows us to show that the equivariant deformation functor can be reformulated as the functor of deformations of the map $X/G_K \to BG_K$ leaving the target fixed. We can then simply use the results obtained in \Cref{section split deformations} and \Cref{section mixed deformations}. 

Notice that this reasoning could be extended to flat group schemes locally of finite presentation by using the fppf topology instead of the etale topology.

\begin{prop}[\cite{nikolaus2015principal}] \label{equivalence NSS}
    Let $G$ be a group object in an $\infty$-topos $\mathcal{C}$.
    Let $$B_G \colon \Act_{\cat}(G) \to \mathcal{C}_{/BG}$$ be the quotient functor and $$\Omega_G \colon \mathcal{C}_{/BG} \to \Act_{\cat}(G)$$ be the functor which sends naturally a map $\Y \to BG$ to the groupoid $(\Y \times_{BG} *) \times_{\Y} (\Y \times_{BG} *)$ or, equivalently, to the groupoid $\Y \times_{BG} G_{\bullet}$.
    Those two functors form an equivalence of categories 
    $$Act_{\cat}(G) \simeq \mathcal{C}_{/BG}.$$

    Let $X \in \cat$. The above equivalence induces an equivalence between fibers $$\Act_{\cat}(G, X) \simeq \mathcal{C}_{/BG} \times_{\mathcal{C}} \{X\}.$$
\end{prop}

\begin{notation}
    Write $\Stk$ for the $\infty$-topos of derived stacks in the etale topology.
\end{notation}

\begin{defi}
    Let $f \colon K \to Z$ be a map in $\Stk$.
    \begin{enumerate}
    \item Denote by $f^* \colon \Grpd_{\Stk_Z} \to \Grpd_{\Stk_K}$ the base change functor defined simplex by simplex.

    \item Let $\X \to Z$ be a map of prestacks.
    The above base change functor restricts to a functor $f^* \colon \Grpd_{\Stk_Z}(\X) \to \Grpd_{\Stk_Z}(\X \times_Z K)$.

    \item Let $G_Z \in \Grp_{\Stk_Z}$. Write the fiber $G_K = G_Z \times_Z K$.
    Base change also induces a functor $$f^* \colon Act_{\Stk_Z}(G_Z) \to Act_{\Stk_K}(G_K).$$
    \end{enumerate}
\end{defi}

\begin{lemma} \label{equivalence NSS pullback}
    Let $f \colon K \to Z$ be a map in $\Stk$. Let $G_Z \in \Grp_{\Stk_Z}$. The equivalence of \Cref{equivalence NSS} is compatible with base change, that is, the following diagram commutes:
\[\begin{tikzcd}
	{\Act_{\Stk_Z}(G_Z)} & {\Stk_{/BG_Z}} \\
	{\Act_{\Stk_K}(G_K)} & {\Stk_{/BG_K}}
	\arrow["B", shift left, from=1-1, to=1-2]
	\arrow["{f^*}"', from=1-1, to=2-1]
	\arrow["\Omega", shift left, from=1-2, to=1-1]
	\arrow["{f^*}", from=1-2, to=2-2]
	\arrow["B", shift left, from=2-1, to=2-2]
	\arrow["\Omega", shift left, from=2-2, to=2-1]
\end{tikzcd}\]
Hence, if $X$ is a derived $K$-scheme endowed with a $G_K$-action encoded in a groupoid $U$, it induces an equivalence of categories on fibers which is natural in $Z$ $$B \colon\begin{tikzcd}
	\Act_{\Stk_Z}(G_Z) \times_{\Act_{\Stk_K}(G_K)} \{U\}
	&
	\Stk_{/BG_Z} \times_{\Stk_{/BG_K}} \{X/G_K \sto BG_K\}
	\arrow[shift left, from=1-1, to=1-2]
	\arrow[shift left, from=1-2, to=1-1]
\end{tikzcd} \colon \Omega$$
\end{lemma}

\begin{proof}
    Colimits are universal in a topos, so the quotient functor $B$ commutes with base change.
\end{proof}

\begin{defi} \label{forget group action}
    Let $K \to W$ be a map of coherent affine schemes. Let $G$ be a derived group scheme over $W$. For all $Z \in \Schaffart_{K/W}$, write the $Z$-fiber $G_Z = G \times_W Z$.
    Let $X$ be a derived scheme over $K$ endowed with a $G_K$-action, encoded in a map $U \sto G_K$ in $\Grpd_{\Sch_{/K}}$.

    Let $\Def_{G_K \acts X}^{\sch} \colon \Schaffartop{K/W} \to \widehat{\Cat}_{\infty}$ be the functor which sends $Z \in \Schaffart_{K/W}$ to 
    
    $$\Def_{G_K \acts X}^{\sch}(Z) = \Act_{\sch_Z}(G_Z) \times_{\Act_{\sch_K}(G_K)} \{U\}.$$
    Write $U_G \colon \Def_{G_K \acts X}^{\sch} \to \Def_X^{\sch}$ for the natural transformation which forgets the $G$-action.
\end{defi}

\begin{rk}
    We will see in \Cref{equivariant deformations as deformations over BG} that the functor $\Def_{G_K \acts X}^{\sch}$ lands in fact in spaces.
\end{rk}

We want to interpret the data of a $G$-action as a map to the etale classifying stack $BG$.

\begin{prop}[{\cite[3.3]{toen2014derivedalgebraicgeometry}}] \label{BG Artin stack}
    Let $G_K$ be a smooth derived group scheme which acts on a derived group scheme $X$. Then $X/G_K$ 
    is an Artin stack. In particular, $BG_K$ is an Artin stack.
\end{prop}

\begin{prop}[{\cite[Chapter 1, 7.4.3]{gaitsgory2019study}}] \label{Artin stacks have deformation theory}
    Artin stacks have deformation theory.
\end{prop}

\begin{coro} \label{equivalence NSS laftdef}
    Let $Z$ be an affine scheme. 
    Let $G_Z$ be a smooth derived group scheme over $Z$, and $\X$ be a \comdef{laft} derived $Z$-scheme. The functors $\Omega$ and $B$ from \Cref{equivalence NSS} induce an equivalence of categories $$B \colon\begin{tikzcd}
	\Act_{Z}(G_Z, \X)
	&
	\prstklaftdef_{/BG_Z} \times_{\prstklaftdef_{/K}} \{\X\}
	\arrow[shift left, from=1-1, to=1-2]
	\arrow[shift left, from=1-2, to=1-1]
\end{tikzcd} \colon \Omega$$
\end{coro}

\begin{proof}
    By \Cref{BG Artin stack} and \Cref{Artin stacks have deformation theory}, and the assumption that  $X$ is laft (\cite[Chapter 2, 4.5.2]{gaitsgory2019studyI}) the quotient functor $B$ lands in fact in $\prstklaftdef_{/BG_Z}$.
    Reciprocally, the functor $\Omega$ sends a prestack map $\Y \to BG_Z$ to the groupoid $\Y \times_{BG_Z} (G_Z)_{\bullet}$, which encodes an action of $G$ on $\X$, in the category of derived schemes.
    Then, \Cref{equivalence NSS} gives that those functors form an equivalence.
\end{proof}

\begin{notation}
    Let $G_K$ be a \comdef{smooth} derived group scheme over $K$. Denote $f \colon BG_K \to K$. Let $M \in \QC{_K}$. We write $$M^{G_K} := f_* M,$$ called the (derived) $G_K$-fixed points of $M$.    
\end{notation}

All Artin stacks, and in particular the quotients by $G$, are encompassed in the setting of prestacks with deformation theory that we used to study deformation theory in the previous sections.
Thus, $G$-equivariant deformations of $X$ can be reformulated as deformations of the classifying map $X/G_K \to BG_K$ among prestacks with deformation theory, allowing us to apply all the machinery of previous sections.

\begin{prop} \label{equivariant deformations as deformations over BG}
    Let $G$ be a smooth group scheme over $W$. For all $Z \in \Schaffart_{K/W}$, write the $Z$-fiber $G_Z = G \times_W Z$.
    Let $X$ be a \comdef{proper} derived scheme locally of finite presentation over $K$ endowed with a $G_K$-action, encoded in a map $U \sto G_K$ in $\Grpd_{\Sch_{/K}}$.  
    Then the inverse maps $\Omega$ and $B$ from \Cref{equivalence NSS laftdef} and \Cref{equivalence NSS pullback} induce an equivalence of functors

    $$\Def_{G_K \acts X}^{\sch} \simeq \Def_{X/G_K \to BG_K}^{K/W}.$$ 
    In particular, $\Def_{G_K \acts X}^{\sch} \in \FMP_{K/W}.$

    Moreover, under this equivalence, the forgetful transformation $U_G \colon \Def_{G_K \acts X}^{\Sch} \to \Def_X^{\sch}$ identifies with the transformation $\fibdef \colon \Def_{X/G_K \to BG_K}^{K/W} \to \Def_{X}^{K/W}$.

    Finally, the action of $G_K$ on $X$ induces an action of $G_K$ on $T_X$ and on $R\Gamma(X, T_X)$, and there is an equivalence in $\QC{K}$ $$T_{\Def_{G_K \acts X}^{\sch}} \simeq R\Gamma(X, T_X)^{G_K}[1].$$

\end{prop}

\begin{rk}
    All the content of this section could be extended to groups which are only flat and locally of finite presentation instead of smooth by considering the fppf topology instead of the etale topology, using \cite{toen2010flatdescentartinnstacks}.
\end{rk}

\begin{proof}
    Let $Z \in \Schaffart_{K/W}$. A map $\Y \to BG_Z$ in $\prstklaftdef_{/BG_Z} \times _{\prstklaftdef_{/BG_K}} \{X/G_K \to BG_K\}$ sits in the following diagram where, by pasting law, all squares are cartesian, and where we denote $\X = \Y \times_{BG_Z} Z$:

\[\begin{tikzcd}
	& X & K \\
	\X & Z & {X/G_K} & {BG_K} \\
	& \Y & {BG_Z}
	\arrow[from=1-2, to=1-3]
	\arrow[from=1-2, to=2-1]
	\arrow[dashed, from=1-2, to=2-3]
	\arrow[from=1-3, to=2-2]
	\arrow[from=1-3, to=2-4]
	\arrow[from=2-1, to=2-2]
	\arrow[from=2-1, to=3-2]
	\arrow[from=2-2, to=3-3]
	\arrow[from=2-3, to=2-4]
	\arrow[dashed, from=2-3, to=3-2]
	\arrow[from=2-4, to=3-3]
	\arrow[from=3-2, to=3-3]
\end{tikzcd}\]
    Remark in particular that $\X$, a priori an object of $\prstklaftdef_{/Z}$, sits in a pullback square

\[\begin{tikzcd}
	X & K \\
	\X & Z
	\arrow[from=1-1, to=1-2]
	\arrow[from=1-1, to=2-1]
	\arrow[from=1-2, to=2-2]
	\arrow[from=2-1, to=2-2]
\end{tikzcd}\]
    where $X$ is a derived scheme. Hence, by \Cref{deformations are the same}, $\X$ is in fact equivalent to a derived scheme over $Z$. Thus, we can use \Cref{equivalence NSS laftdef} and \Cref{equivalence NSS pullback} to get two natural maps $\Omega_Z$ and $B_Z$ which form an equivalence $$\Act_Z(G_Z) \times_{\Act_K(G_K)} \{U\} \simeq \prstklaftdef_{/BG_Z} \times _{\prstklaftdef_{/BG_K}} \{X/G_K \to BG_K\}.$$

    This equivalence is natural in $Z$ so it gives a natural equivalence of functors $$\Def_{G_K \acts X}^{\sch} \simeq \Def_{X/G_K \to BG_K}^{K/W}.$$

    It is clear from the definitions that the forgetful transformation $U_G \colon \Def_{G_K \acts X}^{\sch} \to \Def_X^{\sch}$ identifies with the transformation $\fibdef \colon \Def_{X/G_K \to BG_K}^{K/W} \to \Def_{X}^{K/W}$.
    The computation of the tangent complex comes from this identification and from the computation of the tangent map associated to the transformation $\fibdef$ in \Cref{fiber of deformation tangent mixed}.

\end{proof}

\begin{rk}
    Given that $X$ is locally of finite presentation, the tangent complex $T_X$ is perfect. Hence $R\Gamma(X, T_X)^{G_K}$ coincides with the same usual objects in quasi-coherent sheaves instead of pro-coherent sheaves.
\end{rk}

To conclude this section, we verify that the classicalization of the equivariant deformation functor defined above coincides with the classical equivariant deformation functor.

\begin{defi}
    Let $G$ be a classical group scheme over $W$. Let $X$ be a classical scheme over $K$ endowed with an action of $G_K$.
    Let $\overl[\Def]_{G_K \acts X} \colon \CRing^{\art}_{K/W} \to \Set$ be the classical functor which sends an Artinian extension $A$ of $K$ to the isomorphism class of actions of $G_A$ on flat schemes over $A$ which lift the action of $G_K$ on $X$.
\end{defi}

\begin{prop}
    Let $G$ be a classical group scheme over $W$. Let $X$ be a \comdef{proper} classical scheme locally of finite presentation over $K$ endowed with an action of $G_K$. Then the classical deformation functor $(\Def_{G_K \acts X}^{\sch})^{cl}$ associated to the formal moduli problem $\Def_{G_K \acts X}^{\sch}$ is equivalent to the classical equivariant deformation functor $\overl[\Def]_{G_K \acts X}$.

\end{prop}

\begin{proof}
    This is a corollary of \Cref{classicalisation of deformations}. 
\end{proof}

 \subsection{Unobstructed maps of formal moduli problems} \label{section unobstructed maps}

Let $K = \spec(k)$ be a field, and $W = \Spec(\Lambda)$ a local Noetherian ring with residue field $k$ (the setting of mixed characteristic deformations).
A variety is said to have unobstructed deformations if all deformations can be extended to any Artinian extension of $K$. By extension, any classical deformation functor is said to be unobstructed if it has the analogous property. A (derived) formal moduli problem is said to be unobstructed if its associated classical deformation functor is unobstructed. Unobstructedness can be interpreted as a (classical) formal smoothness conditions at the only $k$-point of the formal moduli problem.

While the definition of unobstructedness remains  purely classical, the use of (derived) formal moduli problems allows us to get a fully functorial tangent complex which contains in particular the classical tangent space and obstruction space. 

There is a natural relative version of unobstructedness, for maps of formal moduli problems. Unobstructed maps of formal moduli problems are very useful to study obstructions because they ``transfer'' obstructions from one side to another. In particular, the source is unobstructed if and only if the target is unobstructed. Hence, a general method to study obstructions of a given formal moduli problem is to link it through unobstructed maps to a formal moduli problem whose obstructions are easier to study.

It is well known that for the case of deformations of a variety, obstructions to deformations lie in $H^2(X, T_X)$. In particular, if $H^2(X, T_X) = 0$, deformations are unobstructed. This principle generalizes in the relative setting, for a map of formal moduli problems: if the ``$H^2$'' (co)homology group of the fiber of the associated map of tangent complexes vanishes, then the map is unobstructed. This is a reformulation in the setting of (derived) formal moduli problem of a result for classical deformation functors, proven in \cite{fantechi1998obstruction} using the theory of obstruction spaces. Here we get it easily in the setting of formal moduli problems which gives full functoriality for obstruction spaces. This condition is very tractable and can be applied in numerous situations.

Note that this notion of unobstructedness is unsatisfying in the derived setting of formal moduli problems because it is a purely classical notion. One could expect a stronger derived notion of unobstructedness, which would make use the whole tangent complex of a formal moduli problem. While the natural condition of derived formal smoothness is too strong to be interesting for the study of deformations, some derived conditions exist: for instance, asking that the induced map of tangent complexes is injective on all homotopy groups, or surjective on all homotopy groups. Such conditions allow to transfer pro-representability of formal moduli problems, as explained in \cite{brantner2025deformationsliftscalabiyauvarieties}. However, for instance in \Cref{section homogeneous varieties}, it can be too hard to prove such derived conditions in practice. 

Let us first remind what is an (un)obstructed formal moduli problem.

\begin{defi}
    A formal moduli problem  $\M \in \FMP_{K/W}$ is said to be unobstructed if for all surjective maps of classical rings $A \to B$ in $\CRing^{\art}_{K/W}$, $\pi_0 \M(A) \to \pi_0 \M(B)$ is surjective. 
    In other words, we say that $\M$ is unobstructed if the associated classical deformation functor $\M^{\cl}$ is unobstructed.

    We say that the formal moduli problem $\M$ has an obstruction if there exists a classical nilpotent extension $A \to B$ such that $\pi_0 \M(A) \to \pi_0 \M(B)$ is not surjective.
\end{defi}

Now this extends to the relative setting, analogously to the usual formal smoothness conditions. Intuitively, it can be thought as asking that all deformations on the right hand side come from deformations on the left hand side.

 \begin{defi} \label{defi unobstructed map}
    Let $\M, \Nr \in \FMP_{K/W}$. 
    A map of formal moduli problems $F \colon \M \to \Nr$ is said to be unobstructed if for all surjective maps of classical rings $A \to B$ in $\CRing^{\art}_{K/W}$, the map $$\pi_0 \M(A) \to \pi_0 \M(B) \times_{\pi_0 \Nr(B)} \pi_0 \Nr(A)$$ is surjective. 
\end{defi}

\begin{rk} \label{surj}
    In particular, if $F \colon \M \to \Nr$ is unobstructed, the map $\M(A) \to \Nr(A)$ is surjective for all $A \in \CRing^{\art}_{K/W}$.

\end{rk}

\begin{rk}
    One can naturally wonder why we do not consider a \emph{derived} smoothness notion, where we would consider all derived Artinian extensions instead of only classical Artinian extensions. The reason is that this notion would be strong for our use: for instance $k[[t]]$ is not derived smooth over $k$, whereas it should be the ``best'' example of unobstructed deformations. 
\end{rk}

The definition of unobstructedness gives immediately that unobstructed map transfers obstructions from one side to another.

\begin{prop} \label{smoothness transfer}
    Let $F \colon \M \to \Nr$ be an unobstructed map in $\FMP_{K/W}$. 

    Let $A \to B$ be a surjective map of classical rings in $\CRing^{\art}_{K/W}$. Let $x \in \pi_0\M(A)$.
    Then $x$ admits a lift to $\pi_0 \M(B)$ if and only if $Fx$ admits a lift to $\pi_0 \Nr(B)$.
    
    Thus, $\M$ is unobstructed if and only if $\Nr$ is unobstructed.
\end{prop}

Unobstructedness is implied by a simple criterion, proved by \cite{fantechi1998obstruction} for deformation functors, that we recall here in a modern setting.

\begin{lemma}[{\cite[3.3, 3.4]{brantner2025deformationsliftscalabiyauvarieties}, \cite[6.2.6]{lurie2004derived}}] \label{composition sqz}
    For every map $A \to B$ in $\calg_{\Lambda}^{\an, \art}$ such that $\pi_0 A \to \pi_0 B$ surjective,  there exists a factorization $$A = A_n \to ... \to A_1 \to A_0 = B$$
    where each of the maps sits in a pullback square 
\[\begin{tikzcd}
	{A_{i+1}} & {A_i} \\
	k & {\sqz_k(V_i[n_i+1])}
	\arrow[from=1-1, to=1-2]
	\arrow[from=1-1, to=2-1]
	\arrow["\lrcorner"{anchor=center, pos=0.125}, draw=none, from=1-1, to=2-2]
	\arrow[from=1-2, to=2-2]
	\arrow[from=2-1, to=2-2]
\end{tikzcd}\]
    with $n_i \geq 0$. 
    If $A$ and $B$ are discrete, the factorization can be chosen such that 
    
\[\begin{tikzcd}
	{A_{i+1}} & {A_i} \\
	k & {\sqz_k(k[1])}
	\arrow[from=1-1, to=1-2]
	\arrow[from=1-1, to=2-1]
	\arrow["\lrcorner"{anchor=center, pos=0.125}, draw=none, from=1-1, to=2-2]
	\arrow[from=1-2, to=2-2]
	\arrow[from=2-1, to=2-2]
\end{tikzcd}\]

\end{lemma}

\begin{lemma} \label{tangent complex homotopy groups}
    Let $\M \in \FMP_{K/W}$. The homotopy groups of the tangent complex of $\M$ can be computed for all $m,n \in \N$ as
    $$\pi_{m-n}(T_{\M}) \simeq \pi_{m}\M(\sqz_k(k[n])).$$
\end{lemma}

\begin{proof}
    It is a consequence of the adjunction between square zero extension and tangent functor of \Cref{properties tangent complex}: $$\Map_{\FMP_{K/W}}(\sqz_k(k[n]), \M) \simeq \Map_{\Mod_k}(k[n], T_{\M}).$$
\end{proof}

\begin{prop} \label{criterion smooth}
    Let $F \colon \M \to \Nr$ be a map in $\FMP_{K/W}$. Let $Q$ be the fiber of the associated map of tangent complexes $T_\M \to T_\Nr$, in $\Mod_k$.
    If $\pi_{-1} Q = 0$, then the map $F$ is unobstructed.
\end{prop}

\begin{rk}
    In the above proposition, the condition $\pi_{-1} Q = 0$ can be reformulated as the following conditions:
    \begin{enumerate}
    \item $\pi_0 T_\M \to \pi_{0} T_\Nr$ is surjective.
    \item $\pi_{-1} T_\M \to \pi_{-1} T_\Nr$ is injective.
    \end{enumerate}
\end{rk}

\begin{proof}
Let $A \to B$ be a surjective map of classical rings in $\CRing^{\art}_{K/W}$. By \Cref{composition sqz}, we can reduce to the case where the map $A \to B$ sits in a pullback square
\[\begin{tikzcd}
	A & B \\
	k & {\sqz_k({k[1]})}
	\arrow[from=1-1, to=1-2]
	\arrow[from=1-1, to=2-1]
	\arrow["\lrcorner"{anchor=center, pos=0.125}, draw=none, from=1-1, to=2-2]
	\arrow[from=1-2, to=2-2]
	\arrow[from=2-1, to=2-2]
\end{tikzcd}\]

This pullback square is sent by $\M$ to the fiber sequence $\M(A) \to \M(B) \xrightarrow{\mathrm{ob}} \M(\sqz_k{k[1]})$, and the same for $\Nr$, hence the diagram 

\[\begin{tikzcd}
	{\pi_1\M(\sqz_k(k[1]))} & {\pi_0\M(A)} & {\pi_0\M(B)} & {\pi_0\M(\sqz_k(k[1]))} \\
	{\pi_1\Nr(\sqz_k(k[1]))} & {\pi_0\Nr(A)} & {\pi_0\Nr(B)} & {\pi_0\Nr(\sqz_k(k[1]))}
	\arrow["\delta", from=1-1, to=1-2]
	\arrow["F_1"', two heads, from=1-1, to=2-1]
	\arrow["f", from=1-2, to=1-3]
	\arrow["F_2"', from=1-2, to=2-2]
	\arrow["\ob", from=1-3, to=1-4]
	\arrow["F_3"', from=1-3, to=2-3]
	\arrow["F_4", hook, from=1-4, to=2-4]
	\arrow["\delta"', from=2-1, to=2-2]
	\arrow["f"', from=2-2, to=2-3]
	\arrow["\ob"', from=2-3, to=2-4]
\end{tikzcd}\] 
where the first vertical map is surjective and the last vertical map is injective, given the assumptions and \Cref{tangent complex homotopy groups}. Let us prove that the map $\pi_0 \M(A) \to \pi_0\Nr(A) \times_{\pi_0\Nr(B)} \pi_0\M(B)$ is surjective.
Let $x \in \pi_0\M(B)$ and $y' \in \pi_0\Nr(A)$ such that $F_3x = fy'$. Denote $$y := F_3x = fy' \in \pi_0 \Nr(B).$$ $y = fy'$ so $\ob(y) = 0$. But $\ob(y) = \ob(F_3x) = F_4(\ob(x))$ so $F_4(\ob(x)) = 0$. By injectivity of $F_4$, $\ob(x) = 0$. Let us make a choice of a representative in $\Nr(B)$ of the point $y$, and a choice of a representative in $\M(B)$ of the point $x$. With those choices, the bottom line and the top line are exact. $fy' = y$ so by exactness, there exists $z' \in \pi_1\Nr(\sqz_k(k[1]))$ such that $\delta z' = y'$. But $F_1$ is surjective so there exists $z \in \pi_1\M(\sqz_k(k[1]))$ such that $F_1 z = z'$. Denote $x' = \delta z \in \pi_0\M(A)$. Then $fx' = x$, and $F_2 x' = F_2 \delta z = \delta F_1 z = \delta z' = y'$.

\end{proof}

\begin{rk}
    In particular, if $\M$ is a formal moduli problem and $\pi_{-1} T_{\M} = 0$, then $\M$ is unobstructed. 
    $\pi_0 T_{\M}$ and $\pi_{-1} T_{\M}$ correspond respectively to what is classically called the ``tangent space'' and the ``obstruction space''.
    For instance, if $\M = \Def^{\sch}_X$ is the formal moduli problem of deformations of a scheme $X$, its tangent complex is $T_{\M} = R\Gamma(X, T_X)[1]$. The condition $\pi_{-1}(T_{\M}) = 0$ corresponds to the cohomological condition $H^2(X, T_X) = 0$. 
    
    Notice that those conditions are sufficient but not necessary in general to get unobstructedness.
\end{rk}

\subsection{Conditions for unobstructed equivariant deformations} \label{section unobstructed equivariant deformations}

Let $K = \spec(k)$ be a field, and $W = \Spec(\Lambda)$ be a local Noetherian ring with residue field $k$. Let $G$ be a group scheme over $W$, let $X$ be a scheme endowed with an action of $G_K$. As we explained in \Cref{section equivariant deformations}, it can be interesting to determine when the $G_K$-action deforms together with deformations of the scheme $X$. There is a genealogy of results around this question. In \cite{rim1980equivariant} and \cite{schroer2018equivariant}, Rim and Schröer-Takayama obtain conditions to get unobstructed equivariant deformations under the form of vanishing of some cohomology groups. However, they do not treat the mixed characteristic situation for group schemes, that is, they don't take into account mixed characteristic lifts of $G_K$.

As we also explained in \Cref{section equivariant deformations}, the situation is analogous to the question of deformations of $G$-equivariant sheaves. In \cite{zdanowicz2018arithmetically}, in the case of equivariant deformations of sheaves, Zdanowicz reformulates the problem in terms of unobstructedness of the map from equivariant deformations of sheaves to usual deformations of sheaves. This way, he gets a criterion in terms of vanishing of $G_K$-cohomology of certain $G_K$-modules in \cite[Theorem 3.3]{zdanowicz2018arithmetically}. 

We use an analogous method in the non-linear case of equivariant deformations of schemes to show that the $G$-action deforms together with $X$ when $H^2(G_K, H^0(X, T_X)) = 0$, $H^1(X, T_X)$ is a trivial $G_K$-module and $H^1(G_K, k) = 0$.

More precisely, the problem of deforming a $G$-action together with deformations of $X$ can be reformulated as the unobstructedness of the forgetful map $\Def_{G_K \acts X} \to \Def_X$. In \Cref{section equivariant deformations}, we observed that this map corresponds to the map $\fibdef^{K/W}_{X/G_K}$. Then, we simply reformulate the unobstructedness conditions from \Cref{criterion smooth} in the specific case of this kind of maps.

\begin{lemma} \label{lemma zdanowicz}
    Let $f \colon K' \to K$ be a map of prestacks with a section $e \colon K \to K'$, and let $M \in \QC{K'}$ be coconnective. Write $L := f_*$ and $U := e^*$, two functors from $\QC{K'}$ to $\QC{K}$. $f_*$ is left $t$-exact. Suppose that $U$ is $t$-exact and faithful. 

    The spectral sequence $E_2^{p,q} = \pi_p L(\pi_q M) \Longrightarrow \pi_{p+q}L(M)$ is a first-quadrant spectral sequence that gives edge maps $\pi_n(LM) \to \pi_0L(\pi_nM)$ for all $n \in \Z^-$, and the same for $U$ (it is trivial). The functoriality of the spectral sequence gives a commutative diagram 
\[\begin{tikzcd}
	{\pi_n(LM)} & {\pi_n(UM)} \\
	{\pi_0L (\pi_n M)} & {\pi_0U (\pi_n M)}
	\arrow[from=1-1, to=1-2]
	\arrow[from=1-1, to=2-1]
	\arrow["\simeq", from=1-2, to=2-2]
	\arrow[from=2-1, to=2-2]
\end{tikzcd}\]
    Then the following assertions hold:
    \begin{enumerate}
    \item The edge map $\pi_n(UM) \to \pi_0U(\pi_nM)$ is an isomorphism.

    \item The map $\pi_0L(\pi_nM) \to \pi_0U(\pi_nM)$ is injective.

    \item The map $\pi_n(LM) \to \pi_n(UM)$ is injective if and only if the map $\pi_n(LM) \to \pi_0L(\pi_nM)$ is injective. 
    
    \item The map $\pi_n(LM) \to \pi_n(UM)$ is surjective if and only if the map $\pi_n(LM) \to \pi_0L(\pi_nM)$ is surjective and the map $\pi_0L(\pi_nM) \to \pi_0U(\pi_nM)$ is an isomorphism.

    \item The following assertions are equivalent: 
    \begin{enumerate}[label=(\roman*)]
    \item the map $\pi_{-1}(LM) \to \pi_{-1}(UM)$ is surjective and the map $\pi_{-2}(LM) \to \pi_{-2}(UM)$ is injective.
    \item $\pi_{-2}L(\pi_0M) = 0$, $\pi_{-1}L(\pi_{-1} M) = 0$ and the map $\pi_0L(\pi_{-1}M) \to \pi_0U(\pi_{-1}M)$ is an isomorphism.
    \end{enumerate}

    \end{enumerate}
\end{lemma}

\begin{proof}
\begin{enumerate}
\item This comes from $t$-exactness of $U$. 
\item Denote $U^R := e_*$ the right adjoint to $U$. $U$ is faithful so the unit map $$\pi_nM \to \pi_0U^R \circ \pi_0U (\pi_nM)$$ is injective, and $L^{\heartsuit}$ is left exact so it preserves injectivity, so the map $$\pi_0 L (\pi_n M) \to \pi_0 L \circ \pi_0U^R \circ \pi_0U (\pi_n M)$$ is injective. But $\pi_0L \circ \pi_0U^R \simeq \id$, so the map $\pi_0L(\pi_nM) \to \pi_0U(\pi_nM)$ is injective.
\item In view of the commutative diagram in the lemma and the isomorphism of (1), the map $\pi_n(LM) \to \pi_n(UM)$ is injective if and only if the map $\pi_n(LM) \to \pi_0L(\pi_nM)$ is injective.
\item In view of the commutative diagram in the lemma and the isomorphism of (1), the map $\pi_n(LM) \to \pi_n(UM)$ is surjective if and only if the map $\pi_n(LM) \to \pi_0L(\pi_nM)$ is surjective and the map $\pi_0L(\pi_nM) \to \pi_0U(\pi_nM)$ is surjective, but it is also injective so it is an isomorphism.
\item Consider the long exact sequence 
\begin{multline*}
    0 \sto \pi_{-1}L(\pi_0M) \sto \pi_{-1}(LM) \sto \pi_0L (\pi_{-1} M) \\ \sto \pi_{-2} L (\pi_0 M) \sto \ker \Big( \pi_{-2}(LM) \rightarrow \pi_0 L (\pi_{-2} M) \Big) \sto \pi_{-1} L (\pi_{-1} M)
\end{multline*}
The conditions of (3) and (4) are verified if and only if $\pi_{-2}L(\pi_0M) = 0$, $\pi_{-1}L(\pi_{-1} M) = 0$ and $\pi_0L(\pi_{-1}M) \iso \pi_0U(\pi_{-1}M)$.
\end{enumerate}
\end{proof}

Let $K = \spec(k)$ be a field, and $W = \Spec(\Lambda)$, a local Noetherian ring with residue field $k$. 
We get a first general result for unobstructedness of maps of the type $F_{\M} \colon f_* \M \to e^* \M$, when the section map $e \colon W \to W'$ is faithfully flat (which is a very strong condition for a section, but corresponds to the case of a quotient map).

\begin{coro} \label{general unobstructedness}
    Let $g \colon W' \to W$ be a laft map of prestacks, with a section $e \colon W \to W'$. Write $K' = K \times_W W'$ and let $\M \in \modstk_{K'/W'}$. Suppose that the pullback $e^*$ is faithful and $t$-exact, and that $T_{\M}[-1]$ is coconnective.
    If the following conditions are verified:
    \begin{enumerate}
        \item $\pi_{-2} g_*(\pi_0 T_{\M}[-1]) = 0$
        \item $\pi_0g_*(\pi_{-1} T_{\M}[-1]) \simeq e^* (\pi_{-1} T_{\M}[-1])$ 
        \item $\pi_{-1}g_*(\pi_{-1}T_{\M}[-1]) = 0,$
    \end{enumerate} 
    then, with \Cref{transformation map from adjunction}, the natural map $F_{\M} \colon g_* \M \to e^* \M$ is unobstructed.

\end{coro}

\begin{proof}
    By \Cref{lemma zdanowicz} applied to $T_{\M}[-1]$, the map $\pi_0(f_* T_{\M}) \to \pi_0(e^* T_{\M})$ is surjective and the map $\pi_{-1}(f_* T_{\M}) \to \pi_{-1}(e^* T_{\M})$ is injective. But the map $f_* T_{\M} \to e^* T_{\M}$ identifies with the tangent map of $F_{\M} \colon f_* \M \to e^* \M$ by \Cref{tangent commutes with adjunction}. Thus, by \Cref{criterion smooth}, the map $F_{\M}$ is unobstructed.
\end{proof}

Using the results of \Cref{section split deformations} and \Cref{section mixed deformations}, \Cref{general unobstructedness} applies in particular to the map $\fibdef_{Y}^{K/W} \colon \Def_{Y \to K'}^{K/W} \to \Def_X^{K/K}$.

\begin{coro} \label{unobstructedness of fiber of deformations}
Let $g \colon W' \to W$ be a laft prestack over $W$, with a section $e \colon W \to W'$. Write $K' = K \times_W W'$, let $Y \in \prstklaftdef_{/K'}$ with the structure map denoted $f \colon Y \to K'$, and write $X = Y \times_{K'} K$.
    Suppose that the pullback $e^*$ is faithful and $t$-exact, and that $T_{Y \sto K'}$ is coconnective.
    If the following conditions are verified:
    \begin{enumerate}
        \item $\pi_{-2}g_*(\pi_0 f_*T_{Y \sto K'}) = 0$
        \item $\pi_0g_*(\pi_{-1} f_*T_{Y \sto K'}) \iso e^*(\pi_{-1} f_*T_{Y \sto K'})$
        \item $\pi_{-1}g_*(\pi_{-1} f_*T_{Y \sto K'}) = 0$,
    \end{enumerate}
    then the map $\fibdef_{Y}^{K/W} \colon \Def_{Y \sto K'}^{K/W} \to \Def_X^{K/K}$ is unobstructed.
\end{coro}

\begin{proof}
    We apply \Cref{general unobstructedness} to $\M = \Def_Y^{K'/W'}$, using the computation of the tangent of \Cref{fiber of deformation tangent mixed}
\end{proof}

Finally, using the correspondence between equivariant deformations and deformations over $BG$ from \Cref{equivariant deformations as deformations over BG}, \Cref{unobstructedness of fiber of deformations} gives criteria for unobstructedness of the forgetful map $\Def_{G_K \acts X}^{\sch} \to \Def_X^{K/W}$, that is, to be able to lift the $G_K$-action to deformations of $X$. We shift to cohomological conventions to stay close to analogous classical results.

\begin{prop} \label{criterion}
    Let $G$ be a \comdef{smooth} group scheme over $W$, let $X$ be a smooth \comdef{proper} scheme over $K$ endowed with an action of the fiber $G_K := G \times_W K$. It induces an action of $G_K$ on $T_X$ and on $R\Gamma(X, T_X)$. If the following conditions are verified:
    \begin{enumerate}
    \item $H^1(G_K, k) = 0$

    \item $H^2(G_K, H^0(X, T_X)) = 0$ 
    \item $H^1(X, T_X)$ is a trivial $G_K$-module
    \end{enumerate}
    then the forgetful map $\Def_{G_K \acts X}^{\sch} \to \Def_X^{\sch}$ (from \Cref{forget group action}) is unobstructed.
\end{prop}

\begin{proof}
    Write $f \colon X/G_K \to BG_K$.
    By \Cref{equivariant deformations as deformations over BG}, up to the equivalence between $\Def_{G_K \acts X}^{\sch}$ and $\Def_{X/G_K \sto BG_K}^{K/W}$, the forgetful map $$\Def_{G_K \acts X}^{\sch} \to \Def_X^{\sch}$$ identifies with the map $F_{\Def_{X/G_K}^{BG_K/BG}}$, and $T_{\Def_{X/G_K}}^{BG_K/BG} \simeq f_*T_{X/G_K \sto BG_K}[1]$. If $X$ is a smooth scheme, $R\Gamma(X, T_X)$ is coconnective. Moreover, if $e \colon W \to BG$ is the natural map, the pullback $e^*$ is faithful and $t$-exact. Hence by \Cref{unobstructedness of fiber of deformations}, this map is unobstructed if \begin{enumerate}
        \item $H^2(G_K, R^0f_*T_{X/G_K \sto BG_K}) = 0$
        \item $H^0(G_K, R^1f_*T_{X/G_K \sto BG_K}) \simeq R^1f_*T_{X/G_K \sto BG_K}$
        \item $H^1(G_K, R^1f_*T_{X/G_K \sto BG_K}) = 0$.
    \end{enumerate}
    But condition (2) is equivalent to the condition that $R_1f_*T_{X/G_K \sto BG_K}$ is a trivial $G$-module. Then, in this case, given that $H^1(G, k) = 0$, it automatically implies that $$H^1(G, R^1f_*T_{X/G_K \sto BG_K}) = 0.$$
    Finally, through the pullback square
\[\begin{tikzcd}
	X & K \\
	{X/G_K} & {BG_K}
	\arrow[from=1-1, to=1-2]
	\arrow[from=1-1, to=2-1]
	\arrow["\lrcorner"{anchor=center, pos=0.125}, draw=none, from=1-1, to=2-2]
	\arrow[from=1-2, to=2-2]
	\arrow[from=2-1, to=2-2]
\end{tikzcd}\]
    we identify $H^i(X, T_X)$ with $R_if_*T_{X/G_K \sto BG_K} \in \qc_{BG_K}$ for all $i \geq 0$. Hence $H^i(X, T_X)$ admits a $G_K$-action and we get the announced result.
\end{proof}

\begin{rk}
    One could extend this result to a flat group $G$ locally of finite presentation using fppf topology instead of etale topology in \Cref{section equivariant deformations}. One could also extend it to a proper classical scheme $X$ which is local complete intersection.
\end{rk}

\begin{rk}
    Let $G$ be a smooth group scheme such that $H^1(G, k) = 0$. We recover the fact that if $H^0(X, T_X) = 0$ and $H^1(X, T_X) = 0$, then an action of $G_K$ on a scheme $X$ lifts together with $X$. If those basic conditions are sometimes sufficient, as in \cite{ekedahl2004nonliftable}, to prove that some deformations are obstructed, they are too strong in general. \Cref{criterion} gives sharper conditions, which will reveal very useful in \Cref{section homogeneous varieties}.
\end{rk}

\subsection{Equivariant deformations for transitive actions} \label{section transitive case}

Let $k$ be an algebraically closed field and $\Lambda$ be a local Noetherian ring with residue field $k$. In this section all schemes are classical. Let $G$ be a group scheme of finite type over $\Lambda$, and write $G_A = G \otimes_{\Lambda} A$ for any $\Lambda$-algebra $A$. We specialize our study of equivariant deformations to the case of a transitive action of $G_k$ on a variety $X$. When $X$ admits a $k$-point, it means that $X$ is $G_k$-homogeneous, that is, $X$ is isomorphic to a quotient $G_k/P$ where $P$ is the stabilizer of the action.
In that case, \textit{pointed} equivariant deformations of $X$ are equivalent to deformations of the stabilizer, as a subgroup scheme of $G_k$. Moreover, when $X$ is smooth, deformations of $X$ come with a deformation of the point $x \colon \Spec(k) \to X$. Finally, when the conditions of \Cref{criterion} are verified, one gets a chain of unobstructed maps of deformation functors:
\[\begin{tikzcd}
	{\Defcl_{P \sto G_k}} & {\Defcl_{G_k \acts X}^*} & {\Defcl_{G_k \acts X}} & {\Defcl_{X}}
	\arrow["\simeq", from=1-1, to=1-2]
	\arrow[from=1-2, to=1-3]
	\arrow[from=1-3, to=1-4]
\end{tikzcd}\]
In that case, it allows us to transfer the study of deformations of the variety $X$ to the study of deformations of $P$ as a subgroup scheme of $G_k$, which are often simpler to study.

\begin{rk}
    Let $k$ be an algebraically closed field, let $G_k$ be a group scheme of finite type acting on a separated scheme $X$ of finite type. The action is called transitive if it is transitive at the level of $k$-points. The scheme $X$ is called $G_k$-homogeneous if the action is transitive and there exists a point $x \in X(k)$ such that the orbit map $G \to X$ is faithfully flat, or equivalently if $X = G/P$, where $P$ is the stabilizer of the action. If $X$ is reduced, a transitive action automatically gives rise to a $G_k$-homogeneous space by \cite[Proposition 7.4]{milne2017algebraic}. 

    Over an Artinian extension $A$ of $k$, we will only use the notion of $G_A$-homogeneous scheme $\X$, in the sense that $\X = G_A/\PP$ where $\PP$ is a flat subgroup scheme of $G_A$. 
    
\end{rk}

Let us begin by defining the functor of pointed deformations.

\begin{defi}
    Let $X$ be a scheme over $k$ and let $x \colon \Spec(k) \to X$ be a $k$-point of $X$. Let $$\Defcl_X^x \colon \CRing_{\Lambda/k}^{\art} \to \Set$$ be the functor which sends an Artinian extension $A \in \CRing_{\Lambda/k}^{\art}$ to the set of isomorphism classes of diagrams
\[\begin{tikzcd}
	{\Spec(k)} & X & {\Spec(k)} \\
	{\Spec(A)} & \X & {\Spec(A)}
	\arrow["x", from=1-1, to=1-2]
	\arrow[from=1-1, to=2-1]
	\arrow[from=1-2, to=1-3]
	\arrow[from=1-2, to=2-2]
	\arrow["\lrcorner"{anchor=center, pos=0.125}, draw=none, from=1-2, to=2-3]
	\arrow[from=1-3, to=2-3]
	\arrow[from=2-1, to=2-2]
	\arrow[from=2-2, to=2-3]
\end{tikzcd}\]
where the map $\X \to \Spec(A)$ is flat.
    The functor $\Defcl_X^x$ is called the functor of pointed deformations of $X$. It comes with a forgetful transformation $\Defcl_X^x \to \Defcl_X$.
\end{defi}

We do not prove here that the functor $\Def_X^x$ verifies Schlessinger's condition because we do not need it. Note that the notion of unobstructedness can be defined for any transformation of functors defined on Artinian extensions.

\begin{lemma} \label{lift of the point}
    Let $X$ be a smooth scheme over $k$, and let $x \colon \Spec(k) \to X$ be a $k$-point of $X$. Then the forgetful map $\Defcl_X^x \to \Defcl_X$ is unobstructed.
\end{lemma}

\begin{proof}
    Let $A \to B$ be a nilpotent extension in $\CRing_{\Lambda/k}^{\art}$. Let $X_B$ be a flat lift of $X$ over $B$ with a lift of the point $x_B \colon \Spec(B) \to X_B$. Let $X_A$ be a a flat lift of $X_B$ over $A$. We want to show that there exists a map $\Spec(A) \to X_A$ such that the following diagram commutes
\[\begin{tikzcd}
	{\Spec(B)} & {X_B} \\
	{\Spec(A)} & {X_A} 
	\arrow[from=1-1, to=1-2]
	\arrow[from=1-1, to=2-1]
	\arrow[from=1-2, to=2-2]
	\arrow[from=2-1, to=2-2]
\end{tikzcd}\]
    But $X$ is smooth, and $X_B$ is flat over $\Spec(B)$, so $X_B$ is automatically smooth over $\Spec(B)$, hence formally smooth. Then, the formal smoothness condition precisely gives that the map $\Spec(B) \to X_B$ lifts to a map $\Spec(A) \to X_A$.
\end{proof}

\begin{rk}
    Let $G$ be a group scheme over $\Lambda$, let $X$ be a scheme over $k$ endowed with a $G_k$-action and a $k$-point $x$. One can also define the functor of pointed equivariant deformations $\Defcl_{G_k \acts X}^x$ which encodes deformations of the $G_k$-action \emph{and} deformations of the point $x$ (independently). If $X$ is smooth, the forgetful transformation $\Defcl_{G_k \acts X}^x \to \Defcl_X^x$ is unobstructed.
\end{rk}

Let us define the functor of subgroup scheme deformations $\Def_{P \sto G_k}^{\Grp}$. We do not prove that it verifies Schlessinger's condition because we do not need it here.

\begin{defi}
    Let $G$ be a group scheme over $\Lambda$.  Let $P$ be a closed subgroup scheme of $G_k$.
    Let $$\overline{\Def}_{P \sto G_k}^{\Grp} \colon \CRing_{\Lambda/k}^{\art} \to \Set$$ be the functor which sends an Artinian extension $A \in \CRing_{\Lambda/k}^{\art}$ to the set of isomorphism classes of flat subgroup schemes $\PP$ of $G_A$ which lifts the inclusion map $P \subset G_k$. The functor $\Def_{P \sto G_k}^{\Grp}$ is called the functor of subgroup scheme deformations of $P$.
\end{defi}

\begin{rk}
    The lift of a closed immersion is automatically a closed immersion by \cite[06AG]{stacks-project}.
\end{rk}

Now we construct the two inverse transformations which give the isomorphism between subgroup scheme deformations $\Defcl_{P \sto G_k}^{\Grp}$ and pointed equivariant deformations $\Defcl_{G_k \acts X}^x$ when $X = G/P$.

\begin{notation}
    If a group scheme $G$ acts on a scheme $X$, and if $x \colon \Spec(k) \to X$ is a $k$-point of $X$, we will denote by $\Stab_x$ the stabilizer of the point $x$ under the action of $G$. The stabilizer sits in the pullback square
\[\begin{tikzcd}
	{\Stab_x} & G \\
	{\Spec(k)} & X
	\arrow[from=1-1, to=1-2]
	\arrow[from=1-1, to=2-1]
	\arrow["\lrcorner"{anchor=center, pos=0.125}, draw=none, from=1-1, to=2-2]
	\arrow[from=1-2, to=2-2]
	\arrow["x", from=2-1, to=2-2]
\end{tikzcd}\]
where the right vertical map is the composite $G \xlongrightarrow{x} G \times_k X \xlongrightarrow{\act} X$.

\end{notation}

\begin{lemma} \label{lift of the stabilizer}
    Let $X$ be a reduced scheme of finite type with a $k$-point $x$, endowed with a transitive action of a group scheme of finite type $G$. Let $A \in \CRing_{\Lambda/k}^{\art}$. Suppose that $\X$ is a flat $A$-scheme which lifts $X$, that it is endowed with a $G_A$-action which lifts the $G_k$-action on $X$, and that it is endowed with an $A$-point $\tilde{x} \colon \Spec(A) \to \X$ which lifts the $k$-point $x \colon \Spec(k) \to X$.
    Then the subgroup scheme $\Stab_{\tilde{x}} \subset G_A$ is a flat lift of the subgroup scheme $\Stab_x \subset G_k$.
\end{lemma}

\begin{proof}
    Consider the following commutative diagram
\[\begin{tikzcd}
	& {\Stab_x} && {G_k} \\
	{\Stab_{\tilde{x}}} && {G_A} \\
	& {\Spec(k)} && X \\
	{\Spec(A)} && \X
	\arrow[from=1-2, to=1-4]
	\arrow[from=1-2, to=2-1]
	\arrow[from=1-2, to=3-2]
	\arrow[from=1-4, to=2-3]
	\arrow[from=1-4, to=3-4]
	\arrow[from=2-1, to=2-3]
	\arrow[from=2-1, to=4-1]
	\arrow[from=2-3, to=4-3]
	\arrow["x"{pos=0.3}, from=3-2, to=3-4]
	\arrow[from=3-2, to=4-1]
	\arrow[from=3-4, to=4-3]
	\arrow["{\tilde{x}}"', from=4-1, to=4-3]
\end{tikzcd}\]
The front square, the bottom square and the back square are cartesian, so by pasting law the top square is also cartesian, so the subgroup scheme $\Stab_{\tilde{x}} \subset G_A$ is a lift of the subgroup scheme $\Stab_x \subset G_k$. 
By \cite[Proposition 7.4]{milne2017algebraic}, the orbit map $G_k \to X$ is flat. By \cite[06AG]{stacks-project}, the map $G_A \to \X$ is also flat. Then, by base change, the group scheme $\Stab_{\tilde{x}}$ is flat over $A$. 
\end{proof}

\begin{defi}
    Let $G$ be a group scheme of finite type over $\Lambda$. Let $X$ be a variety over $k$ endowed with a transitive action of $G_k$, which admits a $k$-point $x \colon \Spec(k) \to X$. Using \Cref{lift of the stabilizer}, we define the natural transformation $$\Stab \colon \overline{\Def}_{G_k \acts X}^x \to \Defcl_{P \sto G_k}$$ which sends, for each $A \in \CRing_{\Lambda/k}^{\art}$, the isomorphism class of an action of $G_A$ on $\X$ to the isomorphism class of its stabilizer at the lift of the point.
\end{defi}

\begin{rk}
    \Cref{lift of the stabilizer} fails when the action is not transitive in general, or when $X$ is not reduced: a lift $\PP$ of the stabilizer is not necessarily flat. Hence, the functor $\Stab$ is not well-defined in general.
\end{rk}

\begin{lemma} \label{lift of the quotient}
    Let $G$ be a flat group scheme of finite type over $\Lambda$, let $P \subset G_k$ be a subgroup scheme, and let $A \in \CRing_{\Lambda/k}^{\art}$. Let $\PP \subset G_A$ be a flat subgroup scheme lift of $P \subset G_k$. Then the quotient $G_A/\PP$ is representable by a scheme and is a flat lift of $G/P$, and the natural action of $G_A$ on $G_A/\PP$ is a lift of the natural action of $G_k$ on $G_k/P$.
\end{lemma}

\begin{proof}
    The quotient $G_A/\PP$ is a scheme by \cite[Chapter 4]{anantharaman1973schemas}, and is flat by \cite[I.5.7(2)]{jantzen2003representations}. Given that $\PP \subset G_A$ is a lift of $P \subset G_k$, the action of $\PP$ on $G_A$ is a lift of the action of $P$ on $G_k$, hence the quotient $G_A/\PP$ is a lift of $G_k/P$, and the induced action of $G_A$ on $G_A/\PP$ is a lift of the action of $G_k$ on $G_k/P$.
\end{proof}

\begin{defi}
    Let $G$ be a flat group scheme of finite type over $\Lambda$. Let $P$ be a subgroup scheme of $G_k$. Let $X = G_k/P$ be the quotient scheme. Using \Cref{lift of the quotient}, we define the natural transformation $$Q \colon \Defcl_{P \sto G_k} \to \Defcl_{G_k \acts X}^x$$ which sends, for each $A \in \CRing_{\Lambda/k}^{\art}$, the isomorphism class of a flat lift $\PP \subset G_A$ of $P \subset G_k$ over $A$ to the isomorphism class of the quotient scheme $G_A/\PP$.
\end{defi}

The transformations $\Stab$ and $Q$ that we just defined are inverse to each other. The only subtle point is to show that the lift of a homogeneous variety is still homogeneous.

\begin{prop} \label{equivalence stab quotient}
    Let $G$ be a flat group scheme of finite type over $\Lambda$. Let $X$ be a variety over $k$, endowed with a transitive action of $G_k$, which admits a $k$-point $x \colon \Spec(k) \to X$. Let $P$ denote the stabilizer of this action at the point $x$.
    Then the natural transformations $$\Stab \colon \Defcl_{G \acts X}^x \to \Defcl_{P \sto G_k}$$ and $$Q \colon \Defcl_{P \sto G_k} \to \Defcl_{G_k \acts X}^x$$ form a natural isomorphism $$\Defcl_{P \sto G_k} \simeq \Defcl_{G_k \acts X}^x.$$
\end{prop}

\begin{proof}
    It is clear that  $\Stab \circ Q = \id$.
    Let us show that $Q \circ \Stab = \id$. Let $A \in \CRing_{\Lambda/k}^{\art}$, and let $\X$ be a flat lift of $X$ over $A$ endowed with a $G_A$-action which lifts the $G_k$-action on $X$, and with an $A$-point $\tilde{x} \colon \Spec(A) \to \X$ which lifts the $k$-point $x \colon \Spec(k) \to X$. Denote by $\PP$ the stabilizer of the action of $G_A$ on $\X$ at the point $\tilde{x}$. By the universal property of quotient, the orbit map $G_A \to \X$ factors through a map of schemes $G_A/\PP \to \X$, The reduction of this map to $k$ gives the isomorphism $G_k/P \iso X$. But the reduction functor is conservative, so the map $G_A/\PP \to \X$ is an isomorphism. Hence, $Q \circ \Stab = \id$.
\end{proof}

\begin{rk}
    There is another perspective on this equivalence using quotient stacks. When $X = G/P$, then $X/G \simeq BP$. We know that the deformation functor $\Def_{G_k \acts X}$ is equivalent to the deformation functor $\Def_{X/G \sto BG_k}$, that is, to the deformation functor $\Def_{BP \sto BG_k}$. Hence, \Cref{equivalence stab quotient} can be rephrased as an equivalence between $\Def_{P \subset G_k}$ and $\Def_{BP \sto BG_k}^*$. Forgetting about the inclusion in $G_k$, it amounts to an equivalence between group scheme deformations $\Def_P$ and pointed stack deformations $\Def_{BP}^*$, which simply corresponds to looping/delooping equivalence. This readily generalizes to the derived setting, considering flat derived group schemes deformations. 
\end{rk}

Finally, in the case of a $G_k$-homogeneous variety $X$, when the conditions of \Cref{criterion} are fulfilled, one can transfer the study of deformations of $X$ to the study of deformations of the stabilizer of that action.

\begin{prop} \label{unobstructedness transitive}
    Let $G$ be a flat group scheme of finite type over $\Lambda$.
    Let $X$ be a smooth variety over $k$ endowed with a transitive action of $G_k$. Suppose that $X$ admits a $k$-point $x$, and denote by $P$ the stabilizer at $x$ of the action of $G_k$.
    Suppose that the following conditions are verified:
    \begin{enumerate}
        \item $H^1(G_k,k) = 0$.
        \item $H^2(G_k, H^0(X, T_X)) = 0$ 
        \item $H^1(X, T_X)$ is a trivial $G_k$-module.
    \end{enumerate}
    Then the natural transformation $\Defcl_{P \sto G_k} \to \Defcl_X$ is unobstructed.
\end{prop}

\begin{proof}
    Consider the following chain of maps of deformation functors
\[\begin{tikzcd}
	{\Defcl_{P \sto G_k}} & {\Defcl_{G_k \acts X}^*} & {\Defcl_{G_k \acts X}} & {\Defcl_{X}}
	\arrow["\simeq", from=1-1, to=1-2]
	\arrow[from=1-2, to=1-3]
	\arrow[from=1-3, to=1-4]
\end{tikzcd}\]
    The action of $G_k$ on $X$ is transitive so, by \Cref{equivalence stab quotient}, the first map is an isomorphism. The scheme $X$ is smooth so, by \Cref{lift of the point}, the second map is unobstructed. Then, by \Cref{criterion}, conditions (1), (2), (3) together imply that the third map is unobstructed. Thus, the composite map is unobstructed.
\end{proof}

\section{A study of homogeneous varieties with non-reduced parabolic stabilizer} \label{section homogeneous varieties}

Let $G$ be a reductive group over an algebraically closed field $k$ of characteristic $p > 0$. The theory of projective homogeneous varieties under $G$, that is, generalized flag varieties, is well known. The stabilizer of the action is a parabolic subgroup $P$ of $G$, and the structure and behaviour of parabolic subgroups is very well understood. However, all this theory focuses on \textit{smooth} parabolic subgroups. In \cite{Wenzel1993}, Wenzel showed that non-reduced parabolic subgroups existed, giving rise to smooth projective homogeneous spaces, often called \textit{unseparated flag varieties}, and he classified them. In \cite{HaboushLauritzen1993}, Haboush and Lauritzen established the fundamental properties of those unseparated flag varieties. Their behaviour can be quite mysterious: for instance, Kempf vanishing is not verified for those varieties. The classification of Wenzel was extended recently to characteristic $2$ and $3$ in \cite{maccan2023parabolicsubgroupscharacteristics}.

We continue this study by elucidating some aspects of the behaviour of unseparated flag varieties. When the stabilizer $P$ is reduced, the cohomology of vector bundles on $G/P$ can be reduced to the cohomology on $G/B$, where $B$ is the Borel subgroup. This fails when $P$ is non-reduced. In \Cref{section general results}, we establish a formula to compute the cohomology of equivariant vector bundles on unseparated flag varieties (\Cref{induction formula}).

Moreover, when $P$ is reduced, due to work of \cite{bott1957homogeneous} in characteristic $0$ and \cite{Demazure1977} in characteristic $p$, the higher cohomology of the tangent bundle of $G/P$ is trivial. This also fails for unseparated flag varieties in general (for instance, we will prove in \Cref{section mixed characteristic deformations} that the second cohomology group of the tangent bundle does not vanish in general). In \Cref{section  general results}, we construct a filtration of the tangent bundle of unseparated flag varieties (\Cref{filtration}) which allows us to get first results on the $G$-module structure of the cohomology of the tangent bundle (\Cref{conditions fulfilled}): in most cases, $H^1(X, T_X)$ is a trivial $G$-module and for all $i \geq 2$, $$H^i(G, H^0(X, T_X)) = 0.$$

Finally, in \Cref{section classification simple modules}, we show that the classification of simple representations of non-reduced parabolic or reductive subgroups is analogous to the classical case: they correspond to (appropriately defined) dominant weights (\Cref{BW theorem non-inf}, \Cref{BW thm parabolic}). We also prove the appropriate generalization of Steinberg theorem, that is, that a $p$-adic expansion of a weight leads to a decomposition of the associated simple representation with Frobenius twists of elementary simple representations (\Cref{steinberg}).

\Cref{section filtration and cohomology} and \Cref{section classification simple modules} are completely independent.

\subsection{General results about non-reduced parabolic subgroups} \label{section general results}

Let $G$ be a reductive group over an algebraically closed field $k$. We call parabolic subgroup of $G$ a subgroup such that $G/P$ is a projective variety. In characteristic $0$, all parabolic subgroups are smooth. However, in positive characteristic, since the work of \cite{Wenzel1993}, we know that there exists non-reduced parabolic subgroups which give rise to interesting smooth projective quotient varieties. 

If $P$ is a reduced parabolic subgroup, the parabolic induction $\ind_P^G$ can be simply computed by passing through the Borel subgroup $$\ind_P^G \simeq \ind_B^G \res_B^P$$
When $P$ is non-reduced, this formula fails, making it much harder to compute parabolic induction. This section culminates in \Cref{induction}, where we establish a formula to compute non-reduced parabolic induction for Frobenius-twisted modules.

\subsubsection{Preliminaries: the structure of non-reduced parabolic subgroups} \label{whl}

In characteristic $0$, all group schemes locally of finite type are smooth, and in particular parabolic subgroups are smooth. In positive characteristic, even if parabolic subgroups are often supposed to be smooth by default, there exist non-smooth parabolic subgroups (but the quotient by those non-smooth parabolic subgroups is still smooth and projective!).

In general, in positive characteristic, by homogeneity of $G$, the obstruction for a group scheme to be smooth can only be non-reducedness. The fundamental point is that, whereas in characteristic $0$ the additive group $\G_a$ has no non-trivial subgroups, in positive characteristic it admits non-trivial subgroups, which are non-reduced and infinitesimal (in the sense that they have only one $k$-point): the kernels of the iterated Frobenius morphisms, denoted $\alpha_{p^r}$. Those infinitesimal subgroups are the building blocks on which are built non-reduced parabolic subgroups.

The simplest interesting example of non-reduced parabolic subgroup is, in $\SL_3$, the subgroup $$\begin{pmatrix}
* & * & * \\
\alpha_p & * & * \\
0 & 0 & * 
\end{pmatrix}.$$ The induced quotient is the twisted incidence variety  $x_0y_0^p + x_1y_1^p + x_2y_2^p = 0$ in $\pro^2 \times \pro^2$, as explained for instance in \cite{maccan2023projectivehomogeneousvarietiespicard}.

In \cite{Wenzel1993}, Wenzel classified non-reduced parabolic subgroups. Then, in \cite{HaboushLauritzen1993}, Haboush and Lauritzen carried out a deeper study. This section is mainly extracted from their work.
A non-reduced parabolic subgroup can be uniquely expressed as the product of Frobenius kernels of reduced parabolic subgroups $P^0_{r_0} ... P^l_{r_l}$, where $P_l \subset ... \subset P^0 = G$, and $r_0 < ... < r_l \leq +\infty$.
Those non-reduced parabolic subgroups also naturally give rise to non-reduced reductive subgroups of $G$, their Levi factor, and to infinitesimal subgroups of $G$, their kernels. 

\begin{defi}
    A parabolic subgroup $P$ of $G$ is a subgroup scheme of $G$ which contains a Borel subgroup, or equivalently such that $G/P$ is a projective variety.    
\end{defi}

\paragraph{Notations.}

Until the end of this chapter, we will consider an algebraically closed field $k$ of characteristic $p \geq 5$, a (connected) reductive group $G$, a maximal torus $T$ of $G$ and a Borel subgroup $B$ containing $T$. We will denote by $S$ the set of simple roots with respect to $T$. We will write $s = |S|$ for the number of simple roots, that is, the semi-simple rank of $G$.

\paragraph{The usual classification of reduced parabolic subgroups.}

To every root $\alpha \in R$ is associated a root subgroup $U_{\alpha}$, isomorphic to the additive group $\G_a$, and whose tangent is isomorphic to $\g_{\alpha}$. In the same way that $\g$ is decomposed with the tangent space of the torus and the $\g_{\alpha}$, the group $G$ is generated by the torus and the root subgroups:  $$G = \<T, (U_{\alpha})_{\alpha \in R}\>.$$ 

Reduced  parabolic subgroups which contains the Borel $B$ are in bijection with subsets of simple roots $I \subset S$. The parabolic subgroup associated to $I$ is the subgroup generated by $B$ and the root subgroups for positive roots generated by $I$: $$P^I := \<B, (U_{\beta})_{\beta \in \N I \cap R}\>.$$
For instance the Borel subgroup $B$, that is, the minimal parabolic subgroup containing $B$, corresponds to $I = \emptyset$.

\begin{notation} \label{notation parabolic}
    We will write $P^I$ for the parabolic subgroup generated by $I \subset S$, and $G^I$ for the Levi subgroup of $P^I$.
\end{notation}

\begin{ex}
    In the case of $\SL_n$, simple roots are given by the $\alpha_i = e_i - e_{i+1}$ for $1 \leq i \leq n-1$. Let $I = S \backslash \{\alpha_{i_1}, ... , \alpha_{i_l}\}$, with $1 \leq i_1 < ... < i_l \leq n-1$. The quotient $SL_n/P^I$ is the (partial) flag variety $Fl(i_1, ... , i_l)$. Let us give some low-dimensional examples:
    \begin{enumerate}

    \item In $SL_2$, up to conjugation, the only parabolic subgroups are the Borel subgroup $B$ and the whole group $SL_2$.

    \item In $SL_3$, the parabolic subgroup  $$\begin{pmatrix}
* & * & * \\
* & * & * \\
0 & 0 & * 
\end{pmatrix}$$ corresponds to the subset of simple roots which is the singleton $\{\alpha_1\}$, and induces the quotient variety $\Grass(2,3) \simeq \pro ^2$.

\item In $\SL_4$, the parabolic subgroups $$\begin{pmatrix}
* & * & * & * \\
0 & * & * & * \\
0 & 0 & * & * \\
0 & 0 & * & *
\end{pmatrix}, \begin{pmatrix}
* & * & * & * \\
* & * & * & * \\
* & * & * & * \\
0 & 0 & 0 & *
\end{pmatrix}, \begin{pmatrix}
* & * & * & * \\
* & * & * & * \\
0 & 0 & * & * \\
0 & 0 & * & *
\end{pmatrix}, \begin{pmatrix}
* & * & * & * \\
0 & * & * & * \\
0 & * & * & * \\
0 & 0 & 0 & *
\end{pmatrix}.$$
induce the quotients $\Fl(1,2,4)$, $\Grass(3,4) \simeq \pro^3$, $\Grass(2,4)$ and $\Fl(1,3,4)$ respectively.
\end{enumerate}
\end{ex}

\paragraph{The Frobenius endomorphism of a group scheme.}

On affine schemes, that is, $\F_p$-algebras, one can consider the ordinary Frobenius endomorphism. The Frobenius endomorphims glue to give an endomorphism of any scheme over $\F_p$.

\begin{defi}
    Let $X$ be a scheme over $\F_p$. The absolute Frobenius endomorphism of $X$ is the morphism of $\F_p$-schemes defined as the identity on the underlying topological spaces and the map of sheaves of rings $O_X \to F_*O_X$ given by $f \mapsto f^p$. 
\end{defi}

However, in general, if the scheme $X$ is defined over a base $S$, the absolute Frobenius is not a morphism of $S$-schemes. To get a morphism of $S$-schemes, one needs to consider the relative Frobenius.

\begin{defi}
    The relative Frobenius morphism $F_{X/S}$ is defined with the following diagram: 
\[\begin{tikzcd}
	X \\
	& {X^{(1)}} & X \\
	& S & S
	\arrow["{F_{X/S}}"{description}, from=1-1, to=2-2]
	\arrow["{F_X}", bend left = 20,from=1-1, to=2-3]
	\arrow[bend right = 20, from=1-1, to=3-2]
	\arrow[from=2-2, to=2-3]
	\arrow[from=2-2, to=3-2]
	\arrow["\lrcorner"{anchor=center, pos=0.125}, draw=none, from=2-2, to=3-3]
	\arrow[from=2-3, to=3-3]
	\arrow["{F_S}"', from=3-2, to=3-3]
\end{tikzcd}\]
\end{defi}

Over the base $S = \F_p$, the Frobenius of $S$ is the identity so $X^{(1)} \simeq X$. This holds in fact as soon as the scheme $X$ has an $\F_p$-structure.

\begin{prop} \label{Fp-structure}
    Let $X$ be a scheme with an $\F_p$-structure, that is, there exists an $\F_p$-scheme $\tilde{X}$ such that $X \simeq \tilde{X} \times_{\F_p} S$. Then $X^{(1)} \simeq X$.
\end{prop}

\begin{proof}
    The following diagram
\[\begin{tikzcd}
	& {X^{(1)}} && S \\
	X && S \\
	\\
	\X && {\Spec(F_p)}
	\arrow[from=1-2, to=1-4]
	\arrow[from=1-2, to=2-1]
	\arrow[dashed, from=1-2, to=4-1]
	\arrow["{F_S}"', from=1-4, to=2-3]
	\arrow[from=1-4, to=4-3]
	\arrow[from=2-1, to=2-3]
	\arrow[from=2-1, to=4-1]
	\arrow[from=2-3, to=4-3]
	\arrow[from=4-1, to=4-3]
\end{tikzcd}\]
shows that by pasting law, $X^{(1)} \simeq \tilde{X}  \times_{F_p} S$, hence that $X^{(1)} \simeq X$. 
\end{proof}

We are interested in the Frobenius of a reductive group, which admits an $\F_p$-structure, and even a $\Z$-structure. Thus, we get an honest Frobenius endomorphism $F_G \colon G \to G$.

\begin{defi} \label{defi Frobenius}
    Let $G$ be a reductive group over $k$.  The relative Frobenius endomorphism of $G$ is a group scheme morphism. Identifying $G^{(1)}$ and $G$ by \Cref{Fp-structure}, we get a group endomorphism $F_G \colon G \to G$, called the Frobenius endomorphism of $G$.
    \begin{enumerate}

    \item Write $F_G^r$ for the $r$-th iterated Frobenius endomorphism of $G$.

    \item Write $G_r$ for the kernel of the iterated Frobenius $F_G^r \colon G \to G$. By convention, we will sometimes use the notation $G_{\infty} := G$.
    \end{enumerate}
    Then the quotient $G/G_r$ is the image of the endomorphism $F_G^r$.
\end{defi}

\begin{rk}
    Le $G$ be a reductive group scheme. Frobenius kernels are infinitesimal schemes, in the sense that they have only one $k$-point (except for $G_{\infty}$ of course). There is an ascending chain $$\{1\} = G_0 \subset G_1 \subset ... \subset G_r \subset ...$$
\end{rk}

\begin{notation}
For the rest of this work, if $H$ is a group scheme, we will write $F_H^r$ for the Frobenius morphism from $H$ to $H/H_r$.
\end{notation}

\begin{lemma}[{\cite[I.9.4]{jantzen2003representations}}] \label{kernel of subgroup}
    Let $K$ be a subgroup scheme of a group scheme $L$. Then the Frobenius endomorphism $F_K^r$ of $K$ is the restriction of the Frobenius endomorphism $F_L^r$ of $L$ to $K$, hence $$K_r = K \cap L_r.$$ 
    In particular,
    \begin{itemize}
    \item if $r \leq s$, then $(G_r)_s = G_r$.
    \item if $r \geq s$, then $(G_r)_s = G_s$.
    \end{itemize}
\end{lemma}

\begin{notation}
    Let $\alpha_{p^r}$ denote the kernel of the $r$-th iterated Frobenius endomorphism of the additive group, $F_{\G_a}^r$. As functor of point, for any $k$-algebra $A$, $$\alpha_{p^r}(A) = \{a \in A, a^{p^r} = 0\}.$$
    As a scheme, $\alpha_{p^r} \simeq \Spec(k[T]/T^{p^r})$. 
\end{notation}

\paragraph{Classification of general parabolic subgroups.}

For the classification that we are going to explain, we refer mainly to \cite{Wenzel1993} and \cite{HaboushLauritzen1993}.
Whereas a reduced parabolic subgroup corresponds to a subset of simple roots, a non-reduced parabolic subgroup corresponds to a filtration of the set of simple roots, together with an integer value for each of the subset of the filtration. In other words, it corresponds to a function $f \colon S \to \N \cup +\infty$ (functions with values only in $\{0, \infty\}$ thus correspond to reduced parabolic subgroups).
To a function can also be associated a non-reduced reductive subgroup of $G$, called the Levi factor of the parabolic subgroup.

We slightly enlarge the usual definition by allowing functions which do not take the value $+\infty$, which give rise to infinitesimal non-reduced reductive subgroups of $G$ (and do not correspond to parabolic subgroups). 

\begin{defi}
    A function $f \colon S \to \N \cup +\infty$ is called a $W$-function.
    Let $r_0 < r_1 <... < r_l$ denote the values of $f$. 
    \begin{enumerate}
    \item We call $l$ the \emph{width} of $f$ and $r_0$ its \emph{minimal altitude}.
    \item The $W$-function $f$ is said to be infinitesimal if it does not take the value $+\infty$.
    \item For $0 \leq k \leq l$, write $I_f^k := f^{-1}(\{r_k, ..., r_l\})$. The filtration of the set of simple roots $S$ $$0 \subset I^l \subset I^{l-1} \subset ... \subset I^1 \subset I^0 = S$$ is called the characteristic filtration of $f$.
    \item If $f$ is not infinitesimal and $l \geq 1$, let $\f$ denote the infinitesimal function identical to $f$, except that it sends $I^{l}$ to $r_{l-1}$ instead of $+\infty$. 
    \end{enumerate}
\end{defi}

\begin{defi} \label{definition parabolic}
    Let $f$ be a $W$-function with values $r_0 < r_1 < ... < r_l$ and characteristic filtration $I_l \subset ... \subset I_1 \subset I_0 = S$.
    \begin{enumerate}
    \item For all $0 \leq k \leq l$, we will write the parabolic subgroup $P^k := P^{I^k}$ and its Levi subgroup $G^k := G^{I^k}$, with \Cref{notation parabolic} (the dependence to $f$ is implicit and will be clear). Notice that $P^0 = G$ and $G^0 = G$.
    \item Let $G_f$ denote the subgroup $G^{0}_{r_0}G^{1}_{r_1}...G^{l}_{r_l}$. Notice that $G_f$ is infinitesimal if and only if $f$ is infinitesimal.
    \item If $f$ is not infinitesimal, let $P_f$ denote the subgroup $P^{0}_{r_0}P^{1}_{r_1}...P^{l}_{r_l}.$ In this case, $G_f$ will be called the Levi factor of $P_f$. By extension, we also call $l$ the width of $P_f$ and $G_f$.
    \item For all $0 \leq k \leq l$, we will write $\overline{P}^k := P^0_{r_0}...P^{k-1}_{r_{k-1}}P^k$, which has width $k$. Those subgroups induce a filtration $$P_f = \overline{P}^l \subset \overline{P}^{l-1} \subset ... \subset \overl^1 \subset \overline{P}^0 = G$$ where the width decreases from $l$ to $0$.
    \end{enumerate} 
\end{defi}

\begin{rk}
    Let $f$ be a non-infinitesimal $W$-function.
    \begin{enumerate}
    \item $P_f = G_fB$.
    \item $G_f = G_{\f}G^l$ and $P_f = G_{\f}P^l$.
    \item The subgroup $P^l$ is the reduced part of $P_f$, denoted $(P_f)_{\red}$.
    \end{enumerate}
\end{rk}

\begin{rk}
    Our definition of width is shifted by $1$ compared to the definition of \cite{HaboushLauritzen1993}.
\end{rk}

\begin{rk} \label{structure parabolic subgroups}
    Let $f$ be a non-infinitesimal $W$-function. There is an alternative description of the parabolic subgroup $P_f$. If $\beta$ is a positive root, its support $\supp(\beta)$ means the subset of simple roots whose coefficients in the decomposition of $\beta$ are positive. The function $f$ can be extended to a function on all positive roots by setting $$f(\beta) = \min_{\alpha \in \supp(\beta)} f(\alpha).$$ 
    Then $P_f$ is the subgroup generated by the Borel subgroup $B$ and the kernels of the root subgroups: $$P_f = \<B, ((U_{\beta})_{f(\beta)})_{\beta \in R^+}\>.$$
    The kernel of a root subgroup $(U_{\beta})_{f(\beta)}$ is isomorphic to the kernel of the additive group $\alpha_{p^{f(\beta)}}$. 
\end{rk}

The construction $f \mapsto P_f$ gives a correspondence between $W$-function and parabolic subgroups up to conjugacy. Moreover, the construction $f \mapsto G_f$ gives non-smooth reductive subgroups, the ``Levi factor'' of $P_f$. However, beware that when $P_f$ is non-reduced, $P_f$ does not admit an analogue to the usual Levi decomposition for reduced parabolic subgroups.

\begin{prop}[{\cite[Part 1, Theorem 6]{HaboushLauritzen1993}}]
    For any non-infinitesimal $W$-function $f$,
    $P_f$ is a parabolic subgroup of $G$. This construction establishes a bijection between the set of parabolic subgroups of $G$ containing $B$ and the set of non-infinitesimal $W$-functions.
\end{prop}

\begin{rk}
    Let $f$ be a $W$-function with values $r_0 < ... < r_l$. The group $G_f$ is reductive in the sense that it contains no non-trivial connected unipotent normal subgroups (that is, we simply drop the usual smoothness assumption). We prove it in \Cref{Gf reductive}. Moreover, the Frobenius kernel $T_{r_l}$ is a maximal torus of $G_f$.
\end{rk}

\begin{ex}
    \begin{enumerate}
    \item Frobenius kernels $G_r$ are (infinitesimal) non-reduced reductive subgroup of $G$ of width $0$. \item In \cite[II.9]{jantzen2003representations}, Jantzen studies the subgroups $G_rT$ and $G_rB$. The subgroup $G_rB$ is a non-reduced parabolic subgroup of $G$ of width $1$ and $G_rT$ is its Levi factor.
    \end{enumerate}
\end{ex}

\begin{ex}
    Let $f$ be a non-infinitesimal $W$-function.
    \begin{enumerate}
    \item Reduced parabolic subgroups correspond to the cases where $f$ has values only $0$ or $+\infty$. In this case, $P_f = P^I$ where $I = f^{-1}(\{+\infty\})$.

    \item $P_f$ has width $0$ if and only if $P_f = G$.

    \item $P_f$ has width $1$ if and only if $P_f = G_rQ$ where $Q$ is a proper reduced parabolic subgroup and $0 \leq r < +\infty$.
    Thus, reduced parabolic subgroups have width $1$, except $G$ which has width $0$.

    \item $P$ has width $2$ if and only if $P = G_rR_sQ$ where $Q \subset R \subset G$ are reduced parabolic subgroups and $0 \leq r < s < +\infty$. 
    \end{enumerate}
\end{ex}

\begin{ex}
    For the case of $\SL_n$, let us give some examples of non-reduced parabolic subgroups in small dimension.
    \begin{enumerate}
    \item In $SL_2$, up to conjugacy, the only non-reduced parabolic subgroups are of the form  $$\begin{pmatrix}
* & * \\
\alpha_{p^r} & *
\end{pmatrix}.$$ They have width $1$. We will see that the induced quotient varieties are all isomorphic to $\pro^1$ (as varieties, not as $SL_2$-varieties).

    \item In $SL_3$, we get the first interesting example of non-reduced parabolic subgroup:  $$\begin{pmatrix}
* & * & * \\
\alpha_p & * & * \\
0 & 0 & * 
\end{pmatrix}$$
It has width $2$. The induced quotient variety is isomorphic to the twisted incidence variety $x_0y_0^p + x_1y_1^p + x_2y_2^p = 0$ in $\pro^2 \times \pro^2$.

    \item For $SL_4$, here are some examples of non-reduced parabolic subgroups of width $2$: 
    $$\begin{pmatrix}
* & * & * & * \\
0 & * & * & * \\
0 & 0 & * & * \\
0 & 0 & \alpha_p & *
\end{pmatrix},
\begin{pmatrix}
* & * & * & * \\
\alpha_p & * & * & * \\
\alpha_p & \alpha_p & * & * \\
0 & 0 & 0 & *
\end{pmatrix}, 
\begin{pmatrix}
* & * & * & * \\
\alpha_p & * & * & * \\
0 & 0 & * & * \\
0 & 0 & \alpha_p & *
\end{pmatrix}, \begin{pmatrix}
* & * & * & * \\
\alpha_p & * & * & * \\
0 & 0 & * & * \\
0 & 0 & * & *
\end{pmatrix}, 
\begin{pmatrix}
* & * & * & * \\
* & * & * & * \\
\alpha_p & \alpha_p & * & * \\
0 & 0 & 0 & *
\end{pmatrix}$$
The first three examples verify $(P_f)_{\red} = B$, but not the last two.

    \end{enumerate}
\end{ex}

\begin{ex} \label{example zdanowicz}
One of the motivation of this work is to generalize the work of Zdanowicz in \cite{zdanowicz2018arithmetically}, where he studies the deformations of varieties defined as the Frobenius twist of the projection from a flag variety to a Grassmannian:

\[\begin{tikzcd}
	X & {\Fl(1,k,n)} \\
	{\Grass(k,n)} & {\Grass(k,n)}
	\arrow[from=1-1, to=1-2]
	\arrow[from=1-1, to=2-1]
	\arrow["\lrcorner"{anchor=center, pos=0.125}, draw=none, from=1-1, to=2-2]
	\arrow["\pi", from=1-2, to=2-2]
	\arrow["F"', from=2-1, to=2-2]
\end{tikzcd}\]
for $2 \leq k \leq n-2$.
As pointed out in \cite{totaro2017failure}, this variety is isomorphic to the quotient variety induced by the non-reduced parabolic subgroups of $\SL_n$ having characteristic filtration $0 \subset S \backslash \{\alpha_1, \alpha_k\} \subset S \backslash \{\alpha_k\}$ with degrees of Frobenius twist $0 < 1 < +\infty$, that is, defined by the function $$
f(\alpha_i) = \left\{
    \begin{array}{ll}
        \infty & \mbox{if } i \neq 1, k \\
        0 & \mbox{if } i = k\\
        1 & \mbox{if } i = 1,
    \end{array}
\right.
$$
This non-reduced parabolic subgroup is ``intermediary'' between the reduced parabolic subgroup corresponding to the subset $S \backslash \{\alpha_1, \alpha_k\}$, which induces the quotient variety $\Grass(k,n)$, and the reduced parabolic subgroup corresponding to the subset $S \backslash \{\alpha_1, \alpha_k\}$ which induces the quotient variety variety $\Fl(1,k,n)$.
For instance for $n = 4$ and $d = 2$, Zdanowicz and Totaro's construction corresponds to the non-reduced parabolic subgroup

 $$\begin{pmatrix}
* & * & * & * \\
\alpha_p & * & * & * \\
0 & 0 & * & * \\
0 & 0 & * & *
\end{pmatrix}.$$ 

\end{ex}

If $G$ is a reductive group, one defines the unipotent subgroup $U$ of $G$ and a Borel $B$. We do the same for non-reduced reductive subgroups.

\begin{defi} \label{decomposition product}
    Let $f$ be a $W$-function of width $l$, with values $r_0 < r_1 < ... < r_l$. Write $G_f = G^0_{r_0}...G^l_{r_l}$.
    For $0 \leq k \leq l$, write $U^k$ and $(U^k)^+$ for the usual negative and positive unipotent part of $G^k$. We define:
    \begin{enumerate}

\item the unipotent subgroup of $G_f$ $$U_f := U^0_{r_0}...U^l_{r_l},$$ and its opposite $$(U_f)^+ := (U_0)^+_{r_0}...(U_l)^+_{r_l}.$$

\item $B_f := U_f \rtimes T_{r_l}$

\item if $f$ is not infinitesimal, $$U^0_f := \prod_{\substack{0 \leq k \leq l-1 \\ \la \in R_k \backslash R_{k+1}}} (U_{\la})_{r_k}$$ and its opposite $$(U^0_f)^+ := \prod_{\substack{0 \leq k \leq l-1 \\ \la \in R_k \backslash R_{k+1}}} (U_{\la})_{r_k}^+,$$ where $R_k$ denotes the subset of positive roots generated by $I^k$.
    \end{enumerate}
\end{defi}

\begin{prop}[{\cite[Proposition 7, (3)]{HaboushLauritzen1993}}] \label{product decomposition}
    \begin{enumerate}
        \item Let $f$ be a non-infinitesimal $W$-function. There is the isomorphism of schemes (not of group schemes)
    $$P_f \simeq (U^0_f)^+ \times P^l.$$ and $$G_f \simeq (U^0_f)^+ \times (U_l^+ B_f).$$
    
    \item Let $f$ be an infinitesimal $W$-function. There is the isomorphism of schemes (not of group schemes) $$G_f \simeq (U_f)^+ \times B_f.$$
    \end{enumerate}
 \end{prop} 

When $G$ is a classical reductive group, $U^+B$ is an open dense subscheme of $G$, called a big cell of $G$. This holds analogously for non-reduced reductive groups.

\begin{prop} \label{density}
    Let $f$ be a $W$-function. Then
    \begin{enumerate}
    \item $(U_f)^+B_f$ is dense in $G_f$.
 \item $(U_f)^+B$ is dense in $P_f$.
    \end{enumerate}
\end{prop}

\begin{proof}
    \begin{enumerate}
    \item Let us show that $(U_f)^+B_f$ is dense in $G_f$. \Cref{product decomposition} gives the isomorphisms of schemes $$G_f \simeq (U^0_f)^+ \times (U_l^+ B_f) \simeq (U^0_f)^+ \times U^0_f \times G^l.$$ 
    But $G^l$ is reduced so, by \cite[II.1.9]{jantzen2003representations}, $U^l \times B^l$ is dense in $G^l$, so $(U^0_f)^+ \times U^0_f \times U^l \times B^l$ is dense in $G_f$. By reordering, we get that $(U_f)^+ \times B_f$ is dense in $G_f$.
    \item The proof is analogous to (1).
    \end{enumerate}
\end{proof}

\subsubsection{Quotients and parabolic induction for non-reduced parabolic subgroups} \label{induction from non-reduced}

If $P \subset Q$ are two reduced parabolic subgroups of $G$, the quotient $Q/P$ is a smooth projective variety. This becomes false in general for non-reduced parabolic subgroups. For instance, if $P$ is a non-reduced parabolic subgroup and $P_{\red}$ is its reduced part, the quotient $P/P_{\red}$ is affine (isomorphic to $U_f^0$). The natural surjection $G/P_{\red} \to G/P$ is purely inseparable. When $P$ has width $1$, the map $G/P_{\red} \to G/P$ is nothing else than the Frobenius. 

However, we show that considering the filtration of \Cref{definition parabolic} $$P = \overl^l \subset \overl^{l-1} \subset ... \subset \overl^1 \subset \overl^0 = G,$$ then the quotients of the form $\overl^k/\overl^j$ are smooth projective, and isomorphic, as varieties (not as homogeneous varieties), to $P^k/P^j$ (\Cref{indk}).

Quotients are related to parabolic induction through the (homotopy) pullback diagram
\[\begin{tikzcd}
	{Q/P} & {*} \\
	{G/P} & {G/Q} \\
	BP & BQ
	\arrow[from=1-1, to=1-2]
	\arrow[from=1-1, to=2-1]
	\arrow["\lrcorner"{anchor=center, pos=0.125}, draw=none, from=1-1, to=2-2]
	\arrow[from=1-2, to=2-2]
	\arrow["\pi", from=2-1, to=2-2]
	\arrow[from=2-1, to=3-1]
	\arrow["\lrcorner"{anchor=center, pos=0.125}, draw=none, from=2-1, to=3-2]
	\arrow[from=2-2, to=3-2]
	\arrow[from=3-1, to=3-2]
\end{tikzcd}\]
where $BP$ and $BQ$ are the classifying stacks of $P$ and $G$ in the fppf topology. A $P$-module $M$ corresponds to a $G$-equivariant vector bundle on $G/P$, and the induction $\ind_P^Q M$ corresponds to the pushforward $R\pi_* M$. 

When $P$ and $Q$ are reduced, projectiveness of $Q/P$ ensures that the induction functor $\ind_P^Q$ preserves the trivial module. It implies the important property that induction of a $P$-module to $G$ can be computed by considering only the $B$-module structure and inducing from $B$ to $G$ (\Cref{pushforward of structure sheaf}).
As we said, this fails when $P$ is non-reduced: to compute the induction of a $P$-module to $G$, one cannot simply consider the $B$-module structure. That makes the computation of induction far more intricate. In \Cref{induction}, we establish a formula to compute induction for Frobenius-twisted modules.

\begin{defi}
    A morphism of group schemes $f \colon H \to G$ induces a base change adjunction for quasi-coherent sheaves over the fppf classifying stacks 
    $$Bf^* \colon \begin{tikzcd}
	D(\qc_{BG}) & D(\qc_{BH})
	\arrow[shift left, from=1-1, to=1-2]
	\arrow[shift left, from=1-2, to=1-1]
\end{tikzcd} \colon Bf_*.$$
In the setting of representation theory, this base change adjunction is usually written $$\res_H^G \colon \begin{tikzcd}
	D(\Mod_G) & D(\Mod_H)
	\arrow[shift left, from=1-1, to=1-2]
	\arrow[shift left, from=1-2, to=1-1]
\end{tikzcd} \colon \ind_H^G.$$
The functor $\res_H^G$ is called the restriction functor. It is an exact functor. The functor $\ind_H^G$ is called the (derived) induction functor.
\end{defi}

\begin{prop} \label{pushforward corresponds to induction}
    A morphism of group schemes $H \to G$ sits into a homotopy pullback diagram 
\[\begin{tikzcd}
	{G/H} & {*} \\
	BH & BG
	\arrow["\pi", from=1-1, to=1-2]
	\arrow["q"', from=1-1, to=2-1]
	\arrow["\lrcorner"{anchor=center, pos=0.125}, draw=none, from=1-1, to=2-2]
	\arrow["e", from=1-2, to=2-2]
	\arrow[from=2-1, to=2-2]
\end{tikzcd}\]
    The pullback $q^*$ induces an equivalence between the category of $H$-modules and the category of $G$-equivariant vector bundles on $G/H$.
    Moreover, if $M$ is a $H$-module, there is the base change property $$\pi_* q^* M \simeq e^* \ind_H^G M.$$ 
\end{prop}

\paragraph{Geometry of quotients.}

If $P \subset Q$ are two reduced parabolic subgroups of $G$, projectiveness and irreducibility of $Q/P$ imply that the cohomology of the structure sheaf of the quotient is trivial. By the projection formula, it gives the important consequence that induction from a parabolic subgroup can simply be computed from the Borel subgroup.

\begin{prop}[{\cite[II.1.8]{jantzen2003representations}}] \label{pushforward of structure sheaf}
    Let $P \subset Q$ be two reduced parabolic subgroups of $G$. 
    \begin{enumerate}
    \item The quotient $Q/P$ is a smooth, projective and irreducible variety.

    \item Let $\pi \colon G/P \to G/Q$. The natural map of $O_{G/Q}$-modules $O_{G/Q} \to R\pi_* O_{G/P}$ is an isomorphism. In other words, the following are verified:
        \begin{enumerate}[label=(\roman*)]
            \item the natural map $k \to R^0\ind_P^Q k$ is an isomorphism 
            \item $R^i \ind_P^Q k = 0$ for $i \geq 1$.
        \end{enumerate}
    \item If $M$ is a $Q$-module, then the natural map $$\ind_Q^G M \simeq \ind_B^G \res_B^Q M$$ is an isomorphism.
    \end{enumerate}
\end{prop}

 Hence, when $P$ is a reduced parabolic subgroup and $M$ is a $P$-module, by abuse of notation we will often simply write $\ind_B^G M$ instead of $\ind_P^G M$. To compute this induction, it is enough to know the $B$-module structure of $M$.
 However, when $P$ is non-reduced, it is not true that the induction $\ind_{P_{\red}}^P$ preserves the trivial weight. Indeed, $P/P_{\red}$ is not projective but affine: it is isomorphic to $U^0_f$. 

 \begin{ex}
    In $\SL_2$, let $B$ be a Borel subgroup, such that $G/B \simeq \pro^1$. Write the simple root $\theta$, which corresponds to the line bundle $O(2)$ on $\pro^1$. Then the induced module $R^0\ind_B^G(p^r\thet)$ identifies with $H^0(G/B, O(2p^r))$, which has dimension $2p^r +1$.
Now let $$P = \begin{pmatrix}
* & * \\
\alpha_{p^r} & *
\end{pmatrix}.$$
We will see that there is the isomorphism of varieties $G/P \simeq \pro^1$, but that $$R^0\ind_P^G (p^r \theta) \simeq (R^0\ind_B^G \theta)^{(r)},$$ which identifies with $H^0(\pro^1, O(2))^{(r)}$.
\end{ex}

From the reduced case one only keeps the partial injectivity result $R^0 \ind_P^G M \hookrightarrow R^0 \ind_B^GM$.

\begin{prop} \label{injectivity R0}
    Let $P$ be a non-reduced parabolic subgroup of $G$. The natural map $$\pi \colon G/P_{\red} \to G/P$$ is a faithfully flat finite purely inseparable morphism. 
    The pushforward $\pi_*$ is exact and the induced map of $O_{G/P}$-modules $O_{G/P} \to R^0\pi_* O_{G/P_{\red}}$ is injective.
    As a consequence, if $M$ is a $P$-module, the natural map $R^0 \ind_P^G M \to R^0\ind_B^G \res_B^P M$ is injective. 
\end{prop}

\begin{proof}
    The map $\pi$ is a faithfully flat finite purely inseparable morphism by \cite[Part 2, Corollary 7]{HaboushLauritzen1993}.
    Faithful flatness implies that the natural map $O_{G/P}$-modules $O_{G/P} \to \pi_* O_{G/P_{\red}}$ is injective. The pushforward $\pi_*$ is exact because $\pi$ is finite, hence affine.

    Tensoring preserves injectivity so the map $M \otimes_{O_{G/P}} O_{G/P} \to M \otimes_{O_{G/P}} R^0\pi_* O_{G/P_{\red}}$ is injective. By the projection formula, the natural map $$M \otimes_{O_{G/P}} \pi_*O_{G/P_{\red}} \to \pi_*(\pi^* M \otimes_{O_{G/P_{\red}}} O_{G/P_{\red}})$$ is an isomorphism. Hence the map $M \to R^0\pi_* \pi^* M$ is injective, that is, the map $$M \to R^0 \ind_{P_{\red}}^P \res_{P_{\red}}^P M$$ is injective. The conclusion then comes from \Cref{pushforward of structure sheaf} (3).
\end{proof}

In \cite{mehta1985frobenius}, Metha and Ramanathan introduced the theory of Frobenius splitting, which studies and uses the splitting of the Frobenius endomorphism of a variety. 
As we will prove in \Cref{r0}, if $P$ is a parabolic subgroup of width $1$, that is, $P = G_rQ$ with $Q$ reduced, then $G/P \simeq G/P_{\red}$ as varieties (not as $G$-varieties) and the map $G/P_{\red} \to G/P$ is the Frobenius.
The notion of splitting generalizes to the case of any non-reduced parabolic subgroup.

\begin{defi}
    Let $P$ be a parabolic subgroup of $G$, and let $\pi \colon G/P_{\red} \to G/P$ be the natural map. The variety $G/P$ is said to be $\pi$-split if the natural injective map $$O_{G/P} \to \pi_* O_{G/P_{\red}}$$ admits a retraction.
\end{defi}

In the case of non-reduced parabolic subgroups of width $1$, the notion of $\pi$-splitting coincides with the notion of Frobenius splitting of \cite{mehta1985frobenius}.
This property is very helpful because it allows to get injectivity for the map between higher derived induction functors of a $P$-module to $G$ and higher induction from its restriction to $B$.

\begin{prop}
    Let $P$ be a parabolic subgroup of $G$. If $G/P$ is $\pi$-split, then for all $i \geq 0$, the map $R^i \ind_P^G M \to R^i \ind_B^G M$ is injective.
\end{prop}

\begin{proof}
    We get an injective map $M \to \pi_* \pi^* M$ which admits a retraction, hence for all $i \geq 0$, $R^i \ind_P^G M \to R^i \ind_P^G \ind_{P_{\red}}^P \res_{P_{\red}}^P M$ is injective, that is, $R^i \ind_P^G M \to R^i \ind_B^G \res_B^P M$ is injective.
\end{proof}

However injectivity for higher derived induction functors fails in general when $G/P$ is not $\pi_P$-split, which is the case as soon as $P$ has width at least $2$.

\begin{prop}[{\cite[Part 4, Theorem 3]{HaboushLauritzen1993} or \cite{lauritzen1993}}] \label{Frob split}
    Let $P$ be a parabolic subgroup of $G$. The variety $G/P$ is $\pi$-split if and only if $P$ has width at most $1$.
\end{prop}

\paragraph{Quotients of non-reduced parabolic subgroups.} 

\begin{prop} [{\cite[I.9.5]{jantzen2003representations}}] \label{quotient reduced}
    If $H$ is reduced group scheme, then 
    \begin{enumerate}
        \item $H/H_r \simeq H^{(r)}$.
        \item for $r \geq s$, $H_r/H_s \simeq (H^{(s)})_{r-s}$.
    \end{enumerate}
    
    In particular, if $H$ is a reduced reductive or parabolic subgroup of $G$, then 
    \begin{enumerate}
        \item $H/H_r \simeq H^{(r)} \simeq H$.
        \item for $r \geq s$, $H_r/H_s \simeq H_{r-s}$.
    \end{enumerate}
\end{prop}

For non-reduced group scheme, it is not true that $H/H_r \simeq H$. In general, the Frobenius can kill some parts of the group.

\begin{lemma} \label{quotient parabolic}
    Let $f$ be a non-infinitesimal $W$-function with values $r_0 < r_1 < \dots  < r_l$. Then for all $0 \leq k \leq l$,
    \begin{enumerate}
    \item $(P_f)_{r_k} = \overline{P}^k_{r_k} = P^0_{r_0}P^1_{r_1}\dots P^k_{r_k}$. 

    \item $P_f/(P_f)_{r_k} \simeq P^{k+1}_{r_{k+1}-r_k}\dots P^{l}_{r_{l}-r_k}$ (with convention $\infty - n = \infty$). It is the parabolic subgroup associated to the $W$-function $max(0, f-r_k)$. 

    \item $\overline{P}^k/(P_f)_{r_k} \simeq \overline{P}^k/\overline{P}^k_{r_k} \simeq P^k$
    \end{enumerate}
\end{lemma}

\begin{proof}
    \begin{enumerate}
    \item Remember that $(P_f)_{r_k} = P_f \cap G_{r_k}$ by \Cref{kernel of subgroup}. But by definition of $P_f$, $P_f \cap G_{r_k} = P^0_{r_0}P^1_{r_1}\dots P^k_{r_k}$.

    \item $P_f/(P_f)_{r_k} \simeq F_G^{r_k}(P_f) \simeq F_G^{r_k}(P^0_{r_0}\dots P^l_{r_l}) \simeq F_G^{r_k}(P^0_{r_0})\dots F_G^{r_k}(P^l_{r_l}) = P^{k+1}_{r_{k+1}-r_k}\dots P^l_{r_l - r_k}$ by \Cref{quotient reduced}.

    \item This is an immediate consequence of (2) applied to the parabolic subgroup $\overl^k$.
    \end{enumerate}
\end{proof}

\begin{lemma}[{\cite[I.6.2]{jantzen2003representations}}] \label{lemma quotients}
    Let $G$ be any group scheme, let $H$ be a subgroup scheme of $G$ and $N$ be a normal subgroup scheme of $H$. There is a natural isomorphism $$(G/N)/(H/N) \iso G/H.$$ 

\end{lemma}

Beware that even if $P$ is a non-reduced parabolic subgroup, $G/P$ can be isomorphic (as a variety, not as a $G$-variety) to a classical flag variety, when $P$ has width $1$. More generally, any quotient $G/P$ can be computed (as a variety) as the quotient by another parabolic subgroup which has minimal altitude $0$. So, when one is interested in quotients purely as varieties, one can always suppose that the minimal altitude is $0$.  

\begin{coro} \label{r0}
    Let $f$ be a non-infinitesimal $W$-function with values $r_0 < r_1 < \dots  < r_l$.
    Then there are isomorphisms of  varieties (not of $G$-varieties): $$G/P_f \simeq G/(P^1_{r_1-r_0}\dots P^l_{r_l-r_0})\simeq G/P_{f-r_0}$$

    In particular, the case of width $1$ gives for all reduced parabolic subgroup $P$ the isomorphism of varieties $$G/G_rP \simeq G/P.$$ 
\end{coro}

\begin{proof}
    By \Cref{quotient parabolic} (1), $(P_f)_{r_0} = G_{r_0}$ so $$(G/G_{r_0})/(P_f/G_{r_0}) = (G/G_{r_0})/(P_f/(P_f)_{r_0}).$$ Now \Cref{quotient parabolic} (2) and \Cref{quotient reduced} (because $G$ is reduced) give the isomorphism $$(G/G_{r_0})/(P_f/(P_f)_{r_0}) \simeq G/P_{f-r_0}.$$ But by \Cref{lemma quotients}, $$(G/G_{r_0})/(P_f/G_{r_0}) \simeq G/P_f,$$ hence the result.
\end{proof}

Even if in general it is not true that the quotient of two parabolic subgroups is projective, there is the following partial result, which extends \Cref{pushforward of structure sheaf}.

\begin{prop} \label{indk}
    Let $f$ be a non-infinitesimal $W$-function with values $r_0 < r_1 < \dots  < r_l$. Write $P_f = P^0_{r_0}\dots P^l_{r_l}$. Then 
    \begin{enumerate}
    \item More generally, for all $0 \leq k \leq l-1$, there are isomorphisms of varieties $$\overl^k/\overl^{k+1} \simeq P^k/(P^k_{r_k}P^{k+1}) \simeq P^k/P^{k+1}.$$
    In particular $\overl^k/\overl^{k+1} $ is a smooth, projective and irreducible variety.

    \item For all $0 \leq j < k \leq l$, there is an isomorphism of varieties $$\overl^j/\overl^k \simeq P^j/(P^j_{r_j-r_{j-1}}\dots P^{k-1}_{r_{k-1}-r_{j-1}}P^k).$$ In particular, $\overl^j/\overl^k$ is a smooth, projective and irreducible variety.

    \item For all $0 \leq j < k \leq l$, $$\ind_{\overl^k}^{\overl^j} k \simeq k$$ and if $M$ is a $\overl^k$-module, then $$\ind_{\overl^k}^{\overl^j} \res_{\overl^k}^{\overl^j} M \simeq M.$$
    \end{enumerate}
\end{prop}

\begin{proof}
    By \Cref{lemma quotients}, $$\overl^j/\overl^k \simeq (\overl^j/P^0_{r_0}...P^{j-1}_{r_{j-1}})/(\overl^k/P^0_{r_0}...P^{j-1}_{r_{j-1}}).$$ But by \Cref{quotient parabolic}, $$\overl^j/P^0_{r_0}...P^{j-1}_{r_{j-1}} \simeq P^j$$ and $$\overl^k/P^0_{r_0}...P^{j-1}_{r_{j-1}} \simeq P^j_{r_j-r_{j-1}}...P^{k-1}_{r_{k-1}-r_{j-1}}P^k.$$ 

    For (3), the proof is the same as in \Cref{injectivity R0}.
\end{proof}

\begin{prop} \label{Gf reductive}
    Let $f$ be a $W$-function, infinitesimal or not. Then $G_f$ is reductive in the sense that it contains no connected unipotent normal subgroup.
\end{prop}

\begin{proof}
    Write $G_f = G^0_{r_0}...G^l_{r_l}$. We reason by induction on the width of $f$.
    Let $U$ be a connected unipotent normal subgroup of $G_f$. Write the quotient map $\pi \colon G_f \to G_f/G^0_{r_0}...G^{l-1}_{r_{l-1}}$. By \Cref{quotient parabolic}, $$G_f/G^0_{r_0}...G^{l-1}_{r_{l-1}} \simeq G^l_{r_l - r_{l-1}}.$$ By \cite[11.8]{humphreys1967algebraic}, we know that a smooth connected affine group is reductive if and only if its Frobenius kernel contains no non-trivial unipotent normal subgroup. The image $\pi(U)$ is unipotent and normal in $G^l_{r_l - r_{l-1}}$ and $G^l$ is reductive, so $\pi(U)$ is trivial. Hence $U \subset G^0_{r_0}...G_{l-1}^{r_{l-1}}$. But $G^0_{r_0}...G^{l-1}_{r_{l-1}}$ has width $l-1$, so by induction hypothesis $U$ is trivial.
\end{proof}

\paragraph{A formula for induction from non-reduced parabolic subgroups.}

\begin{lemma}[{\cite[I.6.11]{jantzen2003representations}}]\label{quotient}
 Let $L$ be a group scheme, $K$ a subgroup scheme of $L$ and $N$ a normal subgroup scheme of $L$ such that $N \subset K \subset L$ and $H/N$ and $G/N$ are affine. Let $M$ be a $H/N$-module. Denote $\pi \colon G \to G/N$. Then there is an isomorphism of $G$-modules $${\ind_H^G (\pi^*M) \simeq \pi^*(\ind_{H/N}^{G/N} M)}.$$
\end{lemma}

For a non-reduced parabolic subgroup of width $1$, there is the formula $$\ind_{G_rP}^G((F^r)^*M) = (F^r)^* \ind_P^G(M).$$ By reducing to the case of pairs of parabolic subgroups of relative width $1$, we express induction from a non-reduced parabolic subgroups using successive twisted inductions from classical parabolic subgroups.

\begin{prop} \label{induction}
    Let $f$ be a non-infinitesimal $W$-function of width $l$. Write $P_f = P^0_{r_0}P^1_{r_1}\dots P^l_{r_l}$ with $r_0 < r_1 < \dots  < r_l$ the values of $f$. 
    Let $M$ be a $P^l$-module. Let $r \geq r_{l-1}$ (with the convention $r_{-1} = 0$). 
    Then there is an isomorphism of $G$-modules $$\ind_{P_f}^{G}  (F^r_{P_f})^* (M) \simeq (F_{P^0}^{r_0})^*  \ind_{P^1}^{P^0}  (F_{P^1}^{r_1-r_0})^*  \ind_{P^2}^{P^1}  \dots  (F_{P^{l-1}}^{r_{l-1} - r_{l-2}})^*  \ind_{P^{l}}^{P^{l-1}}  (F_{P^{l}}^{r - r_{l-1}})^* (M)$$

    In particular, the case $r = r_{l-1}$ gives the isomorphism of $G$-modules $$\ind_{P_f}^{G} ((F^{r_{l-1}}_{P_f})^*M) \simeq (F_{P^0}^{r_0})^* \ind_{P^1}^{P^0} (F_{P^1}^{r_1-r_0})^* \ind_{P^2}^{P^1} \dots  (F_{P^{l-1}}^{r_{l-1} - r_{l-2}})^* \ind_{P^{l}}^{P^{l-1}} (M)$$
\end{prop}

\begin{proof}
    We prove it by induction on the width $l$.
    \begin{itemize}
    \item If $l = 0$, we have to prove the isomorphism $$\ind_G^G((F^r_{G})^* M) \simeq (F_G^{r-r_{-1}})^* M,$$ which is true.

    \item If $l = 1$, we have to prove the isomorphism $$\ind_{\overl^1}^{\overl^0}((F^r_{\overl^0})^* M) \simeq (F^{r_0}_{\overl^0})^* \ind_{P^1}^{P^0}((F_{P^1}^{r-r_0})^*M).$$
    Remark that by \Cref{quotient parabolic}, $\overl^1_{r_0} \simeq P^0_{r_0} \simeq \overl^0_{r_0}$. Thus, by \Cref{quotient}, $$\ind_{\overl^1}^{\overl^0}((F^r_{\overl^0})^* M) \simeq (F^{r_0}_{P^0})^* \ind_{\overl^1/\overl^1_{r_0}}^{\overl^0/\overl^0_{r_0}}((F_{P^1}^{r-r_0})^*M)$$
    
    Then, by \Cref{quotient parabolic} again, $\overl^0/\overl^0_{r_0} \simeq P^0$ and $\overl^1/\overl^1_{r_0} \simeq P^1$. Hence,
    
    $$\ind_{\overl^1/\overl^1_{r_0}}^{\overl^0/\overl^0_{r_0}}((F_{P^1}^{r-r_0})^*M) \simeq (\ind_{P^1}^{P^0}((F_{P^1}^{r-r_0})^*M)$$

    This proves the expected isomorphism.

    \item Now suppose that the formula is verified in length $l-1$.
    The induction from $P_f$ to $G$ decomposes as $$\ind_{P_f}^G = \ind_{\overl^{l-1}}^{G} \circ \ind_{\overl^l}^{\overl^{l-1}}.$$
    By \Cref{quotient parabolic} $$\overl^{k+1}_{r_k} \simeq P^0_{r_0}P^1_{r_1}\dots P^k_{r_k} \simeq \overline{P}^k_k,$$ and moreover $\overl^l/\overl^l_{r_{l-1}} \simeq P^l$ and $\overl^{l-1}/\overl^{l-1}_{r_{l-1}} \simeq P^{l-1}$. Thus, by \Cref{quotient}, we get isomorphisms
    $$\ind_{\overl^{l}}^{\overl^{l-1}} ((F^r_{\overl^l})^* (M)) \simeq (F^{r_{l-1}}_{\overl^{l-1}})^* \ind_{\overl^{l}/\overl^{l}_{r_{l-1}}}^{\overl^{l-1}/\overl^{l-1}_{r_{l-1}}} ((F^{r-r_{l-1}}_{\overl^l})^* (M))$$
    and $$(F^{r_{l-1}}_{\overl^{l-1}})^* \ind_{\overl^{l}/\overl^{l}_{r_{l-1}}}^{\overl^{l-1}/\overl^{l-1}_{r_{l-1}}} ((F^{r-r_{l-1}}_{\overl^l})^* (M)) \simeq (F^{r_{l-1}}_{\overl^{l-1}})^*(\ind_{P^{l}}^{P^{l-1}}((F^{r-r_{l-1}}_{\overl^l})^* (M))).$$

    Denote $N = \ind_{\overl^{l}}^{\overl^{l-1}} ((F^{r-r_{l-1}}_{\overl^l})^* (M))$, such that $\ind_{P_f}^{G}  (F^r_{P_f})^* (M) \simeq (F_{\overl^{l-1}}^{r_{l-1}})^* (N)$.
    The parabolic subgroup $\overl^{l-1}$ has width $l-1$ so by induction hypothesis, $$\ind_{\overl^{l-1}}^G ((F^{r_{l-1}}_{\overl^{l-1}})^* (N)) \simeq (F_{P^0}^{r_0})^* \ind_{P^1}^{P^0} (F_{P^1}^{r_1-r_0})^* \dots (F_{P^{l-2}}^{r_{l-2} - r_{l-3}})^* \ind_{P^{l-1}}^{P^{l-2}} ((F_{P^{l-1}}^{r_{l-1}-r_{l-2}}(N))).$$

    Finally, we get the expected formula
    $$\ind_{P_f}^{G}  (F^r_{P_f})^* (M) \simeq (F_{P^0}^{r_0})^*  \ind_{P^1}^{P^0}  (F_{P^1}^{r_1-r_0})^*  \ind_{P^2}^{P^1}  \dots (F_{P^{l-1}}^{r_{l-1}-r_{l-2}})^* \ind_{\overl^l}^{\overl^{l-1}} F^{r-r_{l-1}}_{P^l}(M).$$
    \end{itemize}
\end{proof}

\begin{rk}
    Given that when $P \subset Q$ are reduced parabolic subgroups, $\ind_P^Q = \ind_B^Q \circ \res_B^P$, the formula of \Cref{induction} can be computed at each step by looking only at the restriction to the $B$-action. 
\end{rk}

\subsection{Filtration and cohomology of the tangent bundle} \label{section filtration and cohomology}

In this section, we focus on the tangent bundle of the variety $X = G/P$. When $P$ is a reduced parabolic subgroup, it has been proven, by Bott in characteristic $0$ (\cite{bott1957homogeneous}) and by Demazure in characteristic $p$ (\cite{Demazure1977}), that for all $i > 0$, $$H^i(X, T_X) = 0.$$
In particular, deformations of $X$ are rigid.
Those results fail in general when $P$ is non-reduced (later in \Cref{section final conclusion}, we notably show that for most non-reduced $P$, deformations of $X$ are completely obstructed, and in particular that $H^2(X, T_X) \neq 0$). However, those spaces $H^i(X, T_X)$ are naturally endowed with a $G$-module structure, and it is possible to get results concerning their $G$-cohomology. Let $\overline{P}$ denote the smallest reduced subgroup of $G$ which contains $P$, and consider the subset $I$ of simple roots which generates $\overl$. In \Cref{conditions fulfilled}, we show that if $|I| \leq |S|-2$ or if $|I| = |S|-1$ and $I$ is not connected as a subgraph of the Dynkin diagram of $G$, then $H^1(X, T_X)$ is a trivial $G$-module, and for all $i \geq 2$, $$H^i(G, H^0(X, T_X)) = 0.$$

Our main tool to study the tangent bundle $T_X$ consists in building a filtration of $T_X$ whose factors are Frobenius twists of tangent bundles of quotients of classical parabolic subgroups (see \Cref{filtration}). 

Moreover, the variety $X = G/P$ sits into the homotopy pullback square
\[\begin{tikzcd}
	X & * \\
	BP & BG
	\arrow["\pi", from=1-1, to=1-2]
	\arrow["q", from=1-1, to=2-1]
	\arrow["\lrcorner"{anchor=center, pos=0.125}, draw=none, from=1-1, to=2-2]
	\arrow["e", from=1-2, to=2-2]
	\arrow[from=2-1, to=2-2]
\end{tikzcd}\]
where $BP$ and $BG$ are the classifying stacks in the fppf topology.
The pullback $q^*$ induces an isomorphism between $P$-modules and $G$-equivariant vector bundles on $X$. In particular, the tangent bundle $T_X$ identifies with the tangent complex $T_{BP \sto BG}$, which is a $P$-module. Moreover, the base change formula $\pi_*q^* \simeq e^* \ind_P^G$ shows that the cohomology groups $H^i(X, T_X)$ lift to the $G$-modules $R^i\ind_P^G T_{BP \sto BG}$. We will constantly use this identification in the rest of this work. In particular, the filtration of $T_X$ that we build can be seen as a $P$-module filtration of $T_{BP \sto BG}$.

\subsubsection{Warm-up: the weights of the tangent bundle} \label{section weight tangentn}

Before studying the tangent bundle of $X$ as a $P_f$-module, as an introduction we study its $T$-module structure, that is, its weights. It will be a consequence of its $P_f$-module decomposition in \Cref{section filtration tangent}, but we prove it here independently for convenience.

Write $r_0 < r_1 < \dots < r_l$ the values of $f$, and $I_l \subset \dots \subset I_0$ the canonical filtration of the set of simple roots.. Write $R_k$ for the subset of positive roots generated by $I_k$. 
The weights of $T_X$ are simply twists of the positive roots: they are the $p^{r_k}\la$ for $\la \in R_k/R_{k+1}$ and $0 \leq k \leq l-1$.

Non-reduced parabolic subgroups are built from Frobenius kernels of root subgroups. Hence, the idea is to begin by computing the tangent complex of the elementary blocks $T_{B\alpha_{p^r} \sto B\G_a}$, while keeping track of the data of the natural action of the multiplicative group $\G_m$. Then, one essentially has to sum up those elementary blocks to get the total tangent complex.

To begin with, remember the following common result, whose proof can be found for instance in \cite{halpern2020moduli}.

\begin{prop} \label{tangent BG}
    Let $G$ be a smooth group scheme. Then the tangent complex of the classifying stack of $G$ is given by
    $$T_{BG} \simeq T_{G,1}[1],$$ endowed with the adjoint action of $G$.
\end{prop}

\begin{lemma} \label{decomposition TP}
    The tangent complex at $1$ of the parabolic subgroup $P_f$ is given, as a $T$-module, by $$T_{P, 1} \simeq T_{P^l,1} \oplus \bigoplus_{\substack{0 \leq k \leq l-1 \\ \la \in R_k \backslash R_{k+1}}} T_{(U_{\la})_{r_k}, 1}.$$
\end{lemma}

\begin{proof}
    As vector spaces, $T_{P,1} \simeq T_{U^0_f,1} \oplus T_{P^l,1}$ because $P \simeq (U^0_f)^+ \times P_{\red}$ by \Cref{product decomposition}. By definition $$(U^0_f)^+ = \prod_{\substack{0 \leq k \leq l-1 \\ \la \in R_k \backslash R_{k+1}}} (U_{\la})^+_{(r_k)},$$ hence the decomposition.
    Moreover, this decomposition is compatible with the $T$-action by conjugation because $T$ stabilizes the root subgroups. 
\end{proof}

Remember that actions of the multiplicative group $\G_m$ on algebras correspond to gradings, that we will call ``weights'', and that $\G_m$ acts on the additive group $\G_a = \spec k[t]$ by giving weight $1$ to $t$ (hence $t^{p^r}$ has weight $p^r$). The following proposition expresses the $\G_m$-action on Frobenius kernels of the additive group. Part of those computations can also be found in \cite[4.2]{antieau2021counterexampleshochschildkostantrosenbergcharacteristicp}.

\begin{coro} \label{tangent alphap}
    The natural action of $\G_m$ on $\G_a$ induces an action of $\G_m$ on $\alpha_{p^r}$, $\T_{\alpha_{p^r},1}$ and on $\T_{\alpha_{p^r} \to \G_a, 1}$. Then the following assertions hold: 
    \begin{enumerate}

    \item $T_{\alpha_{p^r},1} \simeq k[-1] \oplus k[0]$, and the $\G_m$-action has weight $p^r$ in degree $-1$ and weight $1$ in degree $0$.

    \item $T_{\alpha_{p^r} \sto \G_a, 1} = k[-1]$, and the $\G_m$-action has weight $p^r$ in degree $-1$.
    
    \item $T_{B\alpha_{p^r} \sto B\G_a} \simeq k$, and the $\G_m$-action has weight $p^r$.
    \end{enumerate}
\end{coro}

\begin{proof}
The action of $\G_m$ on $\G_a$ restricts to an action of $\G_m$ on $\alpha_{p^r} = (\G_a)_{(r)}$  and thus by differentiating, an action of $\G_m$ on $\T_{(\G_a)_{(r)},1}$. The closed immersion $\alpha_{p^r} \to \G_a$ is $\G_m$-equivariant so it also induces an action on $\T_{\alpha_{p^r} \sto \G_a, 1}$.

The Frobenius kernel $\alpha_{p^r}$ has cotangent complex $L_{\alpha_{p^r}} = [(t^{p^r})/(t^{p^r})^2 \to k[t]/t^{p^r} dt]$ where $(t^{p^r})/(t^{p^r})^2$ is in degree $0$ and $k[t]/t^{p^r} dt$ in degree $1$.
The closed immersion $(\G_a)_{(r)} \to \G_a$, which identifies with the closed immersion $\spec k[t]/t^{p^r} \to \spec k[t]$, has cotangent complex $L_{(\G_a)_{(r)} \rightarrow \G_a} \simeq (t^{p^r})/(t^{p^r})^2[1]$. 
When one restricts at $1$, one gets $(t^{p^r})/((t^{p^r})^2) \otimes_{k[t]/t^{p^r}} k \simeq k$ and $k[t]/t^{p^r}dt \otimes_{k[t]/t^{p^r}} k \simeq k$  and the differential becomes $0$ (the differential is restricted to constants). So, restricting at $1$ gives $L_{\alpha_{p^r},1} \simeq k[0] \oplus k[1]$  and $L_{(\G_a)_{(r)} \rightarrow \G_a, 1} \simeq k[1]$.

Now $t^{p^r}$ has weight $p^r$ under the action of $\G_m$ on $\G_a = k[t]$ so $\G_m$ acts on $\T_{\alpha_{p^r},1} \simeq k[-1] \oplus k[0]$ with weight $p^r$ in degree $-1$ and weight $1$ in degree $0$, and $\G_m$ acts on $T_{\alpha_p \to \G_a, 1} = k[-1]$ with weight $p^r$. 
Hence $T_{\alpha_{p^r},1} \simeq k[-1] \oplus k[0]$ and the $\G_m$-action has weight $p^r$ in degree $-1$ and weight $1$ in degree $0$; and $\T_{\alpha_{p^r} \to \G_a, 1} \simeq k[-1]$ and the $\G_m$-action has weight $p^r$ in degree $-1$.
Finally by \Cref{tangent BG}, $\T_{B\alpha_{p^r} \to B\G_a} \simeq \T_{\alpha_{p^r} \to \G_a, 1}[1] \simeq k[0]$ and the $\G_m$-action has weight $p^r$.
\end{proof}

These computations of tangent for $\alpha_{p^r}$ transfer to analogous results for root subgroups with action of the torus $T$.

\begin{coro} \label{weight non-reduced root subgroups}
    Let $\la$ be a root of $G$. The torus $T$ acts by conjugation on the root subgroup $U_{\la}$.
    Then the following $T$-module isomorphisms hold:
    \begin{enumerate}

    \item $T_{(U_{\la})_{r}, 1} \simeq p^r\la[-1] \oplus \la[0]$ 
    
    \item $T_{B(U_{\la})_{r} \rightarrow BU_{\la}} \simeq p^r\lambda$
    \end{enumerate}
\end{coro}

\begin{proof}
    The torus $T$ acts with weight $\la$ on $U_{\la}$. Through the isomorphisms $T \simeq (\G_m)^n$ and $U_{\la} \simeq \G_a$, it means that the action of $T$ on $U_{\la}$ is the composition of the natural action of $\G_m$ on $\G_a$ and of $\la \colon (\G_m)^n \to \G_m$. 
    The result then follows from \Cref{tangent alphap}.
\end{proof}

\begin{rk}    
    The previous computations admit the following generalization.
    Let $H$ be a reduced group scheme with $\F_p$-structure, endowed with an action of another group scheme $K$ which commutes with the Frobenius. Let $r \in \N$ and write $\frak{h} := T_{H,1}$. 
    The $K$-action on $H$ equips $T_{H_r, 1}$ and $T_{BH_r \sto BG}$ with a $K$-module structure which verifies:
    \begin{enumerate}
\item $T_{H_r, 1} \simeq \frak{h}^{(r)}[-1] \oplus \frak{h}[0]$. 

    \item $T_{BH_r \to BH} \simeq \frak{h}^{(r)}[0]$.
    \end{enumerate}
    It can be proven using the following cartesian square

\[\begin{tikzcd}
	{H_r} & {*} \\
	H & H
	\arrow[from=1-1, to=1-2]
	\arrow["i"', from=1-1, to=2-1]
	\arrow["\lrcorner"{anchor=center, pos=0.125}, draw=none, from=1-1, to=2-2]
	\arrow[from=1-2, to=2-2]
	\arrow["{F_H^r}"', from=2-1, to=2-2]
\end{tikzcd}\]
which shows that $T_{H_r} \simeq i^* T_{F_H^r}$. 
\end{rk}

Finally, the $T$-module structure on the tangent $T_X$ simply decomposes as the sum of the $T$-modules structure on root subgroups.

\begin{prop} \label{weights}
    The set of weights of the tangent $T_X$ is described as $$\{p^{r_k}\la, \la \in R_k \backslash R_{k+1}, 0 \leq k \leq l-1\}.$$ 
\end{prop}

\begin{proof}
    Consider the fiber sequence $$T_{BP \rightarrow BG} \to T_{BP} \to (Bi)^*T_{BG}$$ where $Bi \colon BP \to BG$.
The group $G$ is smooth so $\T_{BG}$ identifies with $T_{G,1}[1]$ endowed with the adjoint action of $G$ by \Cref{tangent BG}. Thus $(Bi)^*\T_{BG}$ identifies with $T_{G,1}[1]$ endowed the restriction of the adjoint action of $G$ to $P$. As a $T$-module, it decomposes as $$T_{P^l, 1}[1] \oplus \bigoplus_{\substack{0 \leq k \leq l-1 \\ \la \in R_k \backslash R_{k+1}}} T_{U_{\la},1}[1] \simeq T_{BP^l} \oplus \bigoplus_{\substack{0 \leq k \leq l-1 \\\la \in R_k \backslash R_{k+1}}} T_{BU_{\la}}.$$
The map $\T_{BP} \to (Bi)^*\T_{BG}$ decomposes in maps $T_{BP_{\red}} \to T_{BP_{\red}}$ and $T_{B(U_{\la})_{(r_k)}} \to T_{BU_{\la}}$ for all $0 \leq k \leq l-1$ and $\la \in R_k \backslash R_{k+1}$.
So, by \Cref{decomposition TP}, as a $T$-module $T_{BP \rightarrow BG}$ admits a decomposition $$T_{BP_f \rightarrow BG} \simeq T_{BP^l \sto BP^l} \oplus \bigoplus_{\substack{0 \leq k \leq l-1 \\ \la \in R_k \backslash R_{k+1}}} T_{B(U_{\la})_{(r_k)} \sto BU_{\la}}$$
Then, by \Cref{weight non-reduced root subgroups}, it implies the isomorphism of $T$-modules:

$$T_{BP \rightarrow BG} \simeq \bigoplus_{\substack{0 \leq k \leq l-1 \\ \la \in R_k \backslash R_{k+1}}} p^{r_k} \la$$
\end{proof}

\subsubsection{A filtration of the tangent bundle} \label{section filtration tangent}

For a classical reduced parabolic subgroup $P$ of $G$, the tangent space $T_{G/P,1}$ of the smooth projective variety $G/P$ is isomorphic, as a $P$-module, to the quotient $\g/\p$. For a non-reduced parabolic subgroup of width $1$, that is, $P = G_rQ$ with $Q$ a reduced parabolic subgroup, then the tangent $T_X \simeq T_{B(G_rQ) \sto BG}$ is isomorphic, as a $P$-module, to the $r$-th Frobenius twist of the tangent $\g/\qfrak$ of $G/Q$.

Now, and for the rest of this section, let $f$ be a non-infinitesimal $W$-function having values $r_0 < \dots < r_l$. Write $P_f = P^0_{r_0} \dots P^l_{r_l}$ and $\overl^k = P^0_{r_0} \dots P^{k-1}_{r_{k-1}}P^k$. Denote by $\p^k$ the tangent Lie algebra of the classical parabolic subgroup $P^k$.
For parabolic subgroups of higher width, the situation becomes more intricate. However, we can use the filtration $$P_f = \overl^l \subset \dots \subset \overl^1 \subset \overl^0 = G$$ (which corresponds to a chain $X = G/P_f \twoheadrightarrow G/\overl^{l-1} \twoheadrightarrow \dots \twoheadrightarrow G/\overl^1 \twoheadrightarrow G/\overl^0 = *$). Each intermediary step can be reduced to a width-$1$ situation. This filtration induces a filtration of the tangent $T_X$, as a $P_f$-module, whose factors are Frobenius twists of classical quotients $\p^k/\p^{k+1}$.

    Our final goal is to study the induction of the tangent $\ind_{P_f}^G T_{BP_f \to BG}$. Using the formula for induction that we established in \Cref{induction}, we deduce a formula for induction for the factors of the filtration of the tangent.

Let us begin by dealing with the intermediary steps of the above filtration, which corresponds to a width-$1$ situation.

\begin{prop} \label{twisted tangent}
     Let $0 \leq k \leq l-1$.
    There is an isomorphism of $\overl^{k+1}$-modules:
    
    $$T_{B\overl^{k+1} \rightarrow B\overl^k} \simeq (F_{\overl^{k+1}}^{r_k})^* T_{BP^{k+1} \rightarrow BP^k} \simeq (F_{\overl^{k+1}}^{r_k})^* (\p^k/\p^{k+1})$$
\end{prop}

\begin{proof}
    The Frobenius of $\overl^{k+1}$ is the restriction of the Frobenius of $\overl^k$, hence the following cartesian diagram:
\[\begin{tikzcd}
	{\overl^{k+1}} & {\overl^{k}} \\
	{P^{k+1}} & {P^k}
	\arrow["i", from=1-1, to=1-2]
	\arrow["{F_{\overl^{k+1}}^{r_k}}"', from=1-1, to=2-1]
	\arrow["\lrcorner"{anchor=center, pos=0.125}, draw=none, from=1-1, to=2-2]
	\arrow["{F_{\overl^k}^{r_k}}", from=1-2, to=2-2]
	\arrow["i"', from=2-1, to=2-2]
\end{tikzcd}\]

    Hence, $T_{B\overl^{k+1} \rightarrow B\overl^k} \simeq (F_{\overl^{k+1}}^{r_k})^* T_{BP^{k+1} \rightarrow BP^k}$.
\end{proof}

\begin{prop} \label{filtration}

    There is a chain of surjections of $P$-modules $$T_{BP_f \sto BG} \twoheadrightarrow \res_{P_f}^{\overl^{l-1}} T_{B\overl^{l-1} \sto BG} \twoheadrightarrow \dots \twoheadrightarrow \res_{P_f}^{\overl^{1}} T_{B\overl^1 \sto BG} \twoheadrightarrow \res_{P_f}^{\overl^{0}} T_{B\overl^0 \sto BG} = 0$$
    or dually a filtration of $P$-modules 
    $$0 = T_{BP_f \sto B\overl^l} \subset T_{BP_f \sto \overl^{l-1} } \subset \dots \subset T_{BP_f \sto B\overl^1} \subset T_{BP_f \sto BG}$$
    
    whose factors are 
    $$\res_P^{\overl^{k+1}} (F^{r_k}_{\overl^{k+1}})^* (\frak{p}^k/\frak{p}^{k+1})$$ for $0 \leq k \leq l-1$.
\end{prop}

\begin{proof}
    Consider the filtration $P = \overl^l \subset \overl^{l-1} \subset \dots \subset \overl^1 \subset \overl^0 = G$.
    For all $0 \leq k \leq l-1$ the chain $P^{k+1} \subset P^k \subset G$ gives the fiber sequence $$T_{B\overl^{k+1} \sto B\overl^k} \to T_{B\overl^{k+1} \sto BG} \to \res_{\overl^{k+1}}^{\overl^k} T_{B\overl^k \to BG}.$$ By \Cref{indk} all those terms are concentrated in degree $0$ so we get a short exact sequence of $P$-modules: $$0 \to T_{B\overl^{k+1} \sto B\overl^k} \to T_{B\overl^{k+1} \sto BG} \to \res_{\overl^{k+1}}^{\overl^k} T_{B\overl^k \sto BG} \to 0.$$
    But $T_{B\overl^{k+1} \sto B\overl^k} \simeq (F^{r_k}_{\overl^{k+1}})^* (\frak{p}^k/\frak{p}^{k+1})$ by \Cref{twisted tangent}.

\end{proof}

\begin{notation}
For $0 \leq k \leq l-1$, we will denote $$T^k_f = \res_P^{\overl^{k+1}} (F^{r_k}_{\overl^{k+1}})^* (\frak{p}^k/\frak{p}^{k+1}),$$ such that $T_{BP_f \rightarrow BG}$ admits a filtration whose factors are the $(T_k)_{0 \leq k \leq l-1}$.
\end{notation}

This filtration is very satisfying because we know how to compute the induction of its factors from $P$ to $G$.

\begin{prop} \label{new induction}
    For all $0 \leq k \leq l-1$, there is an isomorphism of $G$-modules
$$\ind_P^G (T^k_f) \simeq (F_{P^0}^{r_0})^*  \ind_{P^1}^{P^0}  (F_{P^1}^{r_1-r_0})^*  \ind_{P^2}^{P^1}  \dots (F_{P^{k}}^{r_{k} - r_{k-1}})^*  \ind_{P^{k+1}}^{P^{k}} (\frak{p}^k/\frak{p}^{k+1}).$$
\end{prop}

\begin{proof}
    By \Cref{indk}, $$\ind_P^{\overl^{k+1}} \res_P^{\overl^{k+1}} (F^{r_k}_{\overl^{k+1}})^* (\frak{p}^k/\frak{p}^{k+1}) \simeq (F^{r_k}_{\overl^{k+1}})^* (\frak{p}^k/\frak{p}^{k+1}).$$
    Then, by \Cref{induction}, $$\ind_{\overl^{k+1}}^G(F^{r_k}_{\overl^{k+1}})^* (\frak{p}^k/\frak{p}^{k+1})) \simeq (F_{P^0}^{r_0})^*  \ind_{P^1}^{P^0}  (F_{P^1}^{r_1-r_0})^*  \ind_{P^2}^{P^1}  \dots (F_{P^{k}}^{r_{k} - r_{k-1}})^*  \ind_{P^{k+1}}^{P^{k}} (\frak{p}^k/\frak{p}^{k+1}).$$
\end{proof}

\begin{rk}
    Using \Cref{induction formula}, we can also get a filtration of $T_{BP_f \sto BG}$ whose factors are the simple $P_f$-modules of the form $L_{P_f}(p^{r_k} \lambda)$ for $0 \leq k \leq l-1$ and $\lambda \in R_k \backslash R_{k+1}$. It is a finer filtration that the filtration of \Cref{filtration}. We can also compute the induction of these factors from $P_f$ to $G$.
\end{rk}

\subsubsection{Interlude: classical results about induction} \label{section interlude}

The goal of this section is to record some classical results to study the structure of $G$-module on inductions of $B$-modules. It will be useful to study the $G$-modules $H^0(X, T_X)$ and $H^1(X, T_X)$.

In characteristic $0$, the picture is simple thanks to the Borel-Weil-Bott theorem. First, simple $G$-modules corresponds to the $H^0(\la)$ for $\la$ a dominant weight. Second, there is never higher $G$-cohomology: for all $G$-module $M$, for all $i \geq 1$, $H^i(G, M) = 0$. One reason is that all $G$-modules admit a good filtration, that is, a filtration whose factors are simple modules $H^0(\la)$ for $\la$ a dominant weight, and that higher $G$-cohomology of those factors vanishes: $H^i(G, H^0(\la)) = 0$ for all $i \geq 1$.
Third, all derived inductions of weights $H^i(\la)$ either vanish or can be expressed as $H^0(\mu)$ where $\mu$ is a dual weight.

In positive characteristic, all those properties fail in general. It is still true that $H^i(G, H^0(\la))$ vanishes for all dominant weight $\la$ and all $i \geq 1$. However, it is not true in general that a $G$-module $M$ admits a good filtration, because simple $G$-modules do not correspond to the $H^0(\la)$ anymore in general (see \Cref{section classification simple modules} for more details). Hence, higher $G$-cohomology does not vanish in general. Moreover, while Kempf vanishing is still verified, $H^i(\la)$ admits in general no simpler expression. For instance, there are examples of weights $\la$, such that $H^i(\la)$ and $H^j(\la)$ do not vanish for two indices $i \neq j$.

Thus, we will have to proceed to a careful study of the $G$-module structure of induced modules. We begin by showing that $G$-triviality can be deduced from $G$-triviality of the factors of a filtration (\Cref{reconstruct}). Then, when the Dynkin diagram of $G$ is connected, we show that the induction of a weight which is not generated by a sufficient number of simple roots is a trivial $G$-module (\Cref{s-i}).

Let us begin by recalling the important fact that $G$-cohomology of $H^0(\la)$ is trivial.

\begin{lemma}[{\cite[II.4.13]{jantzen2003representations}}] \label{G-cohomology of H0(lambda)}
    If $\la$ is a non-zero dominant weight of $G$, then $$H^{\bullet}(G, H^0(\la)) = 0.$$
\end{lemma}

\begin{lemma}[{\cite[II.4.11]{jantzen2003representations}}] \label{no extension}
    The $G$-cohomology of the trivial $G$-module $k$ is trivial:
    \begin{enumerate}
    \item $H^0(G, k) \simeq k$
    \item for all $i \geq 1$, $H^i(G, k) = 0$.
    \end{enumerate}
\end{lemma}

\begin{rk}
    There exists some very partial results concerning vanishing of higher cohomology of $G$-modules. Let $G$ be a simple, simply connected and connected algebraic group, and let $M$ be a $G$-module such that $\dim M \leq p$. In \cite{jantzen1997low}, Jantzen proves that $H^1(G, M) = 0$. In \cite{mcninch2000secondcohomologysmallirreducible}, McNinch proves that $H^2(G, V) \neq 0$ if and only if $M$ has a composition factor isomorphic to a Frobenius twist of $\g$.
\end{rk}

The following lemma is inspired by \cite[Lemma 6.6]{zdanowicz2018arithmetically}.
If a $P$-module admits a filtration such that the induction of all its factors are trivial $G$-modules, then the total module also induces to a trivial $G$-module.

\begin{lemma} \label{reconstruct}
Let $k \geq 0$. Let $P$ be a parabolic subgroup of $G$. Let $M$ be a $P$-module which admits a filtration $$0 \subset M_l \subset \dots \subset M_1 \subset M_0 = M$$ such that for all $0 \leq i \leq l-1$, $R^k \indp (M_i/M_{i+1})$ is a trivial $G$-module (resp. $0$). Then $R^k \indp M$ is a trivial $G$-module (resp. $0$). 
\end{lemma}

\begin{proof}
    We reason by induction on the filtration.
    The short exact sequence $$0 \to M_{i+1} \to M_i \to M_i/M_{i+1}$$ gives a long exact sequence of $G$-modules which contains $$\dots \to R^k \indp(M_{i+1}) \xlongrightarrow{f} R^k \indp(M_i) \xlongrightarrow{g} R^k \indp (M_i/M_{i+1}) \to \dots.$$ By induction hypothesis $R^k \indp(M_{i+1})$ is a trivial $G$-module, and we made the assumption that $R^k \indp (M_i/M_{i+1})$ is also a trivial $G$-module. 
    But $\Ext^1_G(k, k) = H^1(G,k) = 0$ by \Cref{no extension} so there is no non-trivial extension of trivial $G$-modules (we apply this to the short exact sequence $0 \to \im f \to R^k \indp (M_i/M_{i+1}) \to \im g$ where $\im(f) = G/\kernel f$ and $\im(g)$ are also trivial $G$-module as quotient and submodule of trivial $G$-modules). So $R^k \indp (M_i)$ is a trivial $G$-module.
    The proof is similar (and simpler) for the case of vanishing.
\end{proof}

\begin{lemma} \label{filtration of B-modules}
Let $M$ be a finite-dimensional $B$-module. Then $M$ admits a filtration with $B$-submodules $$0 \subset M_l \subset \dots \subset M_1 \subset M_0 = M$$ such that $M_i/M_{i+1}$ is $1$-dimensional for all $0 \leq i \leq l-1$ (that is, a weight).
\end{lemma}

\begin{proof}
Take a Jordan-Hölder filtration of $M$, such that factors are simple $B$-modules. Now simple $B$-modules are $1$-dimensional, hence the result.
\end{proof}

\begin{coro} \label{induction of module}
    Let $M$ be a finite-dimensional $B$-module such that for every weight $\la$ of $M$, $H_G^i(\la)$ is a trivial $G$-module. Then $H_G^i(M)$ is a trivial $G$-module.

\end{coro}

\begin{proof}
    This is an immediate consequence of \Cref{reconstruct} and \Cref{filtration of B-modules}.
\end{proof}

When the group $G$ has a connected Dynkin diagram, a non-zero dominant weight has to be written on all the simple roots, in the sense of the following lemma.

\begin{lemma} \label{condition dominant}
    Let $\la$ be a dominant weight generated by the simple roots. Then $\la \geq 0$, which means that the decomposition of $\la$ over the simple roots has only non-negative coefficients.
    Moreover, if $G$ has a connected Dynkin diagram and $\la \neq 0$, then all coefficients in the decomposition are positive. In this case, we will say that $\la$ is $S$-positive.
\end{lemma}

\begin{proof}

    We prove the second assertion of the lemma. Denote by $S_1$ the subset of simple roots such that $\la_{\alpha} > 0$ and by $S_2$ the subset of simple roots such that $\la_{\alpha} = 0$. Let $\beta \in S_2$.  We have $$\<\la, \beta^{\vee}\> = \sum_{\alpha \in S_1} \la_{\alpha} \<\alpha, \beta^{\vee}\>.$$ 
    For all $\alpha \in S_1$, $\alpha \neq \beta$ so $\<\alpha, \beta^{\vee}\> \leq 0$. But $\la$ is dominant so $\<\la, \beta\> \geq 0$ so for all $\alpha \in S_1$, $\<\alpha, \beta\> = 0$.
    Thus for all $\alpha \in S_1$ and $\beta \in S_2$, $\<\alpha, \beta\> = 0$. But the Dynkin diagram of $G$ is connected so either $S_1 = \emptyset$, which means that $\la = 0$, either $S_2 = \emptyset$, which means that $\la$ is $S$-positive. 
\end{proof}

\begin{notation}
    Let $\alpha \in S$ and $\la \in X(T)$.
    Write the dot action $$\s \la := \la - (\<\la, \alpha^{\vee}\> +1) \alpha.$$
\end{notation}

By reducing to the case of $\SL_2$, the weights of $H^1_{P_{\alpha}}(\la)$ can be explicitly described.

\begin{lemma}[{\cite[II.5.2]{jantzen2003representations}}] \label{induction minimal parabolic}
    Let $\la \in X(T)$ and $\alpha \in S$. Write $r = \<\la, \alpha^{\vee}\>$. Remember that $P^{\alpha}$ denote the parabolic subgroup generated by the singleton $\{\al\} \subset S$. Then the following assertions hold:
    \begin{enumerate}
        \item if $r \geq 0$, then for all $i \geq 1$, $H^i_{P^{\alpha}}(\la) = 0$.

    \item if $r = -1$, then for all $i \geq 0$, $H^i_{P^{\alpha}}(\la) = 0$.

    \item if $r \leq -2$, then for all $i \neq 1$, $H^i_{P^{\alpha}}(\la) = 0$, and the weights of $H^1_{P^{\alpha}}(\la)$ are $$\s \la, \s \la - \alpha, \dots, \s \la - (-r-2) \alpha.$$
    \end{enumerate}
    In particular, the weights of $H^1_{P^{\alpha}}(\la)$ are generated by $\la$ and $\alpha$.
\end{lemma}

Then, still in the case of a group with connected Dynkin diagram, it allows us to show that when weights are not generated by enough simple roots, their inductions are trivial $G$-modules.

\begin{coro} \label{s-i}
    Suppose $G$ has a connected Dynkin diagram. Let $M$ be a $B$-module. For any $i \geq 1$, if each weight of $M$ is generated by at most $|S|-i$ simple roots, then $H_G^{i-1}(M)$ is a trivial $G$-module.  In particular:
    \begin{enumerate}    
    \item If each weight of $M$ is generated by at most $|S|-1$ simple roots, then $H_G^0(M)$ is a trivial $G$-module.

    \item If each weight of $M$ is generated by at most $|S|-2$ simple roots, then $H_G^1(M)$ is a trivial $G$-module.
    \end{enumerate}
\end{coro}

\begin{proof}
    We reason by induction.
    For $i = 1$, if all the weights of $M$ are generated by at most $s-1$ simple roots, then they are not $S$-positive, hence there are either $0$ or non-dominant, hence $H^0(M)$ is a trivial $G$-module by \Cref{induction of module}.

    Suppose that the result is true for some $i \geq 1$ and that each weight $\la$ of $M$ is generated by at most $s- (i+1)$ simple roots. Let us prove that $H_G^{i}(\la)$ is a trivial $G$-module. If for all $\al \in S$, $\langle \la, \al^{\vee} \rangle \geq 0$, then $\la$ is dominant and by Kempf vanishing $H^{i}_G(\la) = 0$. Otherwise, take $\alpha \in S$ such that $\langle \la, \al^{\vee} \rangle \leq -1$. By \Cref{induction minimal parabolic}, $H^{i}_G(\la) = H^{i-1}_G \circ H^1_{P^{\alpha}}(\la)$ (remember that by abuse of notation we do not write the restriction to $B$, by \Cref{pushforward of structure sheaf}). Also by \Cref{induction minimal parabolic}, all the weights of $H^1_{P^{\alpha}}(\la)$ are generated by $s-i$ simple roots (the generators of $\la$, plus $\alpha$). By induction hypothesis, $H^{i-1}_G \circ H^1_{P^{\alpha}}(\la)$ is a trivial $G$-module, hence $H^{i}_G(\la)$ is a trivial $G$-module.
    Then, by \Cref{induction of module}, $H_G^{i}(M)$ is a trivial $G$-module.
\end{proof}

\begin{ex} \label{counterexample H^1}
    Let us give an important example of a weight $\la$, generated by $|S|-1$ simple roots, such that $H^1(\la)$ is a non-trivial $G$-module. Let $G = SL_n$ and take, for any $r \geq 1$, $$\la := p^r(e_1 - e_{n-1}) = p^r((e_1 - e_2) + \dots + (e_{n-2} - e_{n-1})).$$ Then by \cite[5.15]{jantzen2003representations}, $H^1(\la)$ is not a trivial $G$-module. Its socle is $$L_G(p^r(e_1-e_{n-1}) - (e_{n-1} - e_n)).$$
\end{ex}

\begin{lemma} \label{preserves generation}
    Let $I \subset S$. If $\la$ is a weight generated by $I$, then for all $i \geq 0$, the weights of $H^i_{P^I}(\la)$ are generated by $I$.
\end{lemma}

\begin{proof}
    Let $\la$ be a weight generated by $I$. We reason by induction on $i \geq 0$.
    By \cite[5.21]{jantzen2003representations}, the weights of $H^0_{P^I}(\la)$ are generated by $I$. Let $i \geq 0$ and suppose that the weights of $H^i_{P^I}(\la)$ are generated by $I$. If for all $\alpha \in I$, $\langle \la, \al^{\vee} \rangle \geq 0$, then $\la$ is $P^I$-dominant, and $H^{i+1}_{P^I}(\la) = 0$ by Kempf vanishing. Otherwise, take $\alpha \in I$ such that $\langle \la, \al^{\vee} \rangle \leq -1$. By \Cref{induction minimal parabolic}, $H^{i+1}_{P^I}(\la) = H^{i}_{P^I} \circ H^1_{P_{\alpha}}(\la)$. Again by \Cref{induction minimal parabolic}, given that $\la$ is generated by $I$ and that $\alpha \in I$ all the weights of $H^1_{P^{\alpha}}(\la)$ are generated by $I$. Then, by induction hypothesis, all the weights of $H^{i+1}_{P^I}(\la) = H^{i}_{P^I} \circ H^1_{P^{\alpha}}(\la)$ are generated by $I$.
\end{proof}

\subsubsection{Cohomology of the tangent bundle} \label{section cohomology tangent bundle}

Let $G$ be a (connected) reductive group with connected Dynkin diagram. In characteristic $0$, as an application of his seminal paper \cite{bott1957homogeneous}, Bott shows that higher cohomology of the tangent bundle of a classical homogeneous variety $X = G/P$ vanishes. In \cite{Demazure1977}, Demazure extends this result to characteristic $p$, but only in the case of reduced parabolic subgroups. Those results implies that deformations of $X = G/P$ are rigid.

However, when the stabilizer $P_f$ is a non-reduced parabolic subgroup of $G$, the situation is much more intricate. In particular, we will show in \Cref{section final conclusion} that when the stabilizer $P_f$ has width at least $2$, $H^2(X, T_X)$ does not vanish. In this section, we carry out a detailed study of the $G$-module structure of the first cohomology groups $H^0(X, T_X)$ and $H^1(X, T_X)$. In \Cref{conditions fulfilled}, we show that under a mild technical assumption on $P$, $H^1(X, T_X)$ is a trivial $G$-module, and if moreover the minimum altitude of $f$ is $0$, then higher $G$-cohomology $H^i(G, H^0(X, T_X))$ vanishes for all $i \geq 2$. 

Write $I_l \subset \dots \subset I_1 \subset I_0$ the characteristic filtration of $f$, and write $P_f = P^0_{r_0} \dots P^l_{r_l}$.
To study the cohomology of the tangent bundle $T_X$, we use the filtration of $T_X$ that we constructed in \Cref{section filtration tangent}, with factors $$T_f^k = \res_P^{P^{k+1}} \frak{p}^k/\frak{p}^{k+1},$$ which go with the induction formula of \Cref{new induction}
$$\ind_P^G (T^k_f) \simeq (F_{P^0}^{r_0})^*  \ind_{P^1}^{P^0}  (F_{P^1}^{r_1-r_0})^*  \ind_{P^2}^{P^1}  \dots (F_{P^{k}}^{r_{k} - r_{k-1}})^*  \ind_{P^{k+1}}^{P^{k}} (\frak{p}^k/\frak{p}^{k+1}).$$

The study will be different between the top factor $T^0_f$ and the rest of the filtration. Thus, we will denote $$E := T_{BP \sto B\overl^1},$$ which admits a filtration with factors $(T^k_f)_{1 \leq k \leq l-1}$, and which fits into a short exact sequence $$0 \to E \to T_X \to T_f^0 \to 0.$$

First, as we explained in \Cref{injectivity R0}, by faithful flatness of the map $G/B \to G/P_f$, for a $P_f$-module $M$, vanishing of $R^0 \ind_B^G M$ (where $M$ is seen as a $B$-module) forces vanishing of $R^0 \ind_{P_f}^G M$. It implies the following first result.

\begin{prop} \label{R0 E}
    $R^0 \ind_{P_f}^G E = 0$
\end{prop}

\begin{proof}
    The weights of $E$ are the union of the weights of the factor of its filtration, so there are multiples of roots of $\p^1$, so they are all generated by $|I_1|$ simple roots, and $|I_1| \leq s-1$, so by \Cref{s-i} $R^0\ind_B^G E = 0$. Hence, by \Cref{injectivity R0}, $R^0 \indp E = 0$.
\end{proof}

Now we use Demazure's results from \cite{Demazure1977} to study the cohomology of the top factor $T^0_f$. He restricts to reduced parabolic subgroups which are said to be \emph{non-exceptional}, defined as follows.

\begin{defi}[\cite{Demazure1977}] \label{def exceptional}
    A proper reduced parabolic subgroup $P$ is called exceptional if one of the following conditions is verified:
    \begin{enumerate}
    \item $G$ has type $C_n$ and $P$ has type $C_{n-1}$.
    \item $G$ has type $A_n$ and $P$ has type $A_{n-1}$.
    \item $G$ has type $G_2$ and $P$ corresponds to the stabilizer of extremal lines in the minuscule representation.
    \end{enumerate}
\end{defi}

\begin{notation}
    The group $G$ has a connected Dynkin diagram so it has only one highest root, denoted $\theta$.
\end{notation}

\begin{rk}
    If $G$ has same-length roots, $\thet$ is the only dominant root of $G$.  Otherwise, the only dominant roots are $\thet$ and $\thet'$, the only short dominant root.
\end{rk}

\begin{prop}[\cite{Demazure1977}] \label{lemma Demazure}
    If $P$ is a non-exceptional reduced parabolic subgroup, then
    \begin{enumerate}
        \item $H^k_G(\g/\p) = 0$ for all $k \geq 1$.
        \item there is an isomorphism of $G$-modules $H^0_G(\g/\p) \simeq H^0_G(\theta)$.
    \end{enumerate}
\end{prop}

\begin{proof}
    The first part is exactly Proposition 2, b) of \cite{Demazure1977}.
    Then \cite[Proposition 3]{Demazure1977} shows that there exists a $B$-submodule $F$ of $\g/\p^1$ such that $(\g/\p^1)/F \simeq \thet$, and such that for all $i \geq 0$, $H^i_G(F) = 0$. The vanishing $H^i_G(F) = 0$ in the long exact sequence of cohomology associated to the short exact sequence $$0 \to F \to \g/\p^1 \to \thet \to 0$$ shows that $H^0_G(\g/\p^1) \simeq H^0_G(\theta)$.
\end{proof}

\begin{coro} \label{induction first factor}
    If $P^1$ is not exceptional, then 
    \begin{enumerate}
    \item $R^k \indp (T^0_f) = 0$ for all $k \geq 1$.
    \item there is an isomorphism of $G$-modules $R^0 \indpf (T^0_f) \simeq (F_G^{r_0})^*H^0_G(\theta)$.
    \item if $r_0 = 0$, then $H^i(G, R^0 \indpf (T_f^0)) = 0$ for all $i \geq 0$.
    \end{enumerate}

\end{coro}

\begin{proof}
    By \Cref{new induction}, we have the formula $\indp(T^0_f) \simeq (F_G^{r_0})^*\ind_B^G(\frak{g}/\frak{p}^1)$, so (1) and (2) are direct consequences of \Cref{lemma Demazure} (1) and (2). Then, if $r_0 = 0$, by \Cref{G-cohomology of H0(lambda)} we know that for all $i \geq 0$, $H^i(G, H^0_G(\theta)) = 0$, so for all $i \geq 0$, $H^i(G, R^0 \indp (T_f^0)) = 0$.
\end{proof}

\paragraph{A digression about the \emph{very good} case.} When the prime $p$ is \emph{very good}, previous results about $R^0 \indpf E$ and $R^0 \indpf (T^0_f)$ directly give the structure of $H^0(X, T_X)$, and its $G$-cohomology.

\begin{defi}
    A prime $p \geq 5$ is said to be very good if neither of the following conditions is verified:
    \begin{enumerate}
    \item when $G$ has type $A_n$, $p$ is a divisor of $n+1$.
    \item when $G$ has type $E_8$, $p \neq 5$.
    \end{enumerate}
\end{defi}

\begin{lemma}[{\cite[4.1]{mcninch2000secondcohomologysmallirreducible}}] \label{mcninch simple}
    If $p$ is very good, $H^0(\thet)$ is a simple $G$-module isomorphic to $\g$.
\end{lemma}

\begin{prop} \label{case very good}
    If $p$ is very good and $P^1$ is not exceptional, then there is an isomorphism of $G$-modules $$H^0(X, T_X) \simeq (F_G^{r_0})^* H^0(\theta) \simeq (F_G^{r_0})^* \g.$$

    As a consequence, when $r_0 = 0$, $H^i(G, H^0(X, T_X)) = 0$ for all $i \geq 0$.
\end{prop}

\begin{proof}
    The short exact sequence of $P$-modules  $$0 \to E \to T_X \to T^0_f$$ induces a long exact sequence $$0 \to R^0 \indp E \to H^0(X, T_X) \to R^0 \indp (T^0_f) \to \dots$$
    But by \Cref{R0 E} $R^0 \indp E = 0$, and by \Cref{induction first factor} $R^0 \indp (T^0_f) \simeq (F_G^{r_0})^*H^0_G(\theta)$, so there is an injective map $$H^0(X, T_X) \to (F_G^{r_0})^*H^0_G(\theta).$$ But if $p$ is very good, $H^0(\thet)$ is simple by \Cref{mcninch simple}, so $(F_G^{r_0})^*H^0_G(\theta)$ is simple, so the map $H^0(X, T_X) \to (F_G^{r_0})^*H^0_G(\theta)$ is an isomorphism.
\end{proof}

In the case where $r_0 \geq 1$, we can even conclude that $H^2(G, H^0(X, T_X))$ does not vanish, which prevents from proving unobstructedness of equivariant deformations of $X$.

\begin{lemma}[{\cite[Theorem C]{mcninch2000secondcohomologysmallirreducible}}] \label{second lemma mcninch}
    Let $h$ denote the Coxeter number of $G$. If $p > h$, then for all $r \geq 1$ $$\dim H^2(G, (F_G^r)^* \g) = 1.$$ 
\end{lemma}

\begin{coro} \label{map obstructed}
    If $p > h$, if $P^1$ is not exceptional, and if $r_0 \geq 1$, then $$\dim H^2(G, H^0(X, T_X)) = 1.$$
\end{coro}

\begin{proof}
    This is an immediate consequence of \Cref{case very good} and \Cref{second lemma mcninch}, given that if $p > h$, then $p$ is a fortiori very good.
\end{proof}

\paragraph{Back to the general case.} 

When $p$ is not very good, we cannot directly compare $H^0(X, T_X)$ and $(F_G^{r_0})^* H^0(\theta)$. However, with an additional result on the $G$-structure of $H^1(E)$, we will be able to compare the $G$-cohomology of $H^0(X, T_X)$ and $(F_G^{r_0})^* H^0(\theta)$, thanks to the following purely homological lemma.

\begin{lemma} \label{H^1 trivial}
    Let $P$ be a parabolic subgroup and write here $H^0 := R^0\indp$ and $H^1 := R^1\indp$ to lighten the notation. Consider a short exact sequence of $P$-modules $$0 \to F \to T \to Q \to 0.$$ 
    Suppose that $R^0\indpf (F)$ and $H^1(F)$ are trivial $G$-modules. Then for all $i \geq 2$, there is an isomorphism $$H^i(G, H^0(T)) \iso H^i(G, H^0(Q)).$$ 
\end{lemma}

\begin{proof}
    The short exact sequence induces a long exact sequence $$0 \to H^0(F) \to H^0(T) \to H^0(Q) \xlongrightarrow{\delta} H^1(F) \to \dots$$ 
    Denote by $K$ the kernel of the map $\delta$. There is a short exact sequence $$0 \to H^0(F) \to H^0(T) \to K \to 0,$$ which induces a long exact sequence $$\dots \to H^i(G, H^0(F)) \to H^i(G, H^0(T)) \to H^i(G, K) \to H^{i+1}(G, H^0(F)) \to \dots$$ 
    But for all $i \geq 1$, $H^i(G, H^0(F)) = 0$ because $H^0(F)$ is a trivial $G$-module, so for all $i \geq 1$, the map $$H^i(G, H^0(T)) \to H^i(G, K)$$ is an isomorphism.
    Denote $I = \im(\delta) \subset H^1(F)$. $H^1(F)$ is a trivial $G$-module so the submodule $I$ is also a trivial $G$-module. Moreover the short exact sequence $$0 \to K \to H^0(Q) \to I \to 0$$ induces a long exact sequence $$\dots \to H^{i-1}(G, I) \to H^i(G, K) \to H^{i}(G, H^0(Q)) \to H^i(G, I) \to \dots$$
    But for all $i \geq 1$, $H^i(G, I) = 0$ so for all $i \geq 2$, the map $$H^i(G, K) \to H^i(G, H^0(Q))$$ is an isomorphism. Hence for all $i \geq 2$ there is an isomorphism $$H^i(G, H^0(T))\iso H^i(G, H^0(Q)).$$
\end{proof}

Now all that remains is to study the $G$-module $R^1 \indpf E$.

\begin{notation}
    To simplify the notation, we will write $i_k = (F_{P^k}^{r_k - r_{k-1}})^*\ind_{P^{k+1}}^{P^k}$ for all $0 \leq k \leq l-1$.
\end{notation}

The following lemma is purely homological.

\begin{lemma} \label{mix spectral sequence}
    Let $M$ be a $P_f$-module.
    If for all $0 \leq k \leq l-1$, $R^0i_0 \circ \dots \circ R^1i_k \circ \dots \circ R^0i_{l-1}(M)$ is a trivial $G$-module, then $R^1(i_0 \circ \dots\circ i_{l-1})(M)$ is a trivial $G$-module.
\end{lemma}

\begin{proof}
    We give the proof in general for the composition of two functors $F \circ F'$, with target the category of $G$-modules. Consider the spectral sequence $R^iF \circ R^jF' \Longrightarrow R^{i+j}(F \circ F')$, and its associated $5$-terms exact sequence  $$\dots \to R^1F \circ R^0F' (M) \to R^1(F \circ F')(M) \to R^0F \circ R^1F' (M) \to \dots$$ 
    If $R^1F \circ R^0F' (M)$ and $R^0F \circ R^1F'(M)$ are trivial $G$-modules, we argue as in \Cref{reconstruct} to get that $R^1(F \circ F')(M)$ is also a trivial $G$-module, given that $H^1(G, k) = 0$.
\end{proof}

We would like to use \Cref{s-i} to show that $R^1 i_1 \circ R^0 i_1 \circ \dots \circ R^0 i_l (E)$ is a trivial $G$-module. For this, we need that all the weights of $R^0 i_1 \circ \dots \circ R^0 i_l (E)$ are generated by at most $s-2$ simple roots. This is always true if $|I_1| \leq s-2$, but not necessarily if $|I_1| = s-1$ (see \Cref{counterexample H^1}). However, when $|I_1| = s-1$, it will still become true if we restrict to the case where $I_1$ is disconnected as a subgraph of the Dynkin diagram.

\begin{lemma} \label{disjoint union induction}
    Let $I, J, K \subset S$. If the Dynkin diagram associated to $I$ is the disjoint union of the Dynkin diagrams associated to $J$ and to $K$, and if $\la$ is a weight generated by $J$, then $$\res_{P^J}^{P^I}\ind_B^{P^I} \la \simeq \ind_B^{P^J} \la.$$
\end{lemma}

\begin{proof}
    We reduce to the derived subgroup of the Levi factor $G^I$ of $P^I$. Write $B^I = B \cap G^I$. Then $\ind_B^{P^I} \la \simeq \ind_{B^I}^{G^I} \la$ where the unipotent part $U^I$ acts trivially. Consider $G^I_{\der}$ the derived subgroup of $G^I$, and $B^I_{\der} = B^I \cap G^I_{\der}$. Then $\ind_{B^I}^{G^I} \la \simeq \ind_{B^I_{\der}}^{G^I_{\der}} \la$ where $Z(G^I)$ acts trivially because $\la$ is generated by $J$, hence by $I$.
    But $B^I_{\der} \simeq B^J_{\der} \times B^K_{\der}$ and $G^I_{\der} \simeq G^J_{\der} \times G^K_{\der}$ because $I$ is the disjoint union of $J$ and $K$, so $$\ind_{B^I_{\der}}^{G^I_{\der}} \la \simeq (\ind_{B_J^{\der} \times B_K^{\der}}^{G_J^{\der} \times G_K^{\der}}) \la.$$
    The weight $\la$ is generated by $J$ so it decomposes as $(\la, 0)$ (that is, with trivial action of $B^K$). So we have the isomorphisms $$\ind_{B^I_{\der}}^{G^I_{\der}} \la \simeq \ind_{B^J_{\der} \times B^K_{\der}}^{L^J_{\der} \times L^K_{\der}} \la \simeq \ind_{B^J_{\der}}^{L^J_{\der}} \la \otimes \ind_{B^K_{\der}}^{L^K_{\der}} k \simeq \ind_{B^J_{\der}}^{L^J_{\der}} \la$$ where $L^K_{\der}$ acts trivially.
    Hence $\res_{P^J}^{P^I} \ind_B^{P^I} \la \simeq \ind_{B^J}^{G^J} \la$.
\end{proof}

\begin{lemma} \label{number of generators disjoint union}
    If $|I_1| = s-1$ and $I_1$ is disconnected as a subgraph of the Dynkin diagram, then all the weights of $R^0 i_1 \circ \dots \circ R^0 i_k (\p^k/\p^{k+1})$ are generated by at most $s-2$ simple roots.
\end{lemma}

\begin{proof}
    Decompose the graph of $I_1$ as the disjoint union of $J$ and $K$. All the $I_k$ for $1 \leq k \leq l$ decompose as $I_k \cap J$ and $I_k \cap K$. Write $P^k_J$ and $P^k_K$ the parabolic subgroups associated to $I_k \cap J$ and $I_k \cap K$ respectively. Write $i_k^J = (F^{r_k - r_{k-1}})^*\ind_{P^{k+1}_J}^{P^k_J}$ for all $0 \leq k \leq l-1$.

    Take $\lambda$ a weight of $\p^k/\p^{k+1}$. The root system $R_I$ splits as $R_I = R_J \sqcup R_K$, so $\lambda$ is generated either by $I_k \cap J$ or $I_k \cap K$. Suppose it is generated by $I_k \cap J$. 
    By \Cref{disjoint union induction}, the weights of $R^0 i_1 \circ \dots \circ R^0 i_k (\lambda)$ are the same as the weights of $R^0 i_1^J \circ \dots \circ R^0 i_k^J (\lambda)$. By \Cref{preserves generation}, all the weights of $R^0 i_1^J \circ \dots \circ R^0 i_k^J (\lambda)$ are generated by $J$, so all the weights of $R^0 i_1 \circ \dots \circ R^0 i_k (\lambda)$ are generated by $J$. But $|J| \leq s-2$ so all the weights are generated by at most $s-2$ simple roots.
\end{proof}

\begin{prop} \label{R^1 induction factors trivial}
    Suppose that $|I_1| \leq s-2$ or, if $|I_1| = s-1$, that $I_1$ is disconnected as a subgraph of the Dynkin diagram.
    Then for all $1 \leq k \leq l-1$, $R^1 \indp T^k_f$ is a trivial $G$-module.
\end{prop}

\begin{proof}
    By \Cref{new induction}, there is an isomorphism of $G$-modules $$R^1 \indp T_k \simeq R^1(i_0 \circ \dots \circ i_k) (\p^k/\p^{k+1}).$$
    By \Cref{mix spectral sequence}, we have to show that for all $0 \leq k \leq l$, $R^0 i_0 \circ \dots \circ R^1 i_j \circ \dots \circ R^0 i_k(\p^k/\p^{k+1})$ is a trivial $G$-module.
    Suppose that $k \geq 1$. By \Cref{preserves generation}, the weights of the $P^1$-module $R^0 i_1 \circ \dots \circ R^1 i_j \circ \dots \circ R^0 i_k(\p^k/\p^{k+1})$ are generated by $I_1$ which has cardinality at most $s-1$, so by \Cref{s-i} $R^0 i_0 \circ R^0 i_1 \circ \dots \circ R^1 i_j \circ \dots \circ R^0 i_k(\p^k/\p^{k+1})$ is a trivial $G$-module.
    
    For the case $k = 0$, if $|I_1| \leq s-2$, the weights of $R^0 i_1 \circ \dots \circ R^0 i_k(\p^k/\p^{k+1})$ are generated by $I_1$ by \Cref{preserves generation}, that is, by at most $s-2$ simple roots. If $|I_1| = s-1$ and $I_1$ is disconnected as a subgraph of the Dynkin diagram, then by \Cref{number of generators disjoint union}, $R^0 i_1 \circ \dots \circ R^0 i_k(\p^k/\p^{k+1})$ has weights generated by at most $s-2$ simple roots. Hence, by \Cref{s-i}, the $G$-module $R^1 i_0 \circ R^0 i_1 \circ \dots \circ R^0 i_k(\p^k/\p^{k+1})$ is trivial, which completes the proof.
\end{proof}

Finally, all the partial results of this section assemble to give a result on the cohomology of the whole tangent bundle $T_X$.

\begin{thm} \label{conditions fulfilled}     
    Suppose that $G$ has a connected Dynkin diagram. Suppose that $|I_1| \leq s-2$ or, if $|I_1| = s-1$, that $I_1$ is disconnected as a subgraph of the Dynkin diagram.
Then the following assertions hold:
\begin{enumerate}
\item $H^1(X, T_X)$ is a trivial $G$-module
\item for all $i \geq 2$, $H^i(G, H^0(X, T_X)) \simeq H^i(G, (F_G^{r_0})^*H_G^0(\thet)).$  
\item if $r_0 = 0$, then for all $i \geq 2$, $H^i(G, H^0(X, T_X)) = 0$.
\end{enumerate}
    
\end{thm}

\begin{proof}
    First, notice that if $|I_1| \leq s-2$, or if $|I_1| = s-1$ and $I_1$ is disconnected as a subgraph of the Dynkin diagram, then $P^1$ is in particular non-exceptional in the sense of \Cref{def exceptional}.

    \begin{enumerate}
        \item by \Cref{induction first factor}, $R^1 \indp T^0_f = 0$. By \Cref{R^1 induction factors trivial}, $R^1 \indp T^k_f$ is a trivial $G$-module for all $1 \leq k \leq l-1$. Then, by \Cref{filtration} and \Cref{reconstruct}, $R^1 \ind_P^G T_{BP_f \sto BG}$ is a trivial $G$-module, that is, $H^1(X, T_X)$ is a trivial $G$-module.
        \item Remember that $E = T_{BP \sto B\overl^1}$ fits into an exact sequence $$0 \to E \to T_X \to  T^0_f \to 0.$$ By \Cref{R0 E}, $R^0 \indp E = 0$ and, as we showed in (1), $R^1 \indp E$ is a trivial $G$-module. By \Cref{induction first factor}, $R^0 \indp (T^0_f) \simeq (F_G^{r_0})^*H^0_G(\theta)$. Thus, we can use \Cref{H^1 trivial} to conclude that for all $i \geq 2$, $H^i(G, H^0(X, T_X)) \simeq H^i(G, (F_G^{r_0})^*H_G^0(\thet)).$
        \item by \Cref{G-cohomology of H0(lambda)}, for all $i \geq 1$, $H^i(G, H_G^0(\thet)) = 0$, so if $r_0 = 0$, by (2), it implies that for all $i \geq 2$, $H^i(G, H^0(X, T_X)) = 0$.
    \end{enumerate}
\end{proof}

\begin{rk}
    We will see in \Cref{section final conclusion} that when $P_f$ has width at least $2$, $H^2(X, T_X) \neq 0$, in contrast with the case of characteristic $0$.
\end{rk}

\subsection{Classification of simple modules for non-reduced subgroups} \label{section classification simple modules}

Let $G$ be a (reduced) reductive group over an algebraically closed field $k$. In characteristic $0$, the Borel-Weil-Bott theorem states that simple $G$-modules correspond to global sections of ample line bundles on the flag variety $G/B$, that is, to induced $G$-modules $R^0 \ind_B^G \la$ where $\la$ is a dominant weight. In positive characteristic, the situation is more complicated. Simple $G$-modules correspond to the unique simple submodule of the induced modules $R^0 \ind_B^G \la$, denoted $L_G(\la)$. Moreover, a $p$-adic expansion of a dominant weight $\la$ leads, by the Steinberg theorem, to a decomposition of $L_G(\la)$ with Frobenius twists of $L_G(\la_i)$ where the $\la_i$ are $p$-restricted dominant weights.

We extend those results to non-reduced reductive subgroups $G_f$ of $G$. In \Cref{section BW theorem}, we show that simple $G_f$-modules are of the form $L_{G_f}(\la)$ for $\la$ a $G_f$-dominant weight. $G_f$-dominant weights are nothing else than $(G_f)_{\red}$-dominant weights. In \Cref{section transfer}, we show that when $\la$ is suitably $p$-restricted, restriction of $G_f$-simple modules to another smaller non-reduced reductive group is still simple. It allows us to prove an extension of the Steinberg theorem to non-reduced reductive groups in \Cref{section steinberg}. When $\la$ has a suitable $p$-adic expansion, $L_{G_f}(\la)$ can be decomposed with Frobenius twists of simple modules over reduced reductive groups. Finally, we use Steinberg theorem to propose a filtration of the tangent bundle of $G/P_f$, alternative to the one we constructed in \Cref{section filtration tangent}, which also comes with a formula for induction of its factors. This new filtration could also have been used to study the cohomology of the tangent bundle as in \Cref{section cohomology tangent bundle}.

We begin by some reminders about socle, radical and distribution algebras, which will be useful in what follows.

\begin{defi} \label{socle radical}
    Let $H$ be a group scheme over $k$. Let $M$ be a finite-dimensional $H$-module. \begin{enumerate}

    \item The socle $\soc_H M$ is the largest semi-simple submodule of $M$.

    \item The radical $\rad_H M$ is the smallest submodule of $M$ such that the quotient is semi-simple.
    \end{enumerate}
    Socle and radical are dual in the sense that $\soc_H(M^*) \simeq (M/\rad_H M)^*$.
\end{defi}

\paragraph{Distribution algebras.}

We refer to \cite[I.7]{jantzen2003representations} for a complete introduction to distribution algebras of group schemes.
In characteristic $0$, representations of $G$ are very well reflected by representations of the Lie algebra $\g$ (which are representations of its enveloping algebra $U(\g)$).
In positive characteristic, even when $G$ is smooth, the classical tangent Lie algebra $\g$ does not contain enough information about $G$. One has to consider the distribution algebra (or hyperalgebra) of $G$, denoted $\Dist(G)$. In characteristic $0$, $\Dist(G)$ coincides with $U(\g)$.
A $G$-module induces a $\Dist(G)$-module. For $\al \in R$, denote by $(X_{\al,n})_{n \in \N}$ the usual generators of $\Dist(U_{\al})$. Denote also by $(H_{i,m})_{1 \leq i \leq r, m \in \N}$ the usual generators of $\Dist(T)$.

\begin{prop}[{\cite[II.3.3]{jantzen2003representations}}] \label{basis r}
    Let $\alpha \in R$.
   $(X_{\al, n})_{0 \leq n < p^r}$ is a basis of $\Dist(U_{\alpha,r})$.

\end{prop}

\begin{lemma}[{\cite[II.1.19(5)]{jantzen2003representations}}] \label{shift}
    Let $\alpha \in R$, let $M$ be a $U_{\alpha,r}T$-module. Then for all $n < p^r$, $$X_{\alpha,n} M_{\la} \subset M_{\la + n\alpha}.$$
\end{lemma}

\begin{coro} \label{trivial action}
    Let $M$ be a $U_{\alpha,r}T$-module and let $\la$ be a maximal weight of $M$. Then $U_{\alpha,r}$ stabilizes $M_{\la}$.
\end{coro}

\begin{proof}
    Given that $\la$ is the maximal weight of $M$, \Cref{shift} implies that for all $n \geq 1$, $X_{\alpha,n} M_{\la} = 0$, hence $\Dist(U_{\alpha,r})$ acts trivially on $M_{\la}$. By \Cref{additive group action}, $U_{\alpha,r}$ acts trivially on $M_{\la}$.
\end{proof}

\begin{prop}[{\cite[p. 39]{HaboushLauritzen1993}}] \label{base Dist}
    Let $f$ be a non-infinitesimal $W$-function. The distribution algebra $\Dist(G_f)$ is spanned by $$(X_{\al, n})_{\al \in R^+, 0 \leq n < p^{f(\al)}}, (H_{i,m})_{1 \leq i \leq r, m \in \N}, (X_{-\al, n})_{\al \in R^+, 0 \leq n < p^{f(\al)}}.$$
\end{prop}

\begin{lemma}[{\cite[I.7.15]{jantzen2003representations}} for $\G_a$ and $\G_m$, {\cite[I.8.6]{jantzen2003representations}} for $(\G_a)_r$ and $(\G_m)_r$] \label{additive group action}
    Let $M$ be a $\G_a$-module (resp. $\G_m$-module, resp. $(\G_a)_r$-module, resp. $(\G_m)_r$-module) and let $M'$ be a $k$-submodule of $M$. Then $M'$ is a $\G_a$-submodule (resp. $\G_m$-submodule, resp. $(\G_a)_r$-submodule, resp. $(\G_m)_r$-submodule) of $M$ if and only if it is a $\Dist(\G_a)$-submodule (resp. $\Dist(\G_m)$-submodule, resp. $\Dist((\G_a)_r)$-submodule, resp. $\Dist((\G_m)_r)$-submodule) of $M$.
\end{lemma}

\begin{coro}[from {\cite[II.1.20]{jantzen2003representations}}] \label{Dist submodule}
    Let $f$ be a $W$-function. Let $M$ be a $G_f$-module, let $N$ be a $\Dist(G_f)$-submodule of $M$. Then $N$ is a $G_f$-submodule of $M$.
\end{coro}

\begin{proof}
    $N$ is a $\Dist(T_{r_l})$-module and a $\Dist((U_{\pm \al})_{f(\al)})$-module for all $\al \in R^+$. Using the analogous results for $\G_m$, $\G_a$ and $\alpha_{p, r}$ from \Cref{additive group action}, we deduce that $N$ is a $T$-submodule and a $U_{\al, f(\al)}$-submodule for all $\alpha \in R$. By \Cref{density}, $U_f^+ B_f$ is dense in $G_f$ so as in \cite[II.1.20]{jantzen2003representations} it implies that $N$ is a $G_f$-module.
\end{proof}

\subsubsection{Classification of simple representations of non-reduced groups} \label{section BW theorem}

In \cite[II.2, II.3, II.9]{jantzen2003representations}, Jantzen proves for the groups $G$, $G_r$ and $G_rT$ respectively, that their simple representations correspond to a subset of weights: $G$-dominant weights for $G$, all $X(T)$ for $G_rT$, and all $X(T_r)$ for $G_r$. In contrast to the case of characteristic $0$, the simple representations are not given simply by $H^0(\la)$ but by its socle.
We extend those results to any non-reduced reductive groups $G_f$, infinitesimal or not, by showing that their simple representations correspond to $G_f$-dominant weights, defined similarly as usual $G$-dominant weights, and we establish what are those $G_f$-dominant weights.   
We simply extend the arguments of Jantzen to the setting of non-reduced reductive groups. Moreover, everything works the same for the non-reduced reductive subgroup $G_f$ or for the non-reduced parabolic subgroup $P_f$ (as everything works the same for $G_rT$ and $G_rB$).

\paragraph{Dominant weights.}

As for $G$-dominant weights, we define $G_f$-dominant weights as the weights $\la$ such that $R^0 \ind_{B_f}^{G_f} \la \neq 0$. Write $G_f = G^0_{r_0} \dots G^l_{r_l}$. We show that $G_f$-dominant weights are nothing else than $(G_f)_{\red}$-dominant weights, extending the fact that all weights are $G_rT$-dominant.

\begin{rk} \label{remark T B}
    If $f$ is a $W$-function, one-dimensional representations of $B_f$, denoted $X(B_f)$, are the same as one-dimensional representations of $T_{r_l}$, denoted $X(T_{r_l})$, because $B_f = U_f \rtimes T_{r_l}$ and $U_f$ is unipotent so it has no non-trivial $1$-dimensional representations. Hence, by abuse of notation, if $\la \in X(T_{r_l})$, we also denote by $\la$ the corresponding representation of $B_f$.
\end{rk}

\begin{rk}
   When $r_l \neq +\infty$, by \cite[II.3.7]{jantzen2003representations}, $$X(T_{r_l}) \simeq X(T)/p^{r_l}X(T).$$
\end{rk}

\begin{defi} \label{defi dominant weights}
    Let $I \subset S$
    Let $X_+^I(T)$ be the subset of $I$-dominant weights $$X_+^I(T) := \{\la \in X(T), \forall \alpha \in I, \<\la, \al^{\vee}\> \geq 0 \}.$$
    Let $r \geq 1$. Let $X_r^I(T)$ be the subset of $p^r$-restricted $I$-dominant weights
    $$X_r^I(T) := \{\la \in X(T), \forall \alpha \in I, 0 \leq \<\la, \al^{\vee}\> < p^r \} \subset X_+^I(T).$$
\end{defi}

\begin{defi}
    Let $f$ be a $W$-function with values $r_0 < \dots < r_l$. Let $\la \in X(T_{r_l})$. \begin{enumerate}
    
    \item We will write $$H^0_{G_f}(\la) := R^0 \ind_{B_f}^{G_f} \la.$$
    A weight $\la \in X(T_{r_l})$ is said to be $G_f$-dominant if $H^0_{G_f}(\la) \neq 0$.

    \item If $f$ is not infinitesimal, we will write $$H^0_{P_f}(\la) := R^0 \ind_B^{P_f} \la.$$
    A weight $\la \in X(T)$ is said to be $P_f$-dominant if $H^0_{P_f}(\la) \neq 0$.
    \end{enumerate}
\end{defi}

\begin{lemma} \label{pushforward tensor}
    Let $f$ be a $W$-function. Write $P_f = P^0_{r_0} \dots P^l_{r_l}$. 
    \begin{enumerate}

    \item Suppose that $f$ is non-infinitesimal. Let $M$ be a $P^l$-module. Then there is an isomorphism of vector spaces $$R^0 \ind_{P^l}^{P_f} M \simeq k[(U^0_f)^+] \otimes_k M$$ 
    
    \item Suppose that $f$ is non-infinitesimal. Let $M$ be a $(U^l)^+B_f$-module. Then there is an isomorphism of vector spaces $$R^0 \ind_{(U^l)^+B_f}^{G_f} M \simeq k[(U^0_f)^+] \otimes_k M$$ 

    \item Suppose that $f$ is infinitesimal. Let $M$ be a $B_f$-module.
    Then there is an isomorphism of vector spaces $$H^0_{G_{f}}(M) \simeq k[U_f^+] \otimes_k M$$ compatible with the action of $(U_f)^+$ by regular representation on $k[U_f^+]$ and trivial on $M$, and with the action of $T_{r_l}$ through conjugation on $k[U_f^+]$ and as given on $M$.
    \end{enumerate}
\end{lemma}

\begin{proof}
By \Cref{product decomposition}, when $f$ is not infinitesimal there are the isomorphisms of schemes $P_f \simeq (U^0_f)^+ \times P^l$ and $G_f \simeq (U_0^f)^+ \times (U^l)^+ B_f$, and when $f$ is infinitesimal there is the isomorphism of schemes $G_f \simeq U_f \times B_f$ . Moreover, this last isomorphism is compatible with the action of $(U_f)^+$ by left multiplication on $(U_f)^+$ and on $G_f$, and with the action of $B_f$ by right multiplication on $B_f$ and on $G_f$. As in \cite[II.3.6]{jantzen2003representations}, it implies the expected results.
\end{proof}

\begin{rk}
    Let us mention a point of view using classifying stacks, for $P_f$ for instance. The isomorphism of schemes $P_f \simeq (U^0_f)^+ \times P^l$ 
    induces the two (homotopy) cartesian diagram

\[\begin{tikzcd}
	{(U^0_f)^+} & {*} \\
	{*} & {B(U^0_f)^+} \\
	{BP^l} & BP_f
	\arrow[from=1-1, to=1-2]
	\arrow[from=1-1, to=2-1]
	\arrow["\lrcorner"{anchor=center, pos=0.125}, draw=none, from=1-1, to=2-2]
	\arrow[from=1-2, to=2-2]
	\arrow[from=2-1, to=2-2]
	\arrow[from=2-1, to=3-1]
	\arrow["\lrcorner"{anchor=center, pos=0.125}, draw=none, from=2-1, to=3-2]
	\arrow[from=2-2, to=3-2]
	\arrow[from=3-1, to=3-2]
\end{tikzcd}\]
Hence by base change $\res_{U^0_f}^{P_f} \ind_{P^l}^{P_f} M \simeq  \ind_{*}^{U^0_f} M \simeq k[U^0_f] \otimes_k M$ where $U^0_f$ acts by regular representation on $k[U^0_f]$ and trivially on $M$.
\end{rk}

\begin{rk}
    In the case $P_f = G_rB$, Jantzen proves that $$\res_{G_r}^{G_rB} H^0_{G_rB}(\la) \simeq H^0_{G_r}(\la).$$
    However this does not extend to general non-reduced reductive or parabolic subgroups, because it is already false for $\res_{G_r}^G H^0_G(\la)$.  
    We will see however that this will become true for restricted weights after taking the socles.
\end{rk}

\begin{prop}
    Let $f$ be a non-infinitesimal $W$-function. \begin{enumerate}

    \item A weight $\la$ is $P_f$-dominant if and only if it is $P^l$-dominant, that is, $\la \in X_+^{I_l}(T)$.
    \item A weight $\la$ is $G_f$-dominant if and only if it is $P^l$-dominant, that is, $\la \in X_+^{I_l}(T)$. In particular $P_f$-dominant weights are the same as $G_f$-dominant weights.
    \end{enumerate}
    If $f$ is infinitesimal, all the weights of $X(T_{r_l})$ are $G_f$-dominant.
\end{prop}

\begin{proof}
    \begin{enumerate}
    \item $H^0_{P_f}(\la) = \ind_B^{P_f} \la = \ind_{P^l}^{P_f} \circ \ind_B^{P^l}(\la) = k[(U_f^0)^+] \otimes_k \ind_B^{P^l}\la$ by \Cref{pushforward tensor}. But $k[(U_f^0)^+]$ is a non-zero vector space, so $H^0_{P_f}(\la) = 0$ if and only if $H^0_{P^l}(\la) = 0$.
    \item The proof is analogous to (1).
    \end{enumerate}
\end{proof}

\begin{ex}
    All the weights $X(T)$ are $G_rT$-dominant and $G_rB$-dominant. All $X(T_r)$ are $G_r$-dominant.
\end{ex}

Now let us define the candidates to be the simple $G_f$-modules.

\begin{defi}
    Let $\la \in X(T)$ be a weight. We will write $$L_{G_f}(\la) := \soc_{G_f} H^0_{G_f}(\la).$$
\end{defi}

\paragraph{Classification of simple representations of non-reduced non-infinitesimal reductive subgroups.}

Here let $f$ be a non-infinitesimal $W$-function.
We prove the classification of simple $G_f$-modules by closely following the proof of \cite[II.2]{jantzen2003representations}, simply adapting it to the non-reduced setting.

\begin{lemma}[{from \cite[2.1(1), 2.1 (1), (2)]{jantzen2003representations}}] \label{unipotent fixed points}
    For all non-zero $G_f$-modules $V$, $V^{U_f}$ is a $T$-submodule of $V$ and $$V^{U_f} \neq 0.$$
As a consequence, there is $\la \in X(T)$ such that \[\Hom_{B_f}(\la, V) \neq 0\]
\end{lemma}

\begin{proof}
    $T$ normalizes $U_f$ so $V^{U_f}$ is a $T$-submodule of $V$.
    We use the general result of \cite[I.2.14(8)]{jantzen2003representations} to show that $V^{U_f} \neq 0$ for the first part. Then for all weight $\la \in X(T)$, in view of \Cref{remark T B}, $(V^{U_f})_{\la} = V_{\la}$ but $V^{U_f} \neq 0$ so there exists $\la \in X(T)$ such that $$\Hom_{B_f}(\la, V) \neq 0.$$
\end{proof}

\begin{rk}
    \Cref{unipotent fixed points} works also with $U_f^+$ and $B_f^+$ instead of $U_f$ and $B_f$.
\end{rk}

\begin{coro}[\cite{jantzen2003representations}] \label{embeds}
    For all finite-dimensional non-zero $G_f$-modules $V$, there exists $\la \in X(T)$ such that$$\Hom_{G_f}(V, H^0_{G_f}(\la)) \neq 0.$$ 
\end{coro}

\begin{proof}
    One applies \Cref{unipotent fixed points} to the dual of $V$ and then dualizes: thus there is $\la \in X(T)$ such that $$\Hom_{B_f}(V, \la) \neq 0.$$   By adjunction between induction and restriction, one gets that $$\Hom_{G_f}(V, H^0(\la)) \neq 0.$$
\end{proof}

\begin{lemma} \label{inclusion}

    Let $M$ be a $(B_f)^+$-module.
    If $\la$ is a maximal weight of $M$, then $M_{\la} \subset M^{(U_f)^+}$.
\end{lemma}

\begin{proof}

    Let $\al \in R^+$ such that $f(\al) > 0$. By \Cref{trivial action}, $U_{\alpha,f(\al)}$ stabilizes $M_{\la}$ so $(U_f)^+$ stabilizes $M_{\la}$.    
\end{proof}

\begin{prop}[{from \cite[II.2.2]{jantzen2003representations}}] \label{main prop}
    The following assertions hold:
    \begin{enumerate}
    \item $\dim H^0_{G_f}(\la)^{(U_f)^+} \leq 1$
    \item $H^0_{G_f}(\la)^{(U_f)^+} = (H^0_{G_f}(\la))_{\la}$
    \item $\la$ is the highest weight of $H^0_{G_f}(\la)$
    \item  if $H^0_{G_f}(\la) \neq 0$, there is an isomorphism of $T$-modules $H^0_{G_f}(\la)^{(U_f)^+} \simeq \la$.
\end{enumerate}
\end{prop}

\begin{proof}
    For (1) and (2), this is exactly the same proof as in Jantzen, given that $(U_f)^+ B_f$ is dense in $G_f$ by \Cref{density}. By \Cref{inclusion}, any other weight of $H^0_{G_f}(\la)$ is lower than $\la$. If $H^0_{G_f}(\la) \neq 0$, $H^0_{G_f}(\la)^{(U_f)^+}$ cannot be zero by \Cref{unipotent fixed points}, and $H^0_{G_f}(\la)^{(U_f)^+} = (H^0_{G_f}(\la))_{\la}$ hence $H^0_{G_f}(\la)^{(U_f)^+} \simeq \la$ as a $T$-module.
\end{proof}

\begin{coro} \label{H^0 simple}
    Let $\la$ be a $G_f$-dominant weight.
    \begin{enumerate}
    \item $L_{G_f}(\la)$ is a simple $G_f$-module.
    \item There is an isomorphism of $T$-modules $L_{G_f}(\la)^{(U_f)^+} \simeq \la$, and $\la$ is the highest weight of $L_{G_f}(\la)$.
    \end{enumerate}
\end{coro}

\begin{proof}
    \begin{enumerate}
    \item If there are two distinct simple submodules $L_1$ and $L_2$ , then $L_1 \oplus L_2 \subset H^0_{G_f}(\la)$ so  $L_1^{(U_f)^+} \oplus L_2^{(U_f)^+} \subset H^0_{G_f}(\la)^{(U_f)^+}$ so $\dim H^0_{G_f}(\la)^{(U_f)^+} \geq 2$, which contradicts \Cref{main prop}.
 
    \item We know that $L_{G_f}(\la)^{U^+} \neq 0$ and $L_{G_f}(\la) \subset H^0_{G_f}(\la)$ so $L_{G_f}(\la)^{(U_f)^+}(\la) \subset H_{G_f}^0(\la)^{(U_f)^+}(\la)$. But $H_{G_f}^0(\la)^{(U_f)^+}(\la)$ has dimension $1$ by \Cref{main prop}, so $$L_{G_f}(\la)^{(U_f)^+}(\la) = H_{G_f}^0(\la)^{(U_f)^+}(\la) \simeq \la.$$ So $\la$ is a weight of $L_{G_f}(\la)$, and by \Cref{main prop} we know that all other weights are smaller.
    \end{enumerate}
\end{proof}

\begin{thm} \label{BW theorem non-inf}
    There is a bijection between simple $G_f$-modules (up to isomorphism) and $G_f$-dominant weights, that is, $X_+^{I_l}(T)$. In one direction, it sends a simple $G$-module $M$ to its highest weight, and in the other direction it sends a $G_f$-dominant weight $\la$ to the simple $G_f$-module $L_{G_f}(\la)$.
\end{thm}

\begin{proof}
    If $\la$ is $G_f$-dominant, $L_{G_f}(\la)$ is indeed simple by \Cref{H^0 simple}.
    The map $L_{G_f}(-)$ is injective because by \Cref{H^0 simple} the highest weight of $L_{G_f}(\la)$ is $\la$.
    Now let $V$ be a simple $G_f$-module. By \Cref{embeds}, there is a weight $\la$ (necessarily $G_f$-dominant) and a non-zero map $\phi \colon V \to H^0_{G_f}(\la)$. Given that $V$ is simple, $\phi$ is injective. By \Cref{H^0 simple}, it implies that $V = L_{G_f}(\la)$. So the map $L_{G_f}(-)$ is surjective.
\end{proof}

\paragraph{Classification of simple representations of non-reduced parabolic subgroups.}

Let $f$ be a non-infinitesimal $W$-function. The classification of simple representations of $G_f$ transposes \textit{mutatis mutandis} to simple representations of non-reduced parabolic subgroups $P_f$. The proof is exactly the same as for $G_f$.

\begin{defi}
    Let $\la \in X(T)$ be a weight. We will write $$L_{P_f}(\la) := \soc_{P_f} H^0_{P_f}(\la).$$
\end{defi}

\begin{prop} \label{BW thm parabolic}
    There is a bijection between simple $P_f$-modules (up to isomorphism) and $P_f$-dominant weights, that is, $X_+^{I_l}(T)$. In one direction, it sends a simple $P_f$-module $M$ to its highest weight, and in the other direction it sends a $P_f$-dominant weight $\la$ to the simple $P_f$-module $L_{P_f}(\la)$.
\end{prop}

\paragraph{Classification of simple representations of infinitesimal reductive groups.}

Here let $f$ be an infinitesimal $W$-function. We prove the classification of simple $G_f$-modules by following the arguments of \cite[II.3]{jantzen2003representations}, where Jantzen deals with the case of $G_r$. 

\begin{lemma} \label{lambda socle}
    Let $\la \in X(T)$.
    Then $\la = \soc_{B_f^+} H_{G_f}^0(\la)$.
\end{lemma}

\begin{proof}
    There is the decomposition $B_f^+ = U_f^+ \rtimes T_{r_l}$ where $U_f^+$ is unipotent and $T_{r_l}$ is diagonalizable. Thus, by \cite[II.3.11]{jantzen2003representations}, $\la \simeq \soc_{B_f^+} k[U_f^+] \otimes_k \la$, where $U_f^+$ acts by regular representation on $k[U_f^+]$ and trivially on $\la$, and $T_{r_l}$ acts by the conjugation representation on $k[U_f^+]$ and as given on $\la$. But $H^0_{G_f}(\la)$ is precisely of this form by \Cref{pushforward tensor}, hence $\la = \soc_{B_f^+} H_{G_f}^0(\la)$.
\end{proof}

\begin{lemma} \label{simple module for infinitesimal groups}
    Let $\la \in X(T)$. The following assertions hold:\begin{enumerate}

    \item $H^0_{G_f}(\la)^{U_f^+} \simeq \la$.

    \item $L_{G_f}(\la)$ is simple.

   \item  $L_{G_f}(\la)^{U_f^+} \simeq \lambda$.
    \end{enumerate}
\end{lemma}

\begin{proof}
    $H^0_{G_f}(\la)^{U_f^+} \subset \soc_{B_f^+} H^0_{G_f}(\la) = \la$ and by \Cref{unipotent fixed points}, $H^0_{G_f}(\la)^{U_f^+} \neq 0$ so $H^0_{G_f}(\la)^{U_f^+}$ has dimension $1$.
    The proof of the rest is the same as in the non-infinitesimal case, in \Cref{H^0 simple}.
\end{proof}

\begin{thm} \label{BW theorem inf}
    There is a bijection between simple $G_f$-modules (up to isomorphism) and $X(T_{r_l})$ which, in one direction sends a simple $G_f$-module $M$ to its highest weight, and in the other associates to $\la \in X(T_{r_l})  $, the simple $G_f$-module $L_{G_f}(\la)$.
\end{thm}

\begin{proof}
    The proof is the same as in \Cref{BW theorem non-inf}.
\end{proof}

For the next section, we need an additional result about coinduction for infinitesimal reductive groups, which is also an extension of the analogous result for $G_r$ in \cite[II.3]{jantzen2003representations}.

\begin{defi}[{\cite[I.8.14]{jantzen2003representations}}]
    The restriction functor also admits a right adjoint which is called coinduction. We write

    $$\tilde{H}_{G_f}^0(\la) := R^0\coind_{B_f}^{G_f}(\la) = \Dist(G_f) \otimes_{\Dist(B_f)} \la.$$
\end{defi}

For finite group schemes, coinduction is dual to induction (it can be seen as an instance of the isomorphism between $f_*$ and $f_!$ for proper maps in six-functor formalism, applied to the map of stacks $BL \to BK$).

\begin{lemma}[{\cite[I.8.15]{jantzen2003representations}}] \label{coinduction finite group}
    Let $K$ be a group scheme which is finite as a scheme. Let $L$ be a closed subgroup of $K$ and $M$ be a finite-dimensional $L$-module.
    Then there is an isomorphism of $G$-modules $$\coind_L^K M \simeq (\ind_L^K M^*)^*.$$
\end{lemma}

It allows us to deduce the following lemma which will be useful in the next section.

\begin{lemma} \label{only quotient semi-simple}
    The quotient $\tilde{H}^0_{G_f}(\la)/(\rad_{G_f} \tilde{H}_{G_f}^0(\la))$ is a simple $G_f$-module.
\end{lemma}

\begin{proof}
    There is an isomorphism $\tilde{H}_{G_f}(\la)/(\rad_{G_f} \tilde{H}_{G_f}(\la)) \simeq (\soc_{G_f} \tilde{H}^0_{G_f}(\la^*))^*$ by duality between socle and radical. But by \Cref{coinduction finite group}, $\tilde{H}^0_{G_f}(\la^*) \simeq (H^0_{G_f}(\la))^*$, hence by \Cref{simple module for infinitesimal groups}, $\tilde{H}_{G_f}(\la)/(\rad_{G_f} \tilde{H}_{G_f}(\la))$ is simple.
\end{proof}

\subsubsection{Restriction of simple modules under different non-reduced structures} \label{section transfer}

When $G$ is a (reduced) reductive group, a simple $G$-module $L_G(\la)$ stays simple under restriction to $G_r$ when $\la$ is $p^r$-restricted (\cite[II.3]{jantzen2003representations}). The study of its $G_r$-module structure can bring interesting information, and in particular gives a proof of Steinberg theorem. We show that an analogous result holds when restricting a simple $G_f$-module $L_{G_f}(\la)$ to $(G_f)_r$. Again, we follow and extend the arguments of Jantzen. More generally, we show that the restriction of $L_{G_f}(\la)$ to another non-reduced reductive subgroup contained in $G_f$ remains simple if $\la$ satisfies an appropriate restriction condition. It allows us in particular to show that simple $P_f$-modules always restrict to simple $G_f$-modules, which generalizes the analogous result for $G_rB$ and $G_rT$. It also allows us to prove a Steinberg theorem for non-reduced reductive groups in \Cref{section steinberg}.
Let $f$ be a non-infinitesimal $W$-function with values $r_0 < \dots < r_{l-1} < r_l = +\infty$, and with characteristic filtration $I_l \subset \dots \subset I_1 \subset I_0 = S$.
Let us begin by a series of lemmas.

\begin{lemma}[{\cite[I.6.16]{jantzen2003representations}}] \label{soc is stable}
    Let $G$ be a group scheme. Let $N$ be a normal subgroup of $G$, and let $H$ be a subgroup of $G$ such that $H(k)$ is dense in $H$ (for instance if $H$ is reduced), such that $G = NH$.
    Then for any $G$-module $M$, $\soc_{N} M$ lifts to a $G$-submodule of $M$. Moreover $\soc_G M \subset \soc_N M$.
    In particular, the restriction to $N$ of a semi-simple $G$-module is also semi-simple.
\end{lemma} 

\begin{lemma} \label{lemma Dist}
   Let $\la \in X_+^{I_l}$ and let $v \in L_{G_f}(\la)_{\la}$, $v \neq 0$. For all $r \geq r_{l-1}$, there is an isomorphism of $(G_{f})_r$-modules $$L_{(G_{f})_r}(\la) \simeq \Dist((G_{f})_r)v.$$
\end{lemma}

\begin{proof}
    $\Dist((G_f)_rv$ is a $(G_f)_r$-submodule by \Cref{Dist submodule}. By \Cref{soc is stable}, given that $G_f = (G_{f})_rG^l$ where $G^l$ is reduced, $\res_{(G_{f})_r}^{G_f} L_{G_f}(\la)$ is a semi-simple $(G_{f})_r$-module, hence so is $\Dist((G_{f})_r)v$. As $(B_{f})_r^+$-modules, $\la = kv$ so there is a natural surjective map 

    $$\tilde{H}^0_{(G_{f})_r}(\la) = \Dist((G_{f})_r) \otimes_{\Dist((B_{f})_r^+)} \la \to \Dist((G_{f})_r)v$$
    But $\Dist((G_{f})_r)v$ is semi-simple so by \Cref{only quotient semi-simple}, $\Dist((G_{f})_r)v \simeq L_{(G_{f})_r}(\la)$.
\end{proof}

\begin{lemma} \label{lemma for 3.15}    Let $\alpha \in I^l$, let $r \geq 1$, let $\la \in X^{I_l}(T)$. Then

    $$\Dist(U_{-\al, r})v = \bigoplus_{0 \leq n < p^r}L_{G_f}(\la)_{\la - n\al}$$

\end{lemma}

\begin{proof}
    By \Cref{Dist submodule}, $\Dist(G_f)v$ is a $G_f$-submodule of $L_{G_f}(\la)$. But $L_{G_f}(\la)$ is a simple $G_f$-module, so $\Dist(G_f)v = L_{G_f}(\la)$. By \Cref{base Dist}, $L_{G_f}(\la) = \Dist(U_f)v$. Then, by \Cref{shift}, $$\Dist(U_{-\al, r})v = \bigoplus_{0 \leq n < p^r}L_{G_f}(\la)_{\la - n\al}$$
\end{proof}

Let us denote $W_G$ the Weil group of a reductive group $G$, generated by the reflections $s_{\alpha}$ for $\alpha \in R$. The Weyl group is isomorphic to the quotient of the normalizer of $T$, $N_G(T)(k)/T(k)$. Any $w \in W$ has a representative $\dot{w}$ in $N_G(T)(k)$. 

\begin{lemma}[{\cite[II.1.19(1)]{jantzen2003representations}}] \label{action of W}
    Let $G$ be a reductive group. Let $\dot{w} \in N_G(T)(k)$ be a representative of some $w \in W_G$. Then for all $\la \in X(T)$, $$\dot{w}M = M_{w(\la)}.$$
\end{lemma}

We can then deduce the following core result, which generalizes the case of reduced reductive groups (\cite[II.3.15]{jantzen2003representations}).

\begin{prop} \label{generalization 3.15}
    Let $f$ be a non-infinitesimal $W$-function of width $l$. Let $r \geq r_{l-1}$ (with the convention $r_{-1} = 0$) and let $\lambda \in X_{r}^{I_l}(T)$.
    Then $$\res_{(G_{f})_{r}}^{G_f}L_{G_f}(\la) \simeq L_{(G_{f})_{r}}(\la).$$ 
\end{prop}

\begin{proof}
    Set $L = \Dist((G_{f})_r)v$, a $(G_f)_r$-submodule of $L_{G_f}(\la)$. By \Cref{lemma Dist}, $L \simeq L_{(G_{f})_r}(\la)$. Let us show that $L$ lifts to a $G_f$-submodule of $L_{G_f}(\la)$. Given that $L_{G_f}(\la)$ is simple, it will prove that $L_{(G_{f})_r}(\la) \simeq \res_{(G_{f})_r}^{G_f}L_{G_f}(\la)   $.

    Given that $G_f = (G_{f})_rG^l$, we only have to prove that $L$ is stable under $G^l$ and, given that $G^l$ is reduced, it is enough to prove that it is stable under $G^l(k)$ because $G^l(k)$ is dense in $G^l$.
    By the Bruhat decomposition $$G_l(k) = \bigcup_{w \in W_{G_l}} B_l^+(k)\dot{w}B_l^+(k),$$ it is enough to prove that $L$ is stable under $\dot{w}$ for all $w \in W_{G^l}$, and under $B_l^+(k)$. But $W_{G_l}$ is generated by the $s_{\al}$ for $\al \in I_l$, so it is enough to show that $L$ is stable under $s_{\al}$ for all $\al \in I_l$.
    Moreover, if $g \in G^l(k)$, $gL$ is a simple $(G_f)_{r}$-submodule of $L_{G_f}(\la)$. So either $gL = L$ or $gL \cap L = \{0\}$. Hence, to prove that $L$ is stable under $g$, it is enough to show that $gv \in L$ and $gv \neq 0$.

    First, if $g \in B_l^+(k) \subset B_f^+(k)$, then $gv \in k^{\times}v$ because $v \in L_{G_f}(\la)_{\la}$.

    Now take a simple root $\alpha \in I_l$. By \Cref{action of W}, $\dot{s_{\al}}v \in L_{G_f}(\la)_{s_{\alpha} \la}$, and $\dot{s_{\al}}v \neq 0$. 
    By \Cref{lemma for 3.15}, for all $n < p^{r}$, 
    $$L_{G_f}(\la)_{\la - n\alpha} \subset \Dist(U_{-\al,r})v \subset \Dist((G_{f})_r)v.$$ 
    Now $s_{\alpha} \la = \la - \<\la, \alpha^{\vee}\> \alpha$ and $0 \leq \<\la, \alpha^{\vee}\> < p^{r}$ since $\la \in X_{r}^{I_l}(T)$, so $$L_{G_f}(\la)_{s_{\al}\la} \subset \Dist((G_f)_r)v.$$ But $\dot{s_{\al}}v \in L_{G_f}(\la)_{s_{\alpha} \la}$ so $\dot{s_{\al}}v \in \Dist((G_f)_r)v$.
\end{proof}

\begin{rk} \label{3.15 parabolic}
    \Cref{generalization 3.15} and its proof also work for the parabolic subgroup $P_f$ instead of $G_f$. We have the isomorphism $$\res_{(P_f)_r}^{P_f} L_{P_f}(\la) \simeq L_{(P_f)_r}(\la)$$ for all $\la \in X_{r}^{I_l}(T)$ and all $r \geq r_{l-1}$.
\end{rk}

\Cref{generalization 3.15} can be extended to the kernels of $G_f$ of degree $r \leq r_{l-1}$ by extending the restriction condition on the weight $\la$.

\begin{coro}
    Let $r \in \N$. Choose $1 \leq k \leq l$ such that $r_{k-1} \leq r < r_k$ (with the convention $r_{-1} = 0)$. Then for all $\la \in X_r^{I_k}(T)$, $$\res_{(G_f)_r}^{G_f} L_{G_f}(\la) \simeq L_{(G_f)_r}(\la).$$ 
\end{coro}

\begin{proof}
    Consider the chain of inclusions $(\overl[G]^k)_r = (G_f)_r \subset G_f \subset \overl[G]^k$. By \Cref{generalization 3.15} we know that for all $\la \in X_{r_k}^{I_k}$, $\res^{\overl[G]^k}_{(\overl[G]^k)_r} L_{\overl[G]^k}(\la) = L_{(\overl[G]^k)_r}(\la)$. In particular it is simple. Hence $\res_{G_f}^{\overl[G]^k} L_{\overl[G]^k}(\la)$ is also simple. By comparing the highest weight, it is equal to $L_{(G_f)_r}(\la)$.
\end{proof}

It is also possible to deduce a restriction result between different non-reduced reductive subgroups, under a suitable restriction condition on the weight $\la$.

\begin{defi}
    We define (with convention $I_{l+1} = \emptyset$) $$X_{r_0, \dots, r_l}^{I_0, \dots, I_l}(T) := \{\la \in X(T), \forall k \in \llbracket 0, l \rrbracket, \forall \alpha \in I_k \backslash I_{k+1}, 0 \leq \<\la, \alpha^{\vee}\> < p^{r_k} \}$$ 
    If $g$ is a $W$-function with characteristic filtration $(I_k)_k$ and values $(r_k)_k$, we denote $$X_g(T) := X_{r_0, \dots, r_l}^{I_0, \dots, I_l}(T)$$ (if $g$ is non-infinitesimal, $r_l$ is an integer and if $g$ is infinitesimal we take the natural convention $p^{+\infty} = +\infty$).
\end{defi}

\begin{coro} \label{restriction between parabolic subgroups}
    \begin{enumerate}
    \item Let $0 \leq j < k \leq l$. For all $\la \in X_{r_{j}, \dots, r_{k-1}}^{I_j, \dots, I_{k-1}}(T)$, $$\res_{\overl[G]^k}^{\overl[G]^{j}} L_{\overl[G]^{j}}(\la) \simeq L_{\overl[G]^k}(\la).$$ 
    \item If $g$ is a $W$-function, infinitesimal or not, and if $\lambda \in X_{g}(T)$, then $$\res_{G_{g}}^{G}(L_G(\la)) \simeq L_{G_{g}}(\la).$$
    \end{enumerate}
\end{coro}

\begin{proof}
    \begin{enumerate}
    \item First, let us prove it for $k = j + 1$. Consider the chain of inclusions $$(\overl[G]^j)_{r_{j}} \subset \overl[G]^{j+1} \subset \overl[G]^j.$$ By \Cref{generalization 3.15} we know that for all $\la \in X_{r_{j}}^{I_j}$, $$\res^{\overl[G]^j}_{(\overl[G]^j)_{r_j}} L_{\overl[G]^j}(\la) \simeq L_{(\overl[G]^j)_{r_j}}(\la).$$ In particular it is simple. Hence $\res_{\overl[G]^{j+1}}^{\overl[G]^j} L_{\overl[G]^j}(\la)$ is also simple, hence by \Cref{BW theorem non-inf} and by comparing the highest weight, it is isomorphic to $L_{\overl[G]^{j+1}}(\la)$.
    The first assertion then follows by immediate induction. 
    
    \item If $g$ is non-infinitesimal and if $\la \in X_g(T)$, the second assertion is included in the first. If $g$ is infinitesimal, we have the chain of inclusions $G_g \subset \overl[G]^l \subset G$. By the first assertion, given that $\la \in X_g(T) = X_{r_0, \dots r_{l}}^{I_0, \dots, I_l}(T) \subset X_{r_{0}, \dots, r_{l-1}}^{I_0, \dots, I_{l-1}}(T)$,   $$\res^G_{\overl[G]^l} L_G(\la) \simeq L_{\overl[G]^l}(\la).$$

    Then by \Cref{generalization 3.15}, given that $\la \in X_g(T) \subset X_{r_l}^{I_l}(T)$ and that $G_g = (\overl[G])_{r_l}$, $$\res^{\overl[G]^l}_{G_g} L_{\overl[G]^l}(\la) \simeq L_{G_g}(\la).$$
    The composition of the two last restrictions gives the second assertion.
    \end{enumerate}
\end{proof}

We also deduce that simple $P_f$-modules restrict to simple $G_f$-modules.

\begin{prop}
    For all $\la \in X^{I_l}_+(T)$, $L_{G_f}(\la) \simeq \res_{G_f}^{P_f}(L_{P_f}(\la))$.
\end{prop}

\begin{proof}
    Let $\la \in X^{I_l}_+(\la)$. There exists $r \geq r_{l-1}$ such that $\la \in X^{I_l}_{r}(T)$ (e.g. the maximum between $r_{l-1}$ and the maximum value of $\<\la, \alpha^{\vee}\>$ for $\alpha \in I_l$). By \Cref{3.15 parabolic}, we know that $\res_{(G_f)_r}^{P_f} L_{P_f}(\la) \simeq L_{(G_{f})_r}(\la)$. In particular it is simple. Hence $\res_{G_f}^{P_f} L_{P_f}(\la)$ is also simple. By comparing the highest weight, it is isomorphic to $L_{G_f}(\la)$.
\end{proof}

\subsubsection{Steinberg theorem for non-reduced reductive subgroups} \label{section steinberg}

For a (reduced) reductive group $G$ and a $p$-adic expansion of a weight $\la$, Steinberg theorem gives a decomposition of the simple $G$-modules $L_G(\la)$ with Frobenius twists of simple $G$-module with $p$-restricted highest weight (\cite[II.3.16]{jantzen2003representations}). We extend Steinberg theorem to non-reduced reductive groups. When $\la$ admits a suitable $p$-adic expansion, it allows us to decompose a simple $G_f$-module as Frobenius twists of simple modules over reduced reductive groups. Again we follow and extend the argument of Jantzen.

Let $f$ be a non-infinitesimal $W$-function with values $r_0 < \dots < r_{l-1} < r_l = +\infty$, and with characteristic filtration $I_l \subset \dots \subset I_1 \subset I_0 = S$. Write $G_f = G^0_{r_0} \dots G^{l-1}_{r_{l-1}} G^l$.
We begin by a series of lemmas.

\begin{lemma}[{\cite[2.14 (3)]{jantzen2003representations}}] \label{components of the socle}
    Suppose that $k$ is algebraically closed. Let $H$ be a group scheme, $M$ a $H$-module, $E$ a simple $H$-module. Then there is an isomorphism of $H$-modules
    $$\Hom_H(E, M) \otimes E \iso (\soc_H M)_E,$$
    where $(\soc_H M)_E$ denotes the $E$-isotypic component of $\soc_H M$, that is the sum of all simple $H$-submodules isomorphic to $E$.
\end{lemma}

\begin{lemma}[{\cite[3.16 (1)]{jantzen2003representations}}] \label{decomposition socle}
    Let $r \geq r_{l-1}$. Suppose that there is a system of representatives of $X(T)/p^rX(T)$, denoted $X_r'(T)$, inside $X_r^{I_l}(T)$.
    Let $M$ be a $G_f$-module. The $(G_f)_r$-submodule $\soc_{(G_f)_r} M$ lifts to a $G_f$-submodule of $M$ and there is an isomorphism of  $G_f$-modules $$\soc_{(G_f)_r} M \simeq \bigoplus_{\mu \in X'_{r}(T)} L_{G_f}(\mu) \otimes \Hom_{(G_f)_r}(L_{G_f}(\mu), M)$$
\end{lemma}

\begin{proof}

    By \Cref{BW theorem inf}, simple $(G_f)_r$-modules are given by the $L_{(G_f)_r}(\mu)$ for all $\mu \in X_{r}'(T)$. Let $\mu \in X'_r(T)$. By \Cref{components of the socle}, there is an isomorphism of $(G_f)_r$-submodules
    $$(\soc_{(G_f)_r} M)_{L_{(G_f)_r}(\mu)} \simeq \Hom_{(G_f)_r}(L_{(G_f)_r}(\mu), M) \otimes L_{(G_f)_r}(\mu).$$
    But $r \geq r_{l-1}$ so $G_f = (G_f)_rG^l$ so by \Cref{soc is stable}, the $(G_f)_r$-submodule $\soc_{(G_f)_r}M$ lifts to a $G_f$-submodule of $M$. Then, by \Cref{generalization 3.15}, for all $\mu \in X_{r}'(T)$, which is a subset of $X_{r}^S(T)$ hence of $X_{r}^{I_l}(T)$, the $(G_f)_r$-module $L_{(G_f)_r}(\mu)$ lifts to a $G_f$-module, $L_{G_{f}}(\mu)$, so there is an isomorphism of $G_f$-modules
    
    $$(\soc_{(G_f)_r} M)_{L_{G_{f}}(\mu)} \simeq \Hom_{(G_f)_r}(L_{G_{f}}(\mu), M) \otimes L_{G_{f}}(\mu).$$
\end{proof}

\begin{rk}[{\cite[II.3.15 Remark 2, II.9.7]{jantzen2003representations}}] \label{rk covering}
    If $G$ is semi-simple and simply connected, $X_{r}^S(T)$ is a system of representatives of $X(T)/p^{r}X(T)$. In general, one can replace $G$ by a covering group $\tilde{G}$ such that $X_r^S(T)$ contains a system of representatives of $X(T)/p^{r}X(T)$.  There is a subgroup $\tilde{G^l}$ covering $G^l$. Let $Z$ denote the  kernel of $\tilde{G} \to G$. We have then $\tilde{G}/Z \iso G$ and $\tilde{G}^l/Z \iso G^l$. In general $G_r$ is not a factor group of $\tilde{G}_r$. However, when $f$ is a non-infinitesimal $W$-function, $G_f$ is a factor group of $\tilde{G}_f$. There are isomorphisms of schemes $U^0_f \times G^l \times (U^0_f)^+ \iso G_f$ and $\tilde{U^0_f} \times \tilde{G^l} \times \tilde{(U^0_f)^+} \iso \tilde{G}_f$. The covering $\tilde{G} \to G$ restricts to isomorphisms $\tilde{U} \to U$ and $\tilde{U^+} \to U^+$, hence also $\tilde{U^0_f} \iso U^0_f$ and $\tilde{(U^0_f)^+} \iso (U^0_f)^+$. This implies that $\tilde{G}_f/Z \iso G_f$. It follows that we can identify the category of $G_f$-modules with the full subcategory of $\tilde{G}_f$-modules whose weights belong to $X(T) \subset X(\tilde{T})$.
\end{rk}

\begin{rk}
    Let $r \geq r_{l-1}$. If $M$ is a simple $G_f$-module, then $\soc_{(G_f)_r} M = M$ since the socle is non-zero, hence there is $\mu \in X_r'(T)$ such that $$M \simeq L(\mu) \otimes \Hom_{(G_f)_r}(L(\mu), M).$$
\end{rk}

\begin{lemma} \label{lemma trivial action}
    Let $r \geq r_{l-1}$ and let $M$ be $G_f$-module. Suppose that $(G_f)_r$ acts trivially on $M$. 
Then there is a $G^l$-module $V$ such that $M = (F^r_{G_f})^* V$.
\end{lemma} 

\begin{proof}
    If $(G_f)_r$ acts trivially on $M$, then the action of $G_f$ on $M$ factorizes through $G_f/(G_f)_r$. But $r \geq r_{l-1}$ so $G_f/(G_f)_r \simeq G^l/(G^l)_r \simeq G^l$.
\end{proof}

\begin{lemma} \label{steinberg 0}

    Let $r \geq r_{l-1}$, let $\la \in X_{r}^{I_l}(T)$ and $\mu \in X_{+}^{I_l}(T)$. 
    Then there is an isomorphism of $G_f$-modules $$L_{G_f}(\la + p^{r}\mu) \simeq L_{G_f}(\la) \otimes (F^{r}_{G_f})^* L_{G^l}(\mu).$$
\end{lemma}

\begin{proof}
    We can replace $G$ by a covering group. Thus by \Cref{rk covering} we may a fortiori assume that there exists a system of representatives $X_{r}'(T)$ of $X(T)/p^{r}X(T)$ inside $X_{r}^{I_l}(T)$ and that $\la \in X_{r}'(T)$.
    Denote $M = L_{G_f}(\la + p^r\mu)$, a simple $G_f$-module. Then by \Cref{decomposition socle}, $$M \simeq \Hom_{(G_f)_r}(L_{G_f}(\la'), M) \otimes L_{G_f}(\la')$$ with $\la' \in X_{r}'(T)$.
    $(G_f)_r$ acts trivially on $\Hom_{(G_f)_r}(L(\la'), M)$, so by \Cref{lemma trivial action} $$\Hom_{(G_f)_r}(L(\la'), M) = (F_{G_f}^r)^* V$$ for a $G^l$-module $V$.
    $M \simeq L_{G_f}(\la') \otimes (F_{G_f}^r)^* V$ is a simple $G_f$-module so $V$ is a simple $G^l$-module so $V \simeq L_{G^l}(\mu')$ with $\mu' \in X_+^{I_l}(T)$.
    We compare the highest weights: $\la + p^{r}\mu = \la' + p^{r}\mu'$ so $\la - \la' \in p^rX^{I_l}(T)$ but $\la, \la' \in X_{r}'(T)$ so $\la = \la'$, hence $\mu = \mu'$. 
    Thus $$L_{G_f}(\la + p^{r}\mu) \simeq L_{G_f}(\lambda) \otimes (F^{r}_{G^l})^* L_{G^l}(\mu).$$
\end{proof}

We can now prove the Steinberg theorem for non-reduced reductive subgroups. In this situation, we need to impose additional restrictions on the $p$-adic expansion: not all $G_f$-dominant weights will admit such $p$-adic expansion. In that case, it allows to reduce the computation of a simple $G_f$-module to the computation of simple modules over the reduced reductive groups $G^k$ which enter into the usual decomposition of $G_f$. 

\begin{thm}\label{steinberg}
    Let $\la = \la_0 + \la_1 p + \dots + \la_n p^n$ where for all $r_{k-1} \leq r < r_k$, $\la_r \in X_{1}^{I_k}(T)$ (with usual convention $r_{-1} = 0$).
    Then there is an isomorphism of $G_f$-modules
    $$L_{G_f}(\la_0+\dots +\la_n p^n) \simeq L_{G_f}(\la_0) \otimes (F^1_{G_f})^* L_{G_f/(G_f)_{(1)}}(\la_1) \otimes \dots \otimes (F^n_{G_f})^* L_{G_f/(G_f)_{(n)}}(\la_n).$$ 
    For all $r \in \llbracket 0, n \rrbracket$, denote by $k_r$ the unique integer such that $r_{k_r-1} \leq r < r_{k_r}$ (then $\la_r \in X_1^{I_{k_r}}(T)$). Then the isomorphism can be reformulated as
    $$L_{G_f}(\la_0+\dots +\la_n p^n) \simeq \res_{G_f}^{G^0} L_{G^0}(\la_0) \otimes \res_{G_f}^{\overl[G]^{k_1}} (F^1_{\overl[G]^{k_1}})^* L_{G^{k_1}}(\la_1) \otimes \dots \otimes \res_{G_f}^{\overl[G]^{k_n}} (F^n_{\overl[G]^{k_n}})^* L_{G^{k_n}}(\la_n).$$
\end{thm}

\begin{proof}
    For the first assertion, we reason by induction on the width $l$ of $G_f$.
    If $l = 0$, $G_f = G$ and it is the usual Steinberg formula for reduced reductive groups. Now take $l \geq 1$ and suppose that the assertion is valid in width $l-1$.
    Using \Cref{steinberg 0} successively with $r = n$, $r = n-1$, ..., $r = r_{l-1}$, we get $$L_{G_f}(\la_0 + \dots + \la_n p^n) \simeq L_{G_f}(\la_0 + \dots \la_{r_{l-1}-1} p^{r_{l-1}-1}) \otimes (F^{r_{l-1}}_{G_f})^* L_{G_l}(\la_{r_{l-1}}) \otimes \dots \otimes (F^{n}_{G_f})^* L_{G_l}(\la_{n}).$$
    Now $\la_0 + \dots + \la_{r_{l-1}-1} p^{r_{l-1}-1} \in X_{r_{l-1}}^{I_{l-1}}(T)$ so by \Cref{restriction between parabolic subgroups} $$L_{\overl[G]^l}(\la_0 + \dots \la_{r_{l-1}-1} p^{r_{l-1}-1}) \simeq \res_{\overl[G]^l}^{\overl[G]^{l-1}} L_{\overl[G]^{l-1}}(\la_0 + \dots \la_{r_{l-1}-1} p^{r_{l-1}-1}).$$
    By induction hypothesis, we know that $L_{\overl[G]^{l-1}}(\la_0 + \dots + \la_{r_{l-1}-1} p^{r_{l-1}-1})$ is isomorphic to $$L_{\overl[G]^{l-1}}(\la_0) \otimes (F^1_{\overl[G]^{l-1}})^* L_{\overl[G]^{l-1}/\overl[G]^{l-1}_{(1)}}(\la_1) \otimes \dots \otimes ({F}^{r_{l-1}-1}_{\overl[G]^{l-1}})^* L_{\overl[G]^{l-1}/\overl[G]^{l-1}_{(r_{l-1}-1)}}(\la_{r_{l-1}-1})$$
    Let $0 \leq k \leq r_{l-1}-1$. The commutation of the following diagram 

\[\begin{tikzcd}
	{\overl[G]^{l}} & {\overl[G]^{l}/(\overl[G]^{l})_k } \\
	{\overl[G]^{l-1}} & {\overl[G]^{l-1}/(\overl[G]^{l-1})_k}
	\arrow["{F_{\overl[G]^{l}}^k}", from=1-1, to=1-2]
	\arrow["i"', from=1-1, to=2-1]
	\arrow["j", from=1-2, to=2-2]
	\arrow["{F_{\overl[G]^{l-1}}^k}"', from=2-1, to=2-2]
\end{tikzcd}\]
    shows that $$\res_{\overl[G]^l}^{\overl[G]^{l-1}} (F_{\overl[G]^{l-1}}^k)^* (L_{\overl[G]^{l-1}/(\overl[G]^{l-1})_k}(\la_k)) \simeq (F_{\overl[G]^{l}}^k)^* \res_{\overl[G]^{l}/(\overl[G]^{l})_k}^{\overl[G]^{l-1}/(\overl[G]^{l-1})_k} (L_{\overl[G]^{l-1}/(\overl[G]^{l-1})_k}(\la_k)).$$
    But $\la_k \in X_{1}^{I_{k_i}}(T) \subset X^{I_{l-1}}_{1}(T)$ and $r_{l-1} - k \geq 1$ so by \Cref{restriction between parabolic subgroups} $$\res_{\overl[G]^{l}/(\overl[G]^{l})_k}^{\overl[G]^{l-1}/(\overl[G]^{l-1})_k} (L_{\overl[G]^{l-1}/(\overl[G]^{l-1})_k}(\la_k)) \simeq L_{\overl[G]^{l}/(\overl[G]^{l})_k}(\la_k).$$
    The restriction functor is symmetric monoidal so we get that $L_{\overl[G]^l}(\la_0 + \dots \la_{r_{l-1}-1} p^{r_{l-1}-1})$ is isomorphic to $$L_{\overl[G]^l}(\la_0) \otimes (F^1_{\overl[G]^l})^* L_{\overl[G]^l/(\overl[G]^l)_{(1)}}(\la_1) \otimes \dots \otimes ({F}^{r_{l-1}-1}_{\overl[G]^l})^* L_{\overl[G]^l/(\overl[G]^l)_{(r_{l-1}-1)}}(\la_{r_{l-1}-1}),$$
    hence the first assertion.

    Then let $r \in \llbracket 0, n \rrbracket$, such that $\la_r \in X_1^{I_{k_r}}(T)$. By \Cref{restriction between parabolic subgroups} $$\res_{G_f/(G_f)_r}^{G^{k_r}} (L_{G^{k_r}}(\la_r)) \simeq L_{G_f/(G_f)_r}(\la_r),$$ 
    and the commutation of the following diagram
\[\begin{tikzcd}
	{G_f} & {G_f/(G_f)_{r}} \\
	{\overl[G]^k} & {\overl[G]^{k_r}/(\overl[G]^{k_r})_r \simeq G^{k_r}}
	\arrow["{F_{G_f}^r}", from=1-1, to=1-2]
	\arrow["i"', from=1-1, to=2-1]
	\arrow["j", from=1-2, to=2-2]
	\arrow["{F_{\overl[G]^{k_r}}^r}"', from=2-1, to=2-2]
\end{tikzcd}\]
    implies the isomorphism $$\res_{G_f}^{\overl[G]^{k_r}} (F_{\overl[G]^{k_r}}^r)^* (L_{G^{k_r}}(\la_r)) \simeq (F_{G_f}^r)^* \res_{G_f/(G_f)_r}^{G^{k_r}} (L_{G^{k_r}}(\la_r)).$$
    This proves the second statement.
\end{proof}

\begin{ex} \label{ex steinberg}
    Let $0 \leq k \leq l-1$, let $r_{k-1} \leq r < r_k$. If $\la \in X_1^{I_k}(T)$, then there is an isomorphism of $G_f$-modules $$L_{G_f}(p^r\la) \simeq \res_{G_f}^{\overl[G]^k} (F^r_{\overl[G]^k})^* L_{G^k}(\la).$$
\end{ex}

\begin{rk} \label{steinberg parabolic}
    \Cref{steinberg} works also by replacing everywhere $G_f$ by the non-reduced parabolic subgroup $P_f$.
\end{rk}

\subsubsection{A filtration of the tangent bundle with simple modules} \label{section alternative filtration}

In \Cref{section filtration tangent}, we constructed a filtration of the tangent bundle of $X = G/P_f$ whose factors are Frobenius twists of $\p^k/\p^{k+1}$ ($\overl^{k+1}$-modules). We established a formula for the induction of this factors, and it allowed us to study the first cohomology groups of $T_X$ in \Cref{section cohomology tangent bundle}. We show here that we can refine this filtration into a filtration with factors which are simple $P_f$-modules whose induction can also be computed using Steinberg theorem for non-reduced parabolic subgroups.
The study of cohomology groups of \Cref{section cohomology tangent bundle} would have worked exactly the same way with this new filtration.

\begin{prop}
    The $P_f$-module $T_X$ (which identifies with $T_{BP_f \sto BG}$) admits a filtration whose factors are $P_f$-modules of the form $L_{P_f}(p^{r_{k-1}} \la)$ for $1 \leq k \leq l$ and $\la \in X_1^{I_k}(T)$.
\end{prop}

\begin{proof}
    The Jordan-Hölder filtration (\cite[I.2.14]{jantzen2003representations}) gives a filtration of $T_{B\overl^{k} \sto B\overl^{k-1}}$ whose factors are simple $\overl^k$-modules. By \Cref{BW thm parabolic}, those simple $\overl^k$-modules are of the form $L_{\overl^k}(\mu)$ for a $\overl^k$-dominant weight $\mu$. Moreover by \Cref{restriction between parabolic subgroups}, $\res_{P_f}^{\overl^k} L_{\overl^k}(\mu)$ identifies with the simple $P_f$-module $L_{P_f}(\mu)$. But $\mu$ is also a weight of $T_X$, and in \Cref{section weight tangentn} we showed that the weights of $T_X$ were of the form $p^{r_{k-1}} \la$ for $1 \leq k \leq l$ and $\la \in R^{k-1} \backslash R^{k}$. Finally $\mu$ is $\overl^k$-dominant, that is, $\mu \in X_+^{I_k}(T)$, hence $\la \in X_1^{I_k}(T)$.
\end{proof}

\begin{prop} \label{induction formula}
    Let $0 \leq k \leq l-1$, let $r_{k-1} \leq r < r_{k}$ and let $\la \in X_1^{I_{k}}(T)$. Then there is an isomorphism of $G$-modules
    $$\ind_{P_f}^G(L_{P_f}(p^r \la)) \simeq (F_{P^0}^{r_0})^*  \ind_{P^1}^{P^0}  (F_{P^1}^{r_1-r_0})^*  \ind_{P^2}^{P^1}  \dots (F_{P^{k-1}}^{r_{k-1} - r_{k-2}})^*  \ind_{P^{k}}^{P^{k-1}} (L_{P^{k}}(\la)).$$
\end{prop}

\begin{proof}
    By \Cref{ex steinberg} and \Cref{steinberg parabolic}, $$L_{P_f}(p^r \la) \simeq \res_{P_f}^{\overl[P]^{k}} (F^{r}_{\overl[P]^{k}})^* L_{P^{k}}(\la).$$

    By \Cref{indk}, $$\ind_{P_f}^{\overl^k} \res_{P_f}^{\overl[P]^{k}} (F^{r}_{\overl[P]^{k}})^* L_{P^{k}}(\la) \simeq (F^{r}_{\overl[P]^{k}})^* L_{P^{k}}(\la).$$

    Then, by \Cref{induction}, $$\ind_{\overl^k}^{G} (F^{r}_{\overl[P]^{k}})^* L_{P^{k}}(\la) \simeq (F_{P^0}^{r_0})^*  \ind_{P^1}^{P^0}  (F_{P^1}^{r_1-r_0})^*  \ind_{P^2}^{P^1}  \dots (F_{P^{k-1}}^{r_{k-1} - r_{k-2}})^*  \ind_{P^{k}}^{P^{k-1}} (L_{P^{k}}(\la)).$$

\end{proof}

\section{Deformations of projective homogeneous varieties} \label{section mixed characteristic deformations}

\subsection{Deformations of non-reduced parabolic subgroups are obstructed} \label{def non reduced}

Non-reduced parabolic subgroups of reductive groups are built from Frobenius kernels of root subgroups, which are isomorphic to Frobenius kernels of the additive group $\G_a$. Those subgroups exist only in characteristic $p$: in higher characteristic, the additive group $\G_a$ admits no non-trivial flat subgroup schemes. As a consequence, non-reduced parabolic subgroups are also a purely characteristic $p$ phenomenon: they do not have any flat lift to characteristic $p^2$.

Let $k$ be a field of characteristic $p > 0$. Let $G$ be a reductive group over $\Z$, and write $G_{W_2(k)} = G \otimes_{\Z} W_2(k)$.

\begin{prop}[{\cite[Example A]{oort1968deformations}, \cite[Proposition 4]{schröer2021deligneillusietheoremexceptionalenriques}}]\label{alphap}
    For any $0< r < +\infty$, the Frobenius kernel $\alpha_{p^r}$ does not admit a flat group scheme lift to $W_2(k)$.
\end{prop}

We deduce the analogous result for non-reduced parabolic subgroups. A parabolic subgroup (reduced or not) over $k$ has an obvious lift to $W_2(k)$ as a flat subscheme of $G$, and \Cref{no lift} implies directly that this obvious subscheme lift cannot be equipped with a group scheme structure. However a priori other subscheme lifts could exist, so one has to be slightly careful in the proof. I would like to thank Piotr Achinger for pointing out a mistake in a previous version of the proof, and for helping to correct it.

\begin{prop} \label{no lift}
    Let $P$ be a parabolic subgroup of $G_k$.
    If $P$ admits a lift to $W_2(k)$ as a flat subgroup scheme of $G_{W_2(k)}$, then $P$ is a reduced parabolic subgroup.
\end{prop}

\begin{proof}
    Let $\tilde{P}$ be a flat lift of $P$ to $W_2(k)$ as a subgroup scheme of $G_{W_2(k)}$. Write $X = G/P$, a smooth variety, and $\tilde{X} = G_{W_2(k)}/\tilde{P}$, a flat lift of $X$, hence also smooth by \cite[Tag 01V4]{stacks-project}. Let $T$ be a maximal torus of $G$ contained in $P$, let $\tilde{T}$ be a maximal torus of $G_{W_2(k)}$ which lifts $T$. 
    The subgroup $P$ is the stabilizer of the base point $x_P \in X(k)$ and $T \subset P$ so $x_P \in X^T(k)$.
    By \cite[9.8]{SGA3XII}, the fixed point subscheme $\tilde{X}^{\tilde{T}}$ is smooth. Moreover it is a lift of $X^T$, so the point $x_P \in X^T(k)$ lifts to a point $y \in \tilde{X}^{\tilde{T}}(W_2(k))$. Now $\tilde{P}$ is flat so the quotient map $\pi : G_{W_2(k)} \to G_{W_2(k)}/\tilde{P}$ is faithfully flat, hence the stabilizer $\Stab_{y}$ is a flat subgroup of $G_{W_2(k)}$, and it is a lift of $\Stab_{x} = P$. Moreover, $\tilde{T} \subset \Stab_{y}$. Thus, up to replacing $\tilde{P}$ by $\Stab_y$, we can suppose that $\tilde{T} \subset \tilde{P}$. 

    The choice of a maximal torus $\tilde{T}$ determines root subgroups $(\tilde{U}_{\beta})_{\beta \in R}$ which are lifts of the root subgroups $(U_{\beta})_{\beta \in R}$ by \cite[Theorem 4.1.4]{Conrad2014Reductive}, and we can choose a Borel subgroup $\tilde{B}$ and a positive unipotent subgroup $\tilde{U}^+$ of $G_{W_2(k)}$. Now let us prove that $\tilde{B} \subset \tilde{P}$. We will prove that if $U_{\beta} \subset P$, then $\tilde{U}_{\beta} \subset \tilde{P}$. Given that $\tilde{H} := \tilde{U}_{\beta} \cap \tilde{P}$ is a lift of $U_{\beta} \cap P = U_{\beta}$, it suffices to prove that $\tilde{H}$ is flat. Suppose that $\tilde{H}$ is not flat. The torus $\tilde{T}$ acts by conjugation on $\tilde{P}$, hence on $\tilde{H}$. Thus, $O(\tilde{H}) = W_2(k)[X_{\beta}]/J$, where $J$ is a $T_A$-homogeneous ideal. But $\tilde{H}$ reduces mod $p$ to $U_{\beta}$, that is, $O(H) = k[X_{\beta}]$. Hence, $J$ consists in a relation $pX^n = 0$. Write $Q = P \cap BU^+
    $ and $\tilde{Q} = \tilde{P} \cap \tilde{B}\tilde{U}^+$. $\tilde{B}\tilde{U}^+$ is an open (big cell) of $G_{W_2(k)}$ so $\tilde{Q}$ is also flat. Now $O(\tilde{Q}) = W_2(k)[(X_{\gamma})_{\gamma}]/I$ where $I$ is a $\tilde{T}$-homogeneous ideal. The relation $pX^n = 0$ in $O(\tilde{H})$ lifts to a relation $pX^n + \sum_i c_iM_i = 0$ in $O(\tilde{Q})$ where, for all $i$, $M_i$ is a monomial of weight $n\beta$, hence is a product of several variables. Now, given the structure of the parabolic subgroup $P$ from \Cref{structure parabolic subgroups}, we know that the only relations in $O(Q)$ have the form $X_{\gamma}^m = 0$ (a subgroup $(U_{\gamma})_{r}$ inside $P$ corresponds to a relation $X_{\gamma}^{p^r} = 0$, and a subgroup $U_{\gamma}$ without intersection with $P$ corresponds to a relation $X_{\gamma} = 0$). In particular, there is no relation involving products of several variables in $O(Q)$. By comparing with the reduction of the relation $pX^n + \sum_i c_iM_i = 0$ from $O(\tilde{Q})$ to $O(Q)$, we see that $\sum_i c_iM_i $ necessarily reduces to $0$, hence $\sum_i c_iM_i = pF$ for some $F \in W_2(k)[(X_{\gamma})_{\gamma}]$. Thus, we get a relation $p(X^n + F) = 0$ which implies that $\tilde{Q}$ is not flat, which is a contradiction. Hence, we finally get that $\tilde{B} \subset \tilde{P}$.

    Then, let us show that $\tilde{P} \cap \tilde{U}^+$ is flat. First, let us show that the multiplication induces an isomorphism of schemes $(\tilde{P} \cap \tilde{U}^+) \times \tilde{B} \iso \tilde{P} \cap \tilde{B}\tilde{U}^+$. Indeed, for any $W_2(k)$-algebra $R$, take $u \in \tilde{U}^+(R)$ and $b \in \tilde{B}(R)$ such that $ub \in \tilde{P}(R)$. Then $u = (ub)b^{-1}$ where $ub \in \tilde{P}(R)$ and $b^{-1} \in \tilde{B} \subset \tilde{P}$ so $u \in \tilde{P} \cap \tilde{U}$. But we know that $\tilde{P} \cap \tilde{B}\tilde{U}^+$ is flat, hence $(\tilde{P} \cap \tilde{U}^+) \times \tilde{B}$ is flat. Moreover the Borel subgroup $\tilde{B}$ is faithfully flat over $W_2(k)$, so by \cite[Tag 0584]{stacks-project},  $\tilde{P} \cap \tilde{U}^+$ is flat.

    Conjugation by $\tilde{T}$ on $\tilde{P} \cap \tilde{U}^+$ corresponds to a $\Z^r$-grading on the Hopf algebra $O(\tilde{P} \cap \tilde{U}^+)$. Let $\al$ be a simple root. Write $\tilde{H} = \tilde{P} \cap \tilde{U}^+$ and $\tilde{T}_{\alpha} := \ker(\alpha)$. We have that $\tilde{U}_{\alpha} = (\tilde{U}^+)^{\tilde{T}_{\alpha}}$, hence that $\tilde{H} \cap \tilde{U}_{\al} = \tilde{H}^{\tilde{T}_{\alpha}}$. By the proof of \cite[Theorem 9.1.1]{GilleReductive}, $\tilde{H}^{\tilde{T}_{\alpha}} = O(\tilde{H})/J$ where $J$ is generated by the $O(\tilde{H})_{\chi}$ for $\chi \notin \Z\alpha$. Since $\tilde{H} \subset \tilde{U}^+$, only positive roots occur, and given that $\alpha$ is simple, this gives $$O(\tilde{H}^{\tilde{T}_{\al}}) \simeq \bigoplus_{n \geq 0} O(\tilde{H})_{n\alpha},$$ that is $$O(\tilde{P} \cap \tilde{U}_{\alpha}) \simeq \bigoplus_{n \geq 0} (O(\tilde{P} \cap \tilde{U}^+))_{n\al}.$$ This is a direct summand of $O(\tilde{P} \cap \tilde{U}^+)$, which is flat over $W_2(k)$, so it is also flat over $W_2(k)$, hence the subgroup scheme $\tilde{P} \cap \tilde{U}_{\al}$ is flat over $W_2(k)$. But $\tilde{P} \cap \tilde{U}_{\al}$ is a lift of $P \cap U_{\alpha}$, which is isomorphic to some group $\alpha_{p^r}$ for $0 \leq r \leq +\infty$.
    By \Cref{alphap}, either $P \cap U_{\al} = \{1\}$ or $U_{\al} \subset P$. It means that $P$ is a reduced parabolic subgroup.
\end{proof}

\subsection{Deformations of homogeneous varieties with non-reduced parabolic stabilizer} \label{section final conclusion}

In \cite{zdanowicz2018arithmetically}, Zdanowicz studies a specific class of homogeneous varieties with non-reduced stabilizer: the examples where the acting group is $G = SL_n$, and where parabolic subgroups correspond to $W$-functions $f$ with characteristic filtration $S \backslash \{\al_1, \al_k\} \subset S \backslash \{\al_k\} \subset S$ for $2 \leq k \leq n-2$, and with values $0 < 1 < +\infty$ (see \Cref{example zdanowicz}). In those cases, he proves that the variety does not lift to characteristic $p^2$.
We extend Zdanowicz's result to any reductive group $G$ with connected Dynkin diagram and to any parabolic subgroup $P$ of $G$ of width at least $2$, with only a technical restriction on parabolic subgroups (that we call \emph{special}). It gives a large class of examples of smooth projective varieties which do not lift to characteristic $p^2$. We get a dichotomy result on deformations of varieties which are projective homogeneous under the action of a reductive group with connected Dynkin diagram: if the stabilizer has width at most $1$, then the variety has unobstructed deformations; if the stabilizer has width at least $2$ and is not special, then the variety does not lift to characteristic $p^2$. 
By \Cref{Frob split}, this dichotomy is the same as the dichotomy between homogeneous varieties which are Frobenius split and those which are not. We conjecture that the technical restriction is only technical, in the sense that it could be removed with more work.

The idea of the proof is to show that a deformation of a variety $X = G/P$ induces most of the time a deformation of the $G$-action, which itself induces a deformation of the stabilizer $P$ of the action, and conclude using the fact that a non-reduced parabolic subgroup $P$ does not lift to characteristic $p^2$ as a subgroup of $G$ (\Cref{no lift}). One simply has to be careful with the fact that when the minimum altitude of $P$ is not zero, a deformation of the variety does \textit{not} induce a deformation of the action. However, by \Cref{r0}, to compute the quotient variety, one can always reduce to a case where the parabolic subgroup has minimal altitude $0$.

Showing that a group action lifts when the variety lifts is a common tool to prove obstruction of deformations. For instance, it has been done by Serre in \cite{Serre1961} to construct the first example of smooth projective variety which does not lift to characteristic $0$, or by Ekedahl to prove that Hirokado's threefold does not lift to characteristic $p^2$ \cite[Proof of Theorem A, p. 12]{ekedahl2004nonliftable}. However, what makes the situation more intricate here is that we want to consider deformations of group scheme actions instead of abstract group actions.

\begin{notation}
    Let $f$ be a $W$-function with characteristic filtration $I_l \subset \dots \subset I_1 \subset I_0 = S$.
    The parabolic subgroup $P_f$ is said to be special if $I_1$ is a connected subgraph of the Dynkin diagram of cardinality $|S|-1$.
\end{notation}

\begin{prop} \label{assumptions verified}
    Let $G_k$ be a reductive group with connected Dynkin diagram, let $P$ be a parabolic subgroup of $G_k$, and let $X = G_k/P$.
    Write $\overl$ for the smallest reduced subgroup of $G_k$ which contains $P$, and write $I$ for the subset of the set of simple roots $S$ which generates $\overl$.
    Suppose that $|I| \leq |S|-2$ or, if $|I| = |S|-1$, that $I$ is not connected as a subgraph of $S$.
    Then the forgetful map of deformation functors $\Defcl_{G_k \acts X} \to \Defcl_X$ is unobstructed.
\end{prop}

\begin{proof}
    By \Cref{conditions fulfilled}, $H^1(X, T_X)$ is a trivial $G_k$-module and $H^2(G_k, H^0(X, T_X)) = 0$ (the condition $|I| < |S|$ implies that $r_0 = 0$ in \Cref{conditions fulfilled}). Moreover, by \Cref{no extension}, $H^1(G_k, k) = 0$. So we can apply \Cref{criterion} to conclude.
\end{proof}

\begin{thm} \label{conclusion finale}
    Let $G_k$ be a reductive group with connected Dynkin diagram, let $P$ be a parabolic subgroup of $G_k$, and let $X = G_k/P$. Deformations of $X$ exhibit a dichotomous behaviour:
    \begin{enumerate}
    
    \item If $P$ has width at most $1$, then the variety $X$ has unobstructed deformations.

    \item If $P$ has width at least $2$ and is not special, then the variety $X$ does not lift to characteristic $p^2$.\end{enumerate}
\end{thm}

\begin{proof}
    Write $P = P_f$ for some non-infinitesimal $W$-function $f$. By \Cref{r0}, we know that $G_k/P_f \simeq G_k/P_{f-r_0}$. Denote $g = f-r_0$. \begin{enumerate}

    \item If $P$ has width at most $1$, then $P_g$ is reduced so it is known that $G_k/P_g$ has unobstructed deformations (it is also a side result of this work).

    \item If $P$ has width at least $2$ and is not special, then $\overl[P_g]$ also has width at least $2$ and is also not special. Moreover $\overl[P_g]$ is a proper subgroup of $G_k$. By \Cref{assumptions verified}, the forgetful map of formal moduli problems $\Defcl_{G_k \acts X} \to \Defcl_X$ is unobstructed, and by \Cref{unobstructedness transitive} the composite $\Defcl_{P \sto G_k} \to \Defcl_{G_k \acts X} \to \Defcl_X$ is unobstructed. But, by \Cref{no lift}, $P$ does not lift to characteristic $p^2$ as a subgroup of $G_k$. Hence, by \Cref{smoothness transfer}, the variety $G_k/P$ does not lift to characteristic $p^2$.
\end{enumerate}
\end{proof}

\begin{rk}
    This shows that if $P$ has width $1$ and $r_0 \geq 1$, the variety $G_k/P$ has unobstructed deformations but the subgroup $P$ has obstructed deformations, hence that the forgetful map of deformation functors $\Defcl_{G_k \acts X} \to \Defcl_X$ is obstructed, as hinted by \Cref{map obstructed}.
\end{rk}

\begin{rk}
    Notice that the dichotomy of behaviour of \Cref{conclusion finale} is the same as the dichotomy between Frobenius split and non-Frobenius split homogeneous varieties, by \Cref{Frob split}.
\end{rk}

\begin{conj}
    If $P$ has width at least $2$ and is special, the variety $X = G_k/P$ does not lift to characteristic $p^2$.
\end{conj}

\bibliographystyle{amsalpha}
\bibliography{biblio}

\end{document}